\documentclass[a4paper,10pt]{amsart}

\usepackage{systeme}
\usepackage{graphicx}
\usepackage[]{units}
\usepackage{color}
\usepackage[dvipsnames]{xcolor}
\usepackage{mathrsfs}
\usepackage{bm}
\usepackage[ansinew]{inputenc}
\usepackage[breaklinks,hidelinks]{hyperref}
\usepackage{float}
\usepackage{amsmath}
\usepackage{isomath}
\usepackage{amsthm}
\usepackage{amssymb}
\usepackage{amssymb}

\usepackage[foot]{amsaddr}
\usepackage[margin=2.7cm]{geometry}
\newtheorem{theorem}{Theorem}[section]

\theoremstyle{remark}
\newtheorem{rmk}[theorem]{Remark}

\usepackage{pgfplots}
\usepackage{tikz}
\usepackage{tkz-euclide}
\usepackage{geometry}
\usepackage{epstopdf}
\usepackage{epsfig}
\usepgfplotslibrary{fillbetween}

\newcommand{\RA}[1]{{\color{black}#1}}

\newcommand{\RBC}[1]{{\color{black}#1}}

\usepackage{hyperref}

\newcommand{\blue}{\color{black}}
\newcommand{\orange}{\color{black}}
\usepackage{bm}
\usepackage{graphicx,color, amsfonts,amsmath}

\usepackage{mwe}
\usepackage{subcaption,float}

\DeclareSymbolFont{matha}{OML}{txmi}{m}{it}
\DeclareMathSymbol{\varv}{\mathord}{matha}{118}

\begin{document}

\title[]{A Reduced-Order Shifted Boundary Method for \\ Parametrized incompressible Navier-Stokes equations}

\author{Efthymios N. Karatzas\textsuperscript{1,3}}
\address{\textsuperscript{1}SISSA, International School for Advanced Studies, Mathematics Area, mathLab Trieste, Italy.}
\email{karmakis@math.ntua.gr \& efthymios.karatzas@sissa.it}
%
\author{Giovanni Stabile\textsuperscript{1}}
\email{gstabile@sissa.it}
\author{Leo Nouveau \textsuperscript{2}}
\email{leo.nouveau@duke.edu}
\address{\textsuperscript{2}Civil and Environmental Engineering, Duke University, Durham, NC 27708, United States.}

\address{\textsuperscript{3}Department of Mathematics, School of Applied Mathematical and Physical Sciences, National Technical University of Athens, Zografou 15780, Greece.}

\author{Guglielmo Scovazzi\textsuperscript{2}}
\email{guglielmo.scovazzi@duke.edu}

\author{Gianluigi Rozza\textsuperscript{1}}
\email{grozza@sissa.it}
%
%
%

%

%

\begin{abstract}
We investigate a projection-based reduced order model of the steady incompressible Navier-Stokes equations for moderate Reynolds numbers. 
In particular, we construct an ``embedded'' reduced basis space, by applying proper orthogonal decomposition to the Shifted Boundary Method, a high-fidelity embedded method recently developed. We focus on the geometrical parametrization through level-set geometries, using a fixed Cartesian background geometry and the associated mesh. This approach avoids both remeshing and the development of a reference domain formulation, as typically done in fitted mesh finite element formulations. Two-dimensional computational examples for one and three parameter dimensions are presented to validate the convergence and the efficacy of the proposed approach.
\end{abstract}
\maketitle
\section{Introduction}\label{sec:intro}
We consider a nonlinear system arising from computational fluid dynamics problems. In particular, a stationary Navier--Stokes system is examined and solved by the Shifted Boundary Method (SBM), an embedded boundary finite element method (EBM) that was recently proposed. The geometrical parametrization is described by means of level sets defined over a fixed (undeformed) background mesh.  
Moreover, a reduced order method (ROM), based on a proper orthogonal decomposition approach (POD-Galerkin) is tested, with the purpose of creating an appropriate embedded reduced order basis and decreasing the computational cost in the numerical solution procedure.

Starting with the pioneering work of Peskin \cite{PESKIN1972252}, the scientific community has shown great interest on embedded or immersed methods. Some recent developments are represented by Ghost-Cell Finite {\orange Difference Methods~\cite{wu_faltinsen_chen_2013}, Cut-Cell Finite Volume Methods~\cite{PASQUARIELLO2016670, TUCKER2000591}, Immersed Interface Methods~\cite{Kolahdouz2018AnII,LI2001822}, the Ghost Fluid Method~\cite{BO2014113, Fedkiw2001}, the Volume Penalty Method~\cite{Maury2009,Shirokoff2015, Maury2008}, Cut-FEMs~\cite{Bu_et_all_14}} and the Shifted Boundary Method (SBM)~\cite{MaSco17_2,MaSco17_3,MaSco17_4}. The interested reader can find additional details in~\cite{MiIa05,SOTIROPOULOS20141,KIM2019301,GeSot07NS,ReEtallimmersed_review2019, CaShuJia19_fsiGhostpoint} and references therein. 
Later on, in the context of Navier--Stokes equations, Cut-FEMs~\cite{FRACHON201977, CLAUS2019185} and  unfitted discontinuous Galerkin method~\cite{HeEngIppBa2013_unfittedDG} were employed. 
{
Recently, the authors of~\cite{MaSco17_2} introduced a new unfitted finite method called the Shifted Boundary Method (SBM) for the Poisson and Stokes problems. Their work circumvented the small cut-cell problem by defining a surrogate domain consisting of all uncut elements of the mesh that lie inside the computational domain. The definition of the surrogate domain, however, poses the challenge of devising effective strategies to the imposition of boundary conditions in order to maintain optimal convergence rates. To that end, weakly enforcing a Taylor extrapolation of the numerical solution at the surrogate boundary had the effect of correcting the boundary condition to account for the discrepancy between the surrogate and true boundary. Numerical results for the Poisson and Stokes problems showed optimal convergence rates and adequate conditioning of the matrix. Afterwards, the SBM was extended to the advection/diffusion equation and the laminar and turbulent incompressible Navier-Stokes equations ~\cite{MaSco17_3} in addition to hyperbolic systems such as wave propagation problems in acoustics and shallow water flows ~\cite{MaSco17_4}. The authors of~\cite{nouveau2019high} proposed a number of strategies to increase the accuracy of the method in elliptic problems and the authors of ~\cite{atallah2020analysis} presented a complete numerical analysis of the SB approach for the Stokes problem. Important advances on the original method was presented in~\cite{atallah2020second}, which contained many simplifications and improvement to the original SBM formulation, and~\cite{TheoreticalPoissonAtallahCanutoScovazzi2020}, which addresses numerical convergence and stability in the case of general domains with corners and edges. In~\cite{li2020shifted}, the SBM was extended to problems with internal interfaces.
}

Although immersed and embedded methods show improved features than fitted mesh methods in the case of geometrical design changes, there are still many cases where the approximate solutions of partial differential equations become {\orange computationally unaffordable}, e.g. in real time problems, uncertainty or parametrization of the geometry etc. In these cases {\orange{\it{reduced order modeling}}} techniques appear beneficial, \cite{HeRoSta16,quarteroniRB2016,ChinestaEnc2017,BeOhPaRoUr17}.

The main goal of this paper is to show how the SBM can solve geometrically parametrized nonlinear partial differential systems within the ROM framework. For this purpose, recently developed POD techniques~\cite{KaBaRO18, KaratzasStabileAtallahScovazziRozza2018, KaratzasStabileNouveauScovazziRozza2018,KaratzasRozzaCH2019} will be applied. \RBC{In these seminal works the methodology was developed on simple linear problems such as the Poisson and the Stokes setting. The main novelty of the current article is the extension of the methodology to more complex and nonlinear settings such as the Navier-Stokes one. Such an extension due to the non-linearity of the problem is not trivial. The ROM formulation had in fact to be reshaped and adapted to the underlying iterative approach used by the FOM solver. As it will be more clear in section \ref{sec:HF} the non-linearity is resolved using a Newton's method. Therefore, at each iteration, it is required the computation of the Jacobian of the residual vector and an incremental solution snapshot is obtained by the resolution of the resulting linear system of equations. The same iterative procedure has been implemented also at the ROM level during the online phase. The problem is therefore written in incremental form and the POD basis functions are computed starting from the incremental snapshots and not from the converged ones. The ROM is obtained by the projection of the Jacobian matrix and of the residual vector onto the POD incremental basis functions in order to obtain an incremental reduced solution.} The key feature of our approach is the avoidance of a remeshing stage and/or morphing (i.e., a mapping of all the deformed geometries to a reference domain, see e.g. \cite{HeRoSta16,RoVe07,Rozza2009,ballarin2015supremizer,RoHuMa13,Rozza2008229,BeOhPaRoUr17} for the use of this strategy in traditional body-fitted mesh finite element methods). 
This contribution is organized as follows:
\begin{enumerate}
\item In Section \ref{sec:HF} we define the continuous strong formulation of the mathematical problem and the Nitsche weak formulation. A discrete Shifted Boundary weak formulation is also presented, for the full-order discretization, together with an incremental iterative scheme needed to solve the high fidelity problem during the offline stage.
\item We introduce the reduced order model formulation, the Proper Orthogonal Decomposition and its main features in Section \ref{sec:ROM}.
\item In section \ref{sec:num_exp}, the proposed ROM-SBM technique is tested on a geometrically parametrized problem of the flow around an embedded rectangular domain, and convergence results, errors and execution times are also reported.
\item Finally, in Section \ref{sec:conclusions}, conclusions and perspectives for future improvements and developments are introduced.
\end{enumerate}

\section{The mathematical model and the full-order approximation} \label{sec:HF}
\subsection{Strong formulation of the steady Navier-Stokes problem}
The stationary Navier-Stokes equations for viscous incompressible flow 
describe the flow of a Newtonian, incompressible viscous fluid in a domain when convective forces are not negligible with respect to viscous forces. Consider an open domain $\mathcal D$ in ${\mathbb R}^d$, with $d=2,3$ the number of space dimensions with Lipschitz boundary $\Gamma$, decomposed into two sub-boundaries $\Gamma_D$, $\Gamma_N$.  Let $\mathcal P \subset \mathbb{R} ^k$ be a $k-$dimensional parameter space and $\mu\in \mathcal P$ a parameter vector. 

The strong form of the stationary Navier-Stokes equations with Dirichlet and Neumann boundary conditions on $\Gamma_D$ and $\Gamma_N$, respectively is given by:
\begin{eqnarray}
\rho\nabla\cdot ({\bm{u}}(\mu)\otimes{\bm{u}}(\mu))
+\nabla p(\mu) -\nabla \cdot(2\nu{\bm{\epsilon}} ({\bm{u}}(\mu)) ) - \rho {\bm{g}}(\mu) &=& 0 ,  \text{ in } {\mathcal D}(\mu),
\\
\nabla\cdot {\bm{u}}(\mu)&=& 0,  \text{ in } {\mathcal D}(\mu),
\\
{\bm{u}}(\mu)&=& {\bm{g}}_D(\mu),    \text{ on } \Gamma _D(\mu),
\\
{-(\rho {\bm{u}}(\mu) \otimes {\bm{u}}(\mu)\chi_{\Gamma^-_N(\mu)} +p(\mu){\bm{I}} - 2 \nu{\bm{\epsilon}}({\bm{u}}(\mu)) ) {\bm{n}} }&=& {\bm{g}}_N(\mu),
{  \text{ on }  \Gamma_N(\mu) ,}
\end{eqnarray}
where the variable $\mu$ introduces a geometrical parameterization.
We denote by $\nu$ the viscosity, $\rho$ the density, ${ \bm  \epsilon ( {\bm u}{(\mu)})} = 1/2\left( \nabla   {  {\bm u}(\mu)} + \nabla {{\bm u}({\mu})}^T \right)$ the velocity strain tensor (i.e., the symmetric gradient of the velocity  ${\bm u}(\mu)$), $p(\mu)$ the pressure, $ {\bm g}(\mu)$ a body force, $ {\bm  g_{D}}(\mu)$ the values of the velocity on the Dirichlet boundary and $  {\bm g_{N}}(\mu)$ the normal stress on the Neumann boundary.

We will use for inflow and outflow boundaries  the notation $\Gamma^-_D(\mu)=\{{\bm{x}}\in \Gamma _D(\mu) |{{\bm{g}}}_D(\mu)\cdot {\bm{n}}<0\}$, $\Gamma^-_N(\mu)=\{{\bm{x}}\in\Gamma_N(\mu)|{\bm{u}}(\mu)\cdot {\bm{n}}<0\}$, and $\Gamma^+_D(\mu)={{\RA{\Gamma_D(\mu)\setminus \Gamma^-_D(\mu)}}}$ and $\Gamma^+_N(\mu)={{\RA{\Gamma_D(\mu)\setminus \Gamma^-_N(\mu)}}}$, 
where ${\mathbb{R}}\text{e}[\nu]=\rho UL/\nu$ is the Reynolds number, $U$ and $L$ are the characteristic speed and length of the problem, $h$ is the characteristic mesh
size,  and $\gamma_{N1}$ and $\gamma_{N2}$ are penalty parameters,
with $\chi_{\Gamma^-_D}(\mu)$ and $\chi_{\Gamma^-_N}(\mu)$ the characteristic functions of the boundaries $\Gamma^-_D(\mu)$ and $\Gamma^-_N(\mu)$, respectively.
\subsection{The conformal Nitsche's weak formulation}
For the sake of simplicity, in this subsection, we will omit the parameter dependency with respect {\orange to $\mu$. 
 Let} ${\boldsymbol{V}_h} ({\mathcal D}(\mu))$ and $Q_h ({\mathcal D}(\mu))$ be the spaces of continuous, piecewise-linear, vector- and scalar-valued functions. Namely:
$$
\bm V_h({\mathcal D}(\mu)) = \left\{\bm v \in (C^0 ({{\mathcal{D}}}(\mu)))^{2} : \bm v|_K \in ({P}^1 (K))^{d}, \forall K \in {\mathcal{{{D}}_T}} (\mu)\right\},
$$
$$
Q_h({\mathcal D}(\mu)) = \left\{ {\orange\varv} \in C^0 ({{\mathcal{D}}(\mu)}) : {\orange\varv}|_K \in  {P}^1 (K), \forall K \in {\mathcal{{{D}}_T}}(\mu)\right \},
$$
{\blue{where $\mathcal {D_T}(\mu)$ is the conformal mesh to the true geometry grid.}}
We can now introduce a Nitsche's~\cite{FreSte95,BuFeHa06,BuFe07} variational formulation of the Navier-Stokes equations: \\[.3cm]
{\it find  ${\bm{u}}\in {\boldsymbol{V}_h} ({\mathcal D}(\mu))$ and $p \in Q_h ({\mathcal D}(\mu))$ such that, $\forall {\bm{w}}\in {\boldsymbol{V}_h} ({\mathcal D}(\mu))$ and $\forall q \in Q_h ({\mathcal D}(\mu))$,
 }
\begin{eqnarray}
0 &=& \left( {\bm{w}} , \rho (
{ \bm{u}}\cdot\nabla{ \bm{u}} - {\bm{g}}) \right)_{{\mathcal D}(\mu)} -\left\langle {\bm{w}} , \rho  {{\bm{g}}_D}\cdot{\bm{n}} ({\bm{u}} - {\bm{g}}_D)\right \rangle _{\Gamma^-_D(\mu)}
{-\left({\bm{w}}, {\bm{g}}_N\right)_{\Gamma _N(\mu)}-\left({\bm{w}}, ({\bm{u}}\cdot{\bm{n}}\right)\rho {\bm{u}})_{\Gamma ^-_N(\mu)}}
\nonumber\\
&&- (\nabla\cdot{\bm{w}}, p)_{{\mathcal D}(\mu)} - (q,\nabla\cdot{\bm{u}})_{{{\mathcal D}(\mu)}}  {+\left\langle q ,  ({\bm{u}} - {\bm{g}}_D)\cdot{\bm{n}}\right\rangle _{{\Gamma_D}(\mu)}}
\nonumber\\
&&+\left({\bm{\epsilon(w)}},2\nu{\bm{\epsilon(u)}}\right)_{{\mathcal D}(\mu)}  -\left\langle {\bm{w}}\otimes {\bm{n}} , 2\nu{\bm{\epsilon(u)}}-p{\bm{I}} \right\rangle _{{{\Gamma_D}}(\mu)}
-\left\langle 2\nu{\bm{\epsilon(w)}} , ({\bm{u}} - {\bm{g}}_D)\otimes {\bm{n}}  \right\rangle _{{\Gamma_D}(\mu)}
 \nonumber\\
 &&+\left\langle {\bm{w}} ,
  \gamma \, \nu /h  {\bm{I}} ({\bm{u}} - {\bm{g}}_D)
  \right\rangle _{{\Gamma}_D(\mu)}.
  \nonumber
\end{eqnarray}
Here we use the standard notation $(\cdot, \cdot)_{{\mathcal{D}}}$, $\langle \cdot, \cdot \rangle _{\Gamma_D}$ , $\langle \cdot, \cdot \rangle _{\Gamma_N}$ for the $L^2({\mathcal{D}})$, $L^2(\Gamma_D)$ and $L^2(\Gamma_N)$ inner products over ${\mathcal{D}}$, $\Gamma_D$ and $\Gamma_N$, respectively, and we denote the diameter of element $K\in {\mathcal T}$ as $h_{K}$. The overall size of the mesh is denoted by $h= \max _{K\in{\mathcal T}} h_K$.

\subsection{Shifted boundary variational formulation}
{
In this subsection, we introduce the Shifted Boundary Method for the Navier-Stokes equations~\cite{MaSco17_2,MaSco17_3}. We define a surrogate computational domain $\tilde {\mathcal{D}}$ with surrogate boundary $\tilde \Gamma$ that approximate the true computational domain $\mathcal{D}$ and its boundary $\Gamma$ as seen in Figure~\ref{SurrogateMesh} and Figure~\ref{fig:SBM}. 
Furthermore, $\bm{\tilde{n}}$ indicates the unit outward-pointing normal to the surrogate boundary $\tilde \Gamma$, which differs from the outward-pointing normal  $\bm  n$  to ${\Gamma}$ (see Figure~\ref{fig:ntd}). 
$\tilde  \Gamma$ is composed of the edges/faces of the mesh that are the closest to the true boundary $\Gamma$ in the sense of the closest-point projection, as shown in Figure~\ref{fig:ntd}.

\begin{figure}
\centering
\begin{minipage}{\textwidth}
\centering
\includegraphics[width=0.2\textwidth]{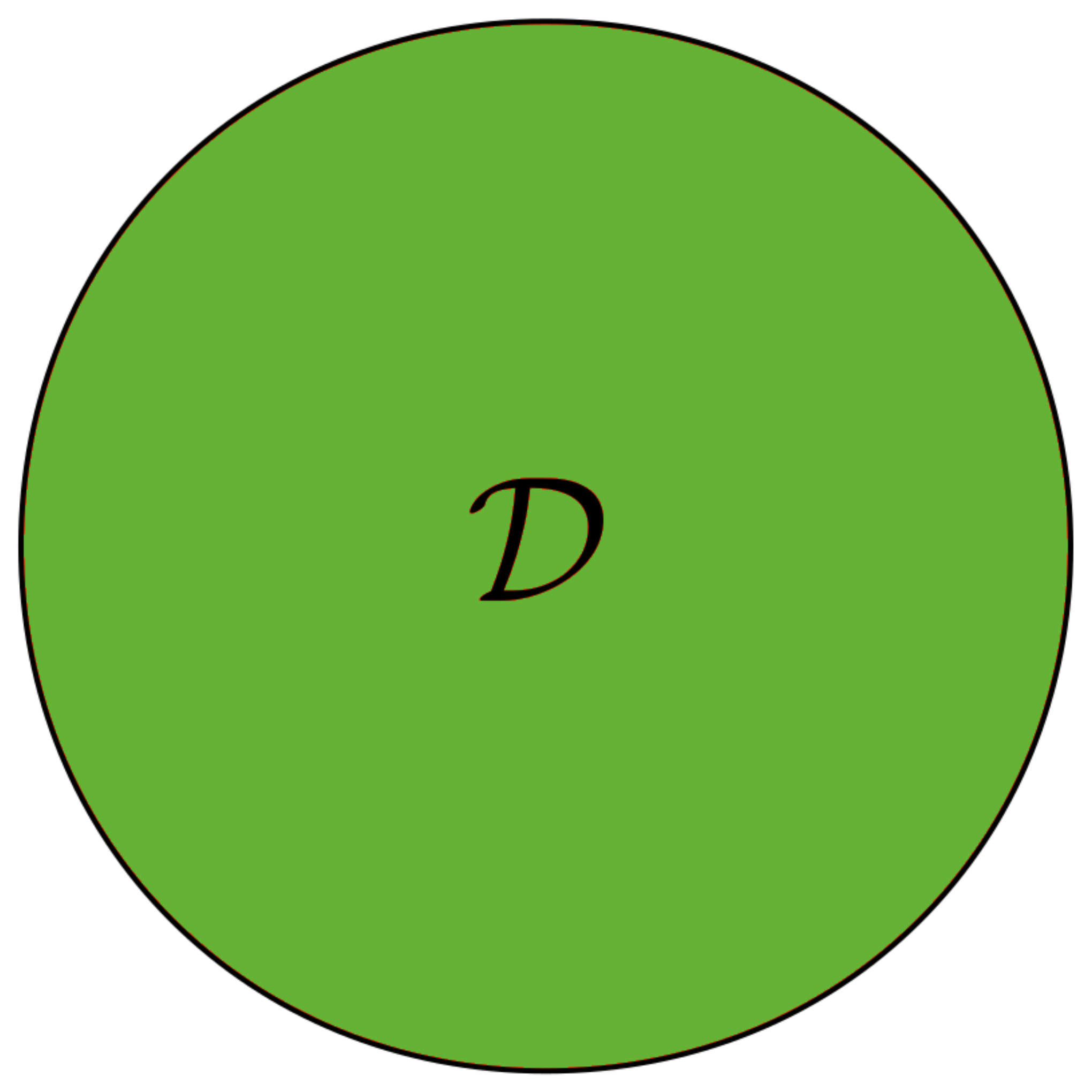}(a)
~~~~~~
\includegraphics[width=0.2\textwidth]{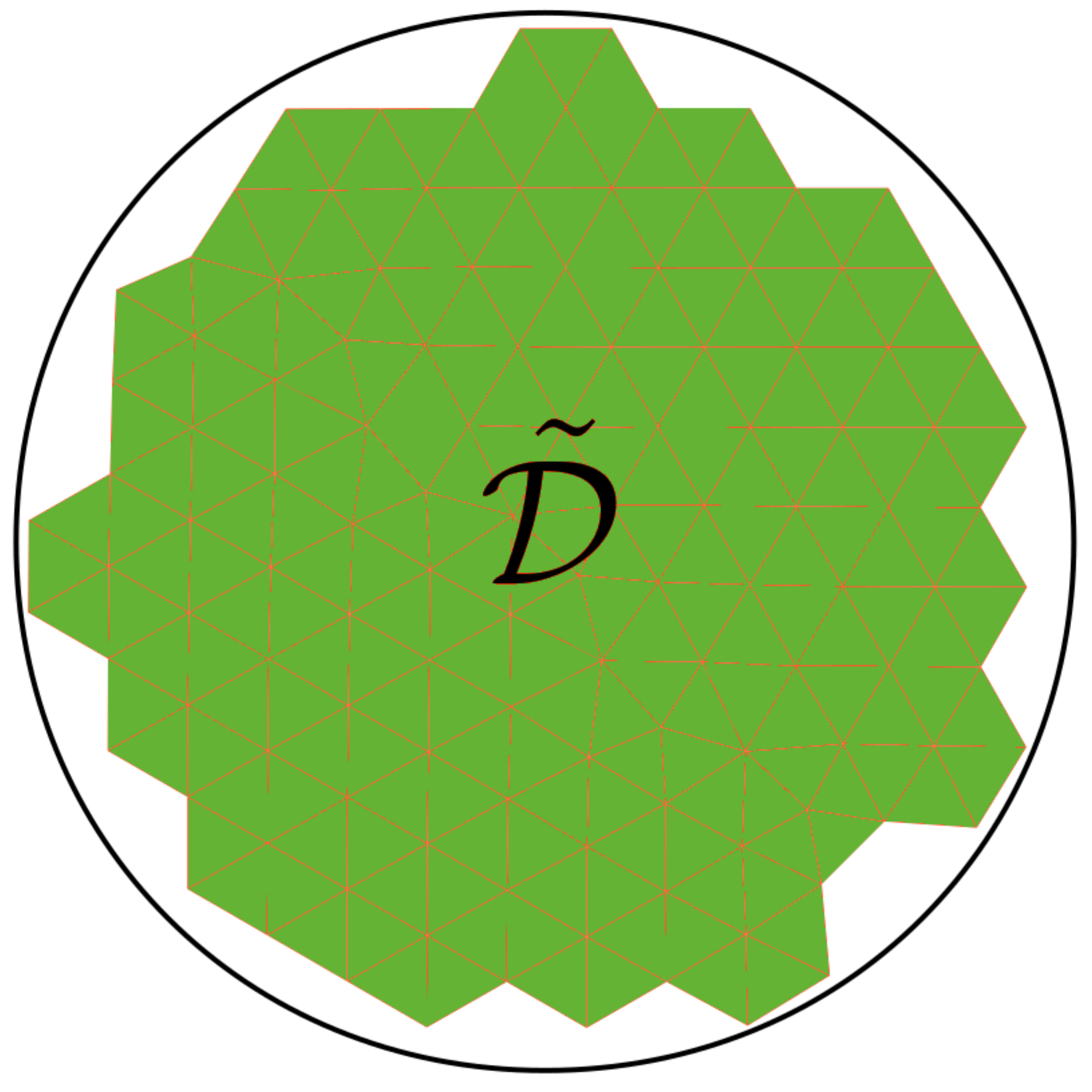}(b)
~~~~~~
\end{minipage}
\caption{
(a) The geometry of a disk and (b) the SBM surrogate geometry.}
\label{SurrogateMesh}
\end{figure}

\newpage

\begin{figure}
\centering
\begin{subfigure}{.4\textwidth}
	\centering
\begin{tikzpicture}[scale=0.63]
\draw [black, line width = 0.95mm, draw=none,name path=surr] plot coordinates { (-2,-3.4641) (-1,-1.73205) (0,-0.5) (1,1.73205) (2.6,2.2) (5,1.73205) (7,2.1) (7.62,0.8) }; 
\draw [blue,  line width = 0.95mm, draw=none,name path=true] plot[smooth] coordinates {(-0.7,-3.4641) (1.75,0.75) (8.25,0.75)}; 
\tikzfillbetween[of=true and surr, split]{gray!15!};
\draw[line width = 0.25mm,densely dashed,gray] (-1,1.73205) -- (0,3.4641);
\draw[line width = 0.25mm,densely dashed,gray] (0,3.4641) -- (2,0);
\draw[line width = 0.25mm,densely dashed,gray] (1,1.73205) -- (1,3.4641);
\draw[line width = 0.25mm,densely dashed,gray] (1,3.4641) -- (0,3.4641);
\draw[line width = 0.25mm,densely dashed,gray] (1,3.4641) -- (2.6,2.2);
\draw[line width = 0.25mm,densely dashed,gray] (1,3.4641) -- (3.25,3.4641);
\draw[line width = 0.25mm,densely dashed,gray] (3.25,3.4641) -- (2.6,2.2);
\draw[line width = 0.25mm,densely dashed,gray] (3.25,3.4641) -- (5,1.73205);
\draw[line width = 0.25mm,densely dashed,gray] (3.25,3.4641) -- (6,3.4641);
\draw[line width = 0.25mm,densely dashed,gray] (6,3.4641) -- (5,1.73205);
\draw[line width = 0.25mm,densely dashed,gray] (6,3.4641) -- (7,2.1);
\draw[line width = 0.25mm,densely dashed,gray] (0,-0.5) -- (-2,0.5);
\draw[line width = 0.25mm,densely dashed,gray] (-2,0.5) -- (-1,1.73205);
\draw[line width = 0.25mm,densely dashed,gray] (-2,0.5) -- (-1,-1.73205);
\draw[line width = 0.25mm,densely dashed,gray] (-2,0.5) -- (-2.5,-2);
\draw[line width = 0.25mm,densely dashed,gray] (-2.5,-2) -- (-1,-1.73205);
\draw[line width = 0.25mm,densely dashed,gray] (-2.5,-2) -- (-2,-3.4641);
\draw[line width = 0.25mm,densely dashed,gray] (0,-0.5) -- (-1,1.73205);
\draw[line width = 0.25mm,densely dashed,gray] (-1,1.73205) -- (1,1.73205);
\draw[line width = 0.25mm,densely dashed,gray] (0,-0.5) -- (2,0);
\draw[line width = 0.25mm,densely dashed,gray] (2,0) -- (1,1.73205);
\draw[line width = 0.25mm,densely dashed,gray] (1,1.73205) -- (0,-0.5);
\draw[line width = 0.25mm,densely dashed,gray] (2,0) -- (2.6,2.2);
\draw[line width = 0.25mm,densely dashed,gray] (2.6,2.2) -- (1,1.73205);
\draw[line width = 0.25mm,densely dashed,gray] (2,0) -- (4,0);
\draw[line width = 0.25mm,densely dashed,gray] (4,0) -- (2.6,2.2);
\draw[line width = 0.25mm,densely dashed,gray] (4,0) -- (2.6,2.2);
\draw[line width = 0.25mm,densely dashed,gray] (2.6,2.2) -- (5,1.73205);
\draw[line width = 0.25mm,densely dashed,gray] (5,1.73205) -- (4,0);
\draw[line width = 0.25mm,densely dashed,gray] (4,0) -- (6,0);
\draw[line width = 0.25mm,densely dashed,gray] (6,0) -- (5,1.73205);
\draw[line width = 0.25mm,densely dashed,gray] (6,0) -- (7,2.1);
\draw[line width = 0.25mm,densely dashed,gray] (7,2.1) -- (5,1.73205);
\draw[line width = 0.25mm,densely dashed,gray] (6,0) -- (8,0);
\draw[line width = 0.25mm,densely dashed,gray] (8,0) -- (7,2.1);
\draw[line width = 0.25mm,densely dashed,gray] (0,-0.5) -- (-1,-1.73205);
\draw[line width = 0.25mm,densely dashed,gray] (-1,-1.73205) -- (1,-1.73205);
\draw[line width = 0.25mm,densely dashed,gray] (2,0) -- (1,-1.73205);
\draw[line width = 0.25mm,densely dashed,gray] (1,-1.73205) -- (0,-0.5);
\draw[line width = 0.25mm,densely dashed,gray] (2,0) -- (3,-1.73205);
\draw[line width = 0.25mm,densely dashed,gray] (3,-1.73205) -- (1,-1.73205);
\draw[line width = 0.25mm,densely dashed,gray] (4,0) -- (3,-1.73205);
\draw[line width = 0.25mm,densely dashed,gray] (2,0) -- (4,0);
\draw[line width = 0.25mm,densely dashed,gray] (4,0) -- (3,-1.73205);
\draw[line width = 0.25mm,densely dashed,gray] (3,-1.73205) -- (5,-1.73205);
\draw[line width = 0.25mm,densely dashed,gray] (5,-1.73205) -- (4,0);
\draw[line width = 0.25mm,densely dashed,gray] (4,0) -- (6,0);
\draw[line width = 0.25mm,densely dashed,gray] (6,0) -- (5,-1.73205);
\draw[line width = 0.25mm,densely dashed,gray] (6,0) -- (7,-1.73205);
\draw[line width = 0.25mm,densely dashed,gray] (7,-1.73205) -- (5,-1.73205);
\draw[line width = 0.25mm,densely dashed,gray] (6,0) -- (8,0);
\draw[line width = 0.25mm,densely dashed,gray] (8,0) -- (7,-1.73205);
\draw[line width = 0.25mm,densely dashed,gray] (0,-3.4641) -- (-2,-3.4641);
\draw[line width = 0.25mm,densely dashed,gray] (-2,-3.4641) -- (-1,-1.73205);
\draw[line width = 0.25mm,densely dashed,gray]  (-1,-1.73205) -- (0,-3.4641);
\draw[line width = 0.25mm,densely dashed,gray] (0,-3.4641) -- (1,-1.73205);
\draw[line width = 0.25mm,densely dashed,gray] (0,-3.4641) -- (2,-3.4641);
\draw[line width = 0.25mm,densely dashed,gray] (2,-3.4641) -- (1,-1.73205);
\draw[line width = 0.25mm,densely dashed,gray] (2,-3.4641) -- (3,-1.73205);
\draw [line width = 0.5mm,blue, name path=true] plot[smooth] coordinates {(-0.75,-3.681818) (1.75,0.75) (7.8,0.75)};
\draw[line width = 0.5mm,red] (1,1.73205) -- (2.6,2.2);
\draw[line width = 0.5mm,red] (2.6,2.2) -- (5,1.73205);
\draw[line width = 0.5mm,red] (5,1.73205) --  (7,2.1);
\draw[line width = 0.5mm,red] (1,1.73205) -- (0,-0.5);
\draw[line width = 0.5mm,red] (0,-0.5) -- (-1,-1.73205);
\draw[line width = 0.5mm,red] (-1,-1.73205) -- (-2,-3.4641);
\node[text width=0.5cm] at (7.6,2.3) {\large${\color{red}\tilde{\Gamma}}$};
\node[text width=3cm] at (2.1,1.25) {\large${\color{red}\tilde{\mathcal{D}}}$};
\node[text width=0.5cm] at (8.35,0.95) {\large${\color{blue}\Gamma}$};
\node[text width=3cm] at (4.5,1.5) {\large${\mathcal{D}} \setminus \tilde{\mathcal{D}} $};
\end{tikzpicture}
\caption{The true domain ${\mathcal{D}}$, the surrogate domain $\tilde{\mathcal{D}} \subset {\mathcal{D}}$ and their boundaries $\tilde{\Gamma}$ and $\Gamma$.}
\label{fig:SBM}
\end{subfigure}
\quad
\begin{subfigure}{.27\textwidth}\centering
	\begin{tikzpicture}[scale=0.82]
\draw[line width = 0.25mm,densely dashed,gray] (0,0.5) -- (-1.5,3);
\draw[line width = 0.25mm,densely dashed,gray] (-1.5,3) -- (0.5,5);
\draw[line width = 0.25mm,densely dashed,gray] (0,0.5) -- (3.5,2);
\draw[line width = 0.25mm,densely dashed,gray] (3.5,2) -- (0.5,5);
\draw [line width = 0.5mm,blue, name path=true] plot[smooth] coordinates {(0.4,-0.5) (2.16,2.55) (1.0,6)};
\draw[line width = 0.5mm,red] (0,0.5) -- (0.5,5);
\node[text width=0.5cm] at (0.5,5.5) {\large${\color{red}\tilde{\Gamma}^h}$};
\node[text width=0.5cm] at (1.75,5.5) {\large${\color{blue}\Gamma}$};
\node[text width=0.5cm] at (1.25,3.25) {\large$\bm{d}$};
\node[text width=0.5cm] at (3,3.5) {\color{blue} \large$\bm{n}$};
\node[text width=0.5cm] at (1.0,2.43) {\color{red} \large$\tilde{\bm{n}}$};
\node[text width=0.5cm] at (2.7,2.25) {\color{blue} \large$\bm{\tau}$};
\draw[->,red, line width = 0.25mm,-latex] (0.25,2.75) -- (1.22,2.63);
\draw[->,line width = 0.25mm,-latex] (0.25,2.75) -- (2.12,3.1);
\draw[->, blue, line width = 0.25mm,-latex] (2.12,3.1) -- (2.28,2.29);
\draw[->, blue, line width = 0.25mm,-latex] (2.12,3.1) -- (2.95,3.25);
\end{tikzpicture}
    \caption{The distance vector $\bm{d}$, the true normal $\bm{n}$ and the true tangent $\bm{\tau}$.}
    \label{fig:ntd}
\end{subfigure}
\quad
\begin{subfigure}{0.27\textwidth}
\centering
\begin{tikzpicture}[scale=0.79]
\draw[line width = 0.5mm,cyan] (0,0) -- (4,0);
\draw[line width = 0.5mm,blue] (4,0) -- (4,-4);
\draw[line width = 0.5mm,red] (3.65,1.75) -- (5.75,-0.5);
\node[text width=0.5cm] at (5,1) {${\color{red} \tilde{E} \subset \tilde{\Gamma}}$};
\draw[->,line width = 0.25mm,-latex] (4.7,0.65) -- (4,0);
\node[text width=0.5cm, rotate=45] at (4.1,0.45) {$\bm{d}(\tilde{\bm{x}})$};
\draw[->,line width = 0.25mm,-latex] (3.65,1.75) -- (3.65,0);
\node[text width=0.5cm, rotate=90] at (3.35,0.85) {$\bm{d}(\tilde{\bm{x}}_{a})$};
\draw[->,line width = 0.25mm,-latex] (5.75,-0.5) -- (4,0);
\node[text width=0.5cm, rotate=-12.5] at (4.6,-0.5) {$\bm{d}(\tilde{\bm{x}}_{b})$};
\draw[->,line width = 0.25mm, cyan, -latex] (2,0) -- (2,-0.75);
\draw[->,line width = 0.25mm, blue , -latex] (4,-2) -- (3.25,-2);
\node[text width=0.5cm] at (0.5,-0.5) {${\color{cyan}\Gamma_{D}}$};
\node[text width=0.5cm] at (1.75,-0.4) {${\color{cyan}\bm{n}_{1}}$};
\node[text width=0.5cm] at (3.6,-3.5) {${\color{blue}\Gamma_{N}}$};
\node[text width=0.5cm] at (3.75,-2.5) {${\color{blue}\bm{n}_{2}}$};
\fill (3.65,1.75) circle (0.5mm);
\node[text width=0.5cm] at (3.85,2.1) {$\tilde{\bm{x}}_{a}$};
\fill (5.75,-0.5) circle (0.5mm);
\node[text width=0.5cm] at (5.95,-0.85) {$\tilde{\bm{x}}_{b}$};
\end{tikzpicture}
\caption{The distance vector $\bm{d}$ may not align with either $\bm{n}_{1}$ or $\bm{n}_{2}$.}
\label{fig:distance_nu}
\end{subfigure}
\caption{
The surrogate domain, its boundary, and the distance vector $\bm{d}$.}
\label{fig:surrogates}
\end{figure}
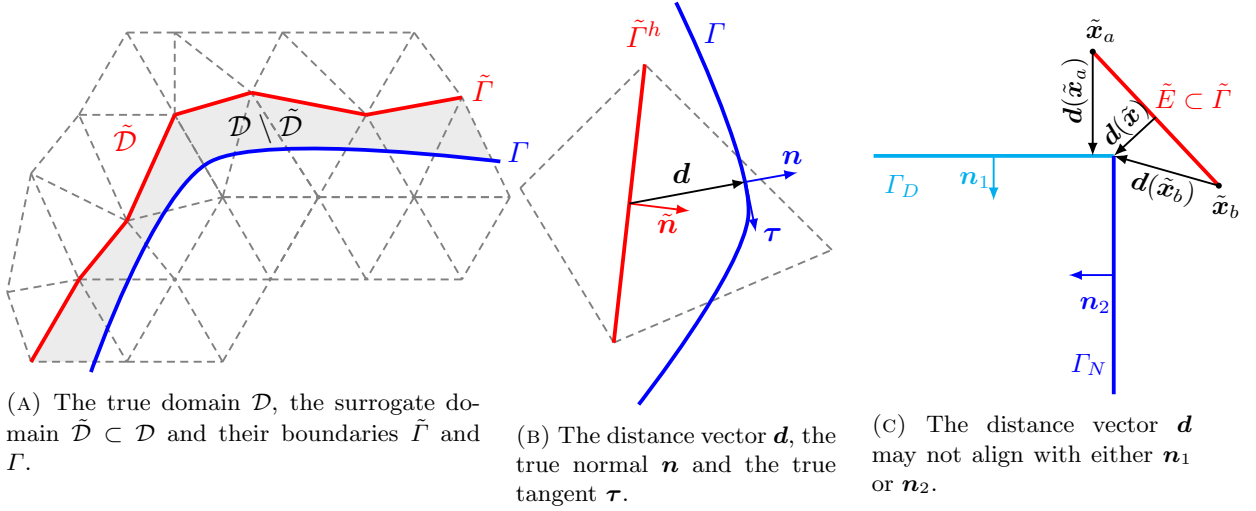
In particular, we consider a mapping
\begin{subequations}\label{eq:defMmap}
\begin{align}
\bm{M}^{h}:&\; \tilde{\Gamma} \to \Gamma \; ,  \\
&\; \tilde{\bm{x}} \mapsto \bm{x}   \; ,
\end{align}
\end{subequations}
as defined in \cite[Section 2.1]{TheoreticalPoissonAtallahCanutoScovazzi2020}, which associates to any point $\tilde{\bm{x}} \in \tilde{\Gamma}$ on the surrogate boundary a point $\bm{x} = \bm{M}^{h}(\tilde{\bm{x}})$ on the physical boundary $\Gamma$.
Through $\bm{M}^{h}$, a distance vector function $\bm{d}_{\bm{M}^{h}}$ can be defined as
\begin{align}
\label{eq:Mmap}
\bm{d}_{\bm{M}^{h}} (\tilde{\bm{x}})
\, = \, 
\bm{x}-\tilde{\bm{x}}
\, = \, 
[ \, \bm{M}^h-\bm{I} \, ] (\tilde{\bm{x}})
\; .
\end{align}
For the sake of simplicity, we set $\bm{d} = \bm{d}_{\bm{M}^{h}} $ where $\bm{d} = \| \, \bm{d} \, \| \bm{\nu}$  and $\bm{\nu}$ is a unit vector.\\
%
%
For smooth surfaces with a single type of boundary condition, $\bm{M}^{h}$ is the closest point projection and therefore $\bm{\nu} = \bm{n}$ (see Figure~\ref{fig:ntd}). On the other hand, when the boundary $\Gamma$ is partitioned into a Dirichlet boundary $\Gamma_{D}$ and a Neumann boundary $\Gamma_{N}$ with $\Gamma = \overline{\Gamma_{D} \cup \Gamma_{N}}$ and $\Gamma_{D} \cap \Gamma_{N} = \emptyset$ (see Figure~\ref{fig:distance_nu}), we need to identify whether a surrogate edge $\tilde{E} \subset \tilde{\Gamma}$ is associated with $\Gamma_{D}$ or $\Gamma_{N}$. Through $\bm{M}^{h}$, we can partition $\tilde{\Gamma}$ as $\overline{\tilde{\Gamma}_D \cup \tilde{\Gamma}_N}$ with $\tilde{\Gamma}_D \cap \tilde{\Gamma}_N = \emptyset$ such that
\begin{align}\label{def:surrogateDir}
\tilde{\Gamma}_D = \{ \tilde{E} \subseteq \tilde{\Gamma} : \bm{M}^{h}(\tilde{E}) \, \subseteq \, \Gamma_{D} \}\end{align}
and $\tilde{\Gamma}_N = \tilde{\Gamma} \setminus \tilde{\Gamma}_D$. In Figure~\ref{fig:distance_nu}, for example, $\bm{M}^{h}$ associates the surrogate edge $\tilde{E}$ with $\Gamma_{D}$ whereby all the boundary information will be transported from $\Gamma_{D}$ to $\tilde{E}$. In turn, the distance $\bm{d}$ is given as the closest point projection from $\tilde{E}$ to $\Gamma_{D}$ which may not align with the true normal to $\Gamma_{D}$. 
For an extensive discussion on how to construct the distance $\bm{d}$ on domains with corners and edges in two and three dimensions, see~\cite{TheoreticalPoissonAtallahCanutoScovazzi2020}.
}

{\RA{Since the true surface is smooth between edges and corners, it is allowed to assume that ${\bm{M}^{h}}$ is continuous, and Lipschitz.}}
Using ${\bm{M}^{h}}$, the unit normal vector $\bm n$ to the boundary $\Gamma$ can easily extend to the boundary $\tilde \Gamma$ as
$\bm \bar {\bm n}(\tilde {\bm x}) \equiv {\bm n}({\bm{M}^{h}}(\tilde {\bm x}))$.
 In the following, whenever we write $\bm n(\tilde {\bm x})$ we actually mean $\bm \bar {\bm n}(\tilde {\bm x})$ at a point ${ \tilde {\bm x}}\in \tilde\Gamma$.
 The above constructions are the basis for the extension of the boundary conditions on $\Gamma$ to the boundary $\tilde  \Gamma$ of the surrogate domain. 

{\RA{We can now introduce the SBM variational formulation~\cite{MaSco17_2,MaSco17_3}. 
Consider the surrogate domain $\tilde{\mathcal{D}}$  and }} a discretization of the continuous boundary value problem with a mesh ${\mathcal{{\tilde{D}}_T}}$ consisting of simplexes $K$ belonging to a tessellation $\mathcal T$. Moreover, we introduce  the discrete spaces ${\bm V}_h(\tilde{\mathcal D}(\mu))$ and $Q_h(\tilde{\mathcal D}(\mu))$, for velocity and pressure and we assume that a stable and convergent base formulation for the  
flow exists for these spaces in the case of conformal grids.
Using the modified continuous piecewise-linear spaces
$$
\bm V_h(\tilde{\mathcal D}(\mu)) = \left\{\bm v \in (C^0 ({\tilde{\mathcal{D}}}(\mu)))^{2} : \bm v|_K \in ({P}^1 (K))^{d}, \forall K \in {\mathcal{{\tilde{D}}_T}} (\mu)\right\},
$$
$$
Q_h(\tilde{\mathcal D}(\mu)) = \left\{ \varv \in C^0 ({\tilde{\mathcal{D}}(\mu)}) : \varv|_K \in  {P}^1 (K), \forall K \in {\mathcal{{\tilde{D}}_T}}(\mu)\right \},
$$
we can now introduce the Shifted boundary variational formulation:
\\[.3cm]
{\it{ find ${\bm{u}}\in {\bm{V}} _h ({\tilde {\mathcal{D}}(\mu)})$ and $p \in Q_h ({\tilde {\mathcal{D}}(\mu)})$ such that, $\forall {\bm{w}} \in {\bm{V}} _h ({\tilde {\mathcal{D}}(\mu)})$ and $\forall q \in Q_h ({\tilde {\mathcal{D}}(\mu)})$,}}
\begin{eqnarray*}
	0 &=& {\mathbb{I_{NS}}}[\nu]({\bm{u}}, p;{\bm{w}},q)
	\nonumber\\
	&=& ({\bm{w}}, 
	\rho({\bm{u}}\cdot {\bm{\nabla u}} - {\bm{b}}))_{\tilde {\mathcal{D}}(\mu)}
	{-\left\langle {\bm{w}}, \rho ({{\bm{g}_D} - (\bm{\nabla u})\bm{d}})\cdot{\bm{\tilde n}} ({\bm{u}} + (\bm{\nabla u})\bm{d} - {\bm{g}_D})\right\rangle_{{\tilde\Gamma^-}_D}}
	{- \left\langle {\bm{w}}, {\bm{g}_N}\right\rangle_{{\tilde\Gamma}_N}}
	%
	\nonumber\\
	&& 
	{- \left\langle {\bm{w}}, {(\bm{u\cdot {\tilde n}})}\rho{\bm{u}}\right\rangle_{{{\tilde\Gamma}^-}_N}}
	- ({\bm{\nabla}}\cdot {\bm{w}}, p)_{{\tilde {\mathcal{D}}(\mu)}}
	-(q,{\bm{\nabla}}\cdot {\bm{u}})_{{\tilde {\mathcal{D}}(\mu)}}
	 {+ \left\langle q,({\bm{u}} + { \chi_{{\tilde{\Gamma}}^+ _D}}
 	({\bm{\nabla u }}){\bm{d}} 
 	{- {\bm{g}}_D})\cdot{\bm{\tilde n}}\right\rangle_{{\tilde\Gamma}_D}
	}
	\nonumber\\
	&&+ ({\bm{\epsilon}}({\bm{w}}),2\nu{\bm{\epsilon}}({\bm{u}}))_{{\tilde {\mathcal{D}}(\mu)}}
	-\left\langle {\bm{w}} \otimes {\bm{\tilde n}},2\nu{\bm{\epsilon}}({\bm{u}}) - p{\bm{I}}
	\right\rangle_{{\tilde\Gamma}_D}
	-\left\langle 2\nu {\bm{\epsilon}}({\bm{w}}),({\bm{u}} + {\chi_{{\tilde{\Gamma}}^+ _D}} ({\bm{\nabla u }}){\bm{d}} {- {\bm{g}}_D}) \otimes {\bm{\tilde n}} \right\rangle_{{\tilde\Gamma}_D}
	\nonumber\\&& 
	+\left\langle {\bm{w}} + { \chi_{{\tilde{\Gamma}}^+ _D}}
	 ({\bm{\nabla w}}){\bm{d}}, \gamma \, \nu/h ({\bm{u}} + {\chi_{{\tilde{\Gamma}}^+ _D}}({\bm{\nabla u}}){\bm{d}} {- {\bm{g}}_D})\right\rangle_{{\tilde\Gamma}_D}
\; ,
\end{eqnarray*}
where we indicated again by $(\cdot, \cdot)_{\tilde{\mathcal{D}}}$, $\langle \cdot, \cdot \rangle _{\tilde\Gamma_D}$ , $\langle \cdot, \cdot \rangle _{\tilde\Gamma_N}$ the $L^2(\tilde{\mathcal{D}})$, $L^2(\tilde\Gamma_D)$ and $L^2(\tilde\Gamma_N)$ inner products over $\tilde{\mathcal{D}}$, $\tilde\Gamma_D$ and $\tilde\Gamma_N$, respectively. 
For more details we refer to \cite{MaSco17_3}.

\subsection{Variational multiscale stabilized finite element formulation}
   Because the proposed variational statement is not numerically stable, we introduce SUPG and PSPG stabilizing operators according to~\cite{HuScoFra04,BroHu_82}, to which the reader can refer for more details. Since this is not the main focus of this paper, we refer the reader to~\cite{MaSco17_3,Stabile2019} for their implementation. The abstract variational form would then read:

   {\it{find ${\bm{u}} \in {\bm{V}}_h (\tilde{\mathcal{D}})$ and $p \in Q_h (\tilde{\mathcal{D}})$ such that, $\forall {\bm{w}} \in {\bm{V}}_h (\tilde{\mathcal{D}})$ and $\forall q \in Q_ h (\tilde{\mathcal{D}})$,
   }}
  \begin{eqnarray*}
0 = {\mathbb {I_{NS}}}[\nu]({\bm{u}}, p;{\bm{w}},q) + {\mathbb {{II}_{STAB}}}[\nu]({\bm{u}}, p;{\bm{w}},q),
\end{eqnarray*}
where ${\mathbb {{II}_{STAB}}}$ $[\nu]({\bm{u}}, p;{\bm{w}},q)$ denotes the SUPG/PSPG stabilization operators.

\bigskip
In what follows, it is useful to define the Navier-Stokes operator 
$$
G(U(\mu))U(\mu) := G\left(  
     \begin{bmatrix} {\bm{u}}(\mu) \\ p(\mu) \end{bmatrix}
  \right)
  =
    \begin{bmatrix} {{\bm{A}}{\bm{u}}(\mu)  + {\bm{\mathcal{C}}}({\bm{u}}(\mu)) +  \bm{B}^Tp(\mu)} \\ ({{\bm{B}} + {\hat{\bm{B}}}}){\bm{u}}(\mu) +{\bm{C}}p(\mu)\end{bmatrix}\\
=
  \begin{bmatrix} \bm{A} + {\bm{\mathcal{C}}}({\bm{u}}(\mu))& \bm{B}^T \\ \bm{B}+ {\hat{\bm{B}}} & {\bm{C}} \end{bmatrix}
  \begin{bmatrix} {\bm{u}}(\mu) \\ p(\mu) \end{bmatrix},
$$
where 
 {\RA{$\bm{A}$ corresponds to the discrete diffusion operator, ${\bm{\mathcal{C}}}$ to the nonlinear convection operator. The rectangular matrix $\bm{B}^T$ represents the discrete gradient operator while $\bm{B}$ represents its adjoint, the divergence operator. $\bm{\hat{B}}$, $\bm{C}$ are instead associated with additional SUPG/PSPG stabilization operators.}}
The right hand side 
$$ 
F(\mu):= \begin{bmatrix} \bm{F}_g(\mu) \\ {\bm F}_q(\mu) \end{bmatrix},
$$
is consisting of forcing and boundary data related to stabilization and Nitsche weak enforcement boundary terms. 
Using these definitions we can express the following residual of  the algebraic system of equations:
$$
R(U(\mu))=G(U(\mu))U(\mu)-F(U(\mu)).
$$
Furthermore, the Jacobian of $G(U(\mu))U(\mu)$ reads
$$
\nabla G(U^{n-1}(\mu)) =  \begin{bmatrix} \frac{\partial(\bm{A}{\bm{u}}(\mu) + {\bm{\mathcal{C}}}({\bm{u}}(\mu)))}{\partial{\bm{u}}}\ & \bm{B}^T \\ {\bm{B}} 
& 0 \end{bmatrix},
$$
and yields the following iterative algebraic system of equations for the increment $\delta U(\mu)=U^{n}(\mu)-U^{n-1}(\mu)$:
\begin{equation}\label{eq:system_linear}
\nabla G(U^{n-1}(\mu))\delta U(\mu) = -R(U^{n-1}(\mu)).
\end{equation}
In the latter system of equations we underline some features that will play important role later in the ROM strategy described in \autoref{sec:ROM}: the discretized differential operators $\bm{A}$, ${\bm{\mathcal{C}}}$, $\bm{B}$ and $\bm{\hat{B}}$ are parameter dependent and in the typical saddle point structure of the problem, the incompressibility equation is partially relaxed adding a stabilization term $\bm{C}$. 
This stabilization term, at full-order level, permits to  the fulfillment of the ``inf-sup'' condition and the use of otherwise unstable pair of finite elements, such as $\mathbb{P}_1-\mathbb{P}_1$ 
it helps to preserve the stability of the reduced order model.

The presented formulation is used to solve the full-order problem during an offline stage and to produce the snapshots necessary for the construction of the ROM of \autoref{sec:ROM}.

\section{Reduced order model with a POD-Galerkin method}\label{sec:ROM}
{{In this section, a POD-Galerkin approach will be analyzed as in \cite{HeRoSta16,Rozza2004ReducedBM}. We will simplify the high fidelity model system to a reduced order one, which preserves its essential properties with the purpose of reducing computational cost in a way adapted to embedded-immersed boundary finite element methods. This approach allows flexibility with geometrical changes and to effectively and efficiently overcome several related issues that appear when using traditional FEM (see for instance~\cite{BaFa2014,KaBaRO18,KaratzasStabileAtallahScovazziRozza2018,KaratzasStabileNouveauScovazziRozza2018}).

Our interest is focused on the nonlinear system of Navier-Stokes equations with parametrized geometries and on the advantages of the SBM. We construct the ROM in a classical way  starting from high dimensional SBM approximations (offline stage which involves the solution of a possibly large number of high fidelity problems). Reduced basis methods have been obtained starting from full-order approximations on fitted mesh methods for non-linear problems~\cite{Veroy2003,Grepl2007,RoQua07} -- while for linear elliptic and linear parabolic equations we refer to~\cite{Rozza2008229,grepl2005}. 

The first step in a reduced order method consists of obtaining a set of high fidelity solutions of the parametrized problem under input parameters variation. The aim of ROMs is to approximate any member of this solution set with a reduced number of basis functions. During the costly {\it{offline stage}}, one works out the solution set and examines its components in order to construct a reduced basis that approximates any member of the solution set to a prescribed accuracy. During a second stage, namely the {\it{online stage}}, after the Galerkin projection of the full-order differential operators describing the governing equations onto the reduced basis spaces, it is possible to solve a reduced problem for any new value of the input parameters  with a reduced computational effort.

In general, POD-Galerkin ROMs for the incompressible Navier-Stokes equations are unstable~\cite{Caiazzo2014598,Gerner2012,RoHuMa13}, due to pressure instabilities, while for dynamic instabilities on transient problems the interested reader could see for instance \cite{Iollo2000,Akhtar2009,Bergmann2009516,Sirisup2005218,taddei2017}. Nevertheless, in the present approach, the SUPG and PSPG stabilization which is applied on the high fidelity solver is strongly propagating through the reduced basis construction procedure to the reduced level, and there is no need for further reduced basis stabilization and supremizers enrichment as in ~\cite{ballarin2015supremizer,RoVe07, stabile_stabilized, StaHiMoLoRo17}, see Appendix A. 
 
The SBM unfitted/surrogate mesh Nitsche finite element method is used to apply parametrization and reduced order techniques considering Dirichlet combined with Neumann boundary conditions. 

We highlight that a parametrized ROM method without the use of the transformation to reference domains will be used taking advantage of the fixed, geometrical parameter independent, background mesh \cite{KaBaRO18, KaratzasStabileNouveauScovazziRozza2018, KaratzasStabileAtallahScovazziRozza2018}.

\subsection{The Proper Orthogonal Decomposition (POD)}\label{subsec_POD_theory}
For the projection of the high fidelity system to the reduced order one, there are several techniques. For more details about the different strategies the interested reader may see \cite{Rozza2008229,ChinestaEnc2017,Kalashnikova_ROMcomprohtua,quarteroniRB2016,Chinesta2011,Dumon20111387}. In the present work, the POD is applied to the parameter-dependent full snapshots matrices. The full-order model is solved for each $\mu^k \in \mathcal{K}=\{ \mu^1, \dots, \mu^{N_k}\} \subset \mathcal{P}$ where $\mathcal{K}$ is a finite dimensional training set of parameters chosen inside the parameter space $\mathcal{P}$. 
The  number of snapshots is denoted by $N_s $ and the number of degrees of freedom for the discrete full-order solution by $N_u^h$,  $N_p^h$ for the velocity and pressure, respectively.  The snapshots matrices $\bm{\mathcal{S}_u}$ and $\bm{\mathcal{S}_p}$, for velocity and pressure, are then given by $N_s$ full-order snapshots:
\begin{gather}
\bm{\mathcal{S}_u} = [\bm{u}(\mu^1),\dots,\bm{u}(\mu^{N_s})] \in \mathbb{R}^{N_u^h\times N_s},\quad\bm{\mathcal{S}_p} = [p(\mu^1),\dots,p(\mu^{N_s})] \in \mathbb{R}^{N_p^h\times N_s}.
\end{gather}
Given a general scalar or vectorial function $\bm{u}:{\mathcal D} \to \mathbb{R}^d$, with a certain number of realizations $\bm{u}_1,\dots, \bm{u}_{N_s}$, the POD problem consists in finding, for each value of the dimension of POD space $N_{POD} = 1,\dots,N_s$, the scalar coefficients $a_1^1,\dots,a_1^{N_s},\dots,a_{N_s}^1,\dots,a_{N_s}^{N_s}$ and functions $\bm{\varphi}_1,\dots,\bm{\varphi}_{N_s}$ that minimize the quantity:
\begin{eqnarray}\label{eq:pod_energy}
E_{N_{POD}} = \sum_{i=1}^{N_s}||{\bm{u}}_i-\sum_{k=1}^{N_{POD}}a_i^k {\bm{\varphi}}_k||^2_{L^2({\mathcal D})}, \,&&\forall
N_{POD} = 1,\dots,N,\\\nonumber
&& \mbox{ with } ({\bm{\varphi}_i,\bm{\varphi}_j}) _{{\mathcal D }} = \delta_{ij}, \mbox{\hspace{0.5cm}}  \mbox{ }\forall
 i,j = 1,\dots,N_s.
\end{eqnarray}
In this case the velocity field $\bm{u}$ is used as an example. It can be shown \cite{Kunisch2002492} that the minimization problem of equation~\eqref{eq:pod_energy} is equivalent to solving the following eigenvalue problem:
\begin{gather}
{\bm {\mathcal{X}}^u}\bm{Q}^u = \bm{Q^u}\bm{\lambda^u} ,\quad \mbox{\hspace{0.5cm} for }\mathcal{X}^u_{ij} = ({\bm{u}_i,\bm{u}_j}) _{{\mathcal D }} \mbox{,\, } i,j = 1,\dots,N_s ,\nonumber
\end{gather}
where $\bm{\mathcal{X}^u}$ is the correlation matrix obtained from the parameter dependent snapshots $\bm{\mathcal{S}_u}$, $\bm{Q^u}$ is a square matrix of eigenvectors and $\bm{\lambda^u}$ is a diagonal matrix of eigenvalues. 

The basis functions can then be obtained with: 
\begin{equation}
\bm{\varphi_i} = \frac{1}{N_s{\lambda_{ii}^u}^{1/2}}\sum_{j=1}^{N_s} \bm{u}_j Q^u_{ij}.
\end{equation}
The same procedure is applied  for the pressure field considering the snapshots matrix consisting of the snapshots $
p_1,p_2,\dots,p_{N_s}$. The correlation matrix of the pressure field snapshots $\bm{C^p}$ is assembled and we solve a similar eigenvalue problem ${\bm{C}^p}{\bm{Q}^p}={\bm{Q}^p}{\bm{\lambda}^p}$. The POD modes ${\chi_i}$ for the pressure field can be computed as
\begin{equation}
 {\chi_i}=\frac{1}{N_s{{\lambda^p _{ii}}^{1/2}}}\sum_{j=1}^{N_s} p_j {Q^p_{ij}}.
\end{equation}
The POD spaces are constructed for both velocity and pressure using the above methodology resulting in the spaces:
\begin{equation}
\begin{split}
\bm{L}_u = [{\bm{\varphi}}_1, \dots , {\bm{\varphi}}_{N_u^r}] \in \mathbb{R}^{N_{u}^h \times N_u^r},\quad
\bm{L}_p = [{\chi_1}, \dots , {\chi_{N_p^r}}] \in \mathbb{R}^{N_{p}^h \times N_p^r}.
\end{split}
\end{equation}
where $N_u^r$, $N_p^r < N_s$ are chosen according to the eigenvalue decay of $\bm{\lambda}^u$ and $\bm{\lambda}^p$, \cite{Rozza2008229,BeOhPaRoUr17}.

\begin{figure} \centering
  \includegraphics[width=0.75\textwidth]{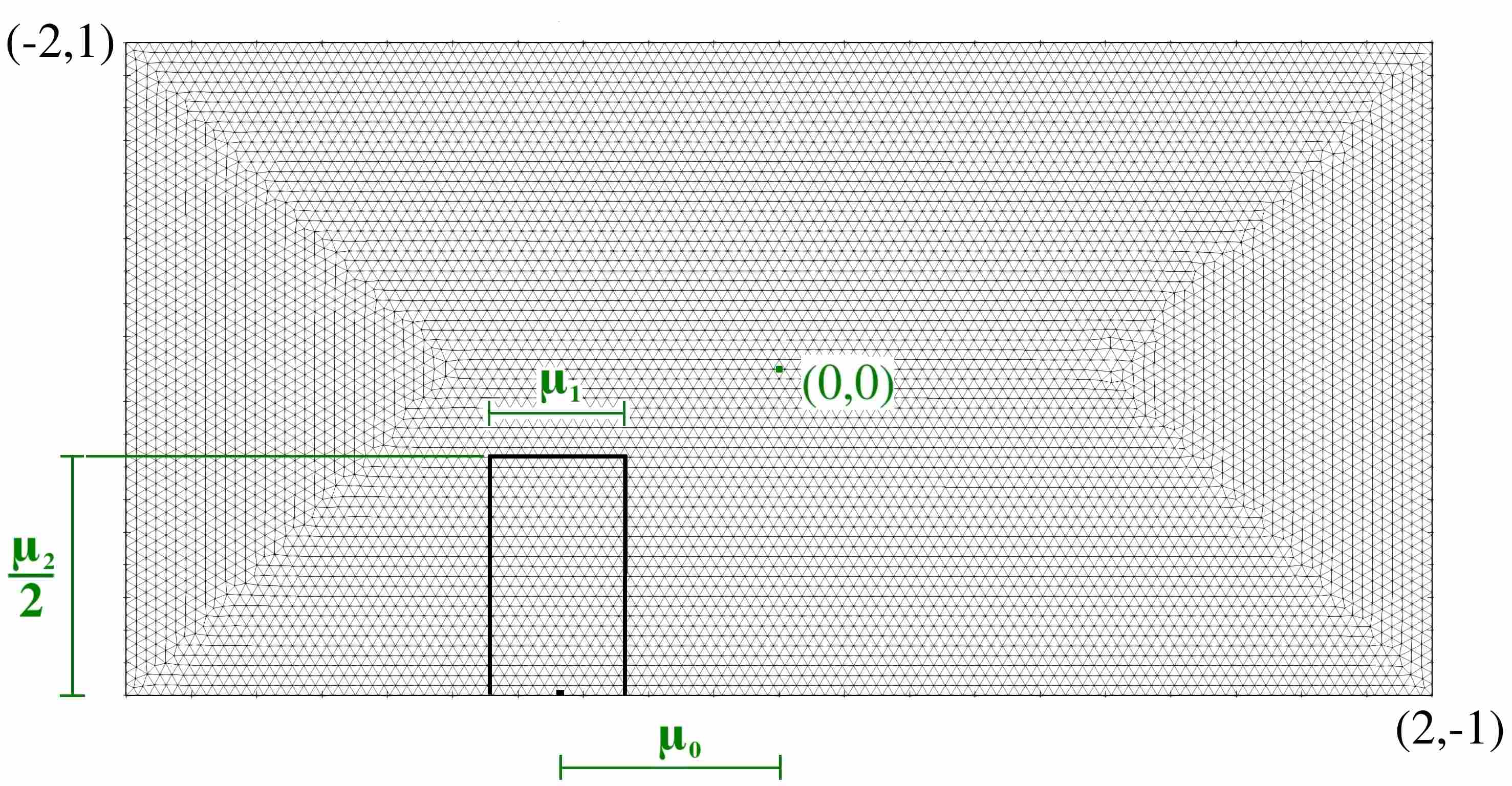}
  \caption{A sketch of the embedded domain and associated background mesh, together with the parameters considered in the numerical examples, {\RA{namely the $x$ position of its center ($\mu_0$) and its aspect ratio ($\mu_2/2\mu_1)$ for both the $x$ and $y$ direction, 
width and height respectively}}.}
  \label{background_mesh}
\end{figure}

\begin{rmk}
Following the ideas of \cite{KaBaRO18}, in the out-of-interest area (i.e., the region outside  the true geometry), we prefer to use the solution values that are computed using the shifted boundary method and in particular the smooth mapping ${\mathcal M}$ from the true to the surrogate domain. This approach allows a smooth extension of the boundary solution to the neighboring ghost elements with values which are decreasing smoothly to zero, see for instance the zoomed image in \autoref{fig:poisson_zoom}. This approach guarantees a regular solution in the background domain and permits therefore the construction of an effective reduced basis.
\begin{figure} 
\centering
\includegraphics[width=0.4\textwidth]{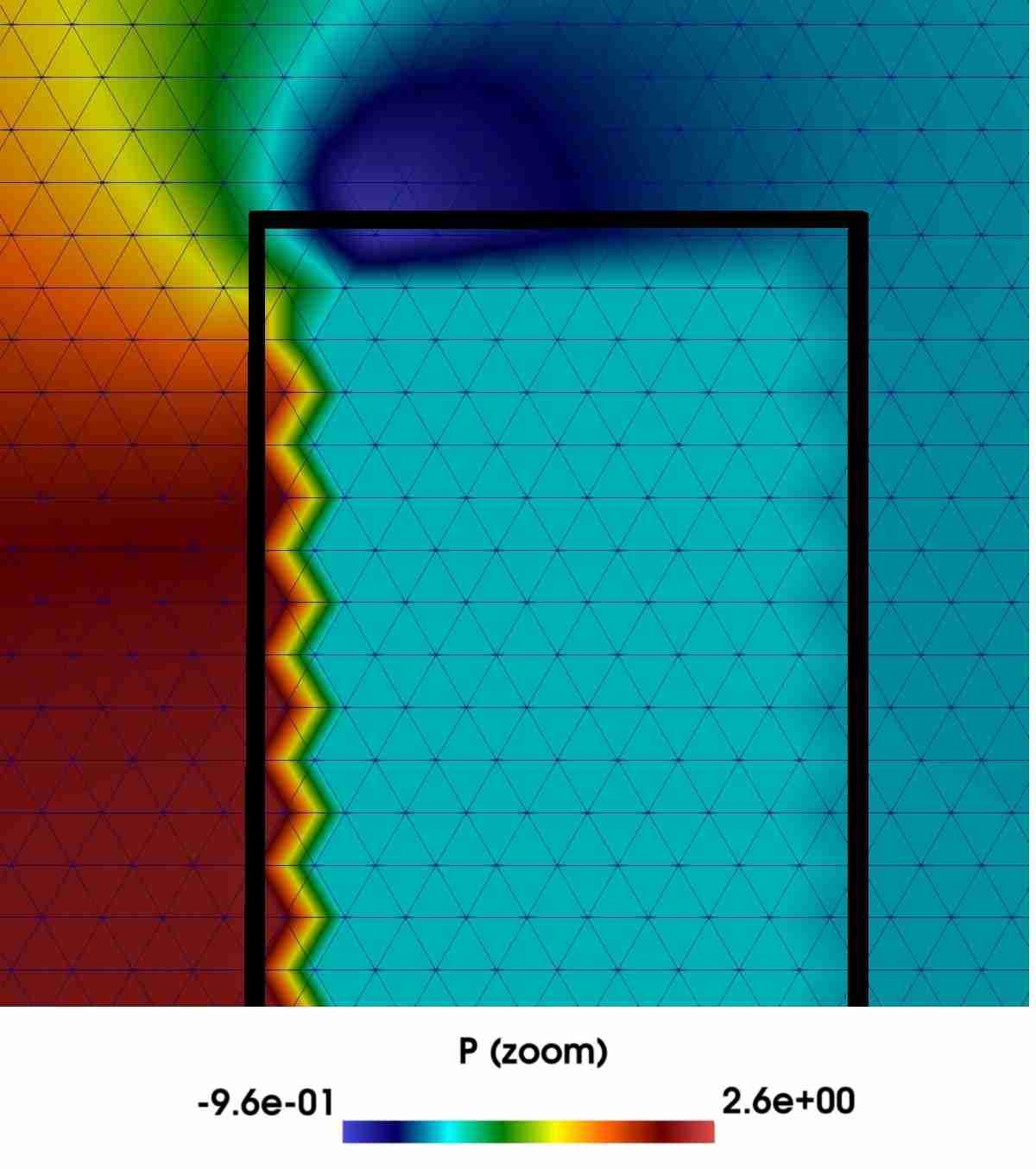}
\caption{A zoom onto the embedded box in order to show the smoothing procedure employed by the SBM method inside the ghost area.}
\label{fig:poisson_zoom}
\end{figure}
\end{rmk}
\subsection{The projection stage and the generation of the ROM}
Once the POD functional spaces are constructed, the reduced velocity and pressure fields can be approximated with: 
\begin{equation}\label{eq:aprox_fields}
\bm{u^r} \approx \sum_{i=1}^{N_u^r} a_i(\mu) \bm{\varphi}_i(\bm{x}) = \bm{L}_u \bm{a}(\mu), \mbox{\hspace{0.5 cm}}
p^r\approx \sum_{i=1}^{N_p^r} b_i(\mu)\chi_i(\bm{x}) = \bm{L}_p \bm{b}(\mu).
\end{equation} 
The reduced solution vectors $\bm{a} \in \mathbb{R}^{N^r_u\times1}$ and $\bm{b}\in \mathbb{R}^{N^r_p\times1}$ depend only on the parameter values and the basis functions $\bm{\varphi}_i$ and ${\chi}_i$ depend only on the physical space. Denoting by ${\bm{L}}=\begin{bmatrix} \bm{L}_u & \bm{0} \\ \bm{0} & \bm{L}_p \end{bmatrix}$ and by ${\bm{L}}^T =\begin{bmatrix} {\bm{L}}^T_u & {\bm{0}} \\ {\bm{0}} & {\bm{L}}^T_p \end{bmatrix}$, the unknown vectors of coefficients $V=\begin{bmatrix} \bm{a} \\ \bm{b} \end{bmatrix}$ then can be obtained through a Galerkin projection of the full-order system of equations onto the POD reduced basis spaces.  The subsequent solution of a reduced iterative algebraic system of equations for the increment $\delta V(\mu)$ then becomes,
\begin{equation}\label{eq:system_linear}
{\bm{L}}\nabla G(U^{n-1}(\mu))\bm{L}^T\delta {V}(\mu) = - \bm{L}^T R(U^{n-1}(\mu)),
\end{equation}
which leads to the following algebraic reduced system:
\begin{equation}\label{eq:system_linear_reduced}
\nabla G^r(V^{n-1}(\mu))\delta {V}(\mu) = - R^r(V^{n-1}(\mu)).
\end{equation}

We remark here that at the reduced order level, we need to assemble the FOM problem in order to compute the reduced differential operator, but this expensive operation could be avoided, for example, using hyper reduction techniques~\cite{stabile_geo_,Xiao20141,BARRAULT2004667,Carlberg2013623}. 
Moreover, during the online stage, also the stabilization term ${\bm{C}}$ (as well as the nonlinear term ${\bm{\mathcal C}}$)  is projected onto the reduced basis space, which allows the inf-sup stabilization condition to propagate in an efficient way into the reduced model, as we will see in the following numerical examples.}}

\begin{rmk}
We remark here that the reduced basis spaces have been generated using the iterative solution snapshots and not only the final solutions. This is justified by the fact that also at the reduced order level an iterative procedure is used to solve the non-linear problem and to obtain the reduced basis solutions. 
\end{rmk}

\section{Numerical experiments}\label{sec:num_exp}
We consider two different test cases, based on the setup shown in Figure \ref{background_mesh}. 
The first one consists of a geometrical parametrization using a one-dimensional parameter space where both parameters $\mu_0$, $\mu_1$ are fixed. The second one consists of a geometrical parametrization with a three-dimensional parameter space where all the parameters $\mu_0$ $\mu_1$ and $\mu_2$ are left free. The problem domain is the rectangle $ [-2, 2] \times [-1, 1]$, in which an embedded rectangular is immersed. The viscosity $\nu$ is set to $1$ and the Reynolds number is set to $100$. A constant velocity in the $x$ direction, $u_{\text{in}} = 1$ is applied at the left side of the domain, and an open boundary condition with {\RA{$\nabla \bm u \cdot \bm n 
= 0$}}  on the right. In addition, a slip (no penetration) boundary condition is applied on the top and bottom edges while on the boundary of the embedded rectangle a no slip boundary condition is applied. 
The results for the test problems have been obtained with a mesh size of $h = 0.0350$ for the background mesh, 
using $15022$ triangles for the discretization and $\mathbb{P}1/\mathbb{P}1$ polynomials.


\subsection{Geometrical parametrization with one-dimensional parameter space}
\begin{figure} \centering
\begin{minipage}{\textwidth}
\centering
\makebox[\textwidth][l]{%
\begin{minipage}{0.24\textwidth}
  \includegraphics[width=\textwidth]{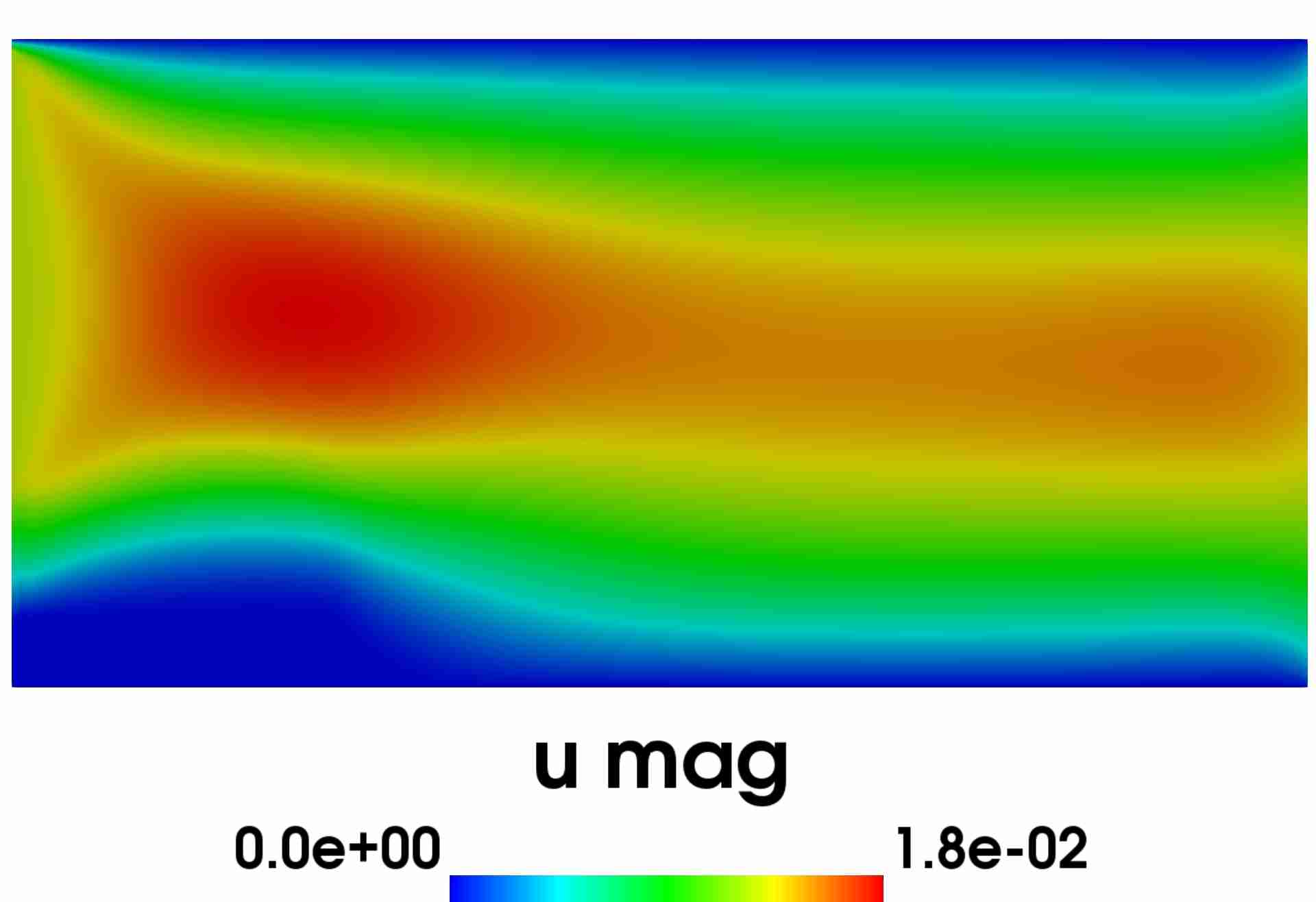} 
\end{minipage}
\begin{minipage}{0.24\textwidth}
  \includegraphics[width=\textwidth]{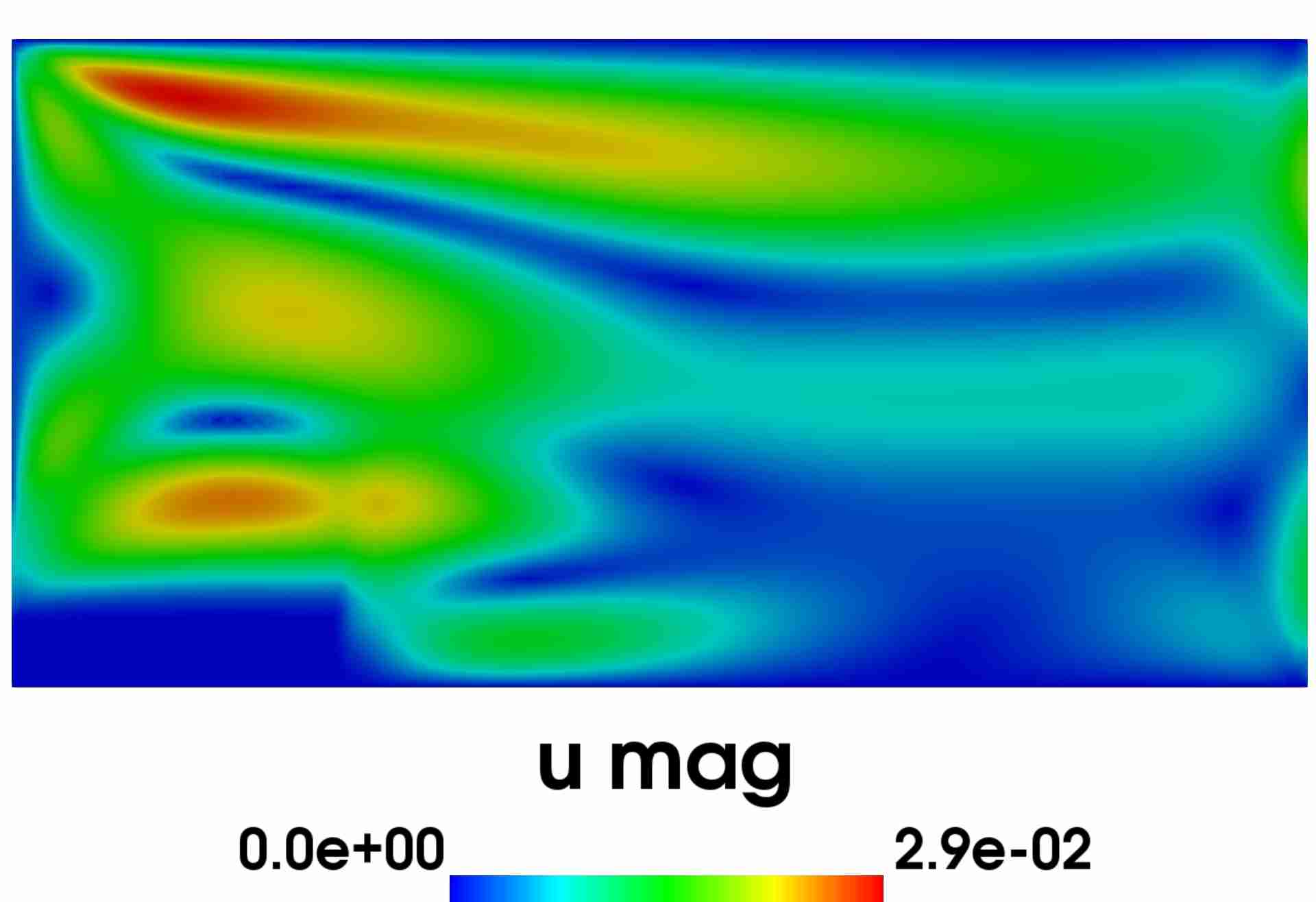}
\end{minipage}
\begin{minipage}{0.24\textwidth}
  \includegraphics[width=\textwidth]{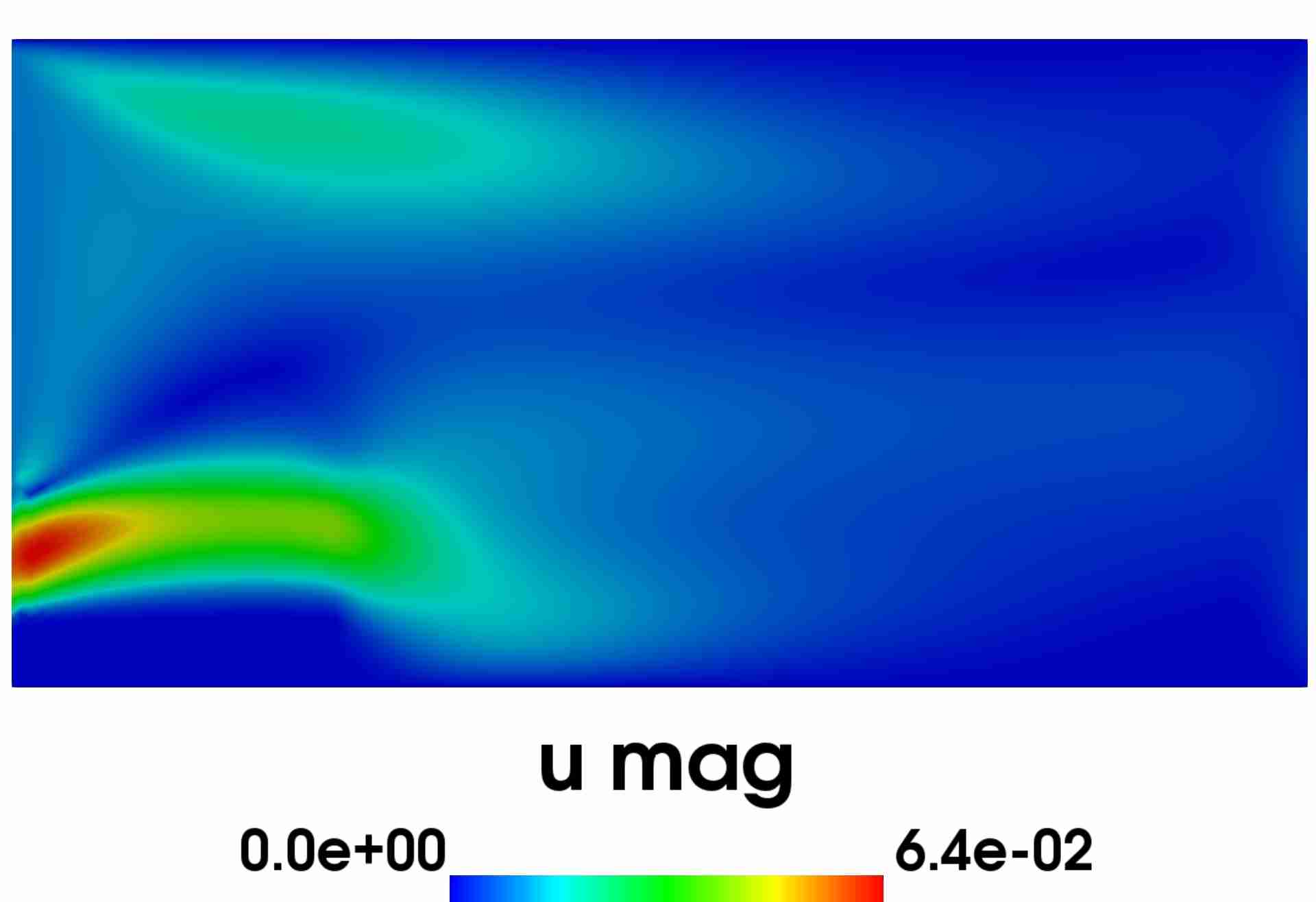}
\end{minipage}
\begin{minipage}{0.24\textwidth}
  \includegraphics[width=\textwidth]{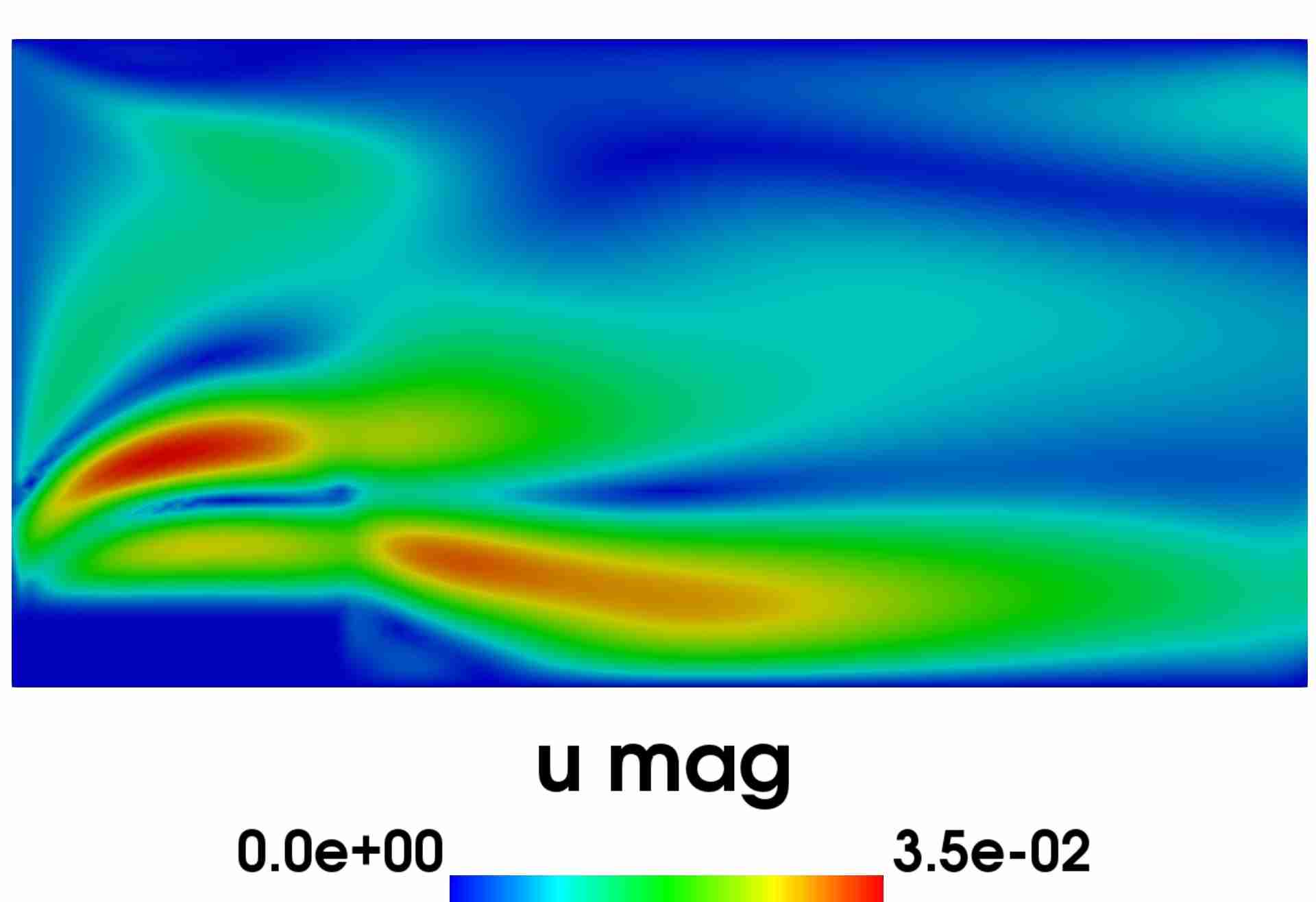}
\end{minipage}
}
\makebox[\textwidth][l]{%
\begin{minipage}{0.24\textwidth}
  \includegraphics[width=\textwidth]{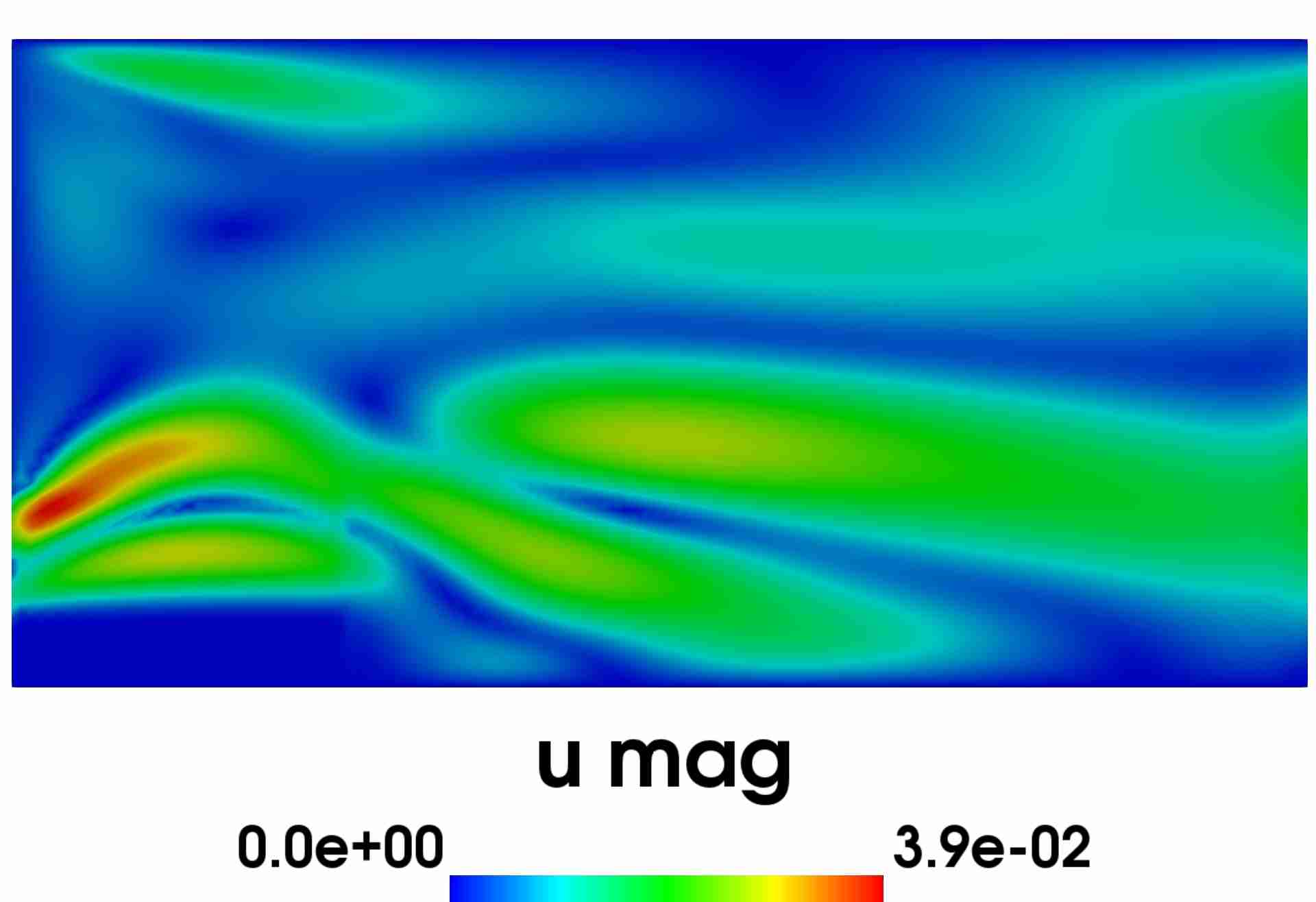}
\end{minipage}
\begin{minipage}{0.24\textwidth}
  \includegraphics[width=\textwidth]{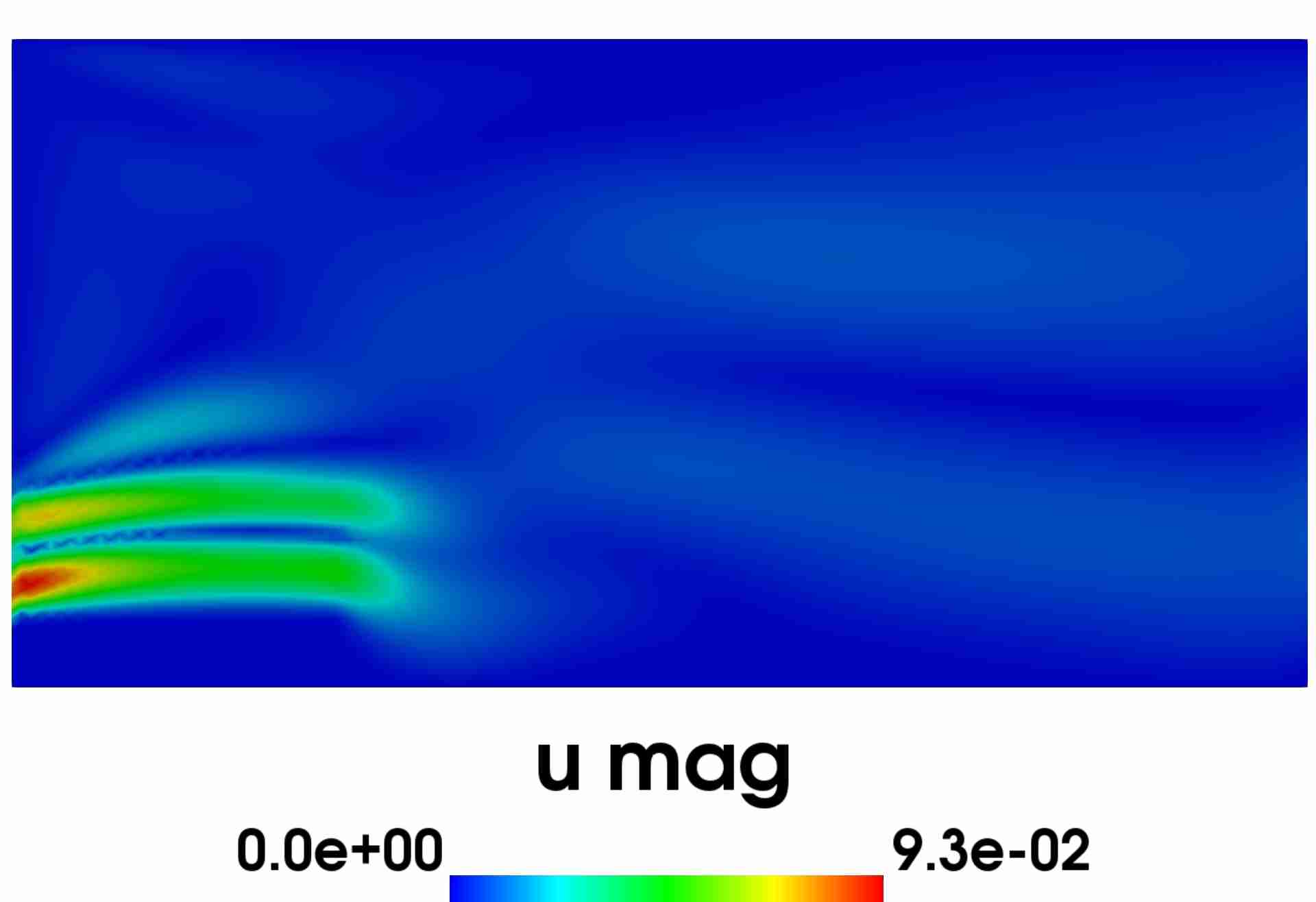}
\end{minipage}
}
\newline
\makebox[\textwidth][l]{%
\begin{minipage}{0.24\textwidth}
  \includegraphics[width=\textwidth]{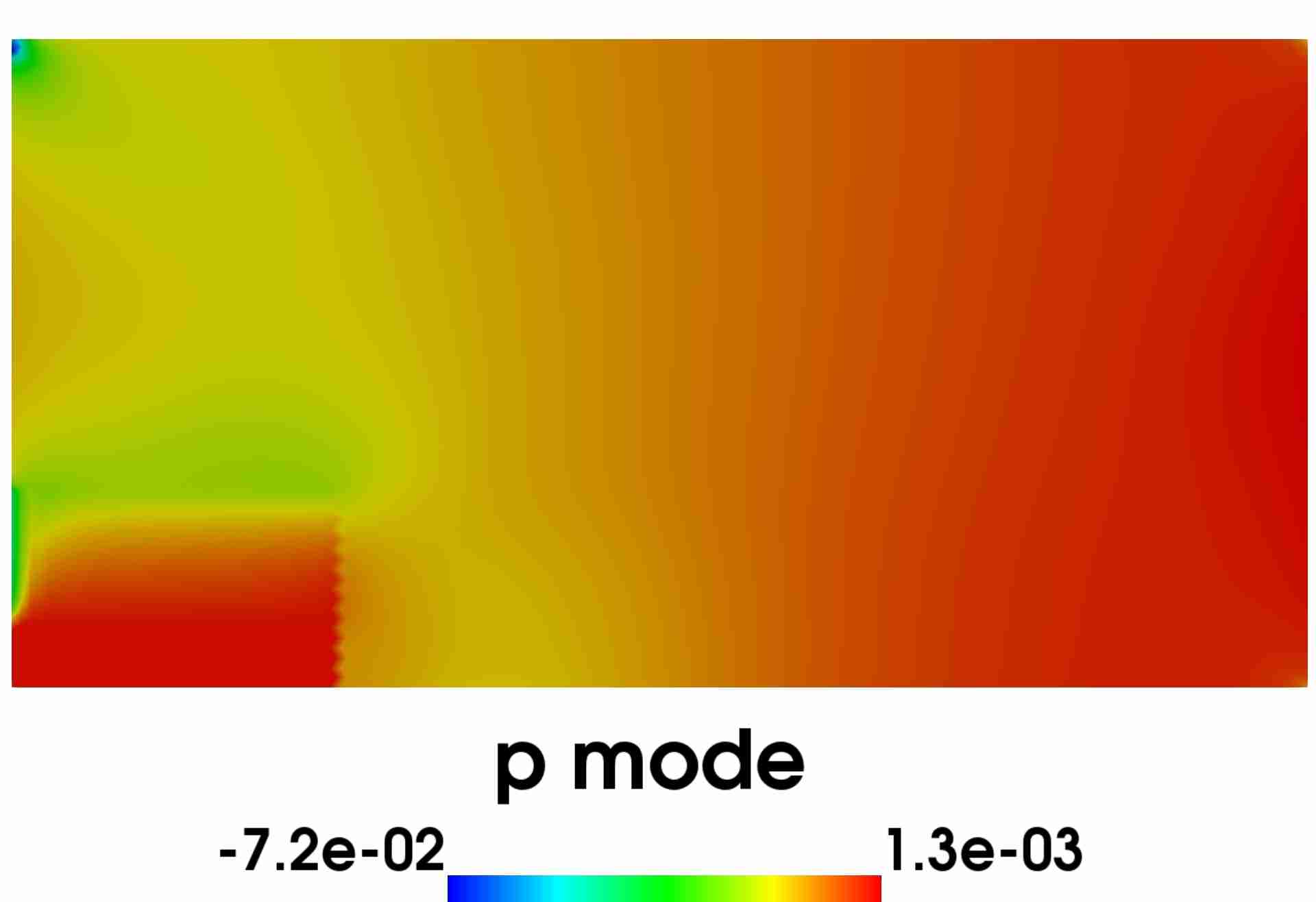} 
\end{minipage}
\begin{minipage}{0.24\textwidth}
  \includegraphics[width=\textwidth]{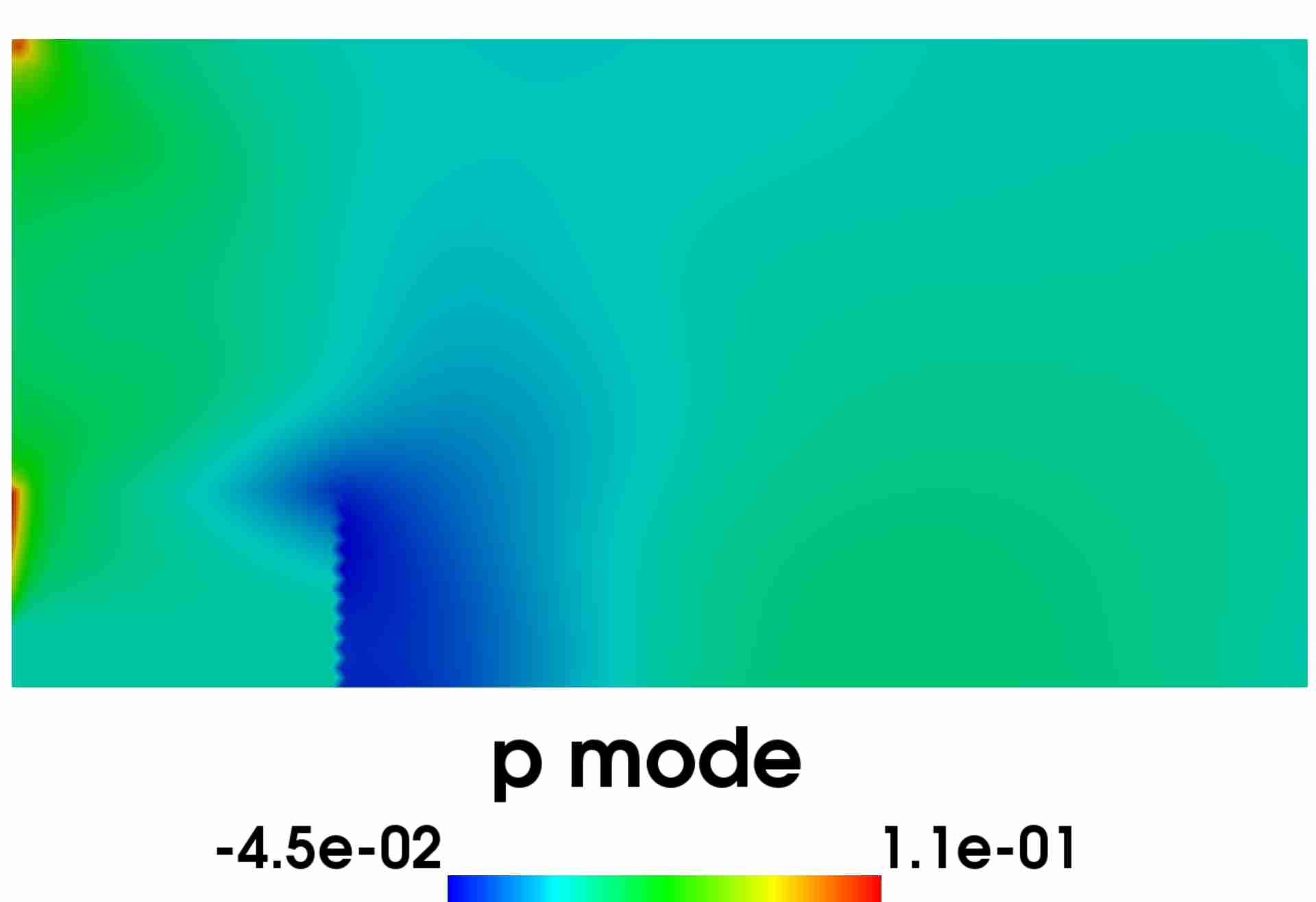}
\end{minipage}
\begin{minipage}{0.24\textwidth}
  \includegraphics[width=\textwidth]{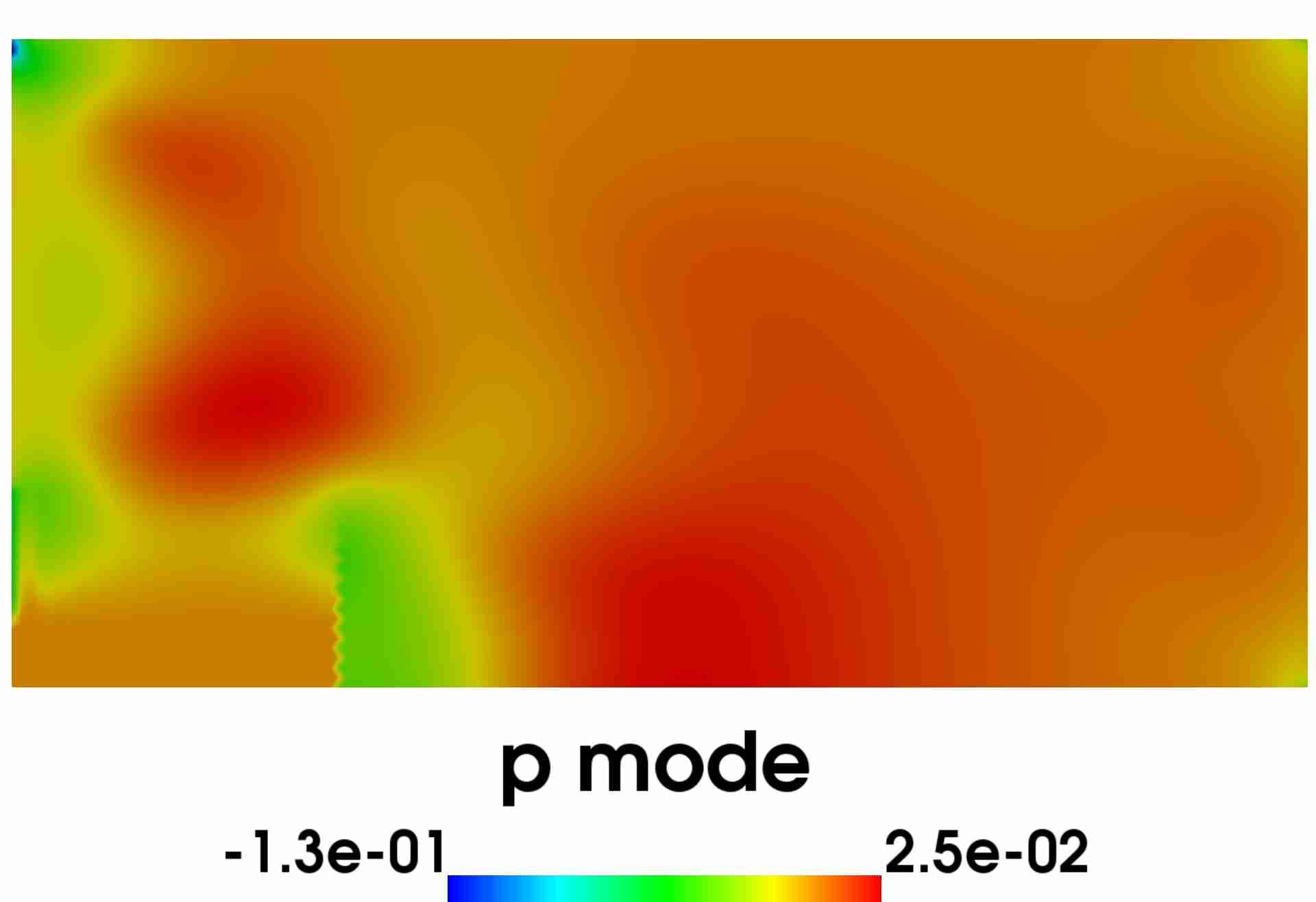}
\end{minipage}
\begin{minipage}{0.24\textwidth}
  \includegraphics[width=\textwidth]{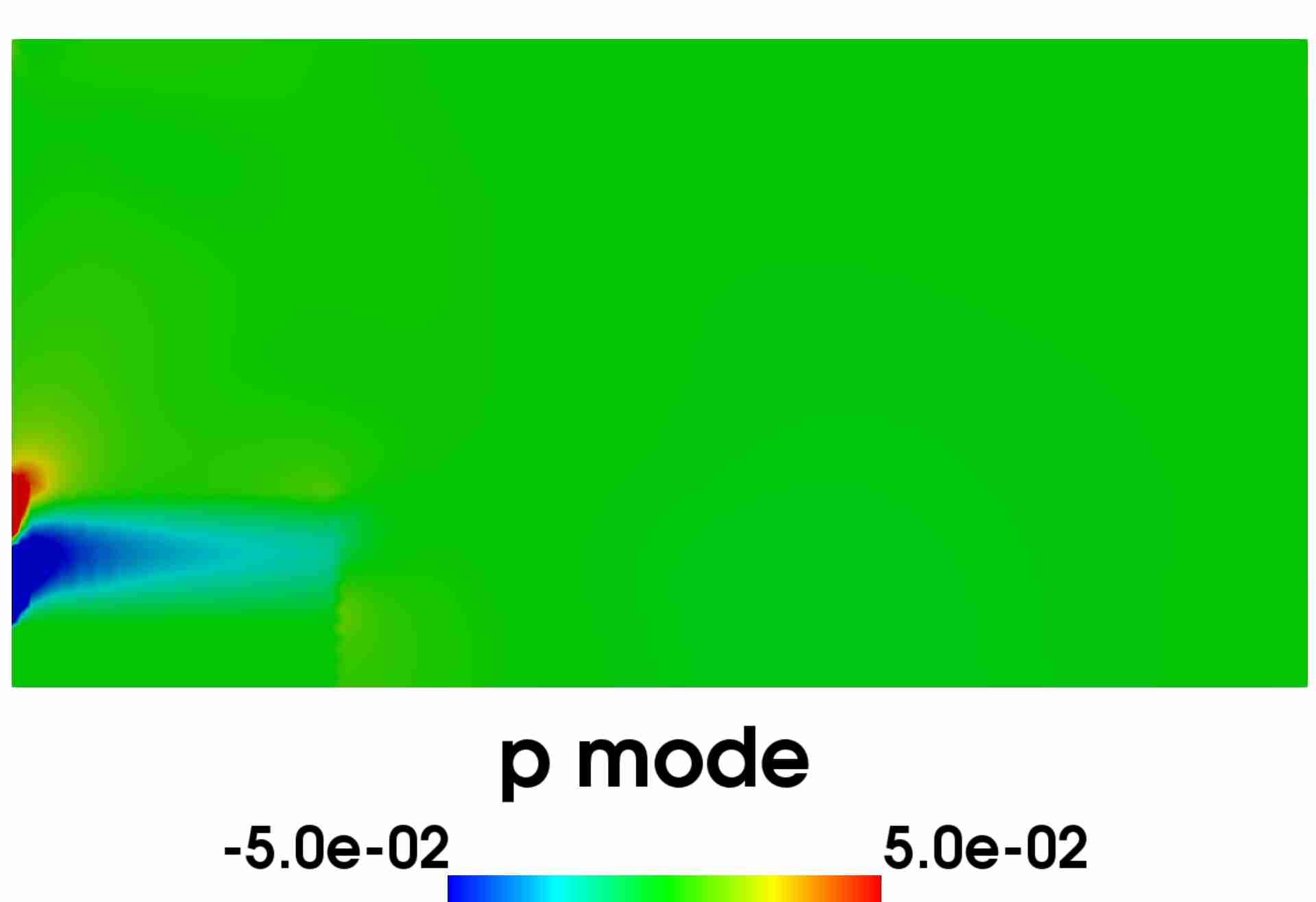}
\end{minipage}
}
\newline
\makebox[\textwidth][l]{%
\begin{minipage}{0.24\textwidth}
  \includegraphics[width=\textwidth]{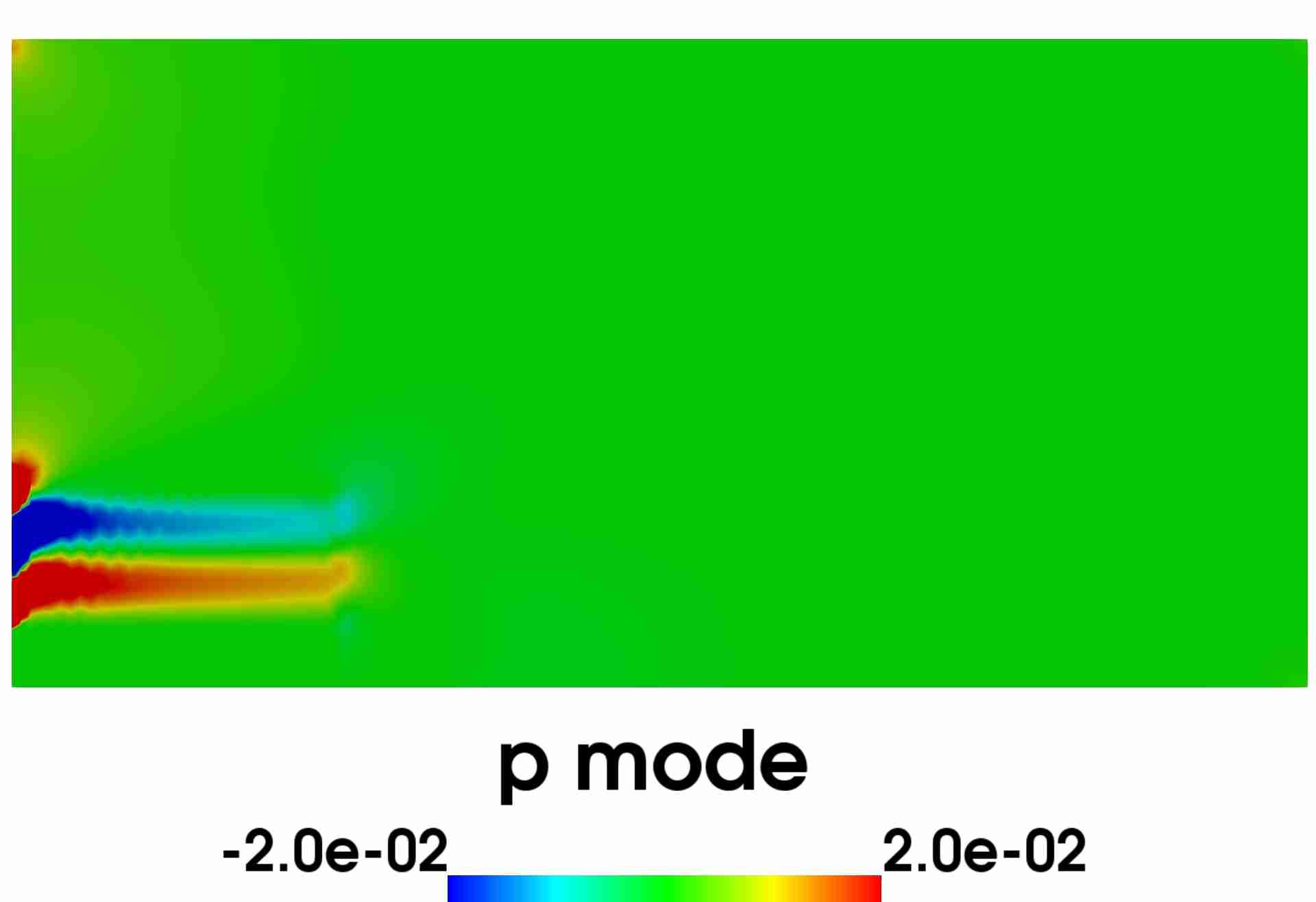}
\end{minipage}
\begin{minipage}{0.24\textwidth}
  \includegraphics[width=\textwidth]{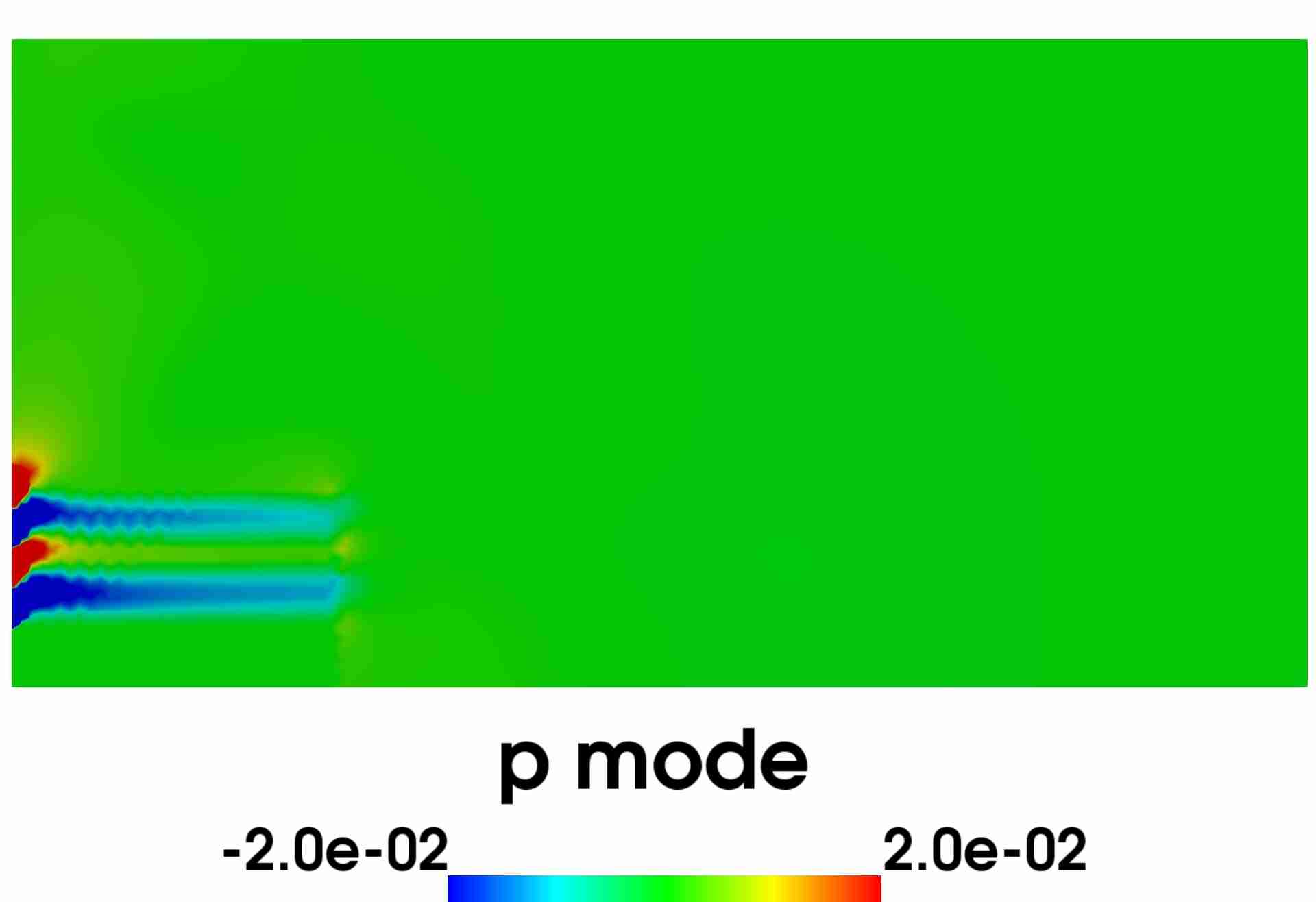}
\end{minipage}
}
\end{minipage}
  \caption{The first $6$ components of the velocity basis (first row) and pressure basis (second row) for the geometrical parametrization with a one-dimensional parameter space case with $\mu \in [1, 2]$.
}  \label{1DStokes_Pressure_and_Velocity_Components_Modes}
\begin{minipage}{\textwidth}
\centering
\makebox[\textwidth][l]{%
\begin{minipage}{0.24\textwidth}
  \includegraphics[width=\textwidth]{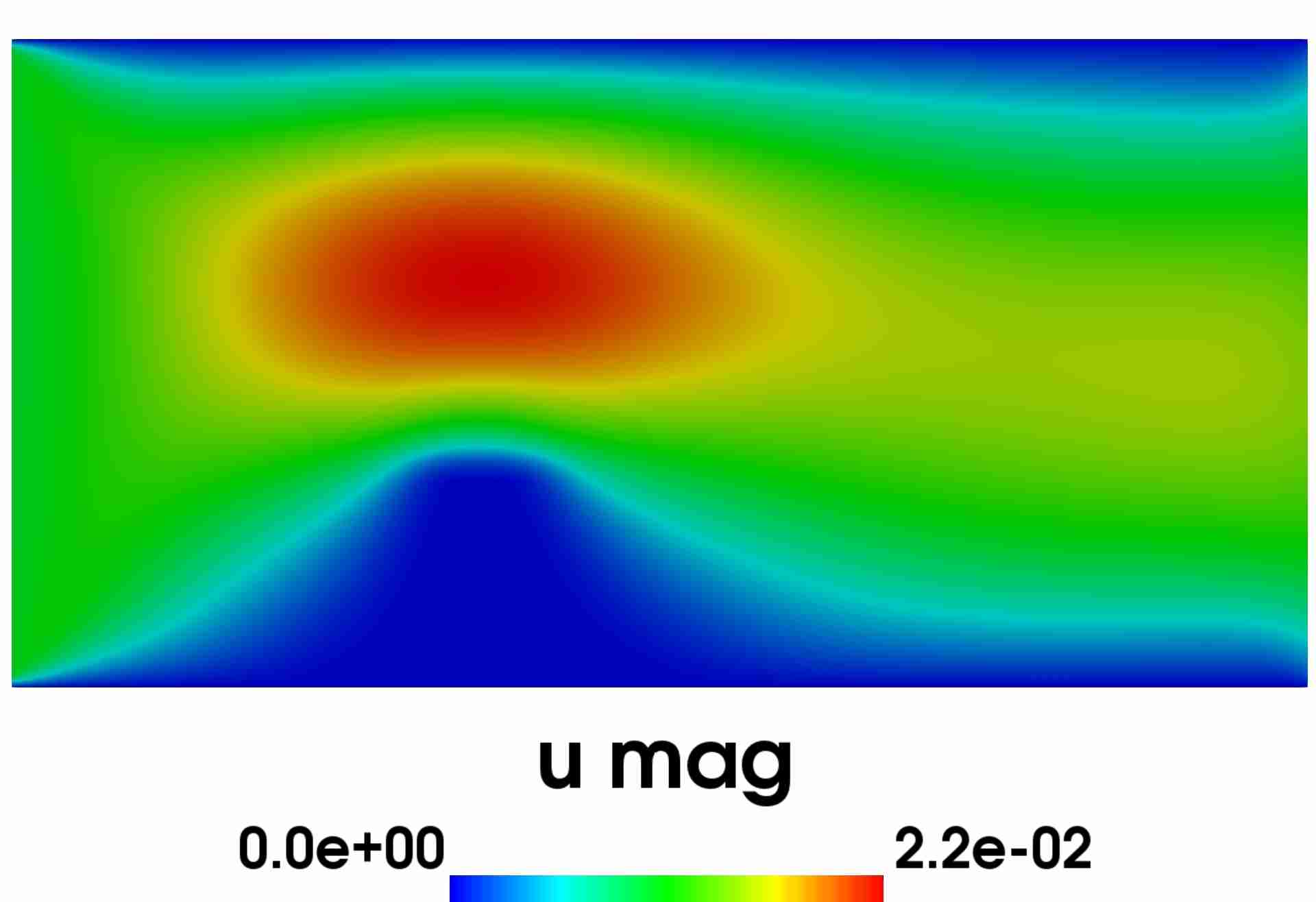} 
\end{minipage}
\begin{minipage}{0.24\textwidth}
  \includegraphics[width=\textwidth]{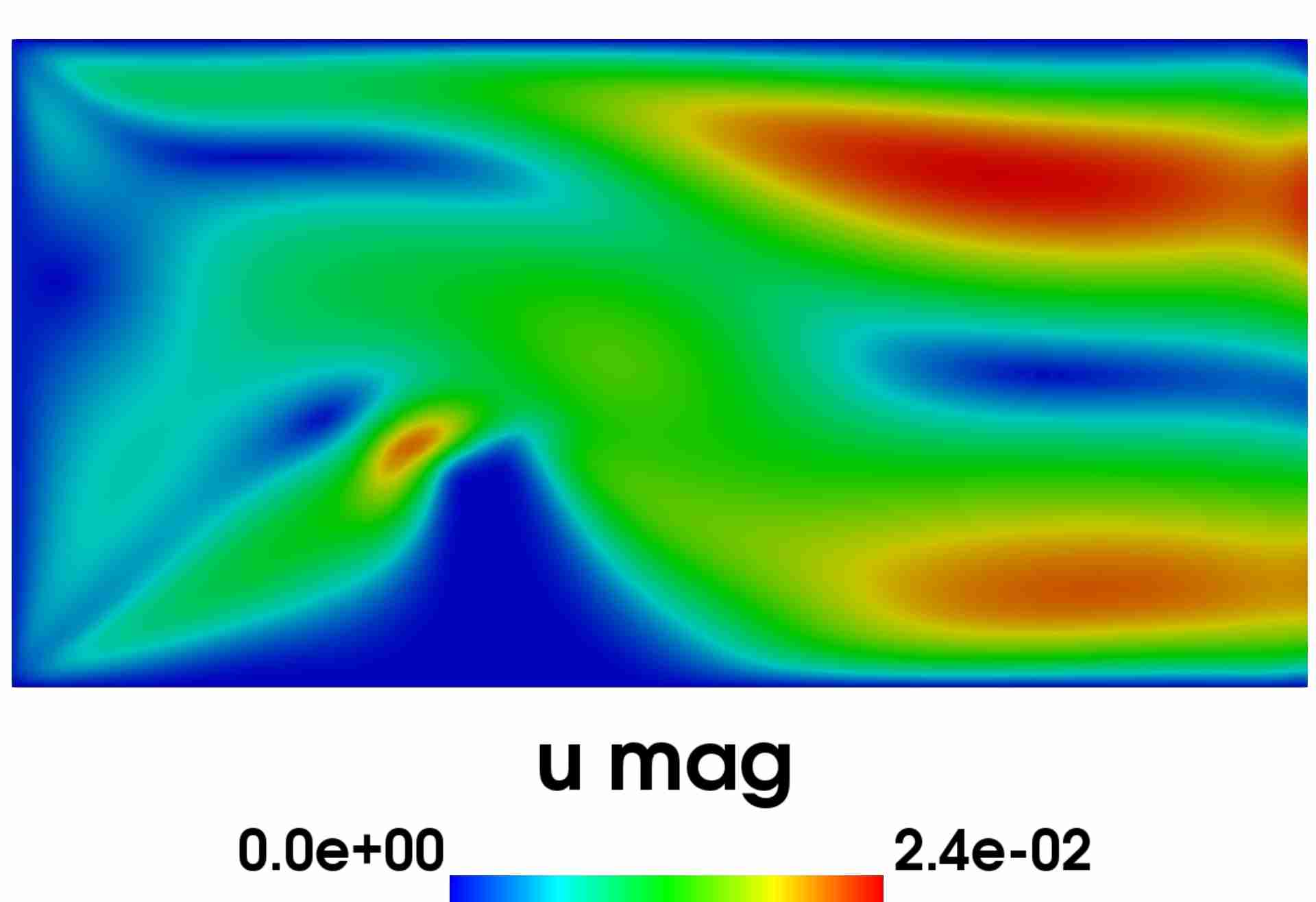}
\end{minipage}
\begin{minipage}{0.24\textwidth}
  \includegraphics[width=\textwidth]{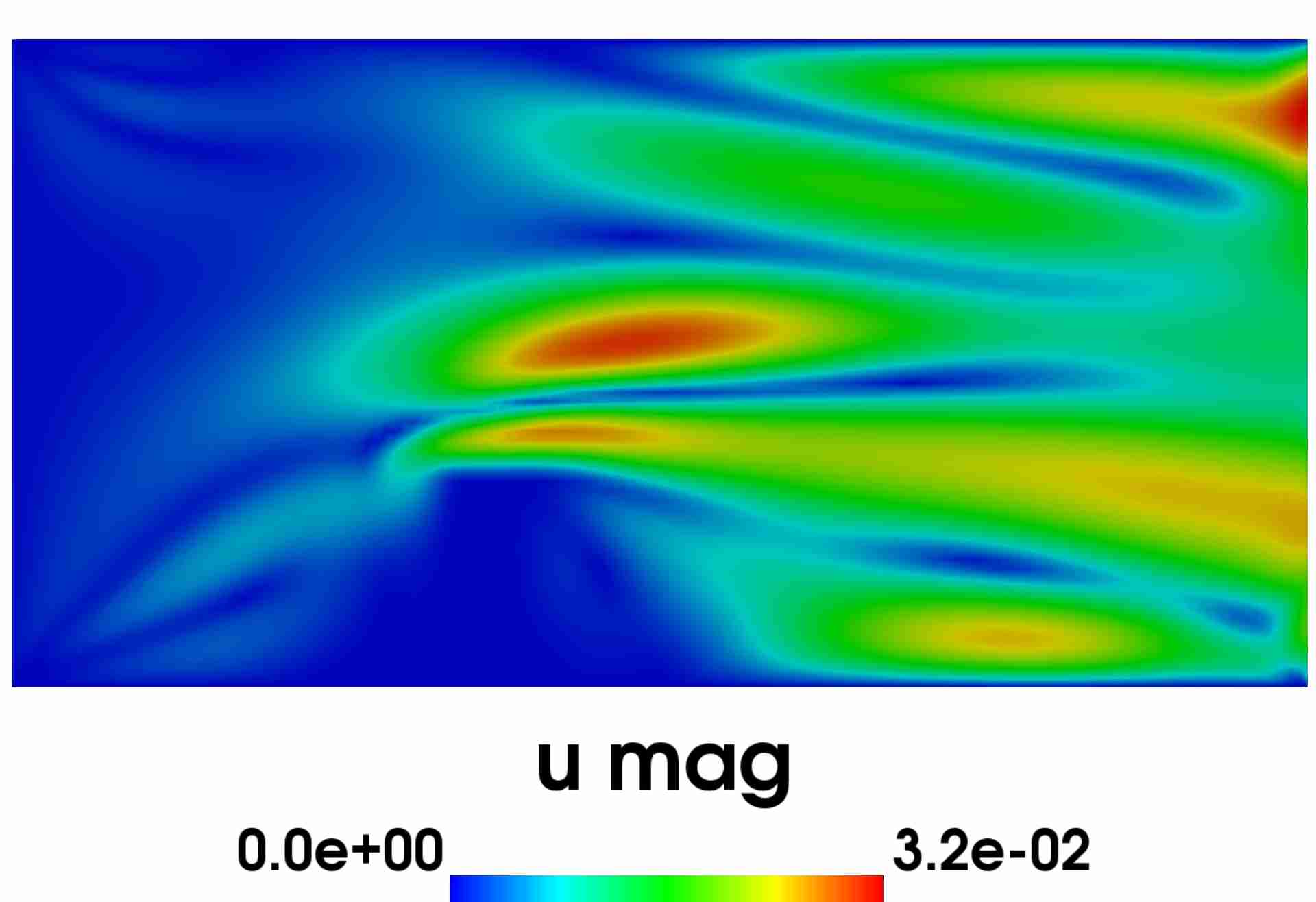}
\end{minipage}
\begin{minipage}{0.24\textwidth}
  \includegraphics[width=\textwidth]{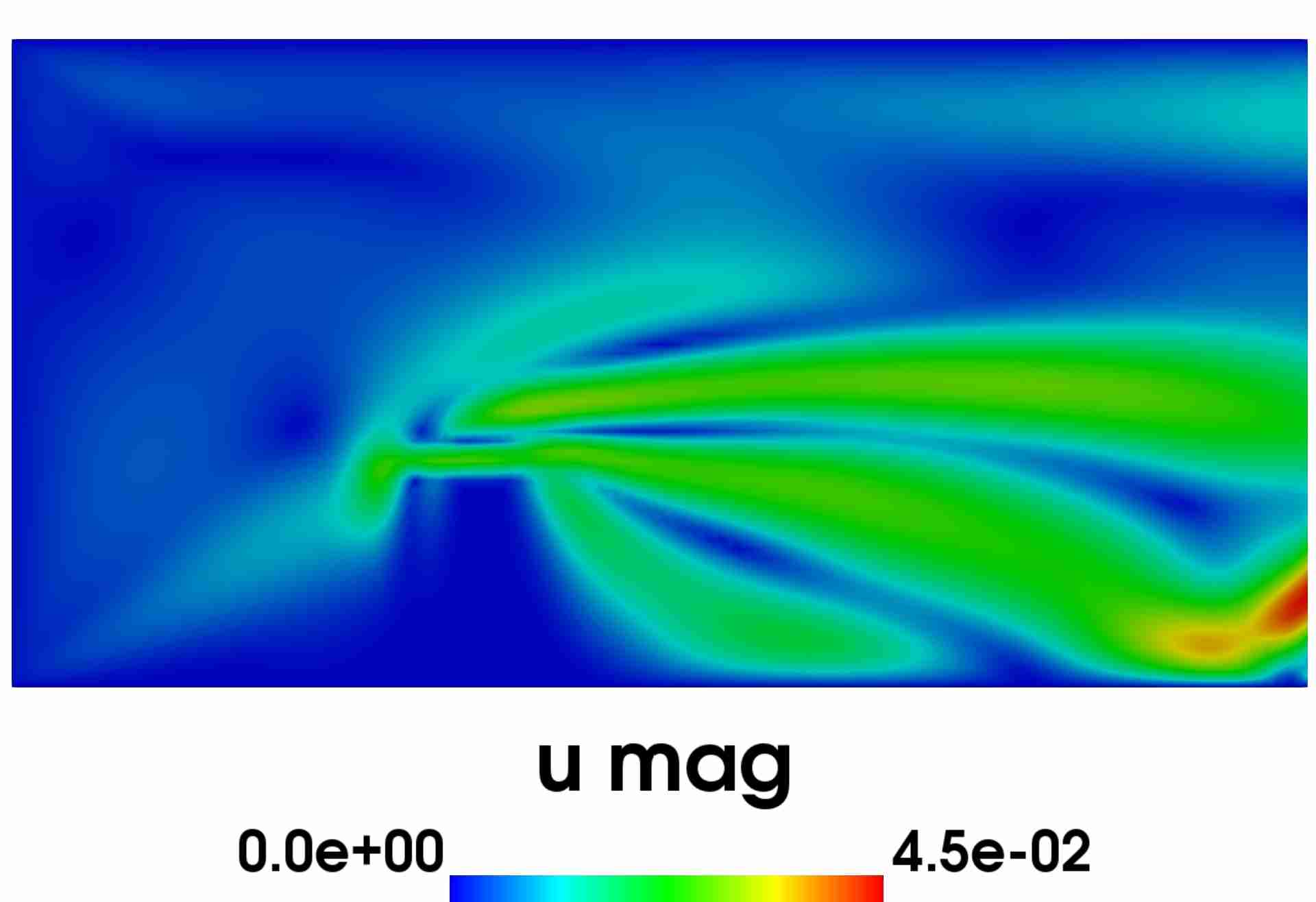}
\end{minipage}
}
\newline
\makebox[\textwidth][l]{%
\begin{minipage}{0.24\textwidth}
  \includegraphics[width=\textwidth]{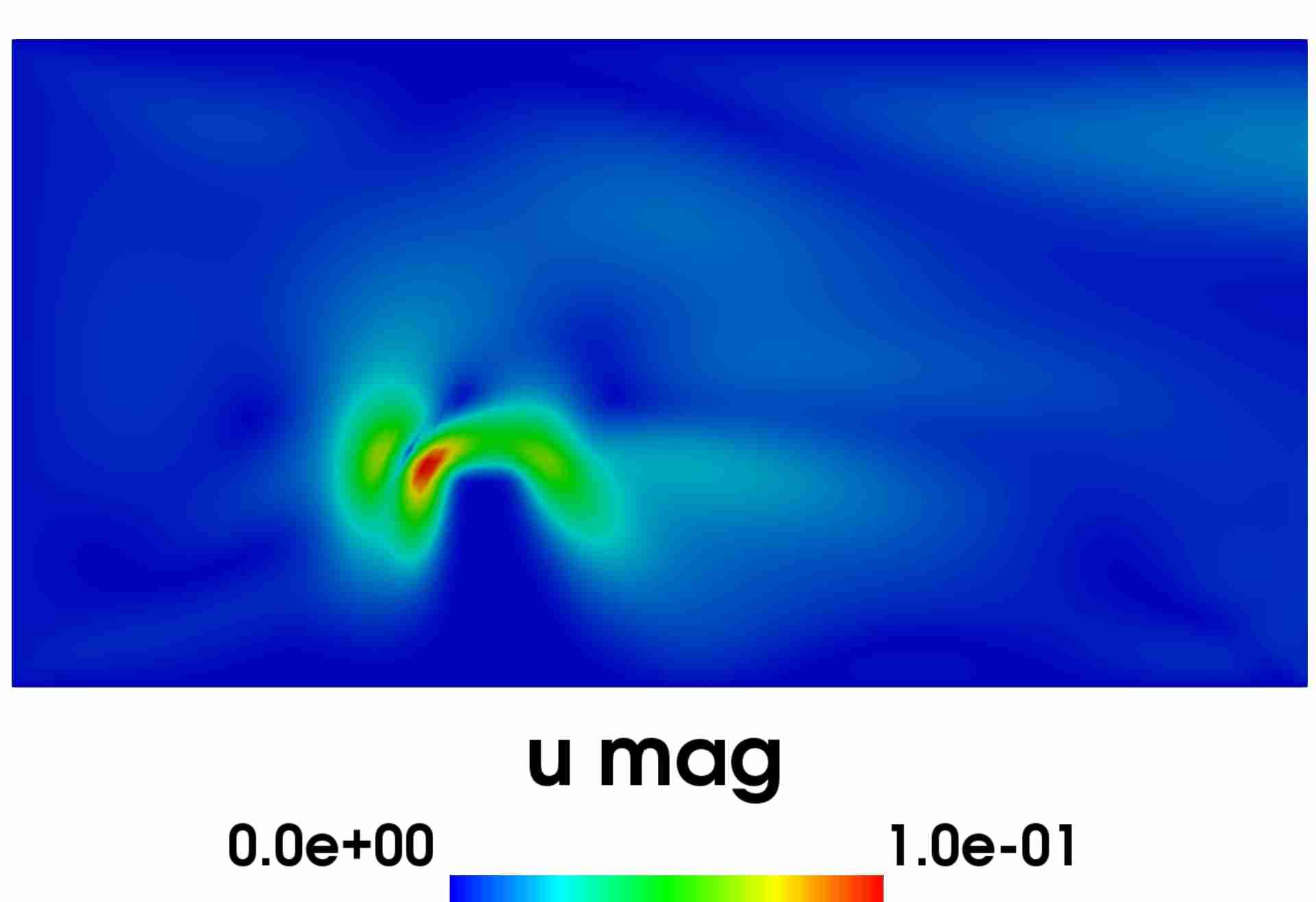}
\end{minipage}
\begin{minipage}{0.24\textwidth}
  \includegraphics[width=\textwidth]{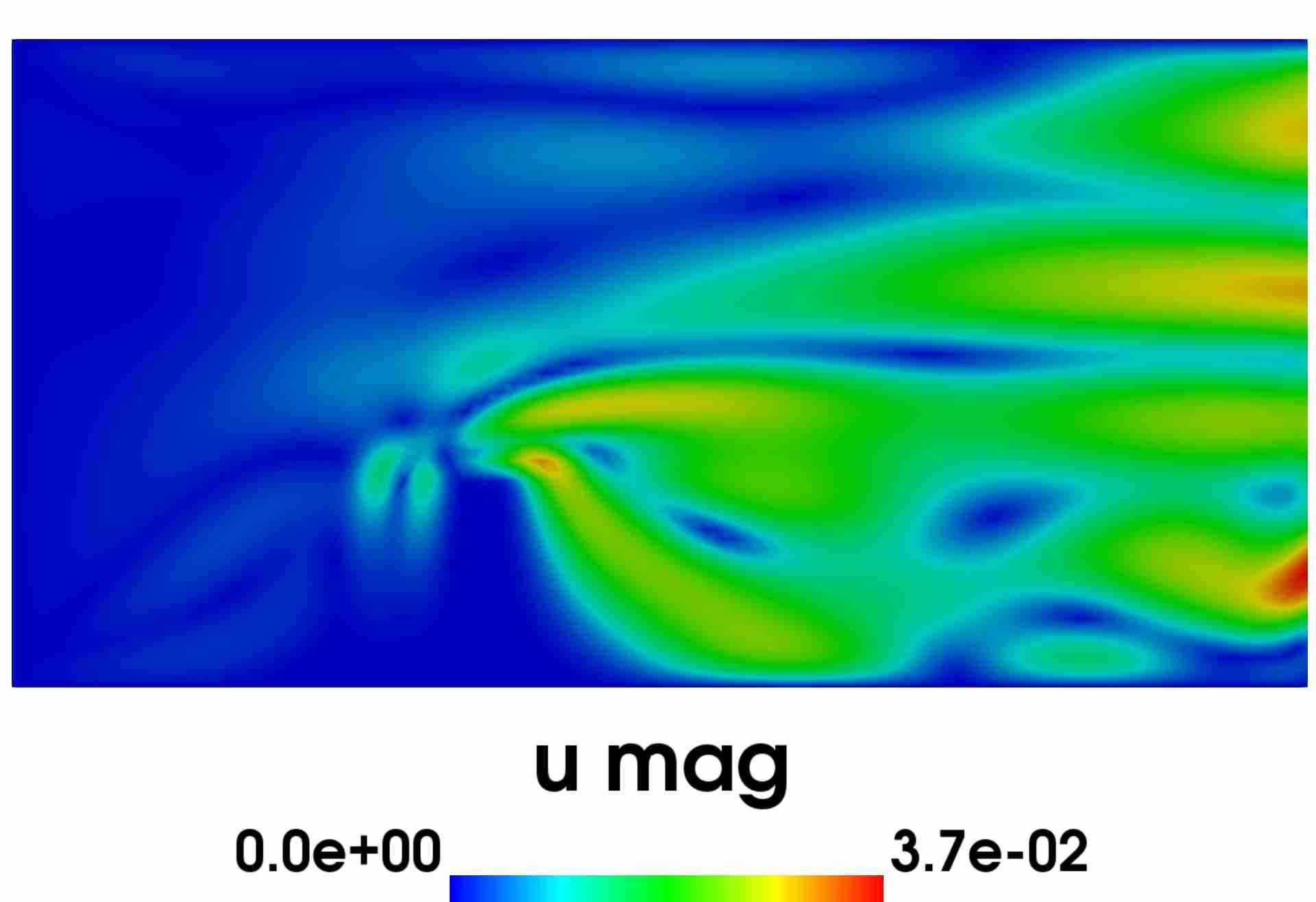}
\end{minipage}
}
\newline
\makebox[\textwidth][l]{%
\begin{minipage}{0.24\textwidth}
  \includegraphics[width=\textwidth]{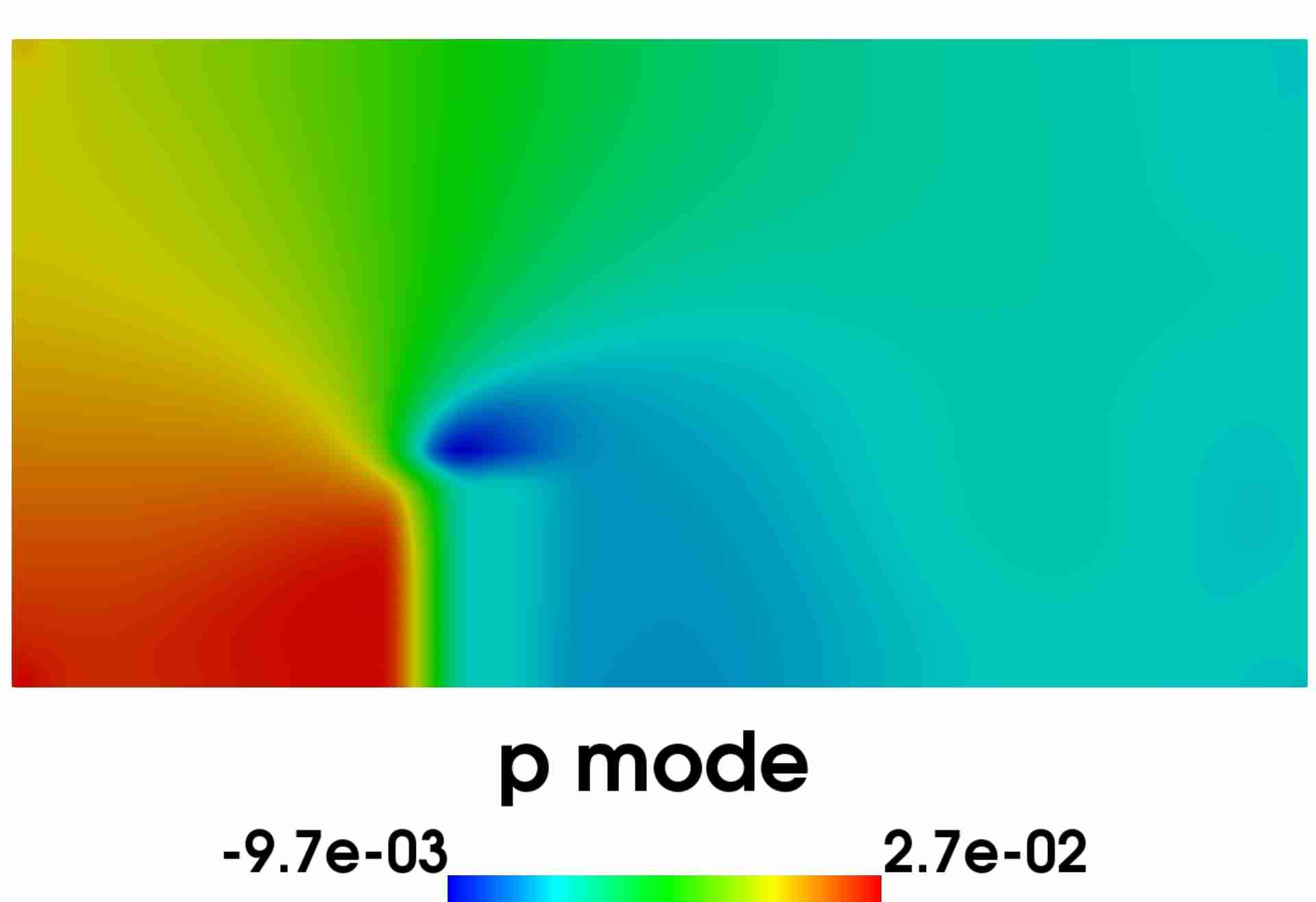} 
\end{minipage}
\begin{minipage}{0.24\textwidth}
  \includegraphics[width=\textwidth]{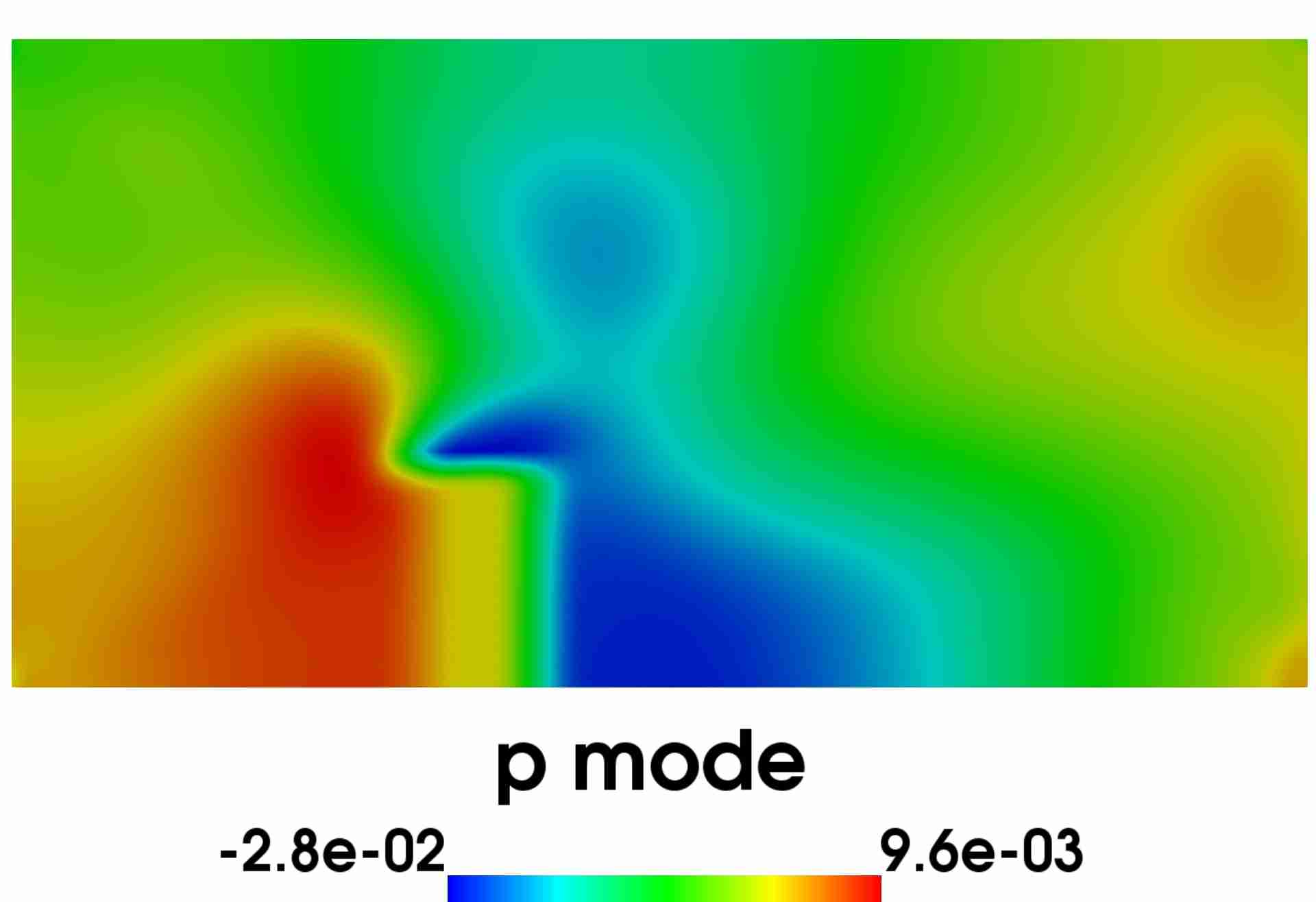}
\end{minipage}
\begin{minipage}{0.24\textwidth}
  \includegraphics[width=\textwidth]{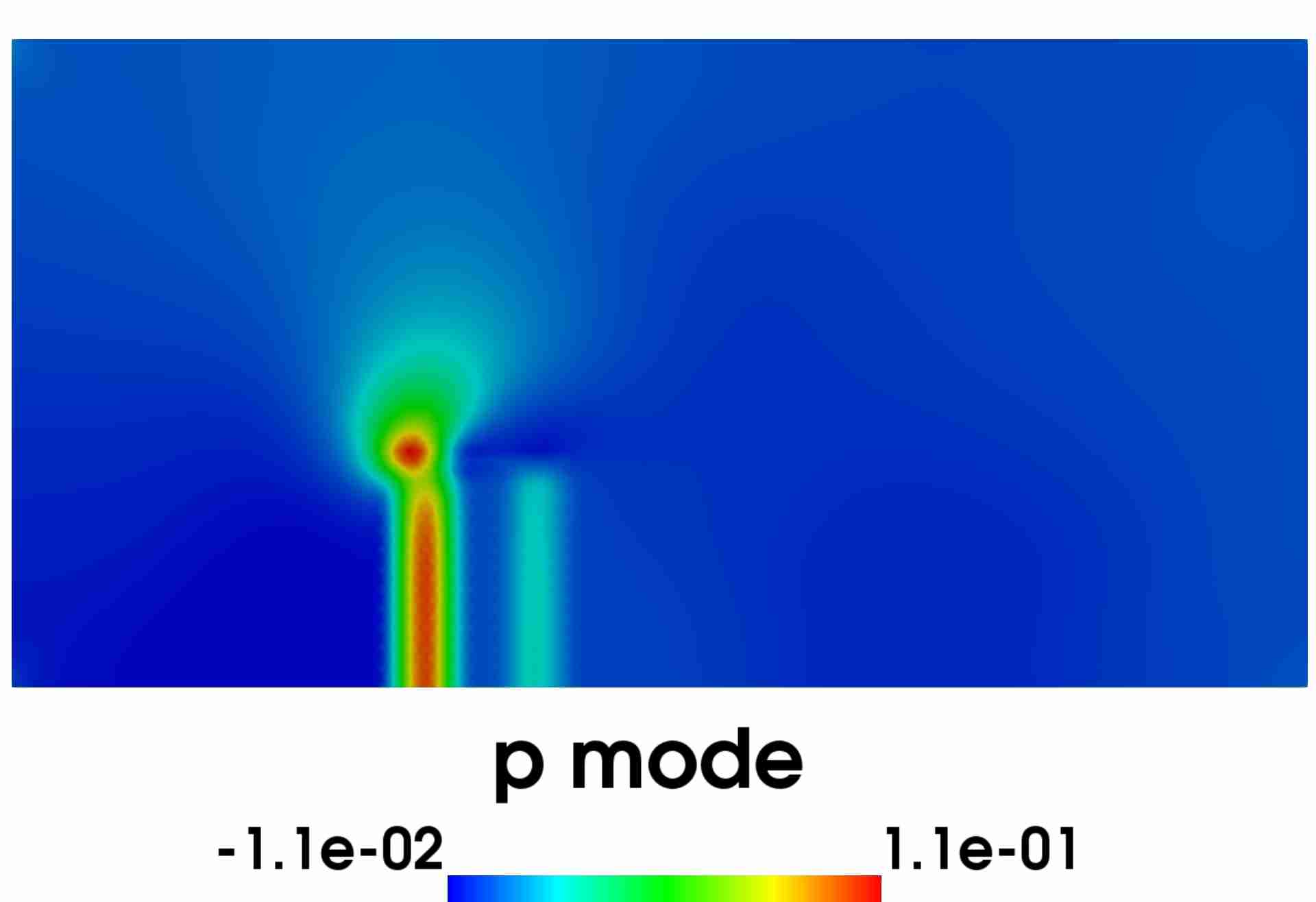}
\end{minipage}
\begin{minipage}{0.24\textwidth}
  \includegraphics[width=\textwidth]{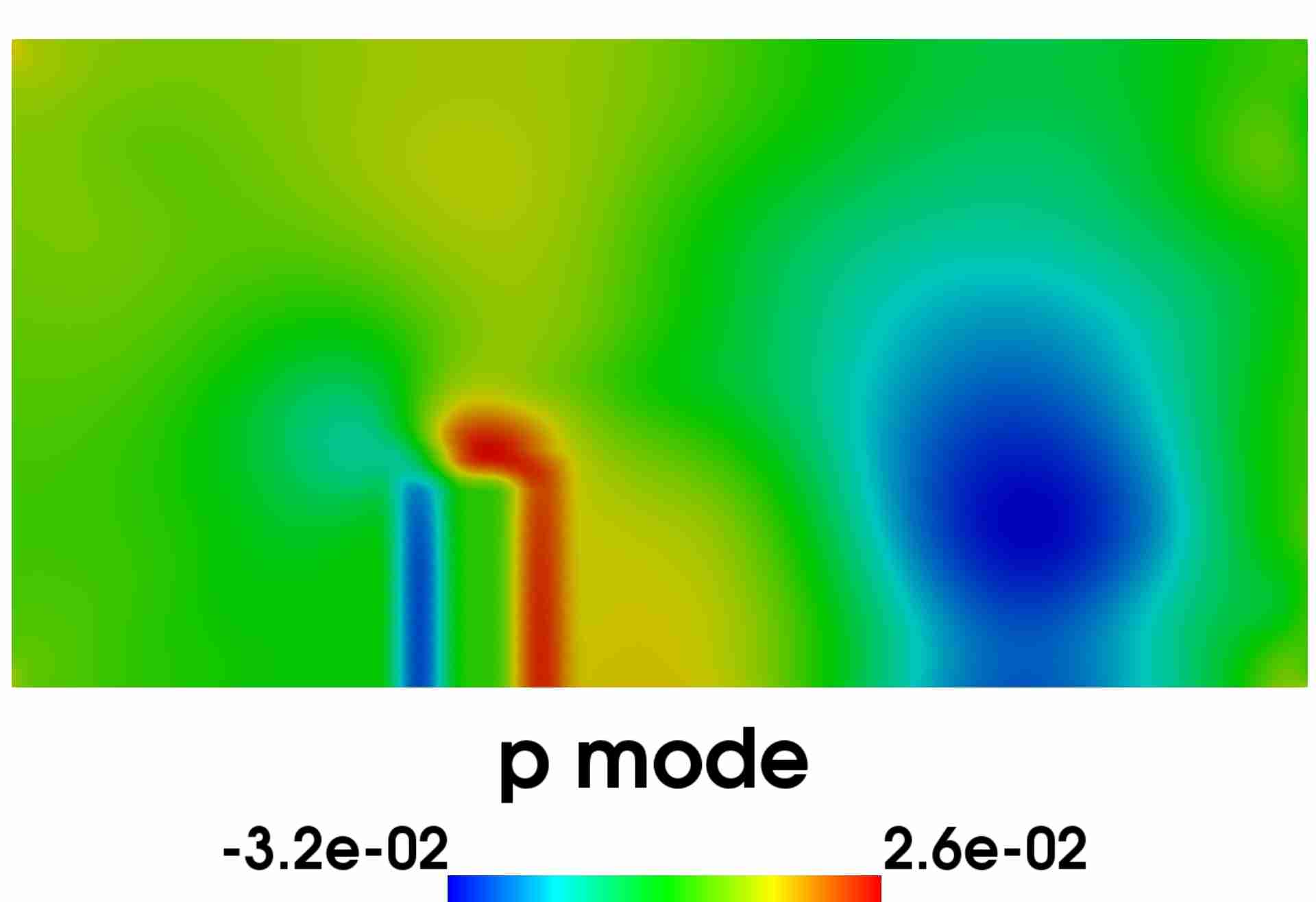}
\end{minipage}
}
\newline
\makebox[\textwidth][l]{%
\begin{minipage}{0.24\textwidth}
  \includegraphics[width=\textwidth]{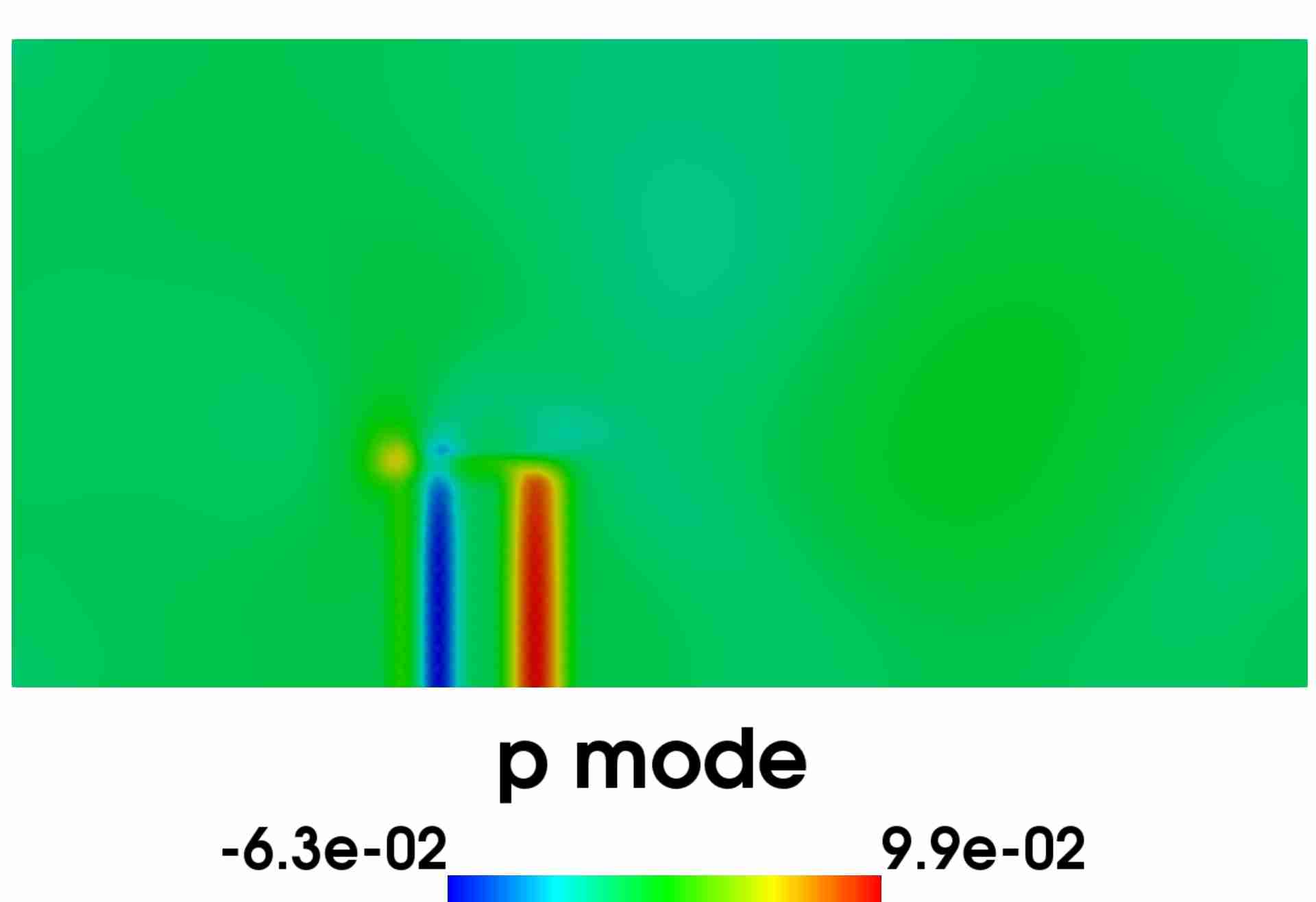}
\end{minipage}
\begin{minipage}{0.24\textwidth}
  \includegraphics[width=\textwidth]{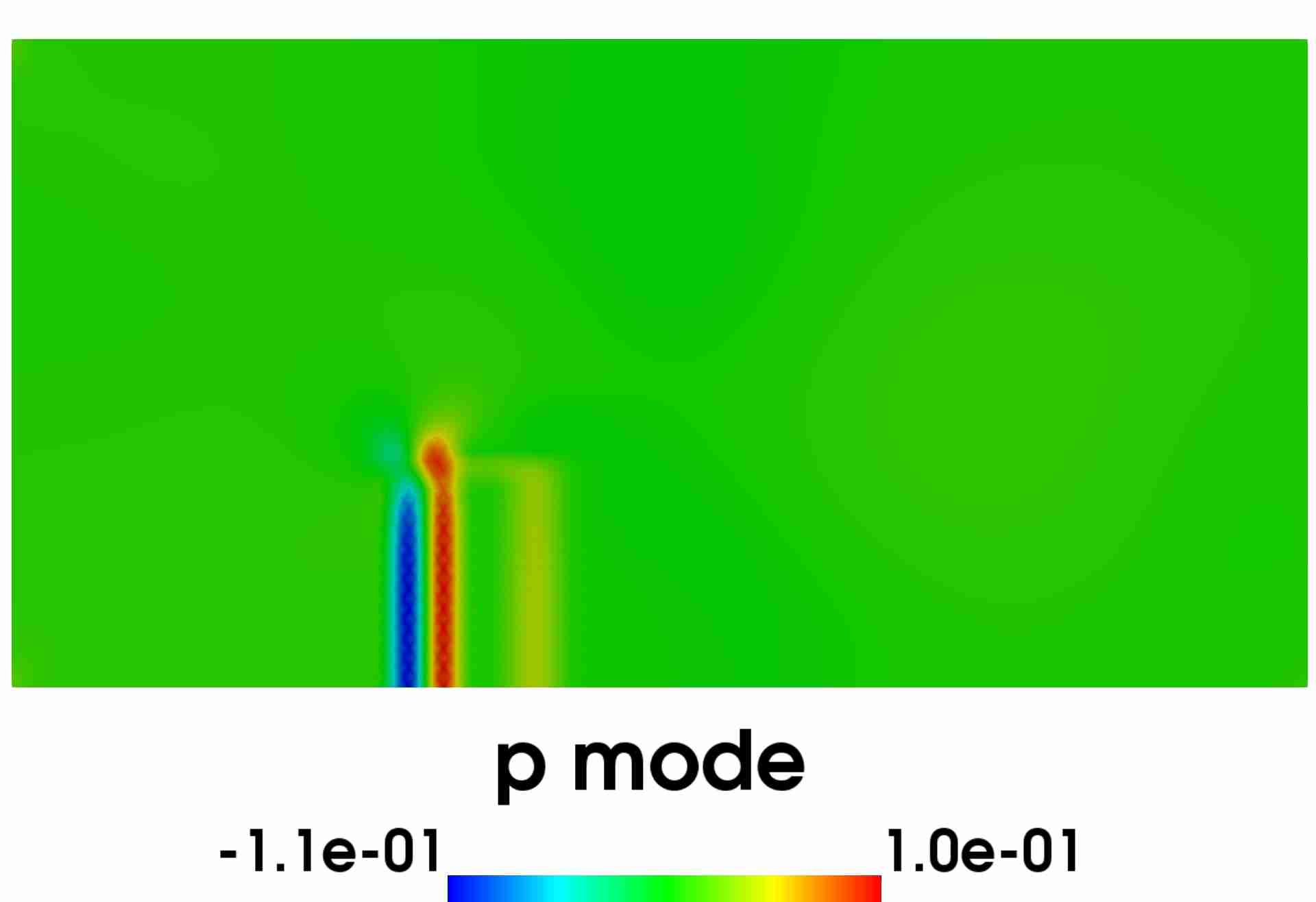}
\end{minipage}
}
\end{minipage}
  \caption{The first $6$ components of the velocity basis (first row) and pressure basis (second row) for the geometrical parametrization with a three-dimensional parameter space case with $\mu \in [-0.65, -0.45]\times[0.25, 0.45]\times[1.25, 1.45]$.}
  \label{3DStokes_Pressure_and_Velocity_Components_Modes}
\end{figure} 
\begin{figure} 
\centering
\begin{minipage}{\textwidth}
\centering
\begin{minipage}{0.24\textwidth}
  \includegraphics[width=\textwidth]{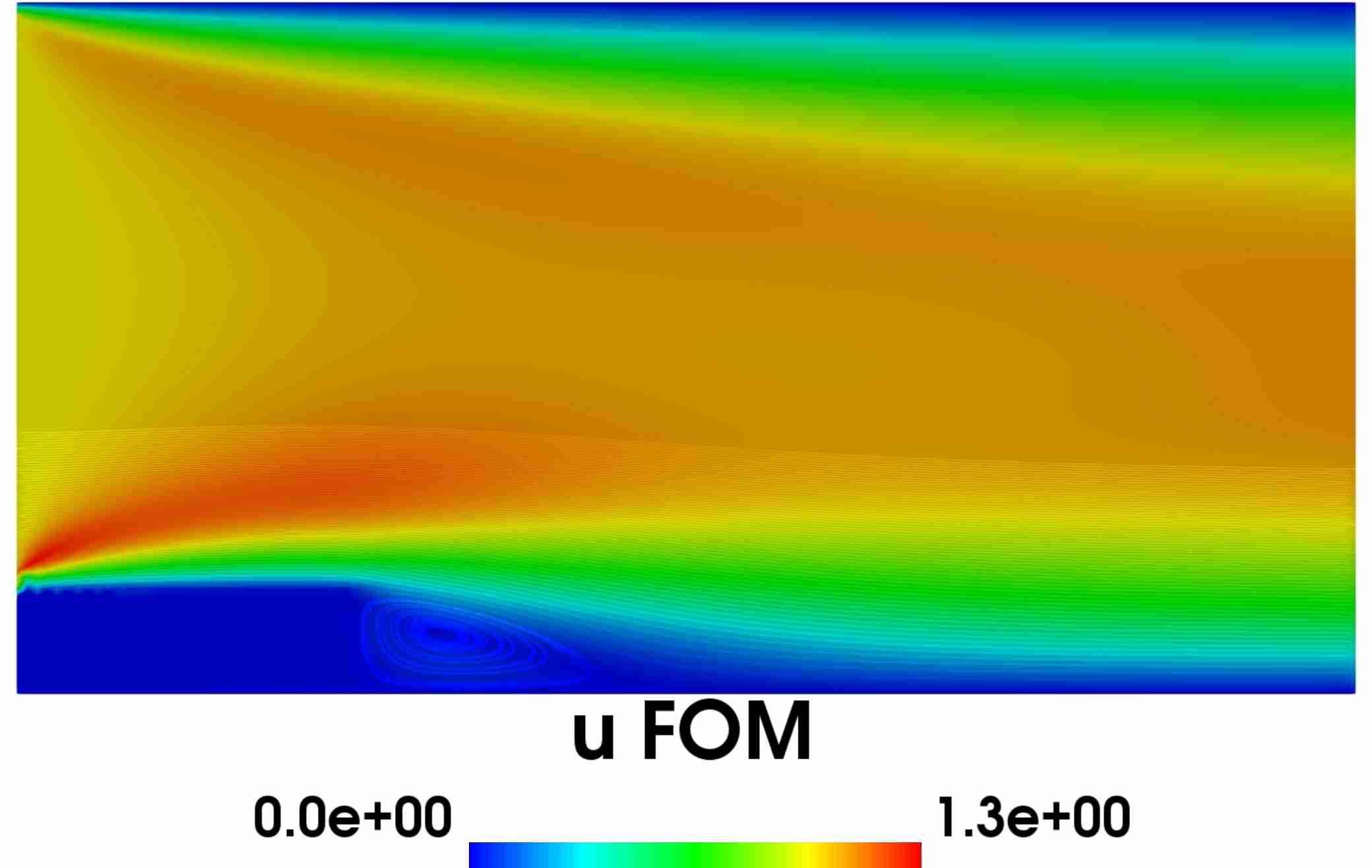} 
\end{minipage}
\begin{minipage}{0.24\textwidth}
  \includegraphics[width=\textwidth]{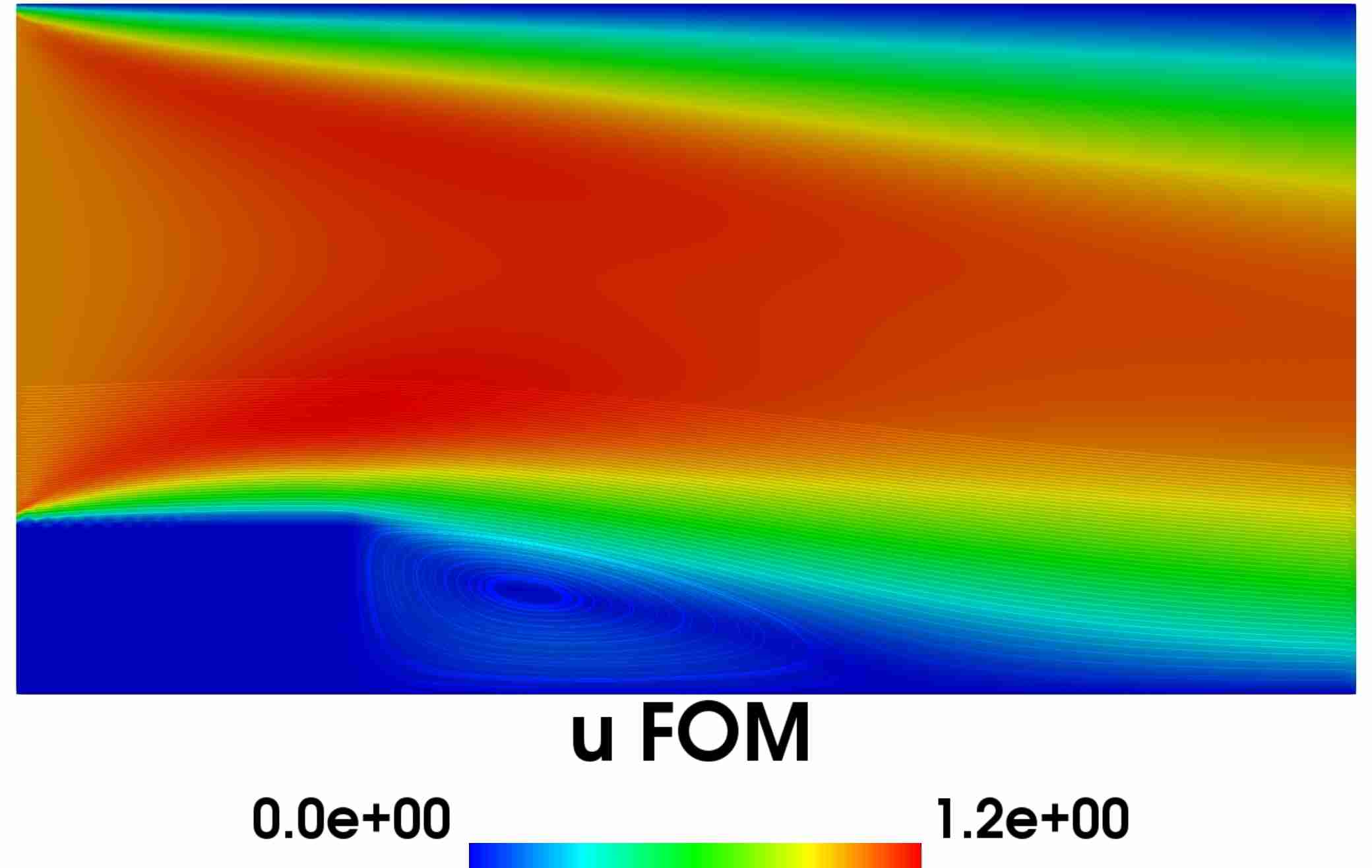}
\end{minipage}
\begin{minipage}{0.24\textwidth}
  \includegraphics[width=\textwidth]{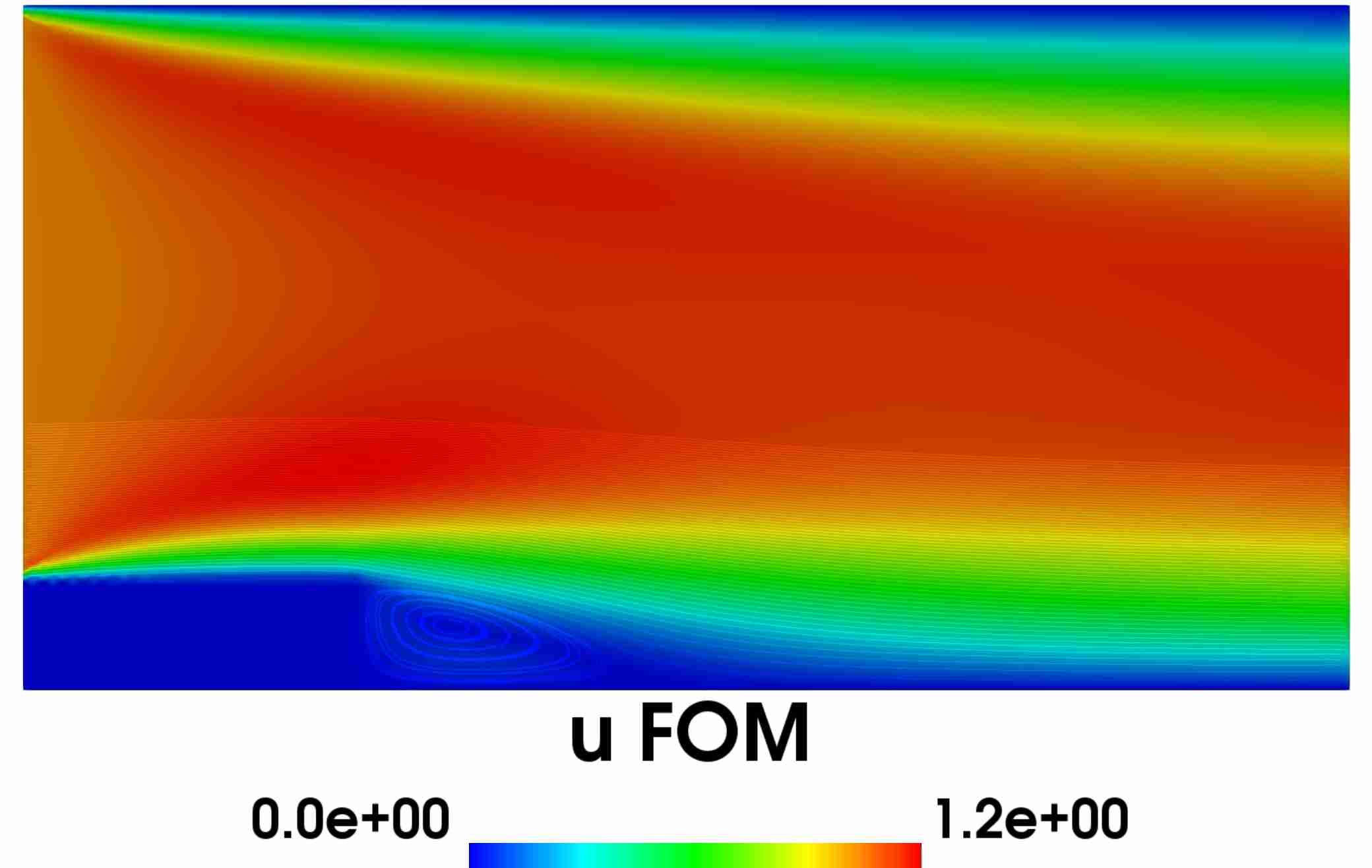}
\end{minipage}
\begin{minipage}{0.24\textwidth}
  \includegraphics[width=\textwidth]{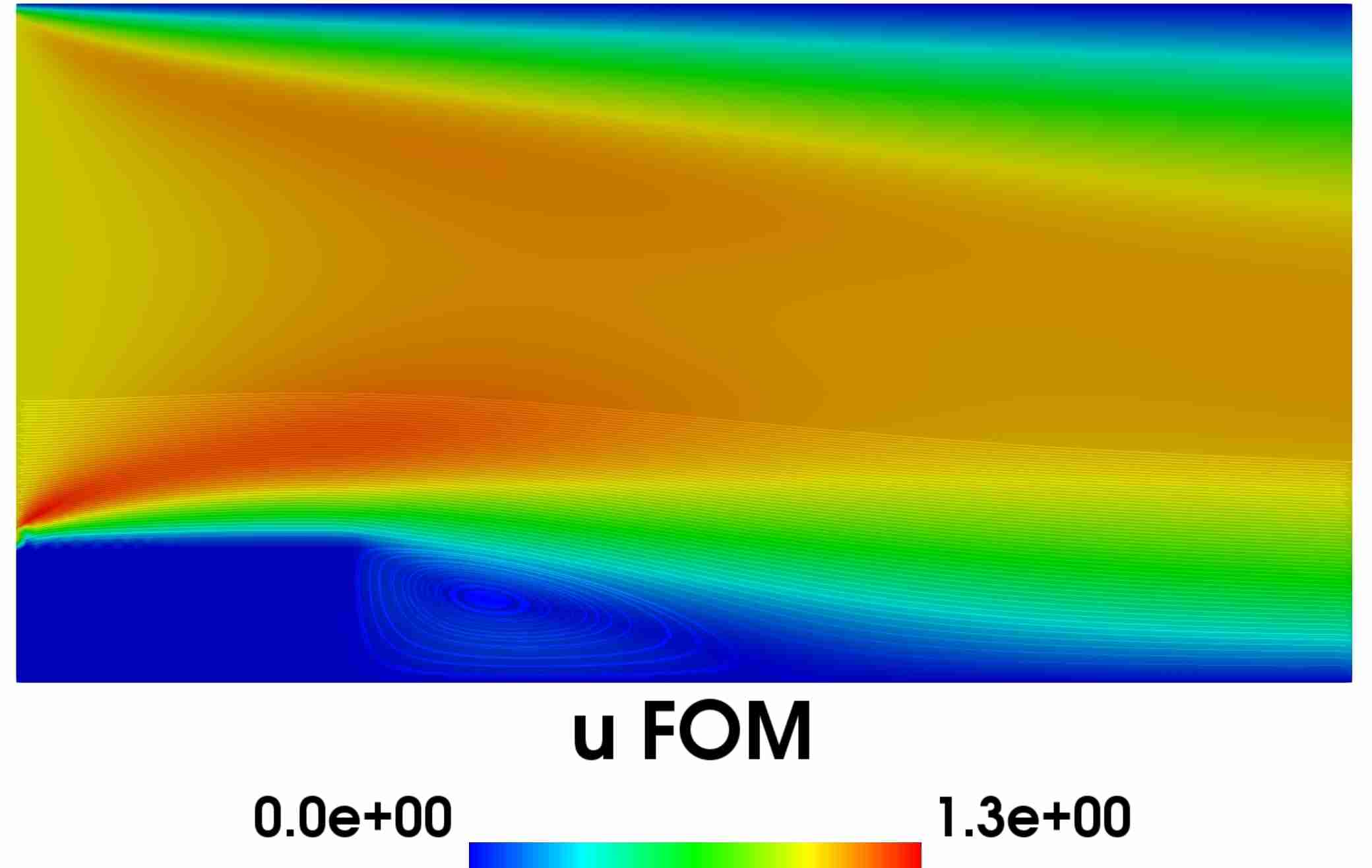}
\end{minipage}
\begin{minipage}{0.24\textwidth}
  \includegraphics[width=\textwidth]{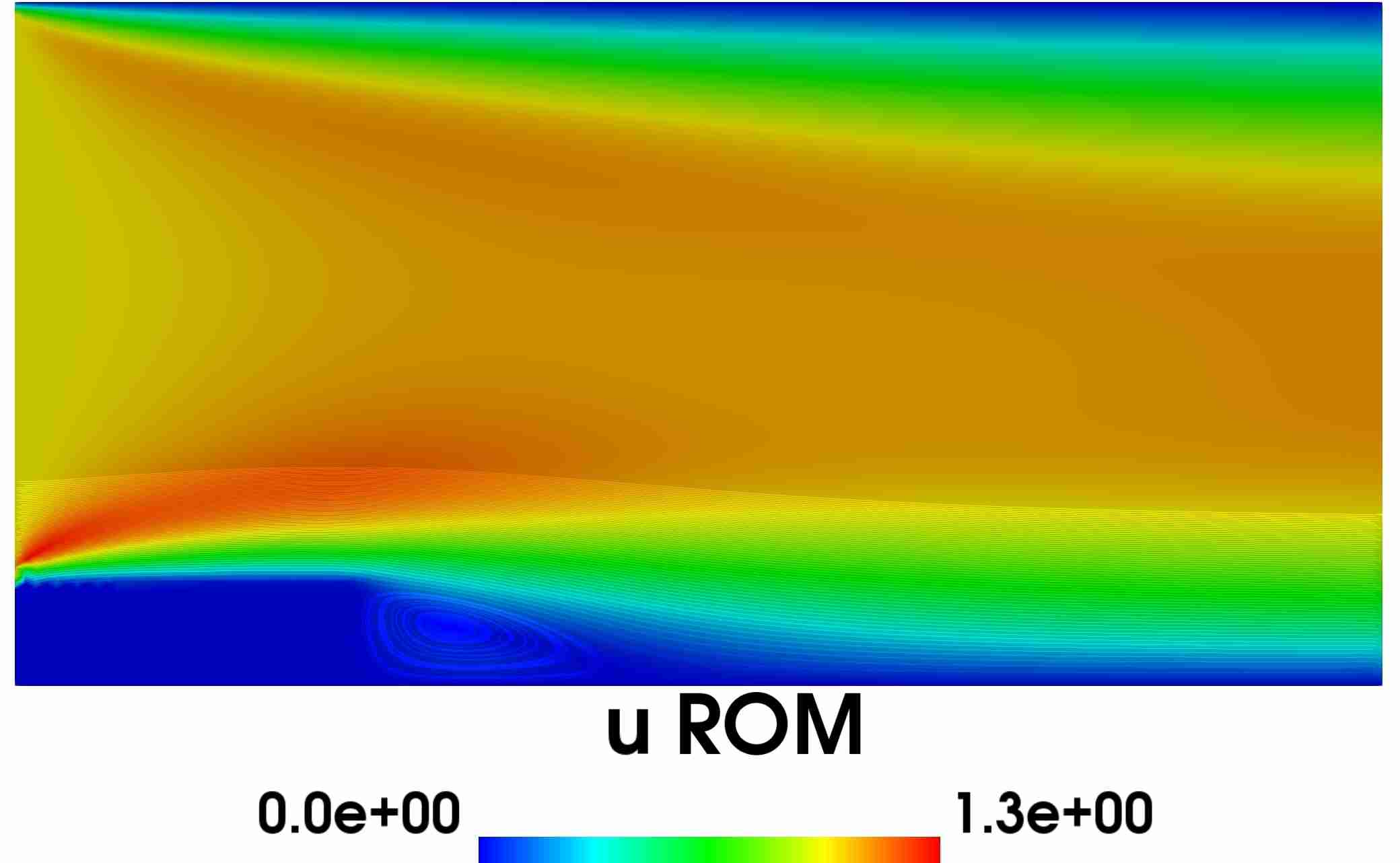} 
\end{minipage}
\begin{minipage}{0.24\textwidth}
  \includegraphics[width=\textwidth]{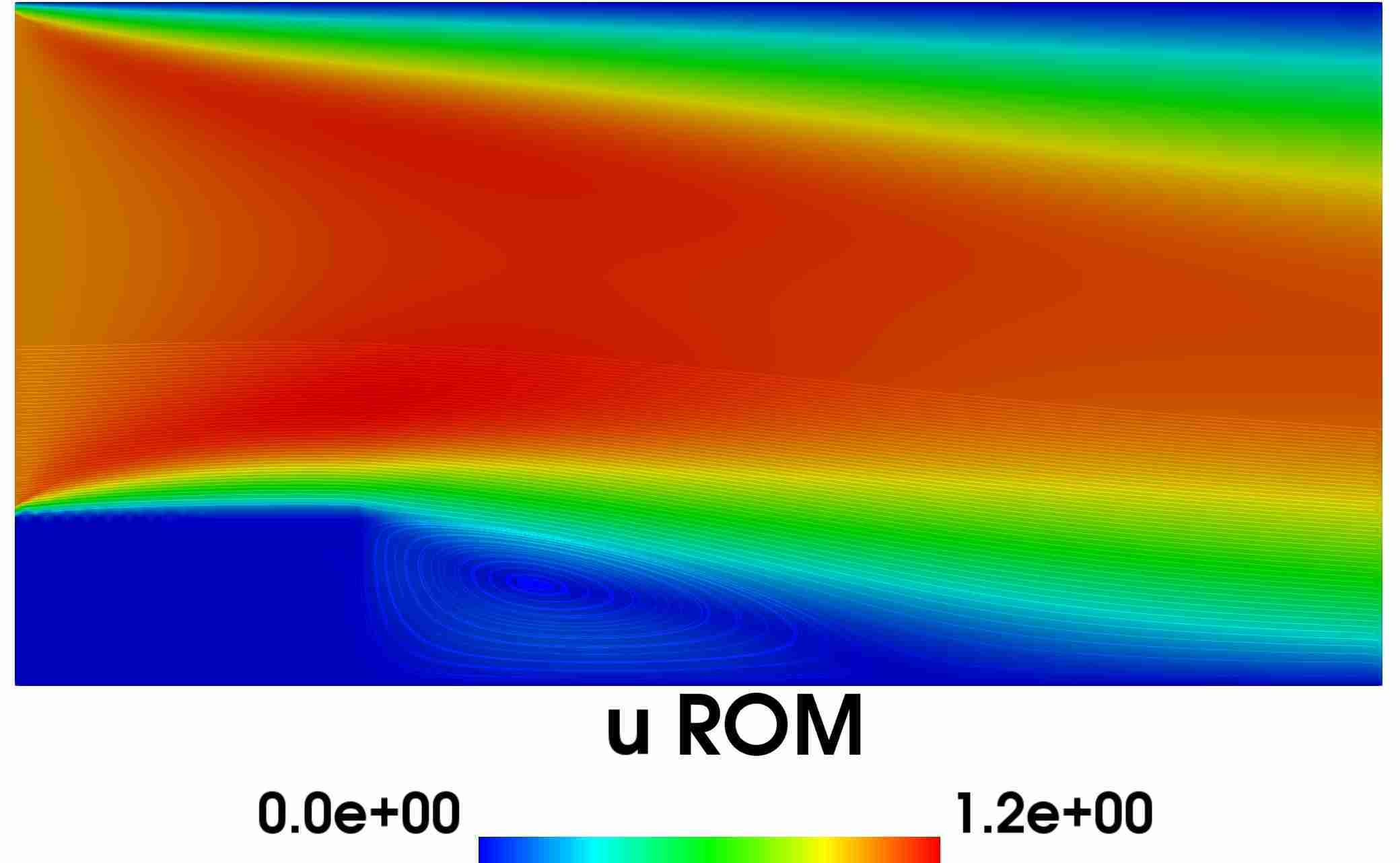}
\end{minipage}
\begin{minipage}{0.24\textwidth}
  \includegraphics[width=\textwidth]{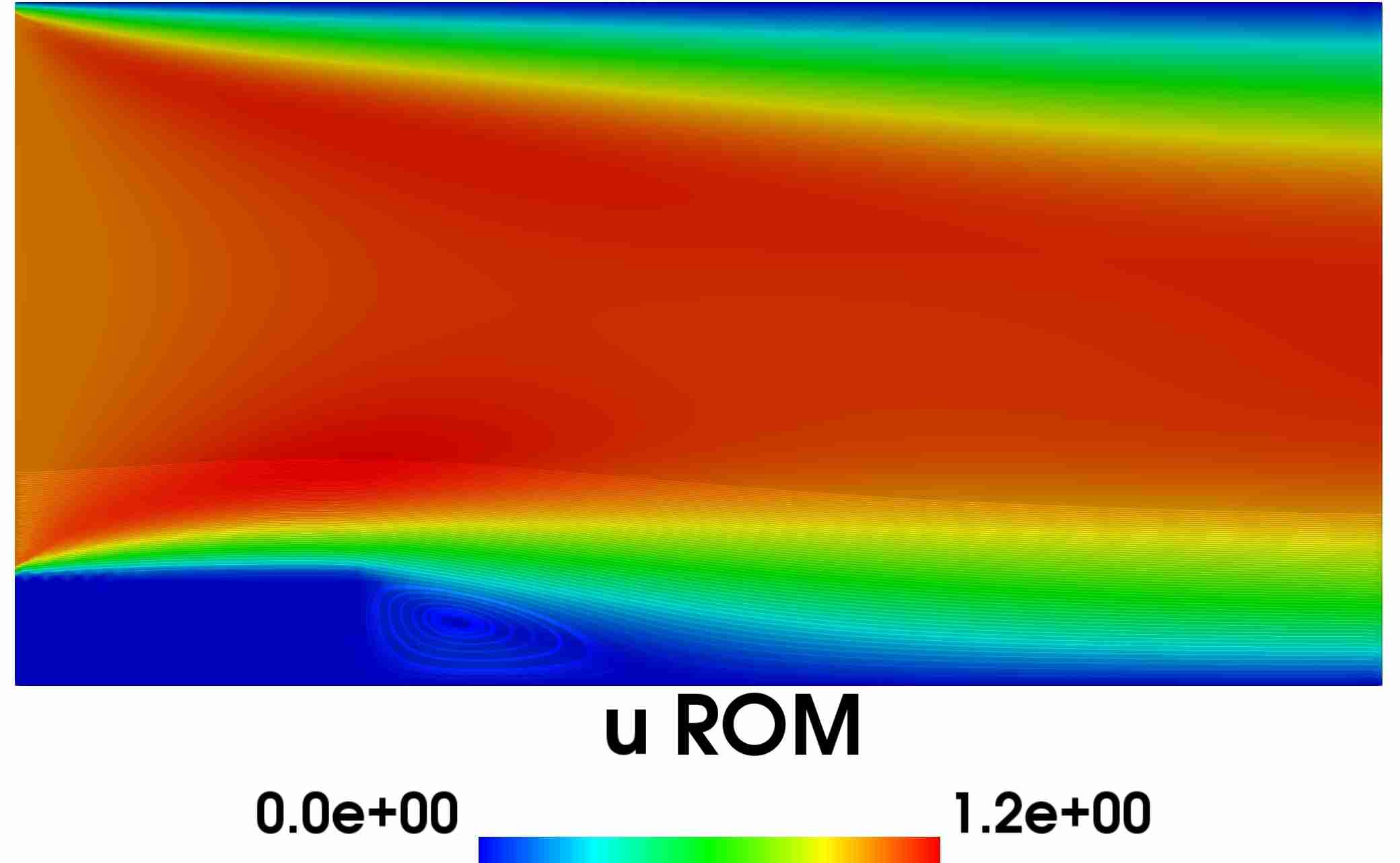}
\end{minipage}
\begin{minipage}{0.24\textwidth}
  \includegraphics[width=\textwidth]{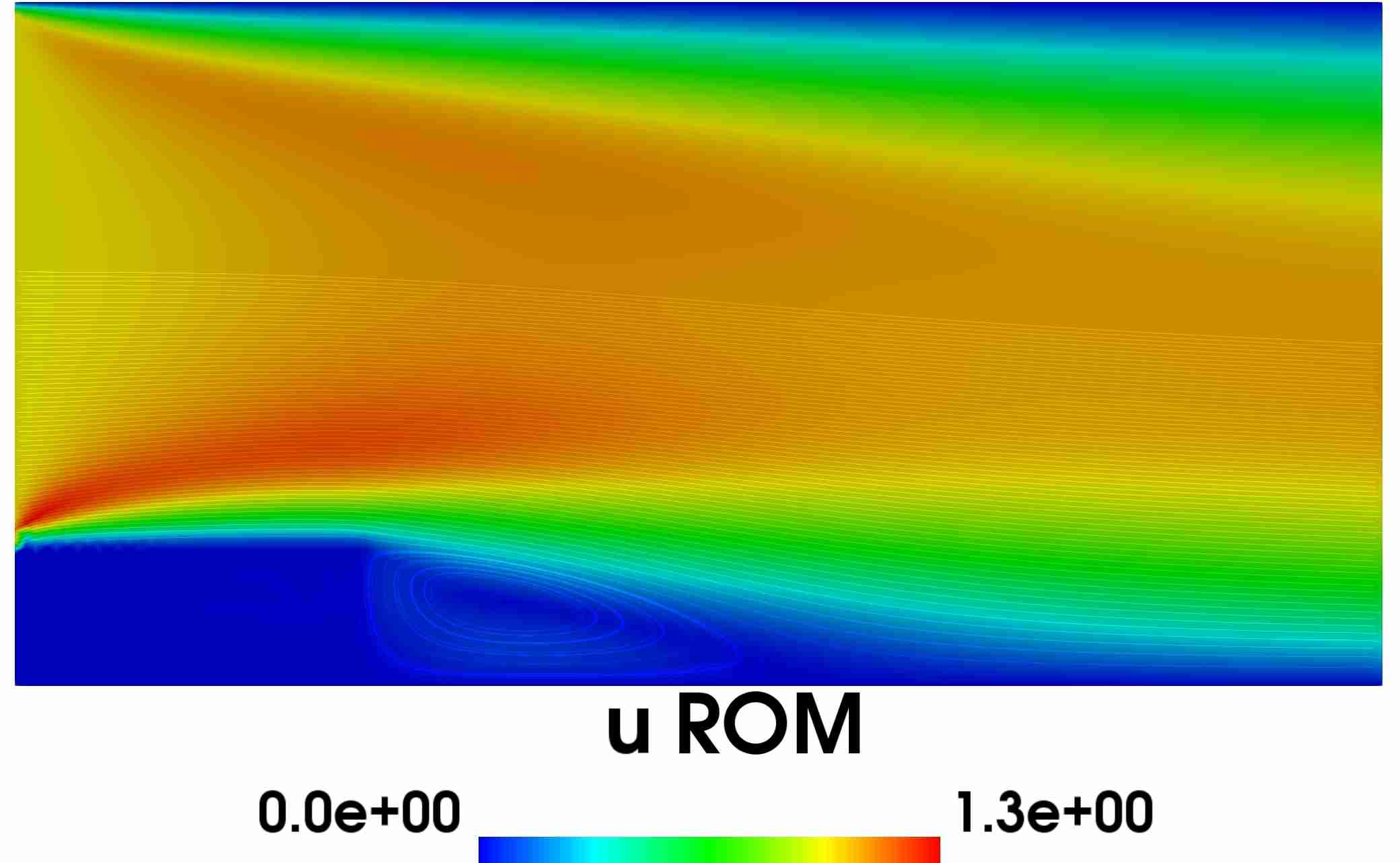}
\end{minipage}
\begin{minipage}{0.24\textwidth}
  \includegraphics[width=\textwidth]{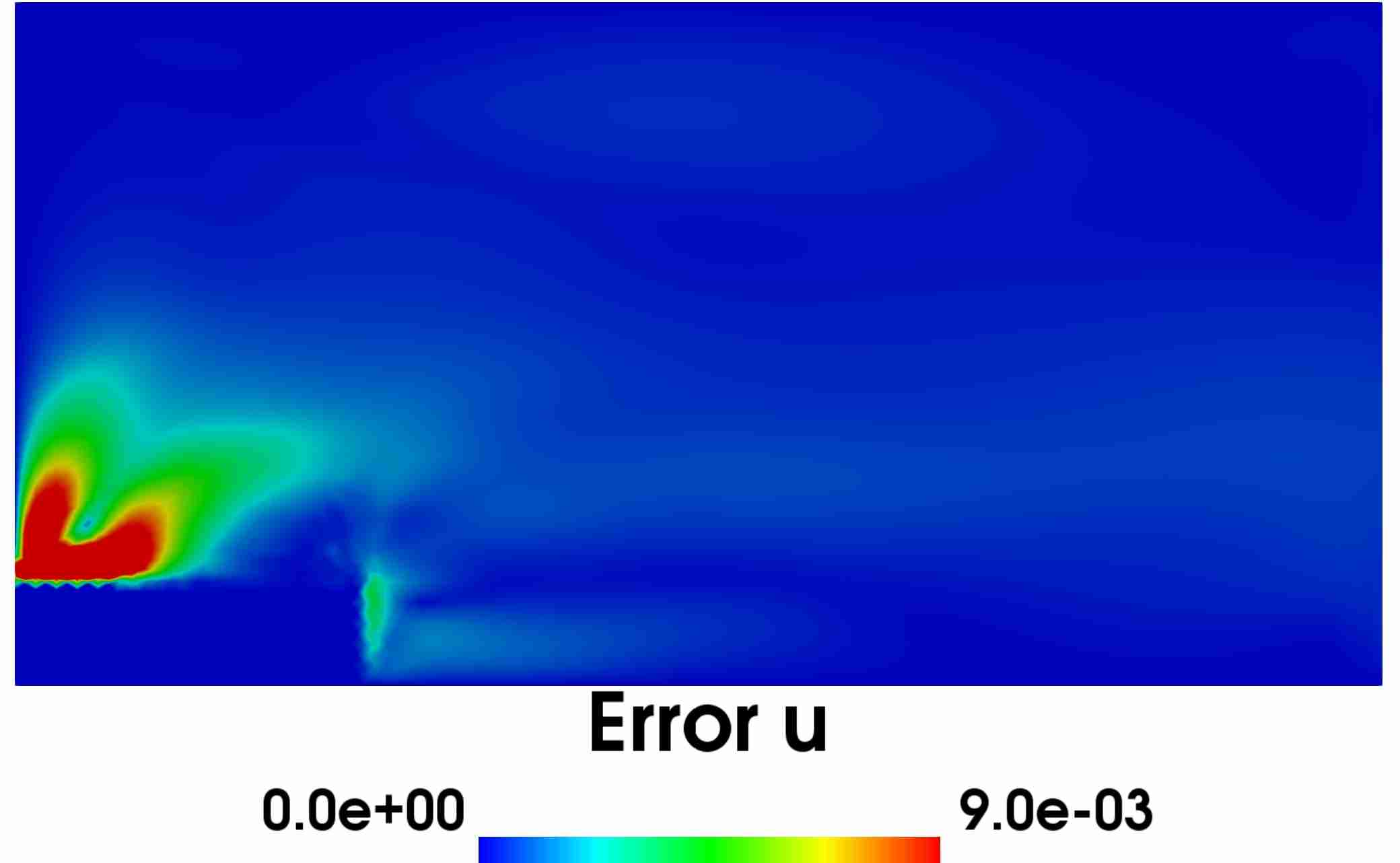} 
\end{minipage}
\begin{minipage}{0.24\textwidth}
  \includegraphics[width=\textwidth]{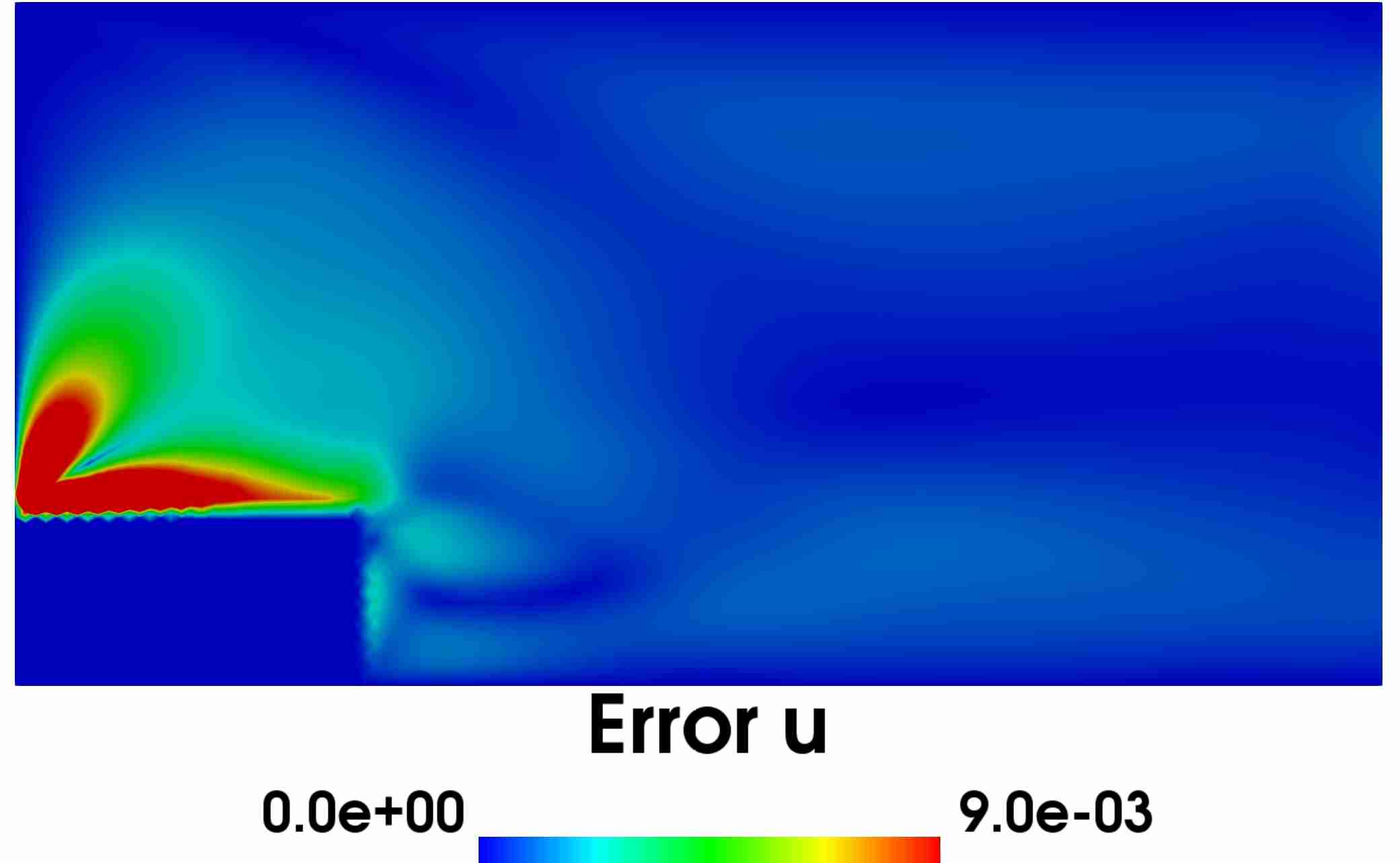}
\end{minipage}
\begin{minipage}{0.24\textwidth}
  \includegraphics[width=\textwidth]{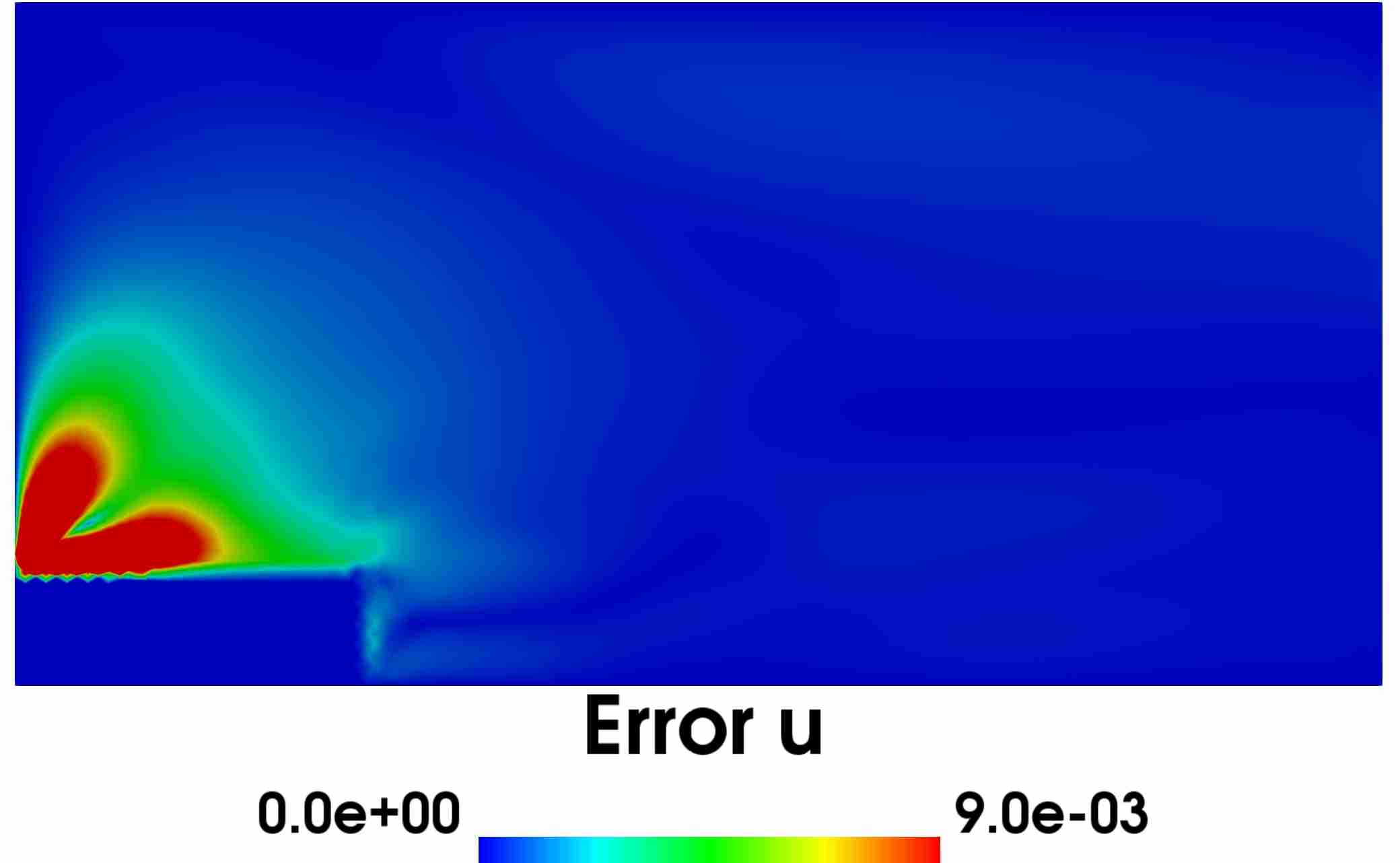}
\end{minipage}
\begin{minipage}{0.24\textwidth}
  \includegraphics[width=\textwidth]{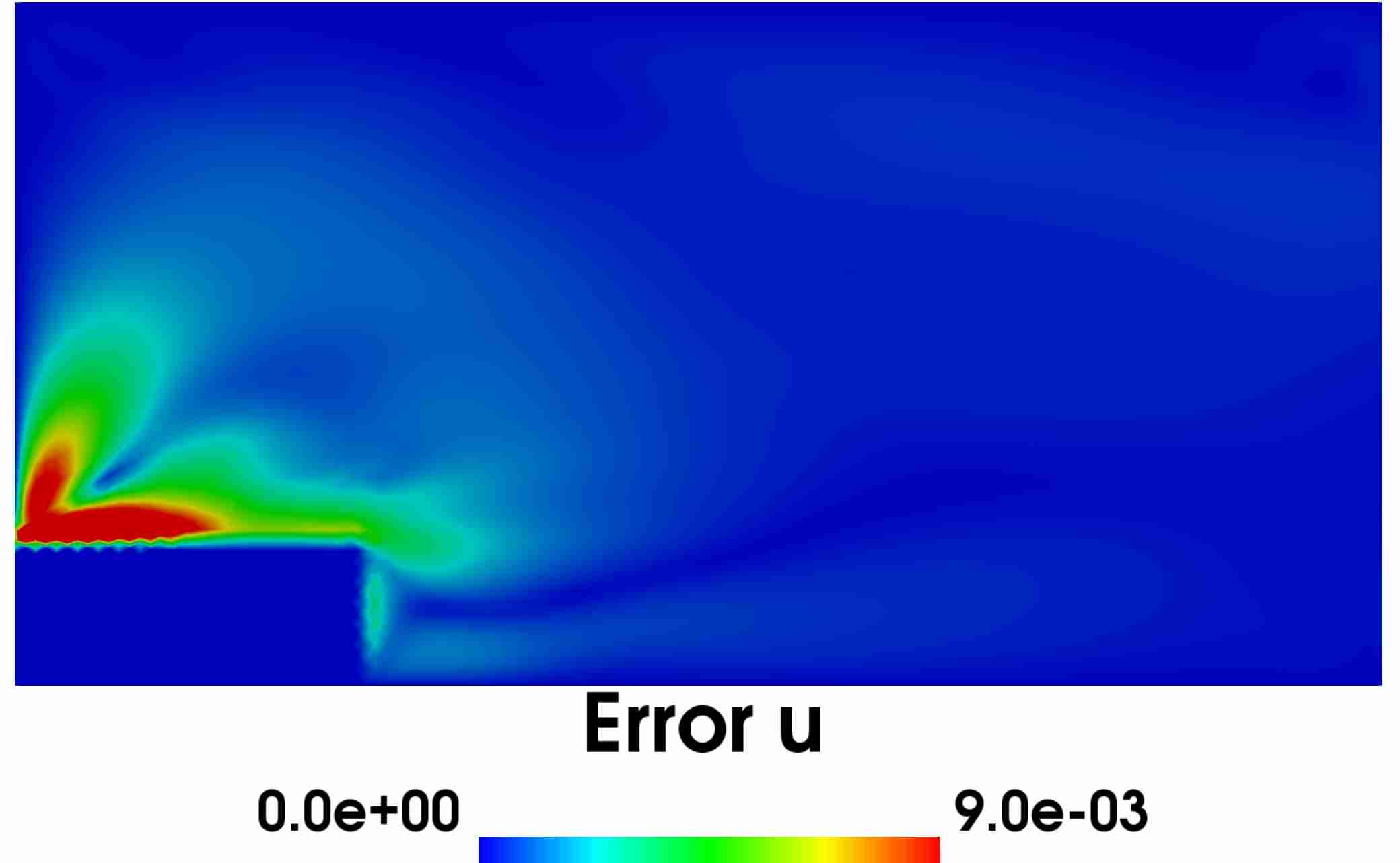}
\end{minipage}
\begin{minipage}{0.24\textwidth}
  \includegraphics[width=\textwidth]{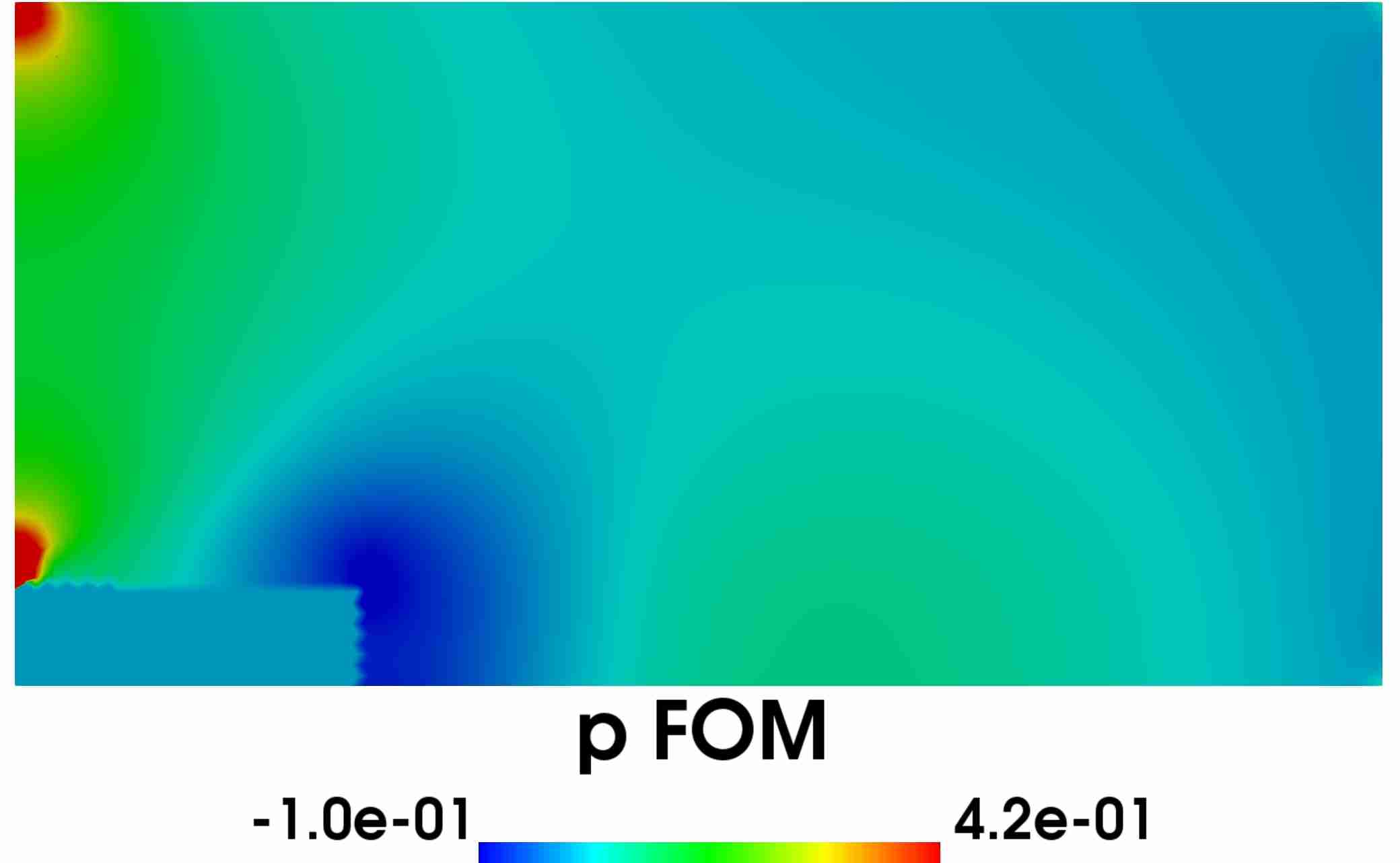} 
\end{minipage}
\begin{minipage}{0.24\textwidth}
  \includegraphics[width=\textwidth]{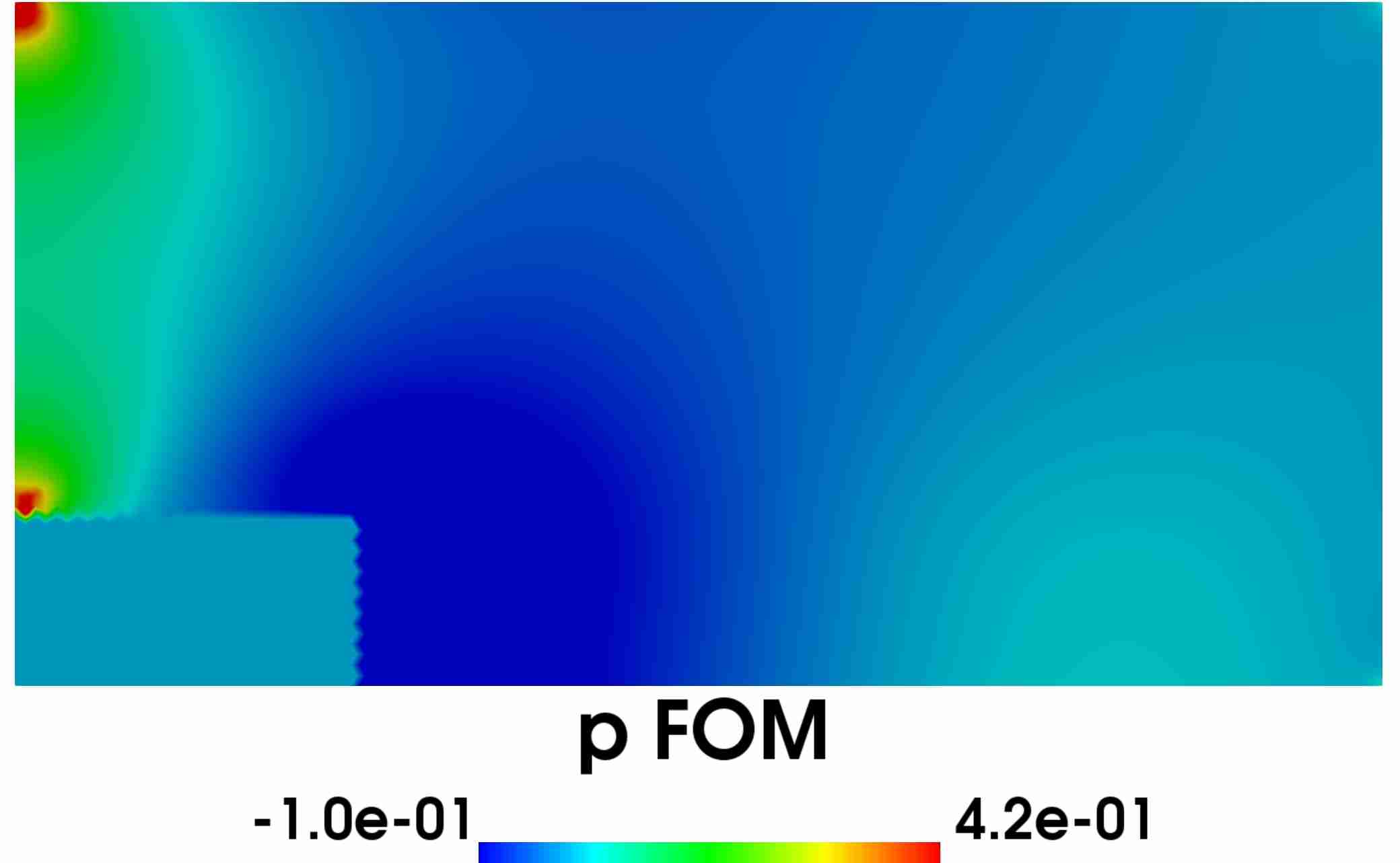}
\end{minipage}
\begin{minipage}{0.24\textwidth}
  \includegraphics[width=\textwidth]{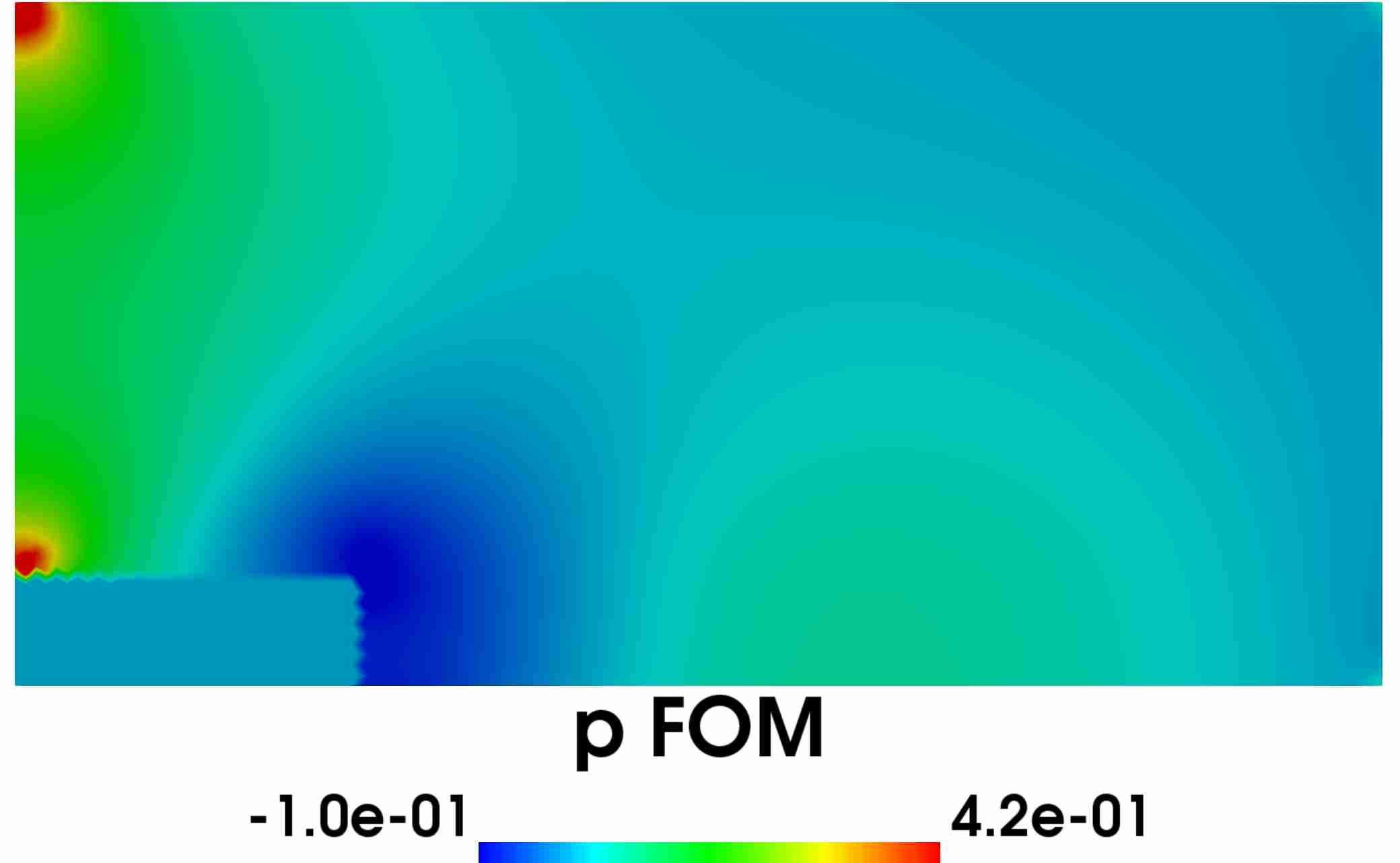}
\end{minipage}
\begin{minipage}{0.24\textwidth}
  \includegraphics[width=\textwidth]{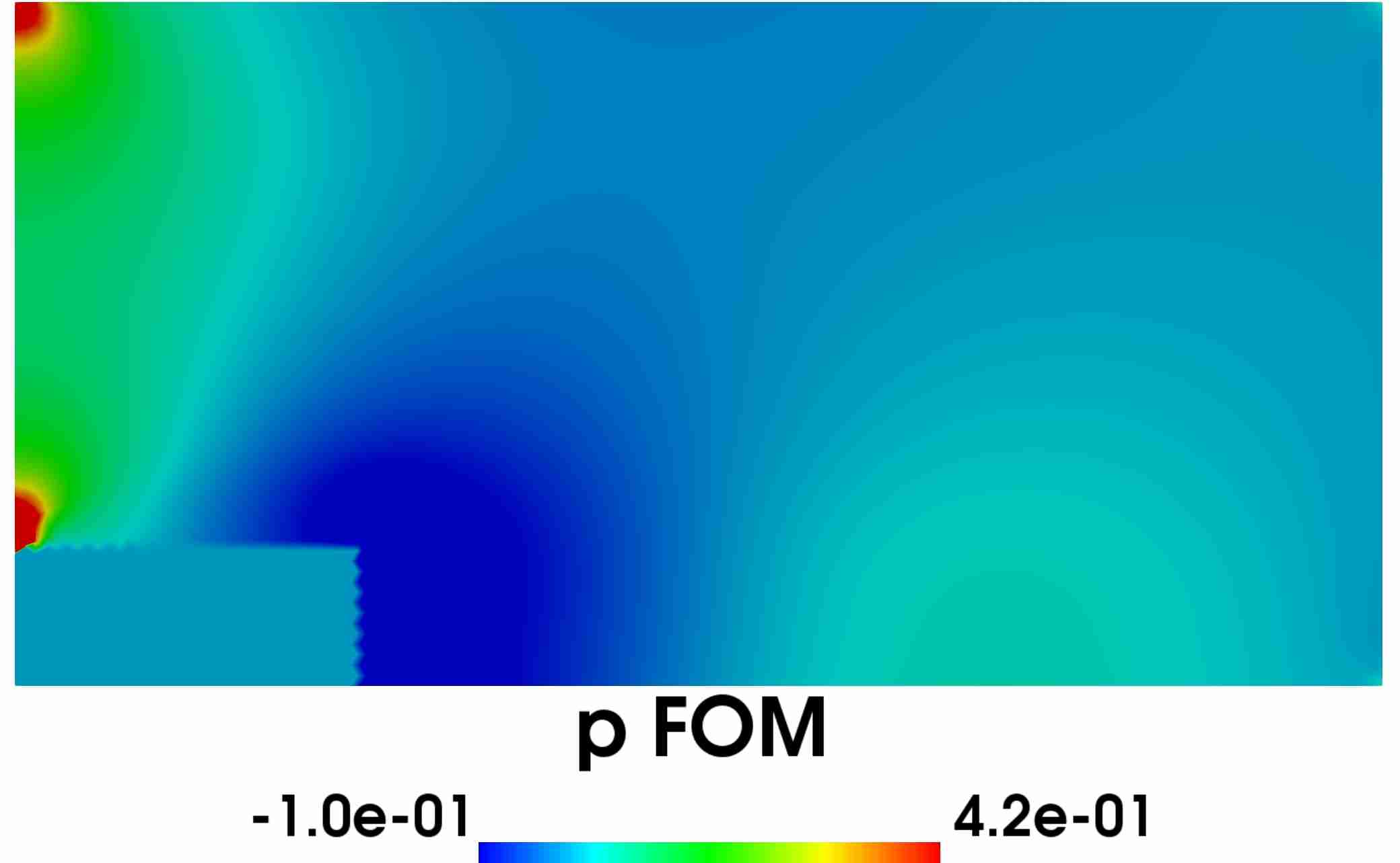}
\end{minipage}
\begin{minipage}{0.24\textwidth}
  \includegraphics[width=\textwidth]{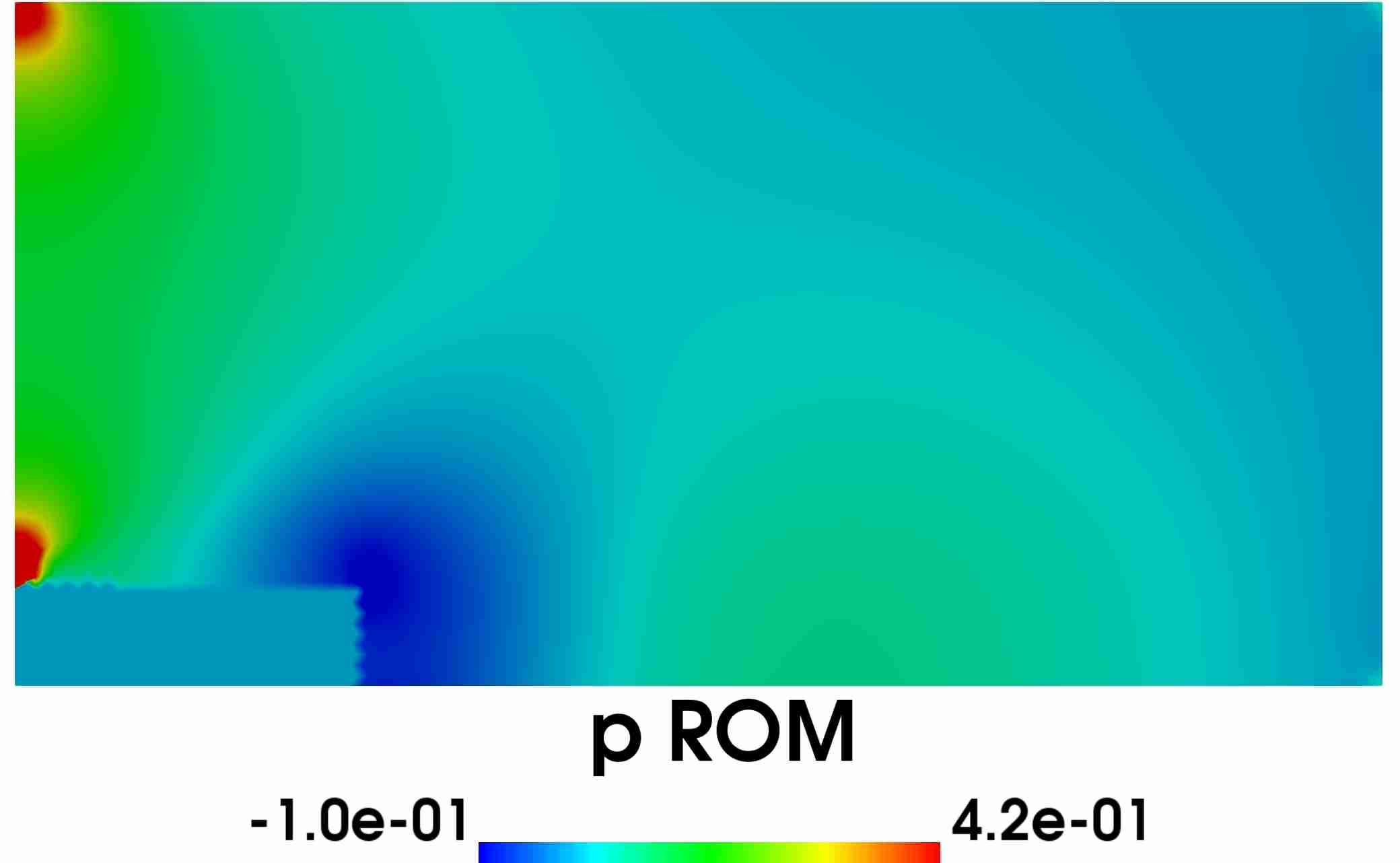} 
\end{minipage}
\begin{minipage}{0.24\textwidth}
  \includegraphics[width=\textwidth]{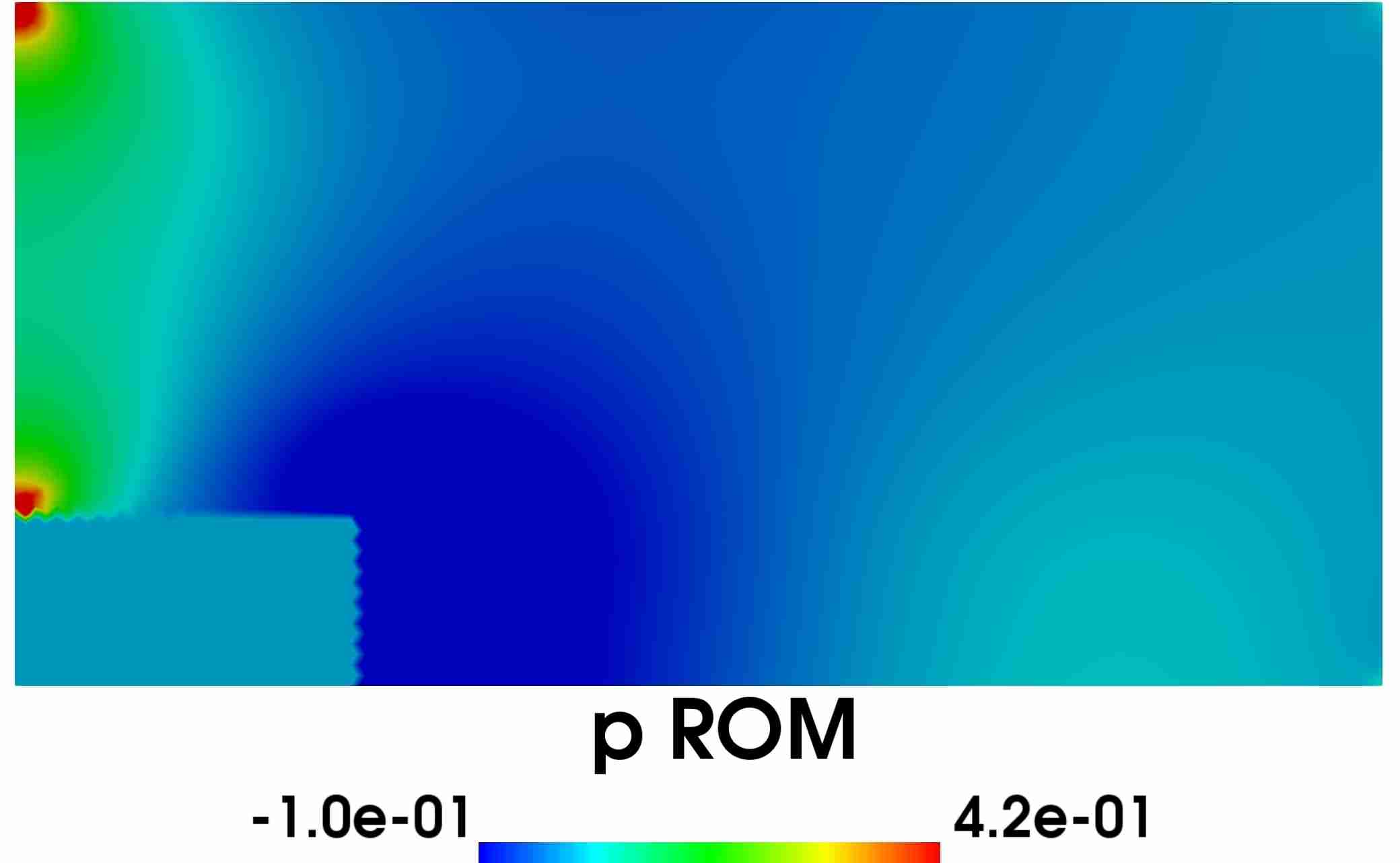}
\end{minipage}
\begin{minipage}{0.24\textwidth}
  \includegraphics[width=\textwidth]{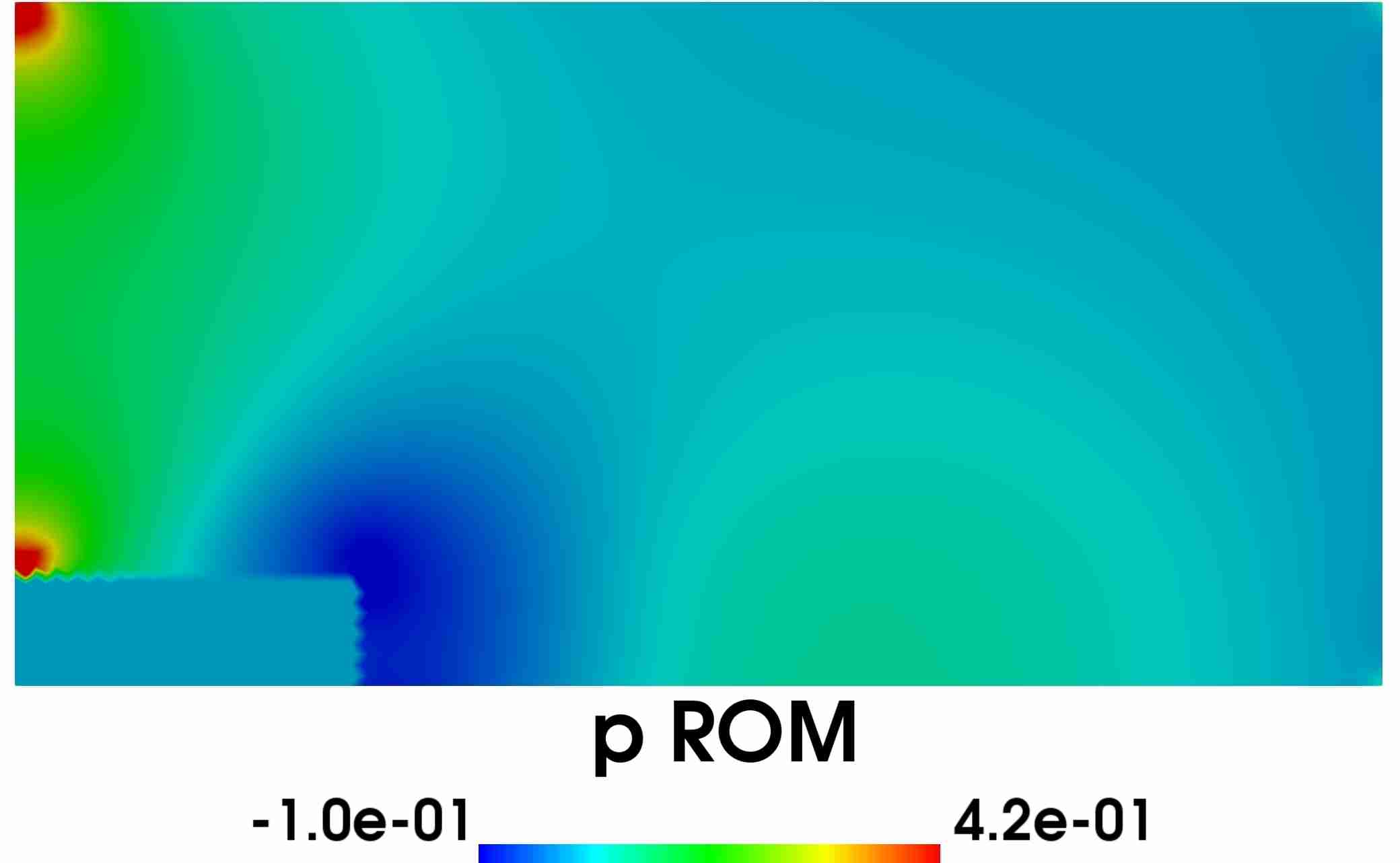}
\end{minipage}
\begin{minipage}{0.24\textwidth}
  \includegraphics[width=\textwidth]{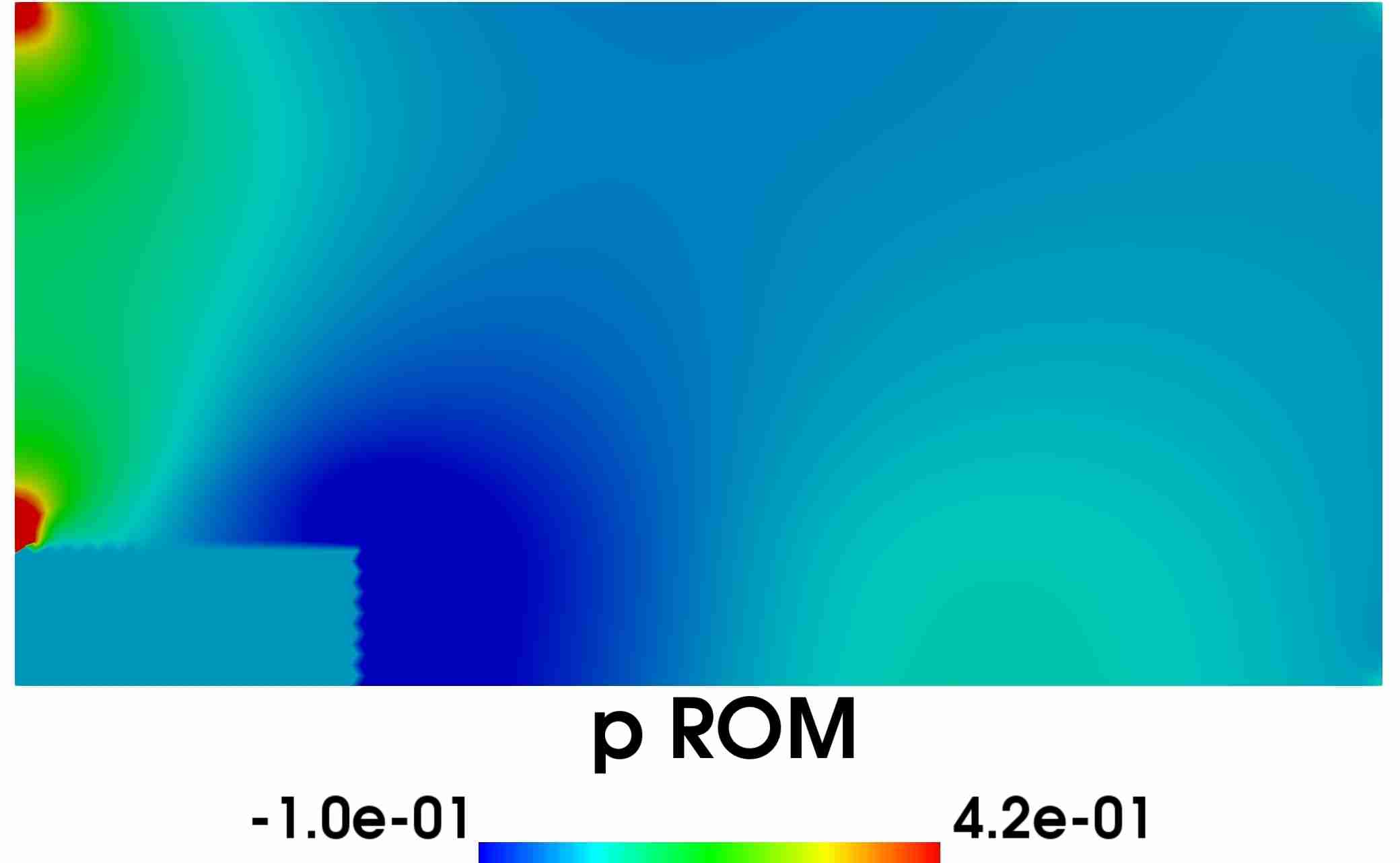}
\end{minipage}
\begin{minipage}{0.24\textwidth}
  \includegraphics[width=\textwidth]{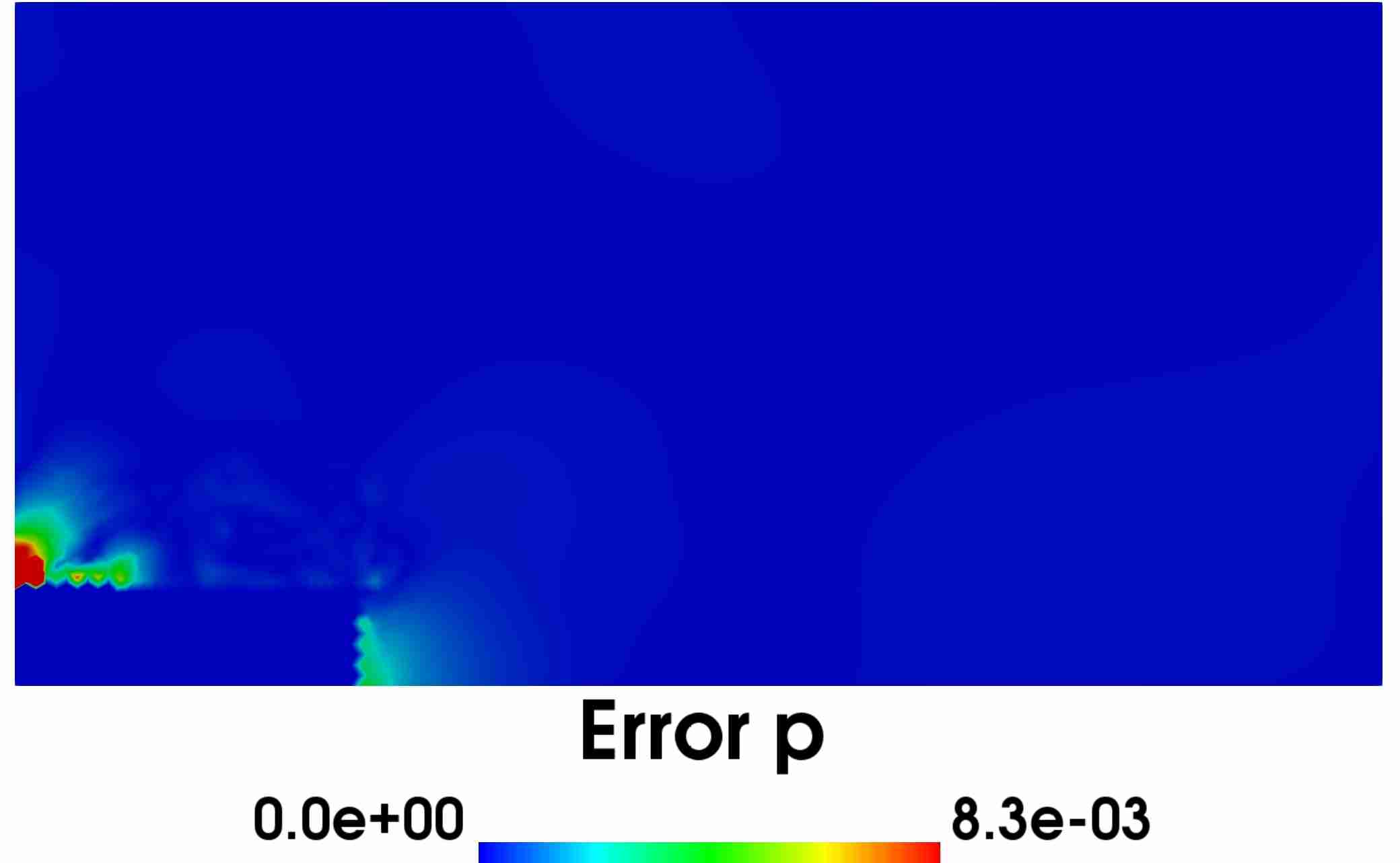} 
\end{minipage}
\begin{minipage}{0.24\textwidth}
  \includegraphics[width=\textwidth]{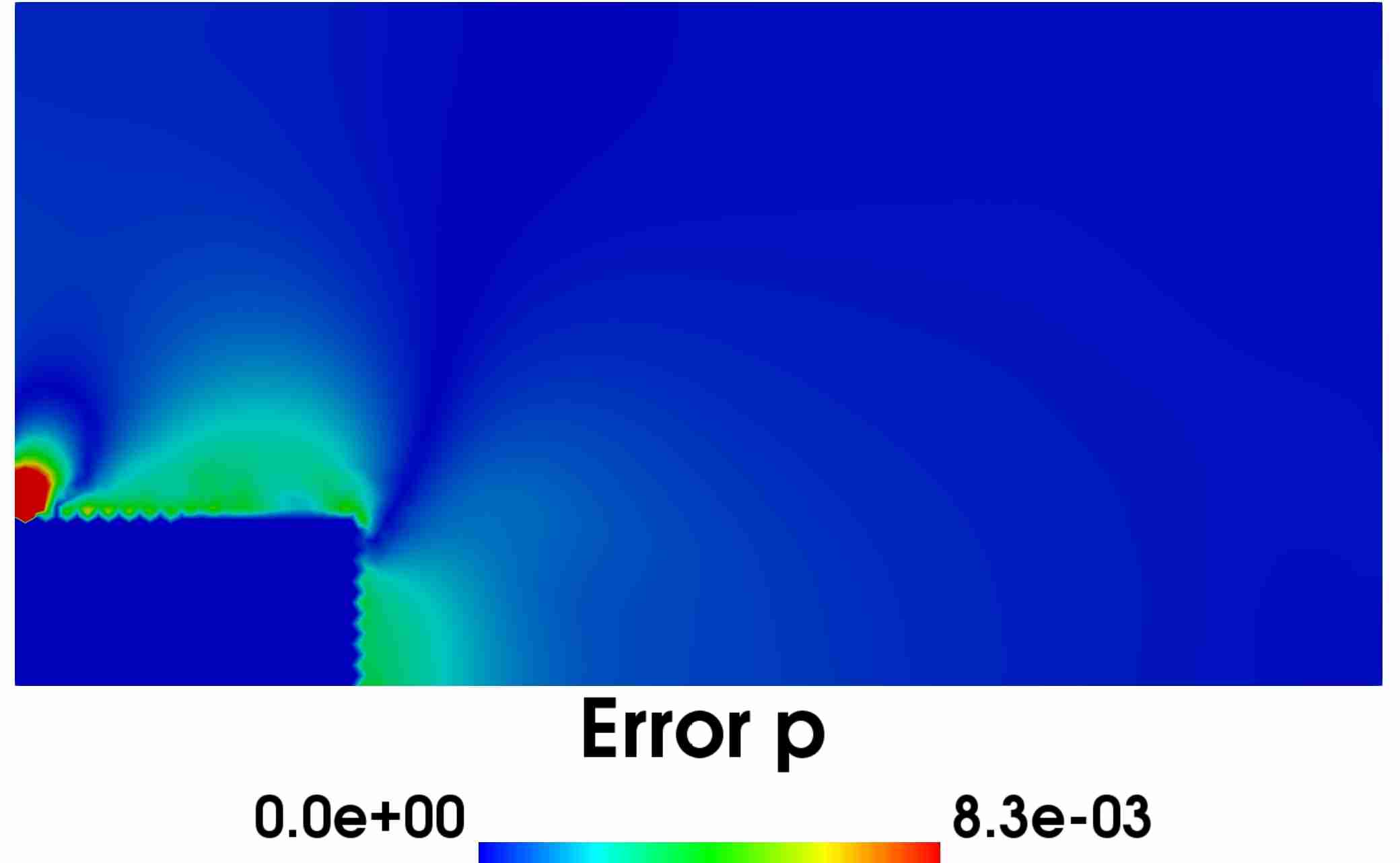}
\end{minipage}
\begin{minipage}{0.24\textwidth}
  \includegraphics[width=\textwidth]{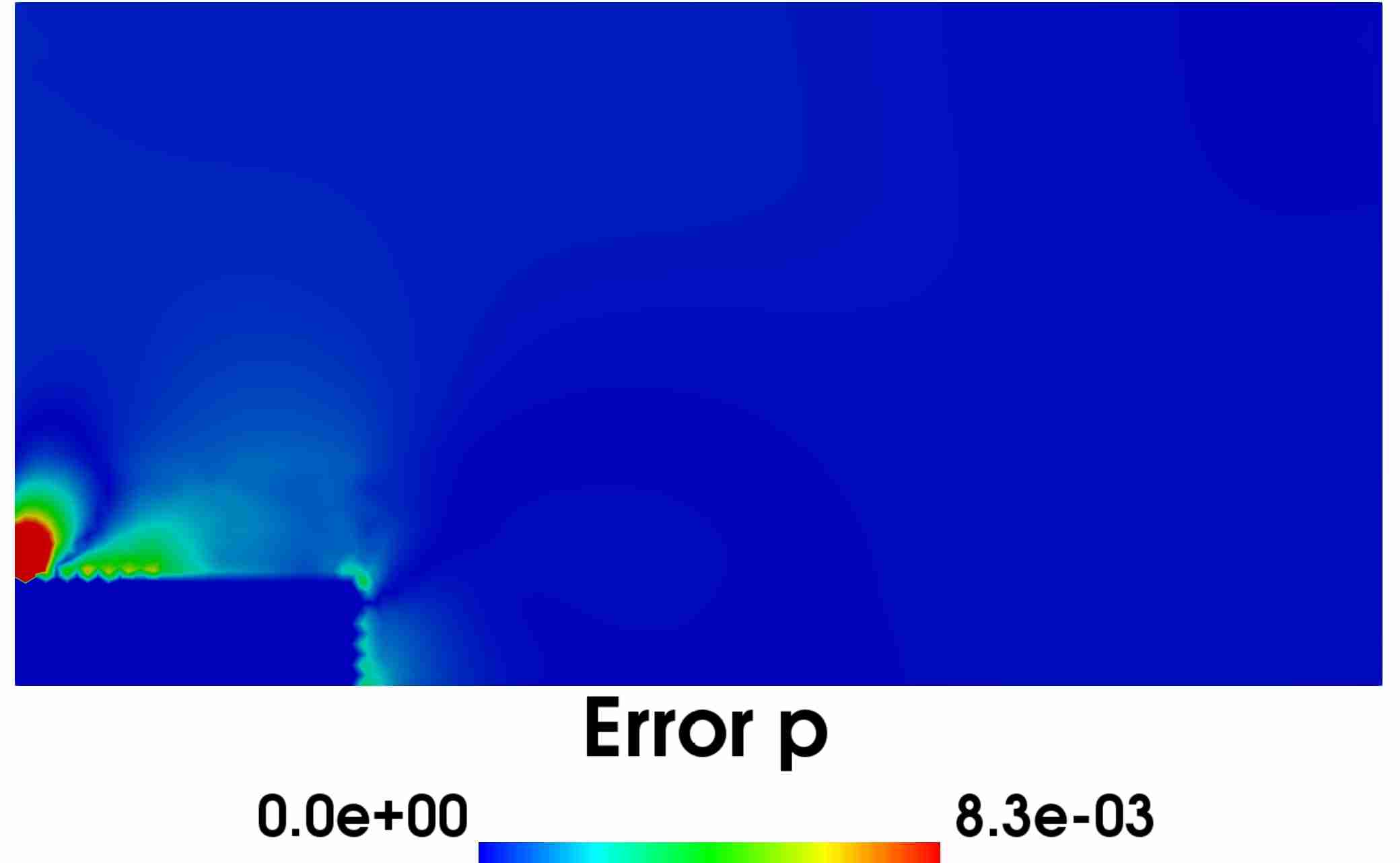}
\end{minipage}
\begin{minipage}{0.24\textwidth}
  \includegraphics[width=\textwidth]{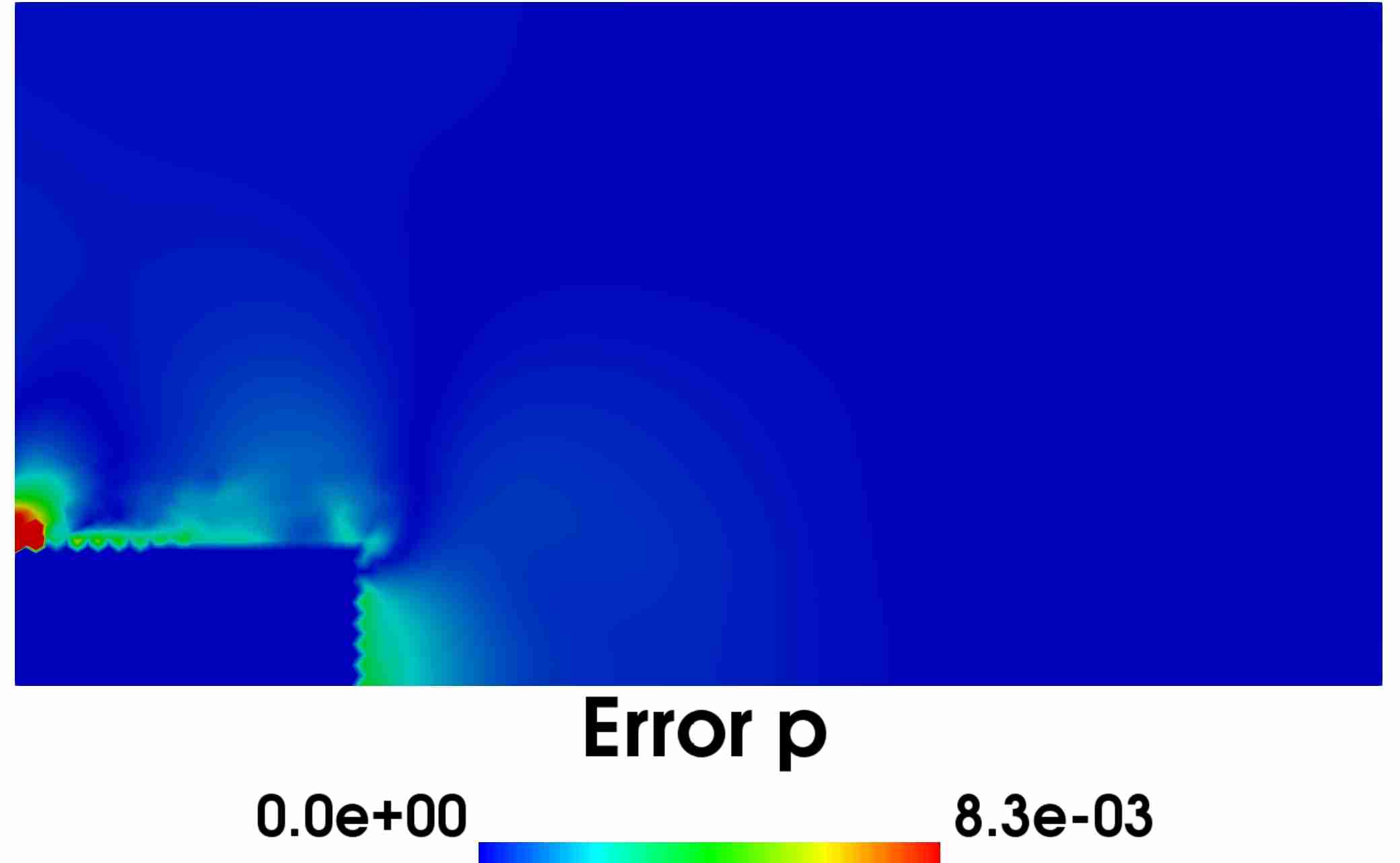}
\end{minipage}
\begin{minipage}{0.24\textwidth}
  \includegraphics[width=\textwidth]{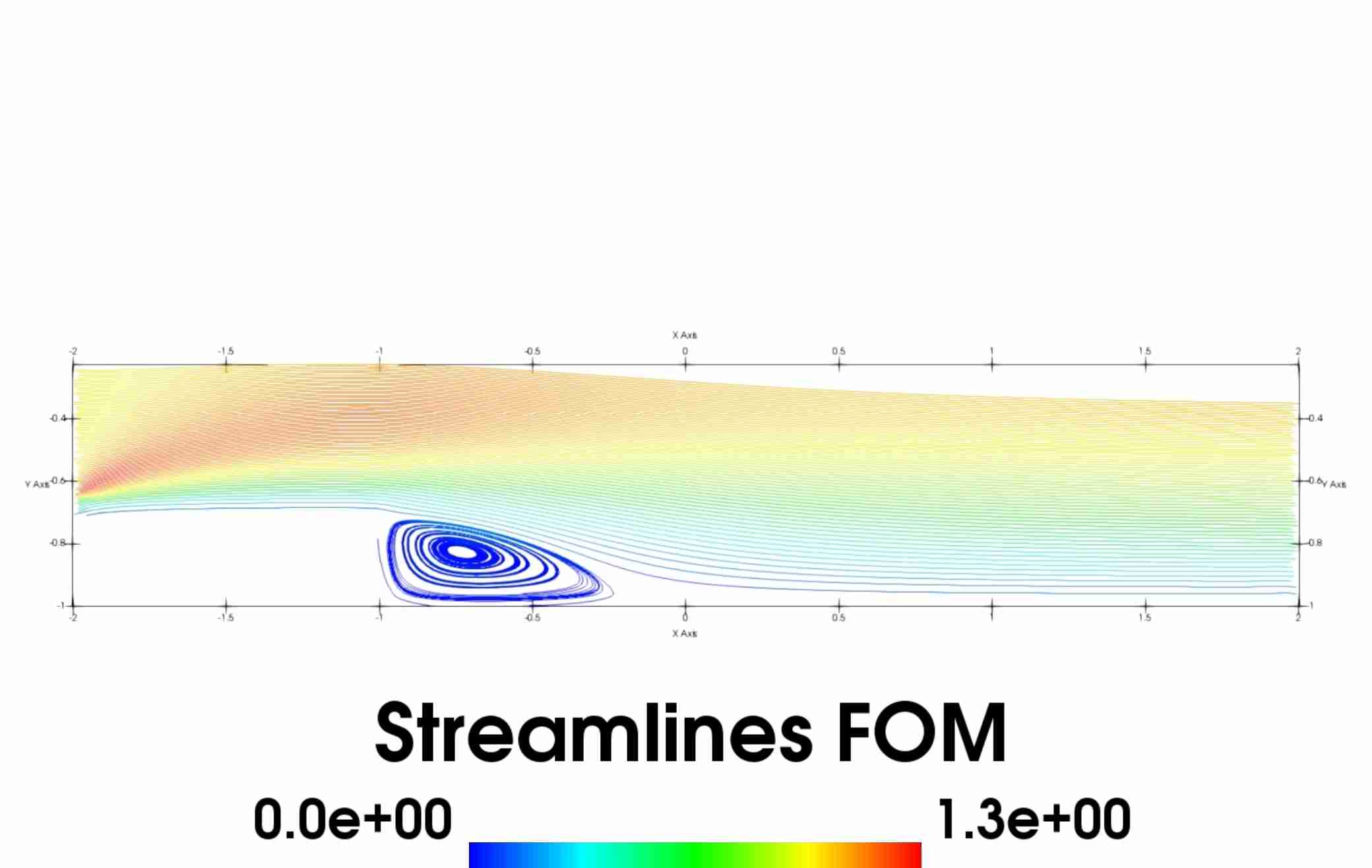} 
\end{minipage}
\begin{minipage}{0.24\textwidth}
  \includegraphics[width=\textwidth]{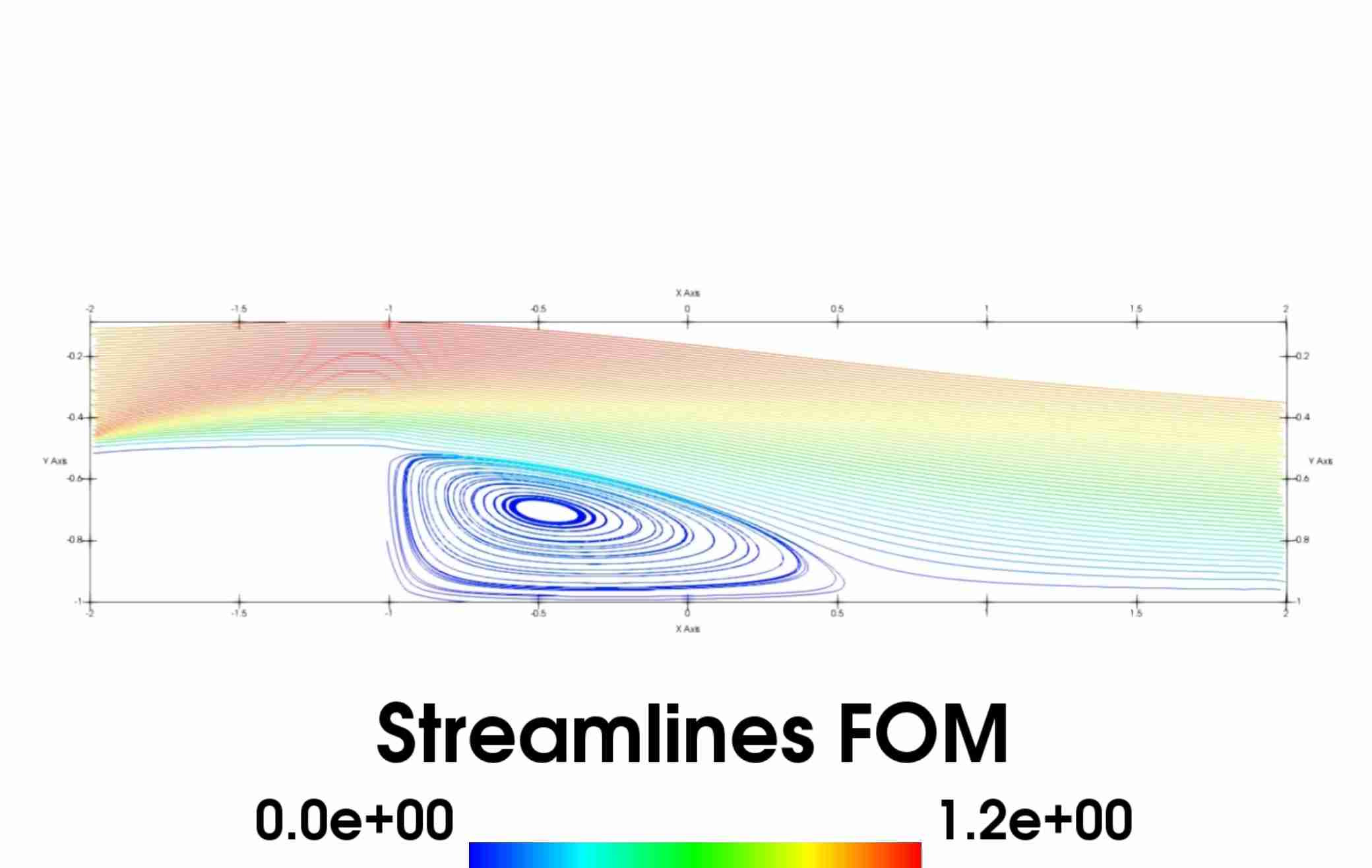}
\end{minipage}
\begin{minipage}{0.24\textwidth}
  \includegraphics[width=\textwidth]{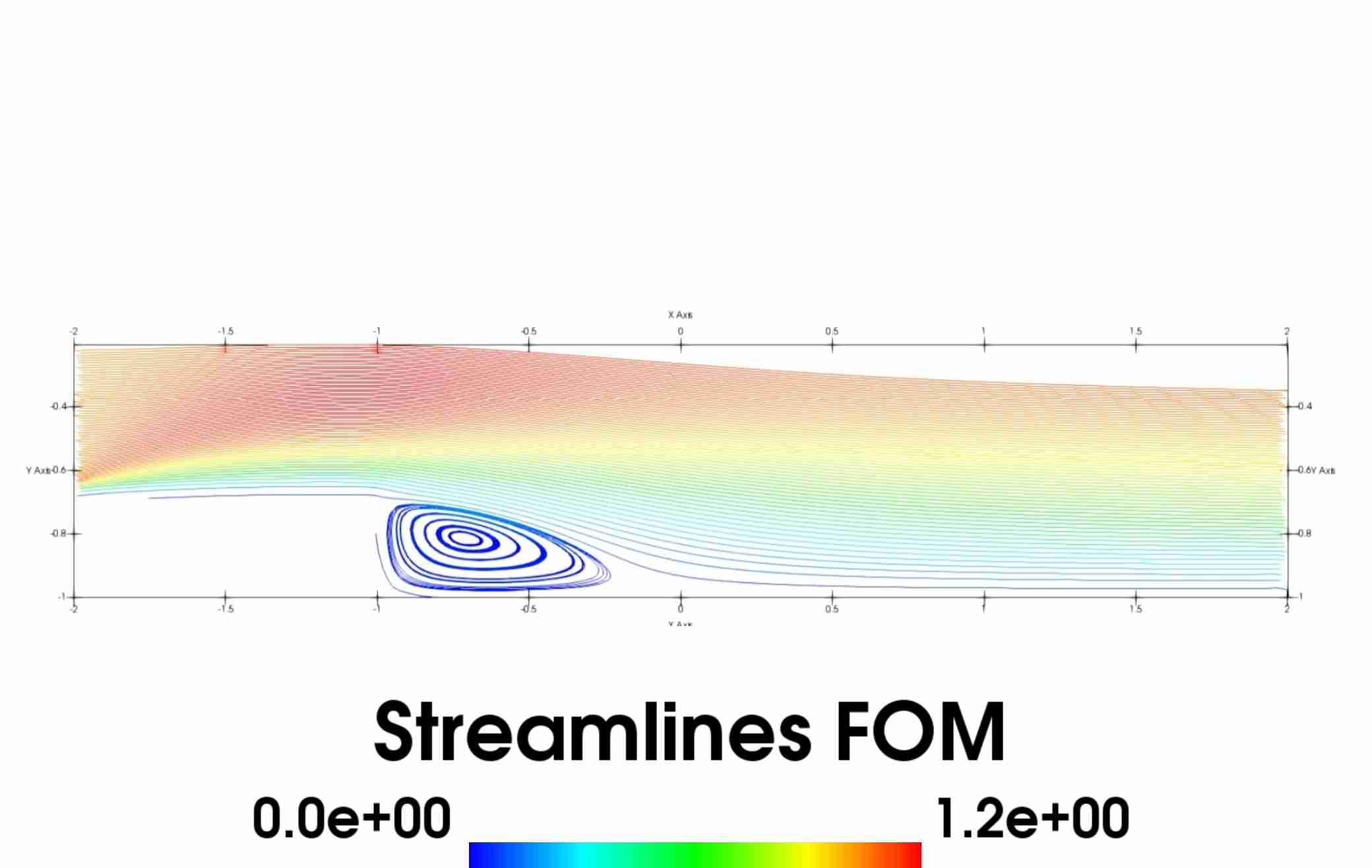}
\end{minipage}
\begin{minipage}{0.24\textwidth}
  \includegraphics[width=\textwidth]{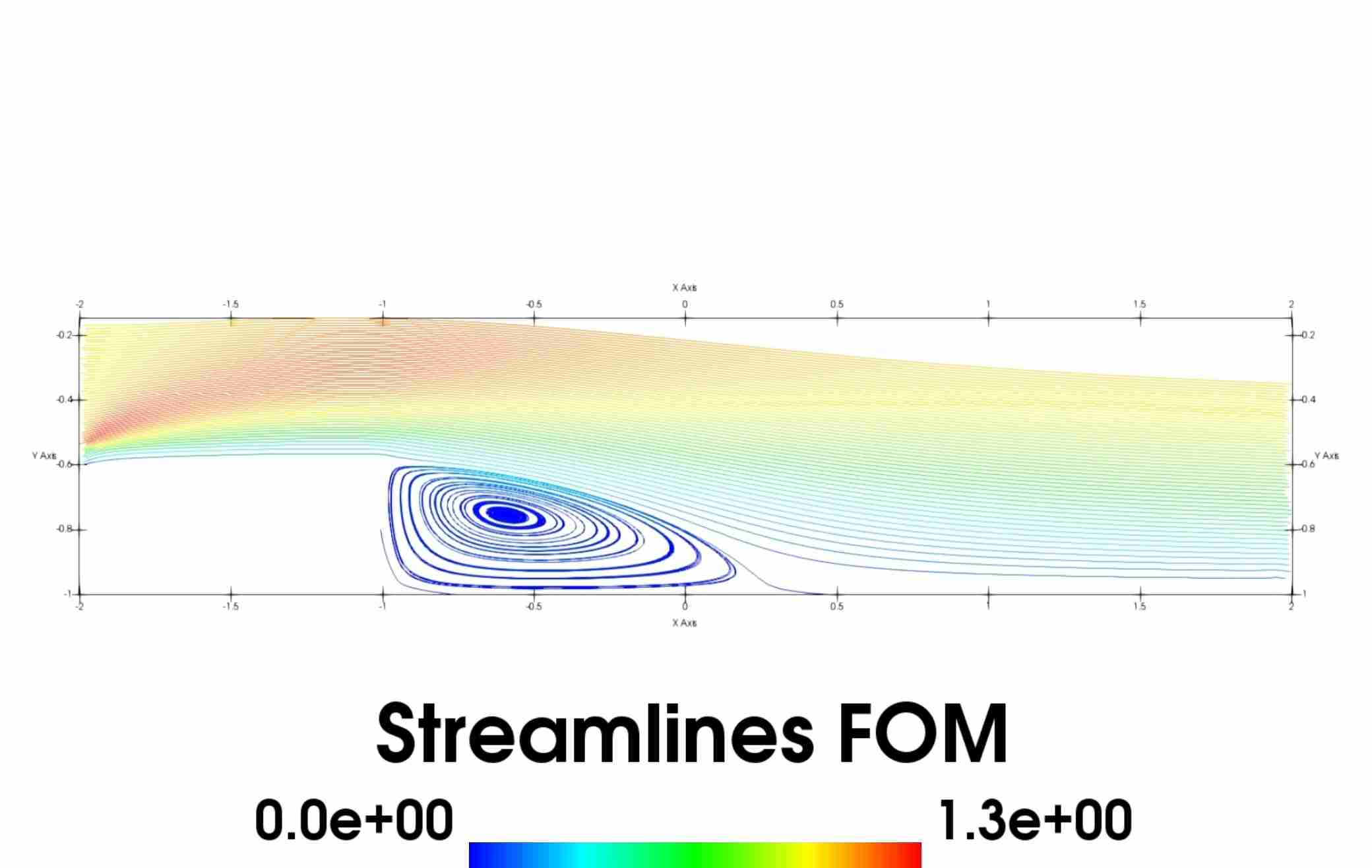}
\end{minipage}
\begin{minipage}{0.24\textwidth}
  \includegraphics[width=\textwidth]{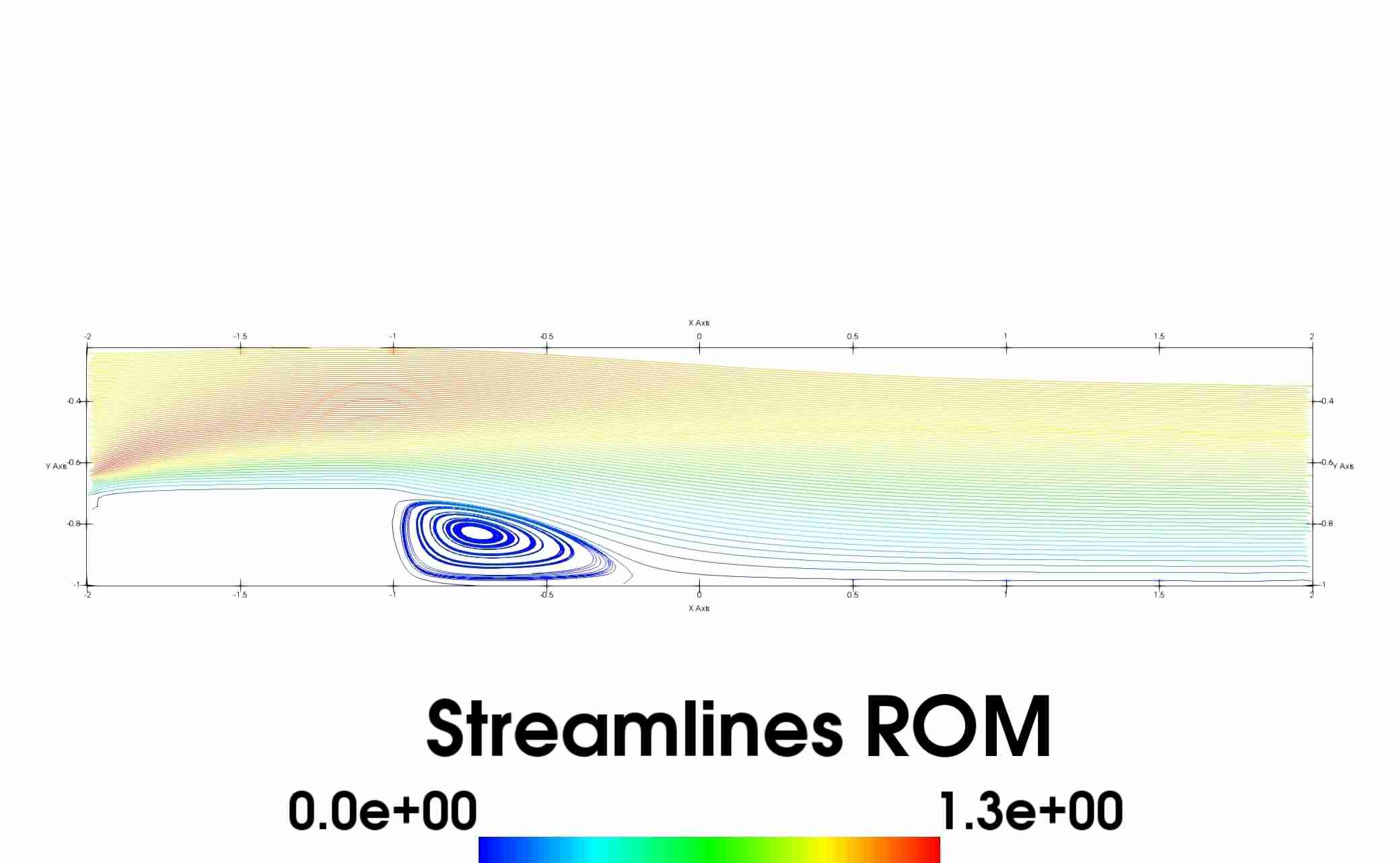} 
\end{minipage}
\begin{minipage}{0.24\textwidth}
  \includegraphics[width=\textwidth]{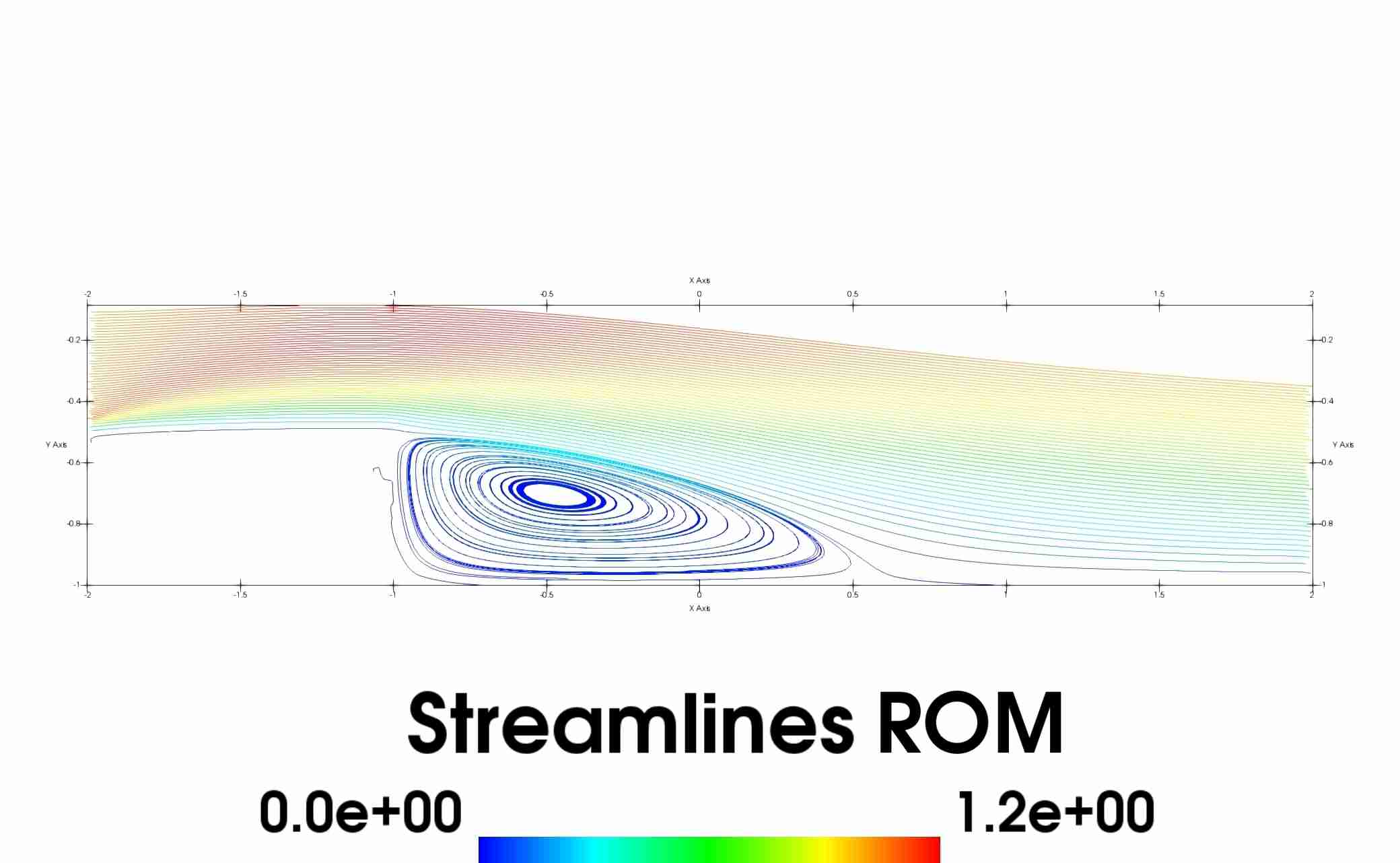}
\end{minipage}
\begin{minipage}{0.24\textwidth}
  \includegraphics[width=\textwidth]{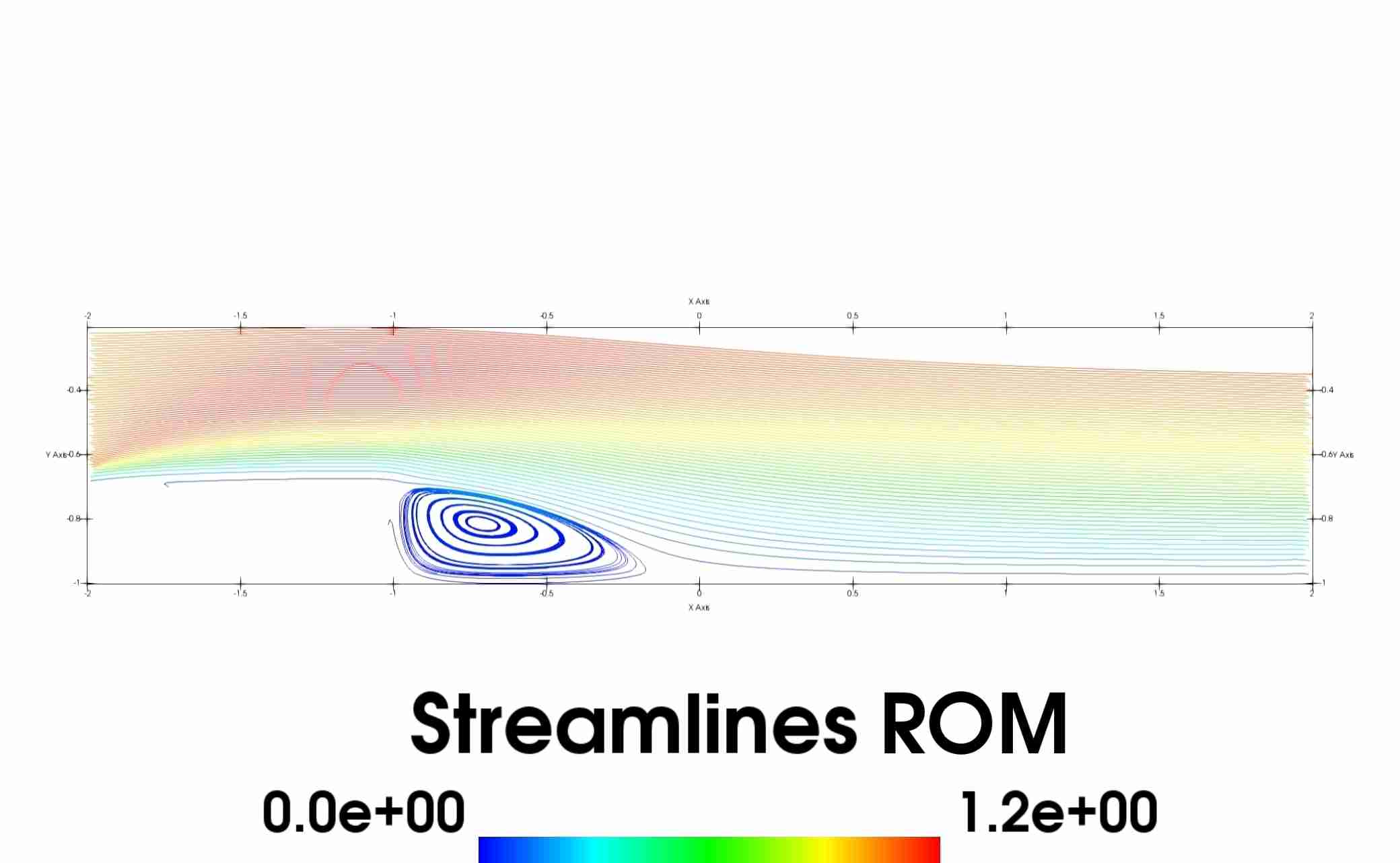}
\end{minipage}
\begin{minipage}{0.24\textwidth}
  \includegraphics[width=\textwidth]{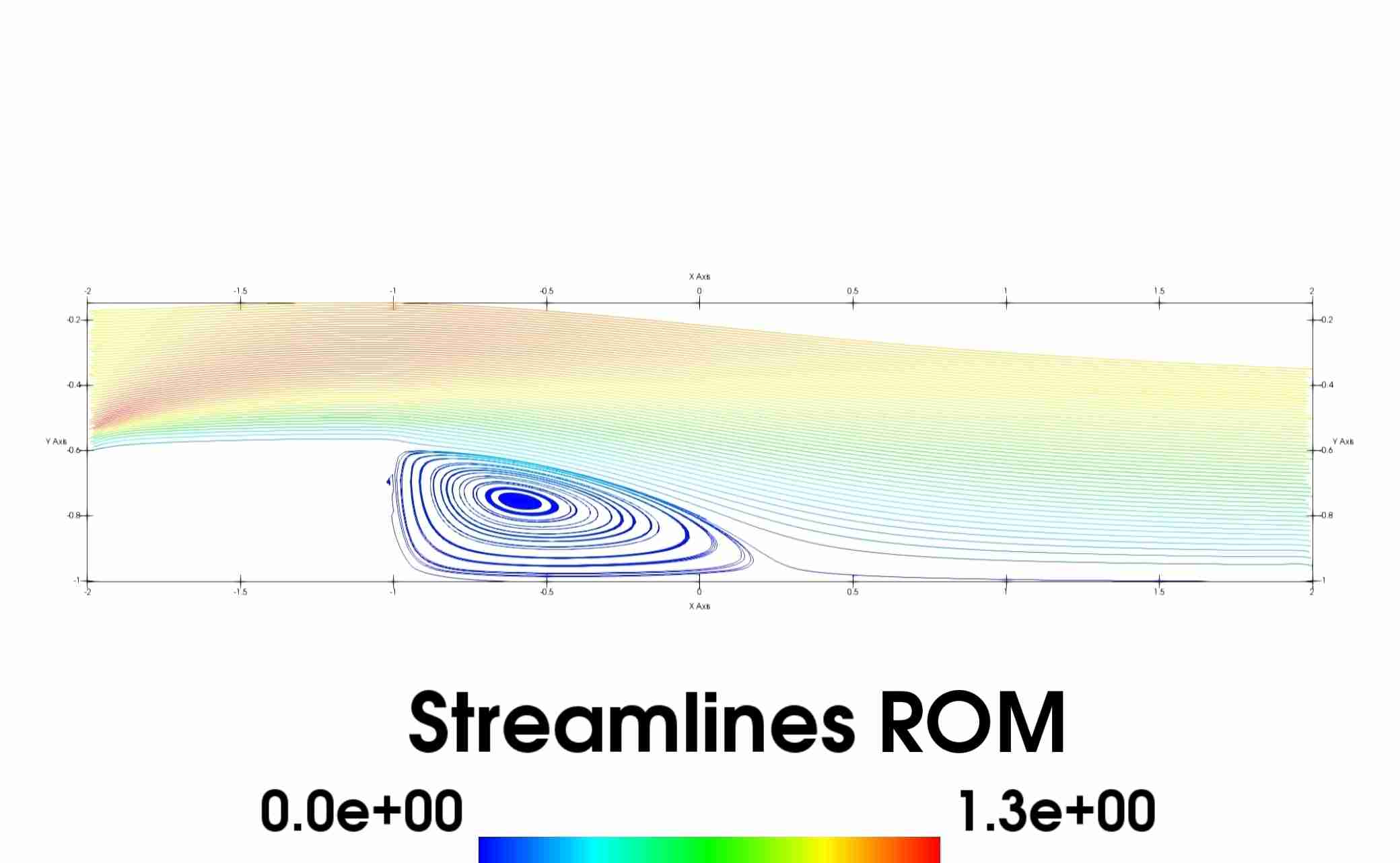}
\end{minipage}
\end{minipage}
\caption{Results for the geometrical parametrization using a one-dimensional parameter space with $\mu_1 \in [0.5, 1.1]$.
In rows $1-3$ we report the full-order solution, the reduced order solution and the absolute error plots for the velocity field, in rows $4-6$ we report the same quantities for the pressure field while in rows $7-8$ the FOM and ROM velocity streamlines. The different columns are for four different values of the input parameter, $\mu_1 = [0.6263,1.0967,1.0169,0.98941569]$.}
\label{fig:FOM_ROM_ERROR_1D}
\end{figure} 
\begin{figure} 
\centering
\begin{minipage}{\textwidth}
\centering
\begin{minipage}{0.24\textwidth}
  \includegraphics[width=\textwidth]{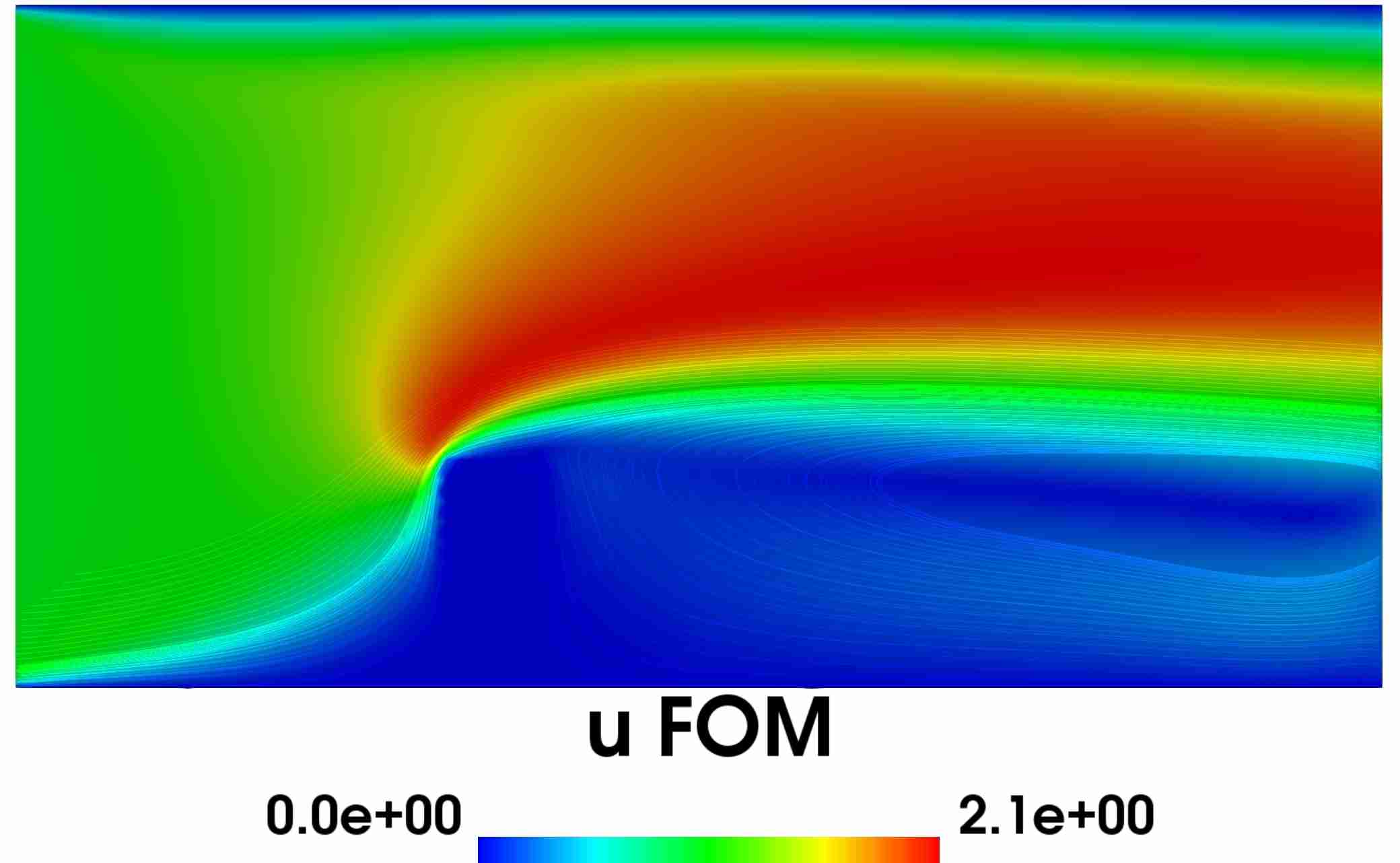} 
\end{minipage}
\begin{minipage}{0.24\textwidth}
  \includegraphics[width=\textwidth]{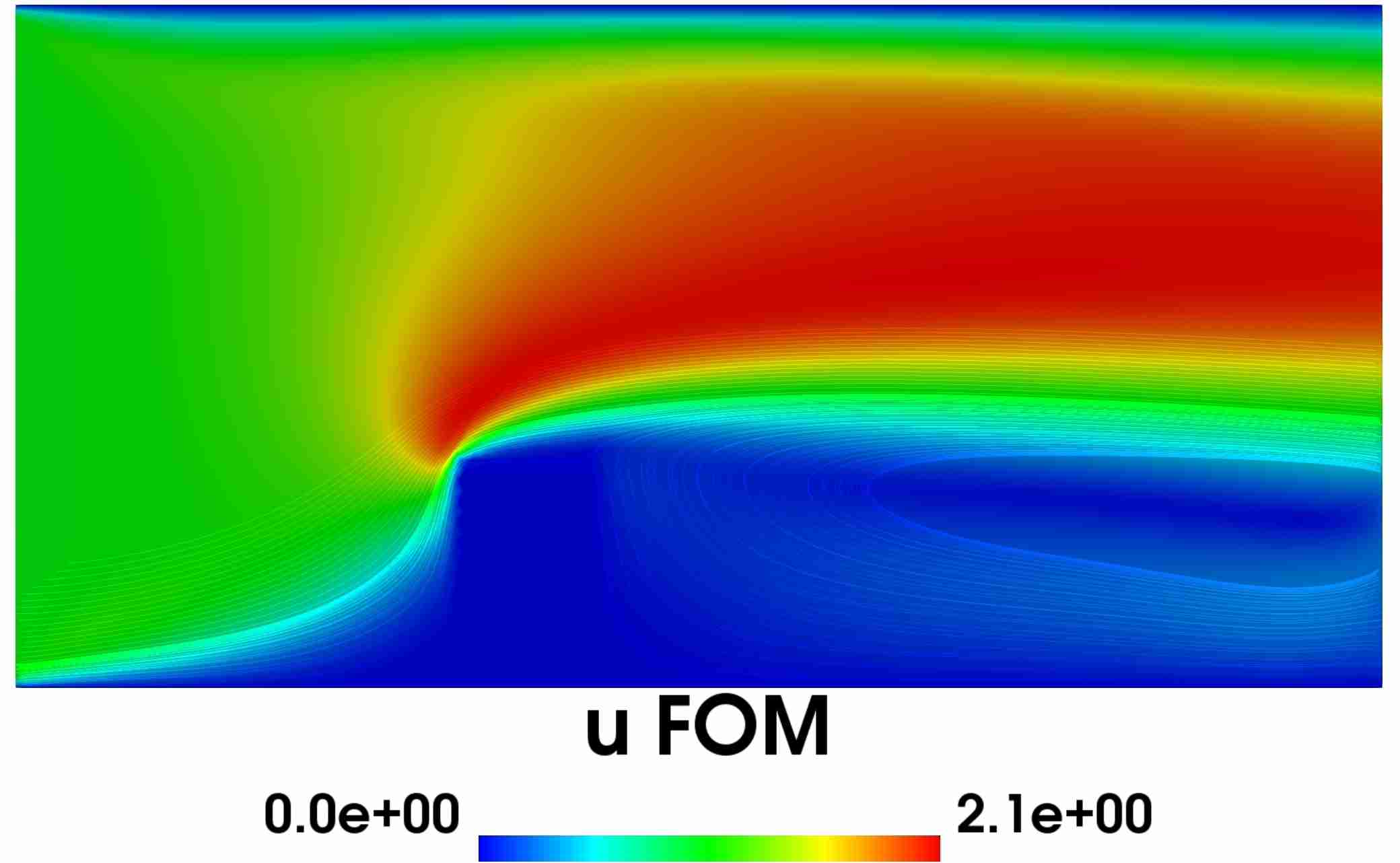}
\end{minipage}
\begin{minipage}{0.24\textwidth}
  \includegraphics[width=\textwidth]{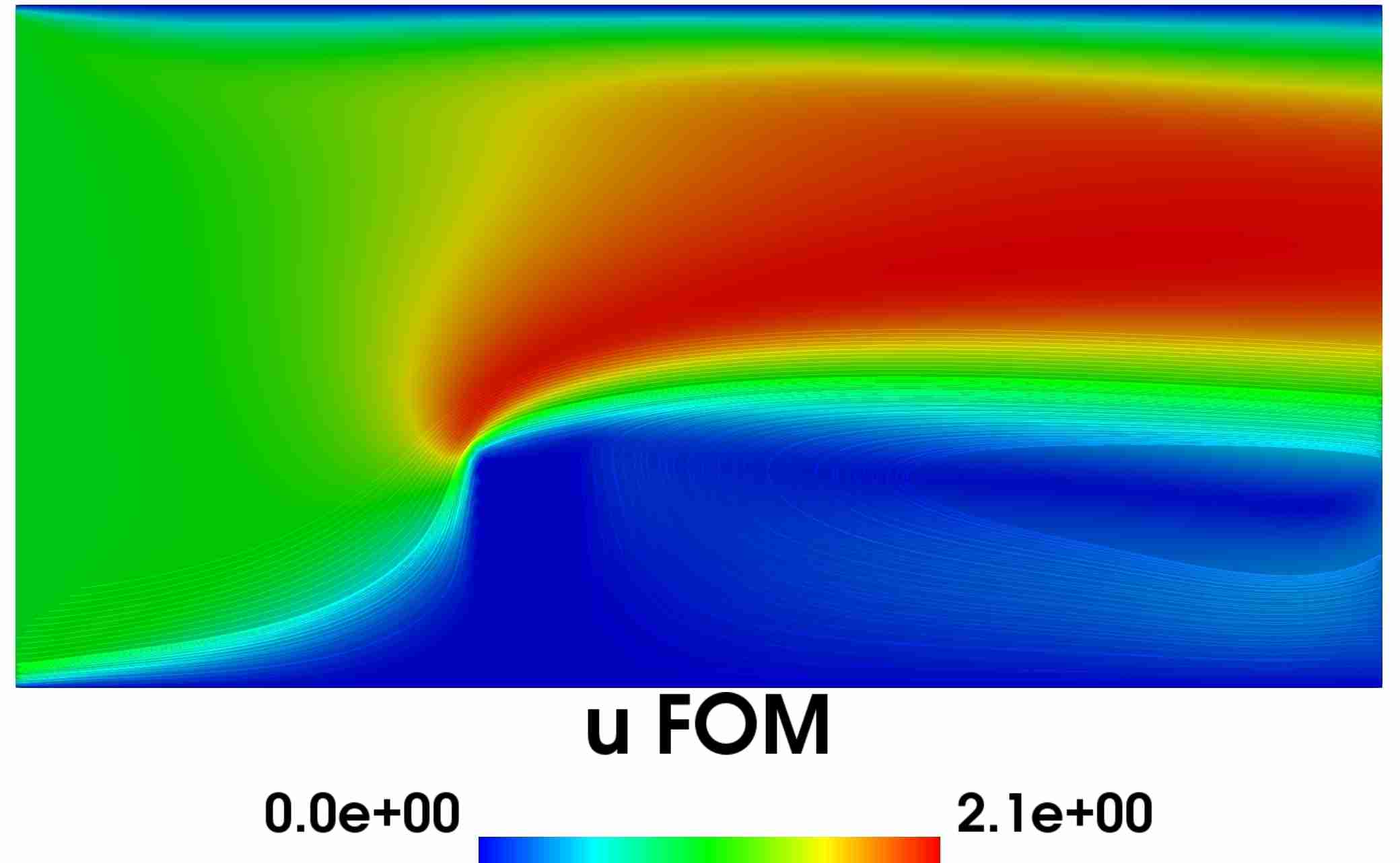}
\end{minipage}
\begin{minipage}{0.24\textwidth}
  \includegraphics[width=\textwidth]{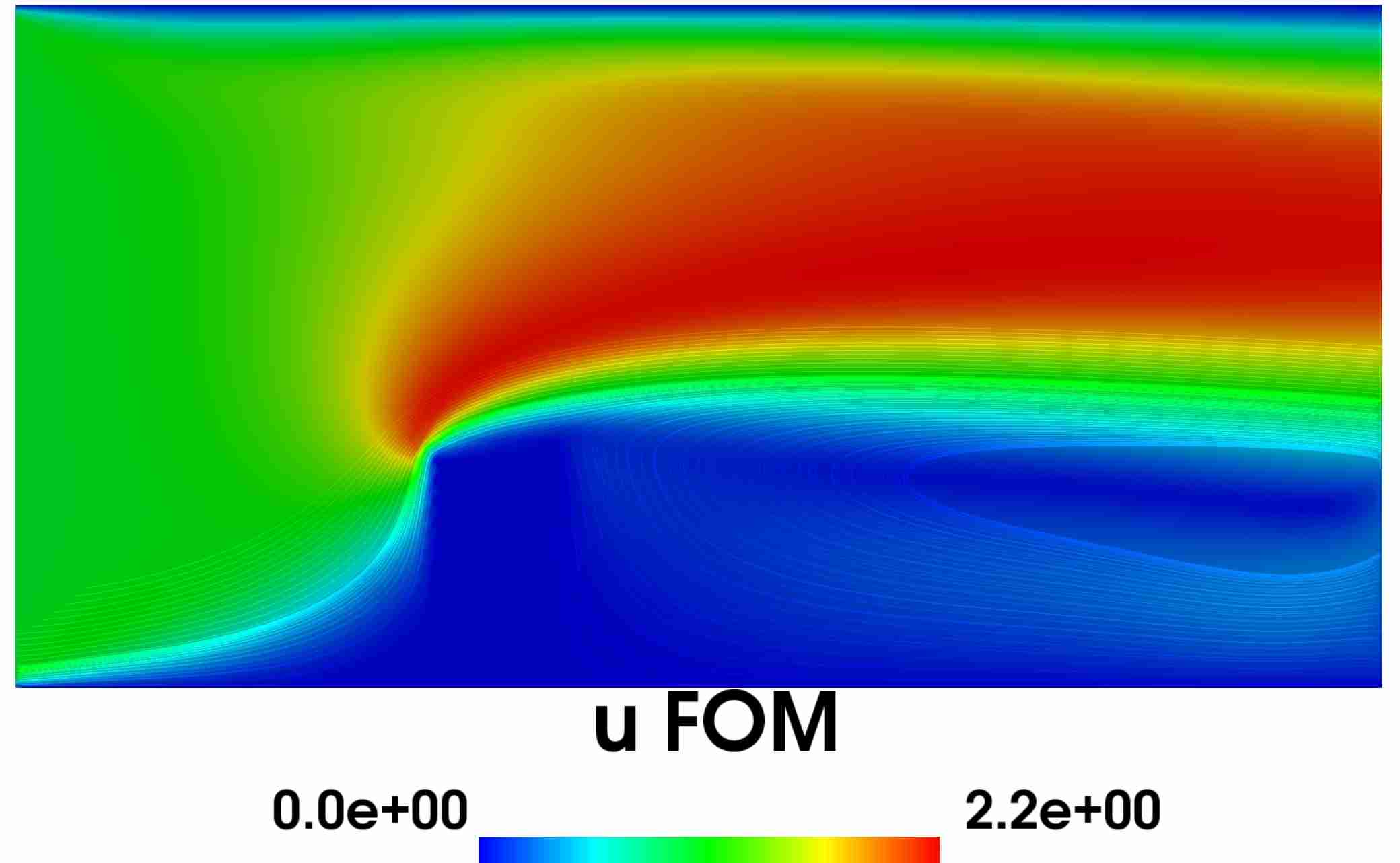}
\end{minipage}
\begin{minipage}{0.24\textwidth}
  \includegraphics[width=\textwidth]{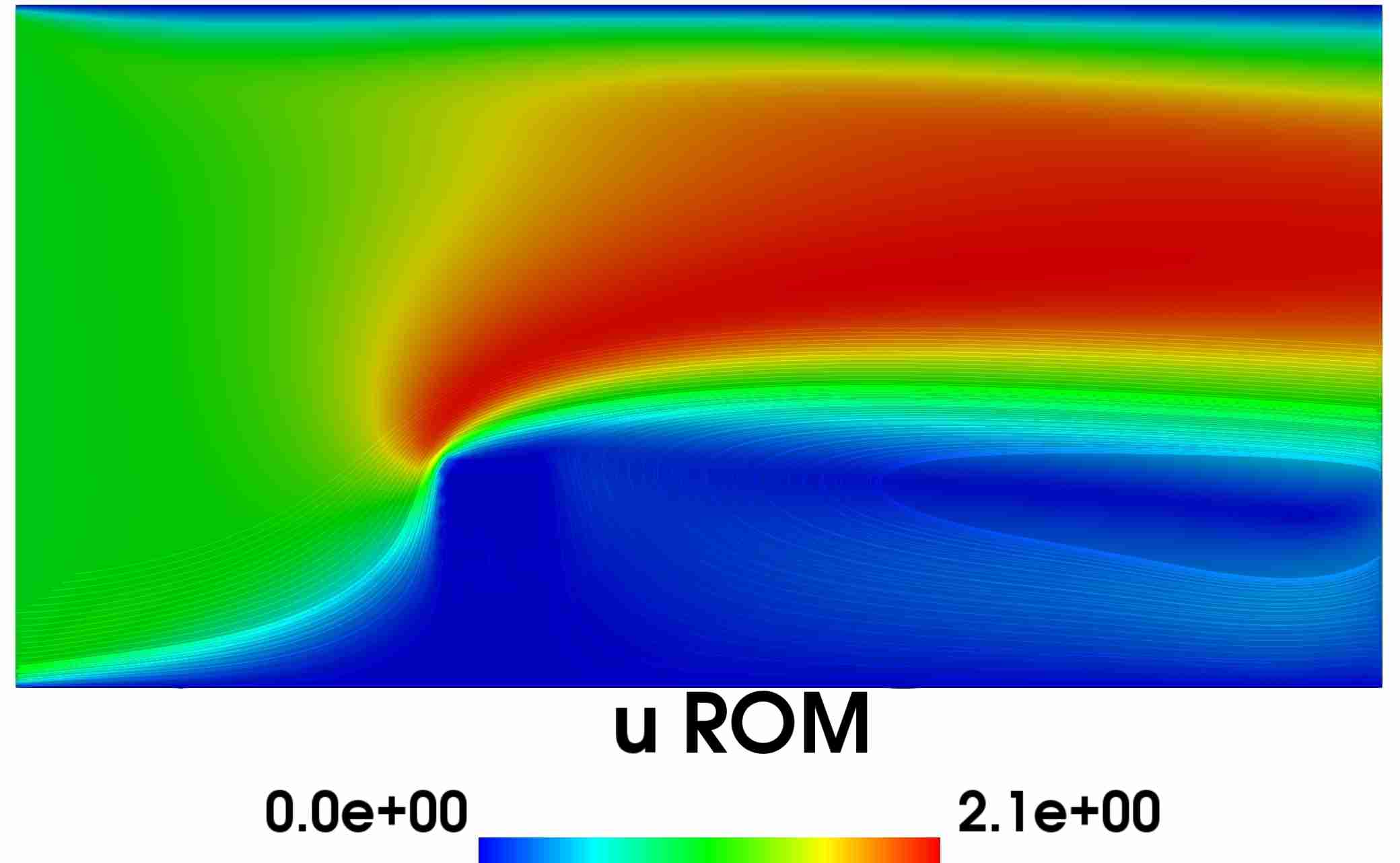} 
\end{minipage}
\begin{minipage}{0.24\textwidth}
  \includegraphics[width=\textwidth]{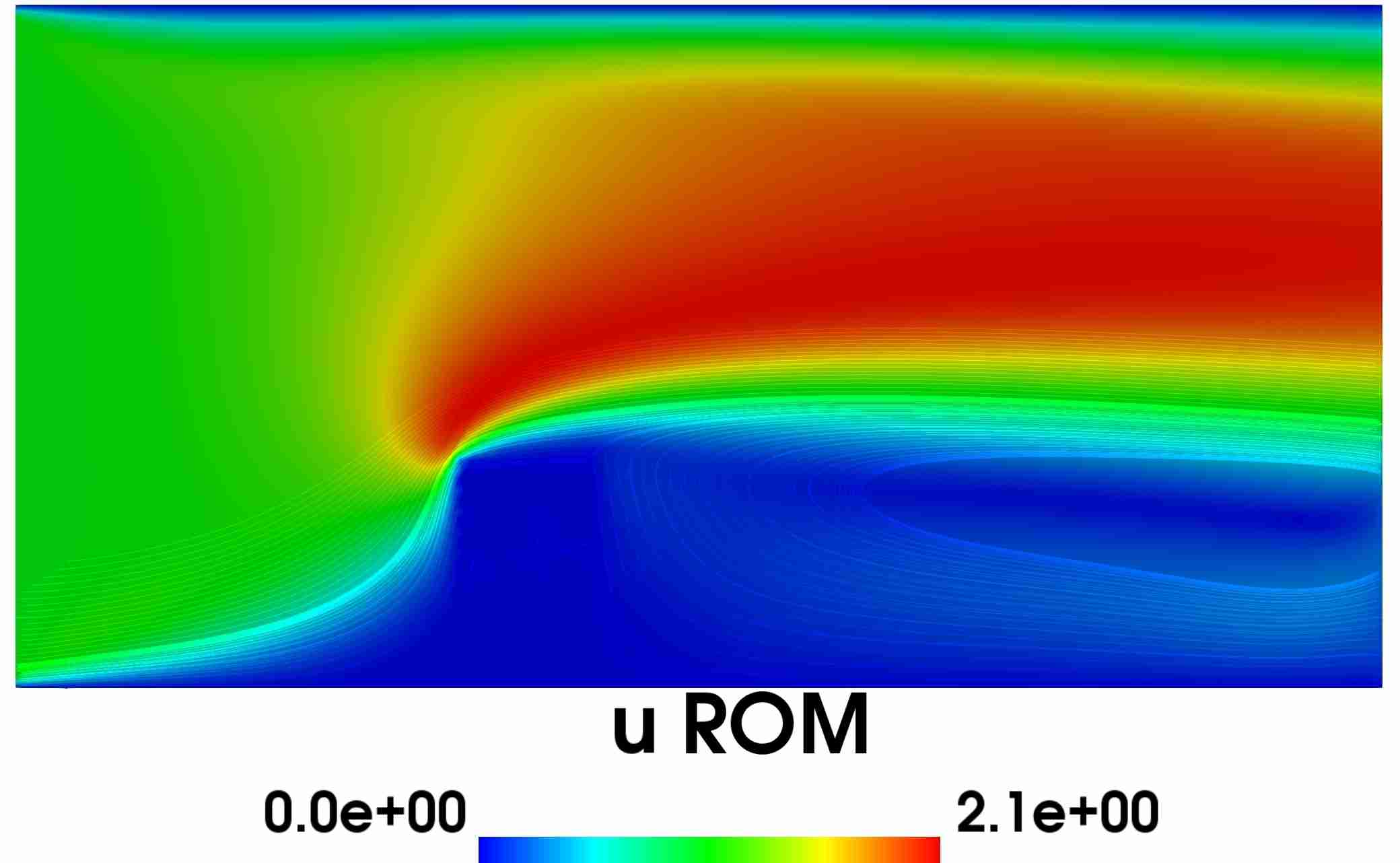}
\end{minipage}
\begin{minipage}{0.24\textwidth}
  \includegraphics[width=\textwidth]{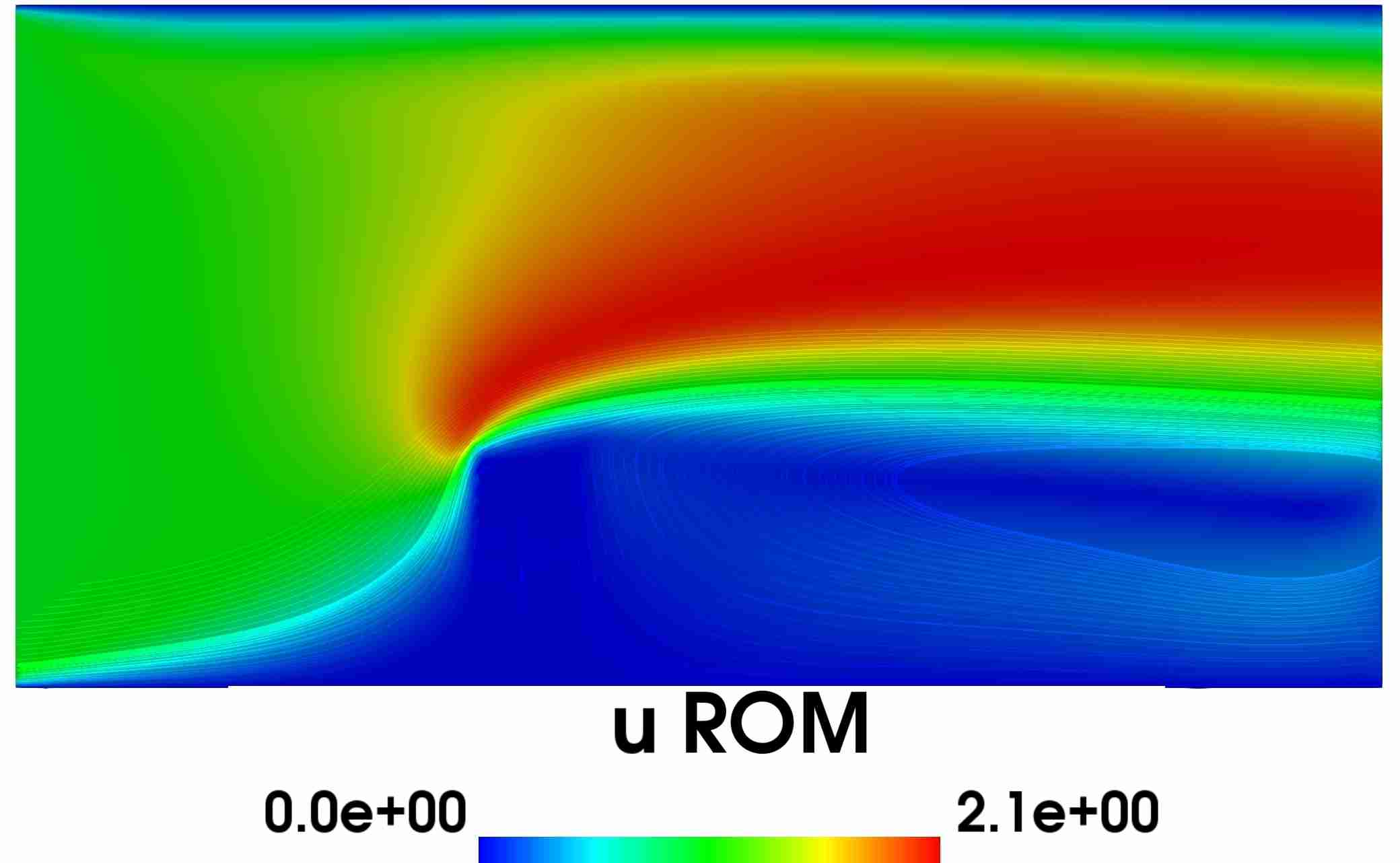}
\end{minipage}
\begin{minipage}{0.24\textwidth}
  \includegraphics[width=\textwidth]{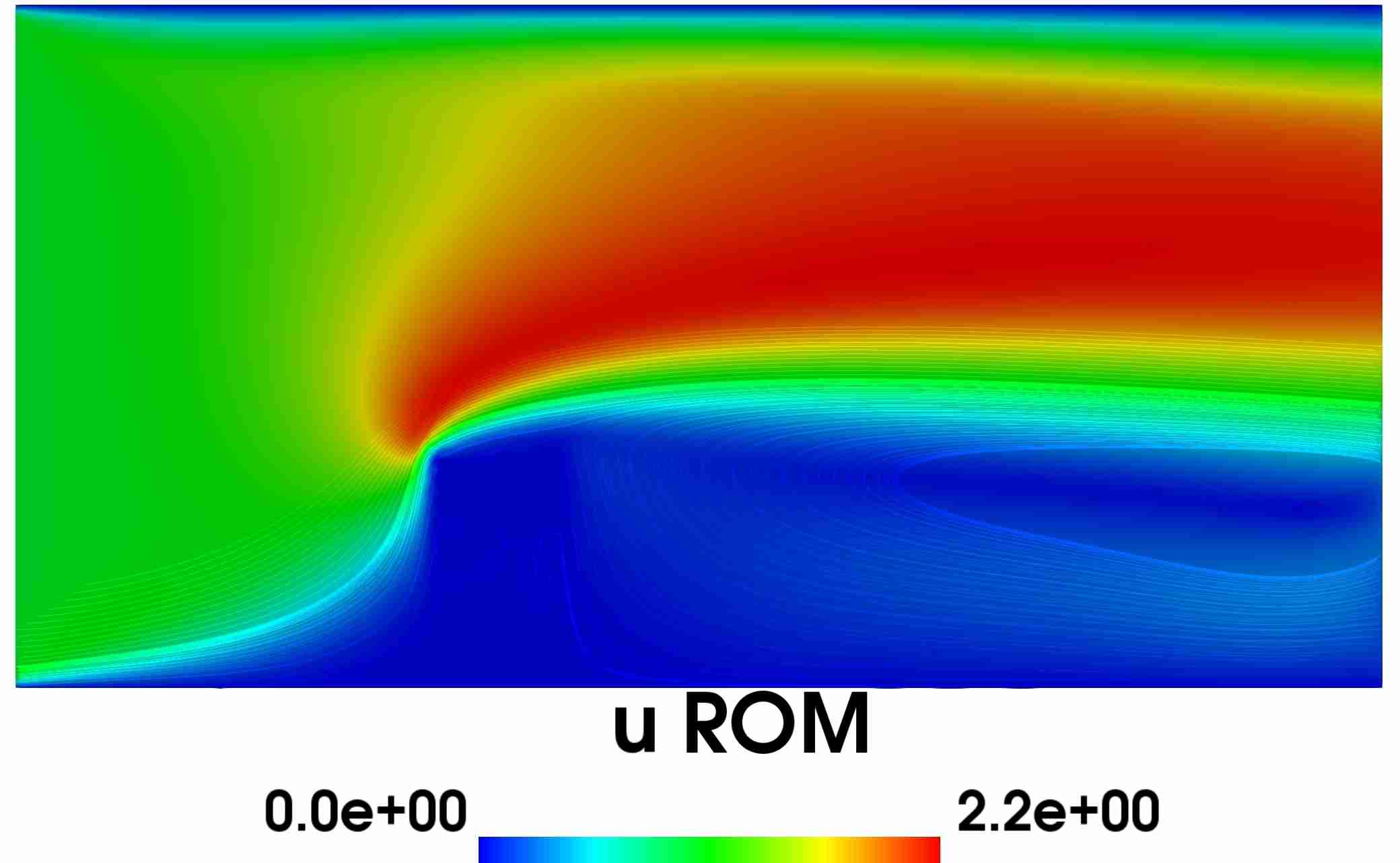}
\end{minipage}
\begin{minipage}{0.24\textwidth}
  \includegraphics[width=\textwidth]{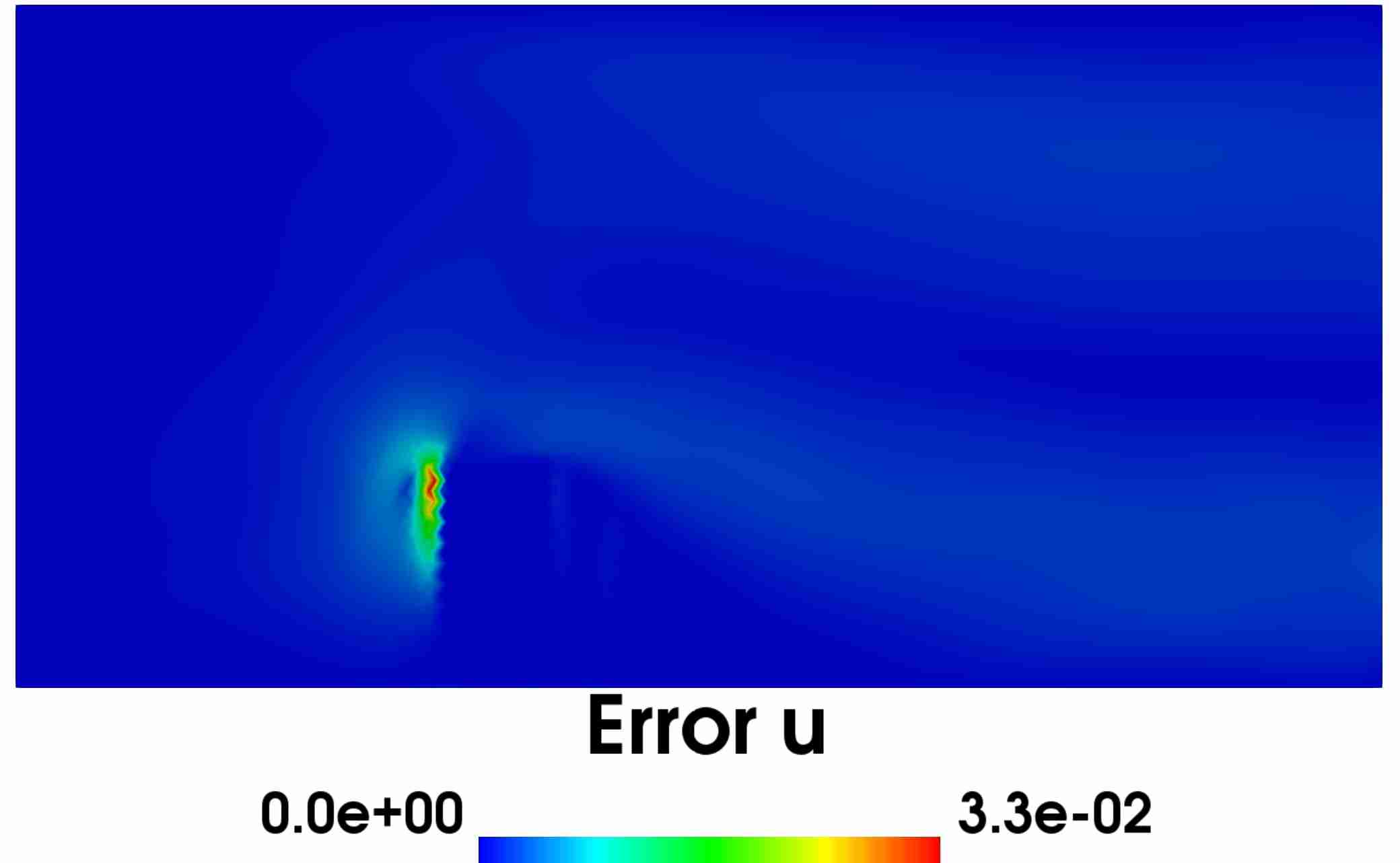} 
\end{minipage}
\begin{minipage}{0.24\textwidth}
  \includegraphics[width=\textwidth]{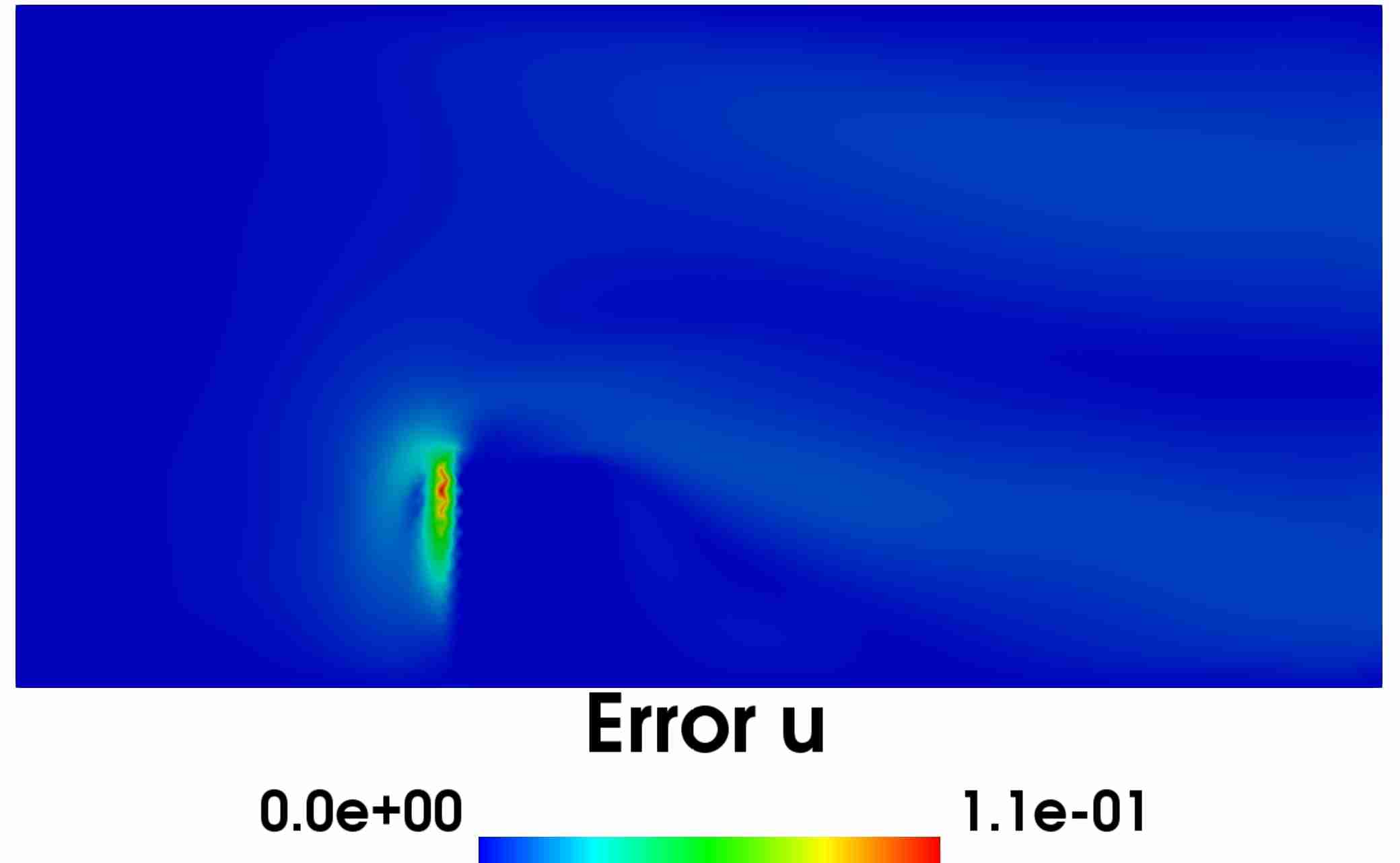}
\end{minipage}
\begin{minipage}{0.24\textwidth}
  \includegraphics[width=\textwidth]{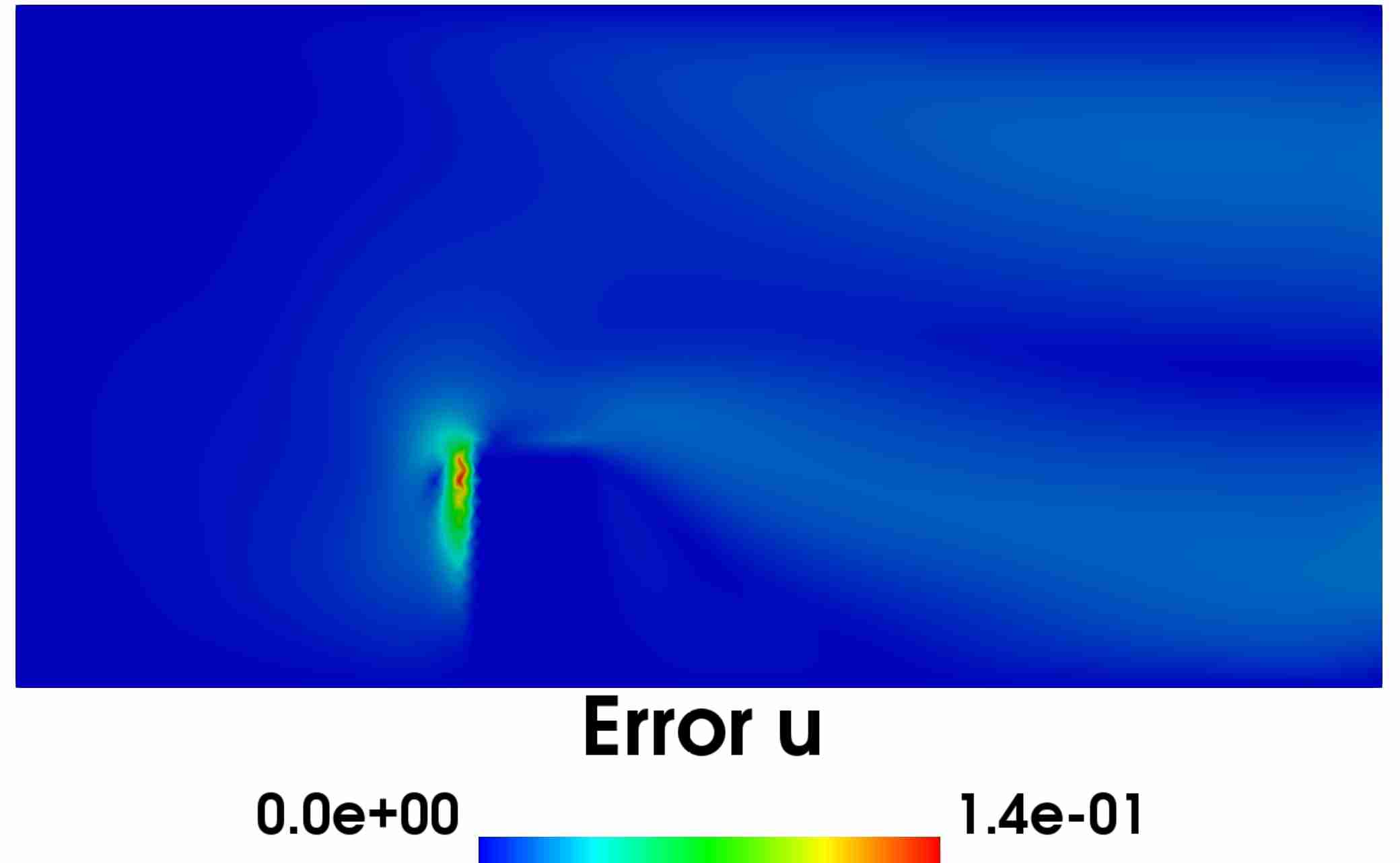}
\end{minipage}
\begin{minipage}{0.24\textwidth}
  \includegraphics[width=\textwidth]{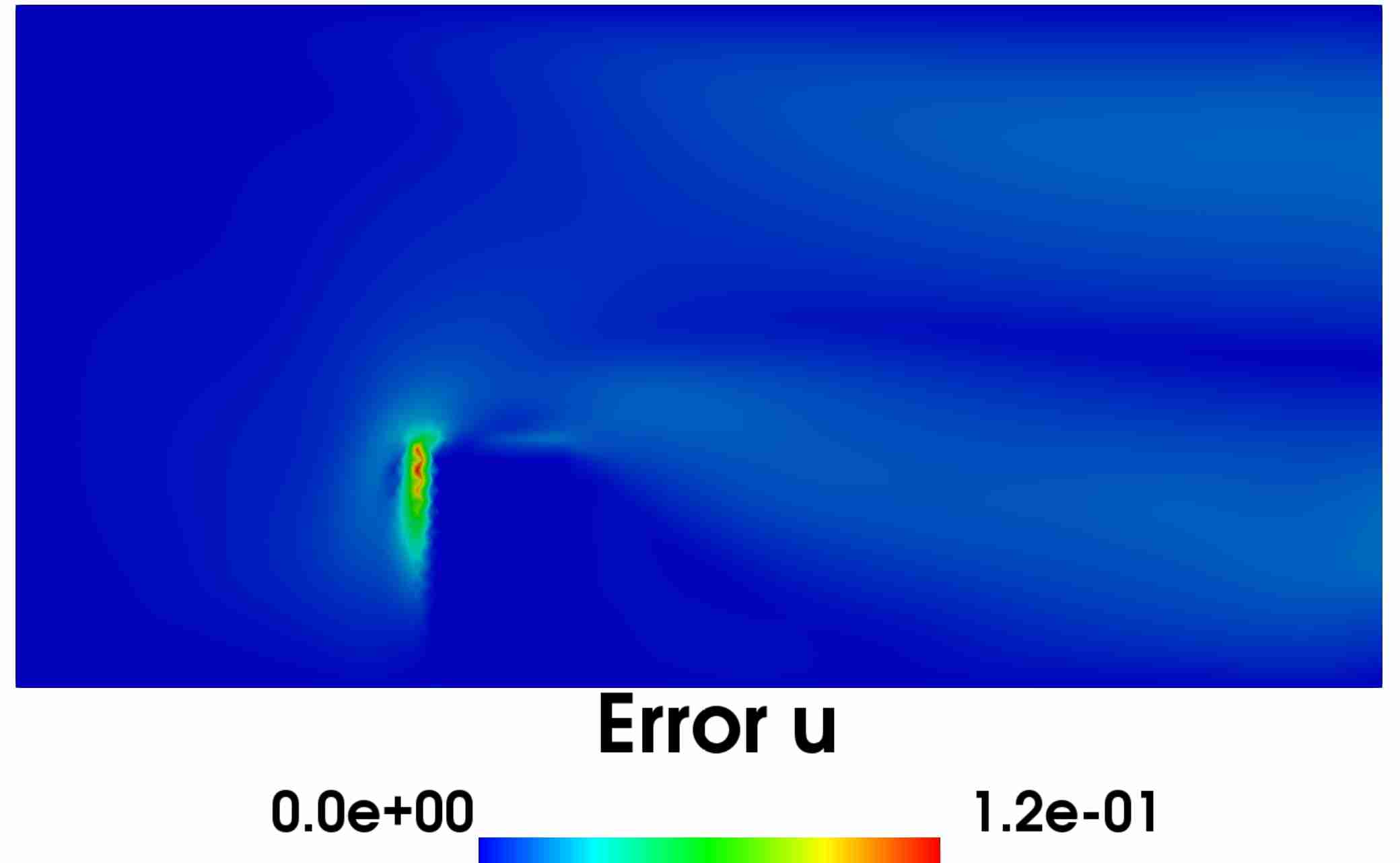}
\end{minipage}
\begin{minipage}{0.24\textwidth}
  \includegraphics[width=\textwidth]{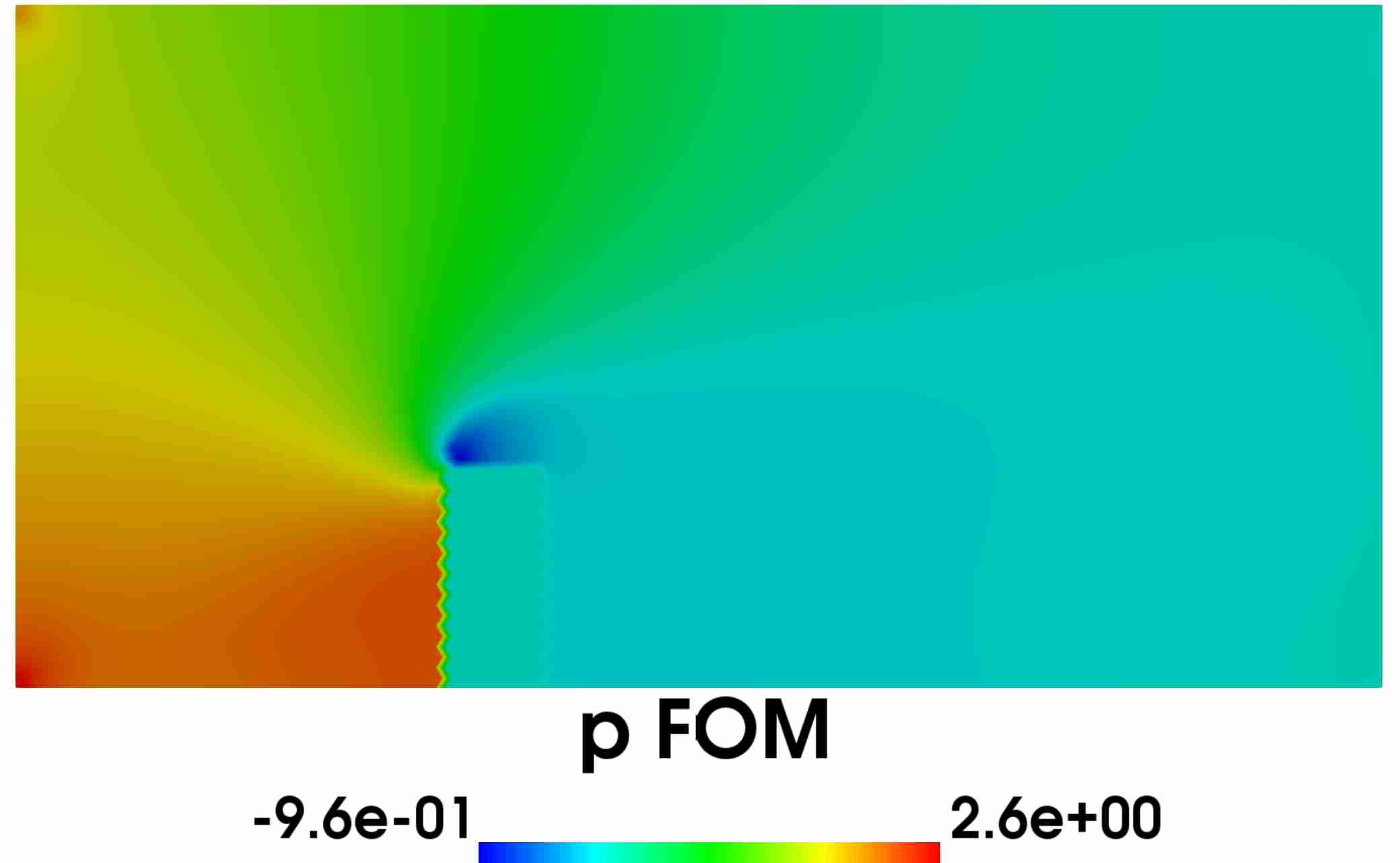} 
\end{minipage}
\begin{minipage}{0.24\textwidth}
  \includegraphics[width=\textwidth]{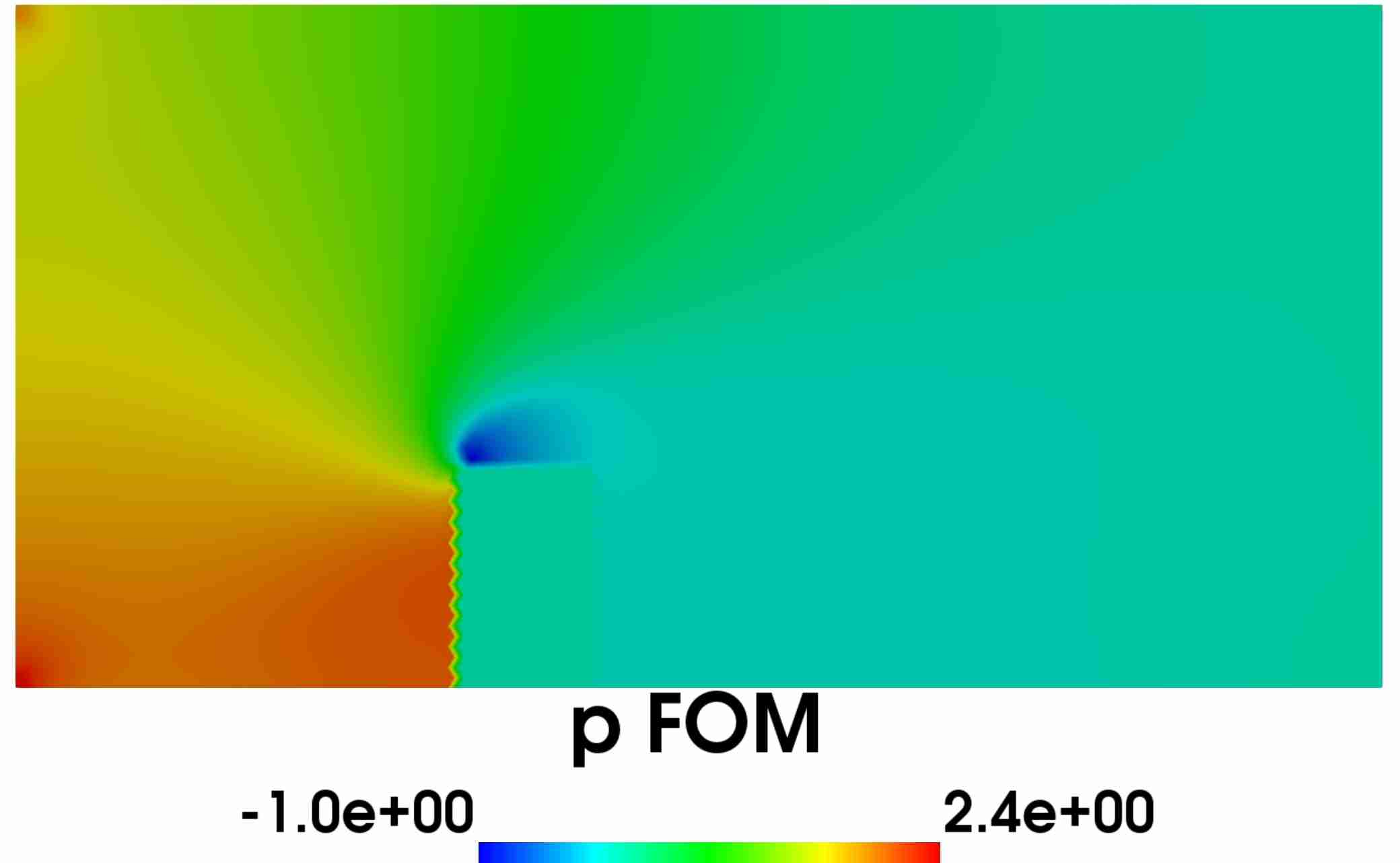}
\end{minipage}
\begin{minipage}{0.24\textwidth}
  \includegraphics[width=\textwidth]{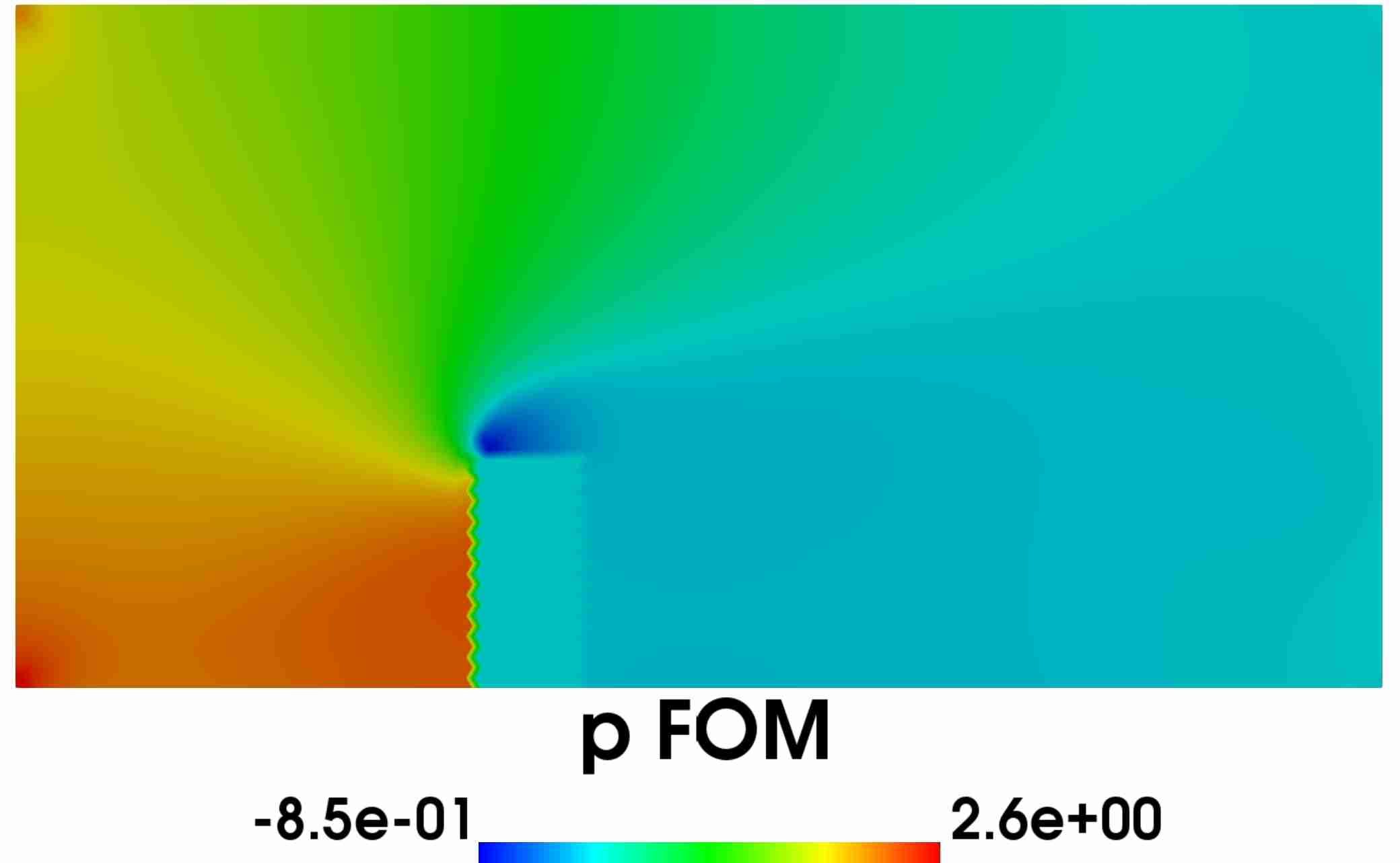}
\end{minipage}
\begin{minipage}{0.24\textwidth}
  \includegraphics[width=\textwidth]{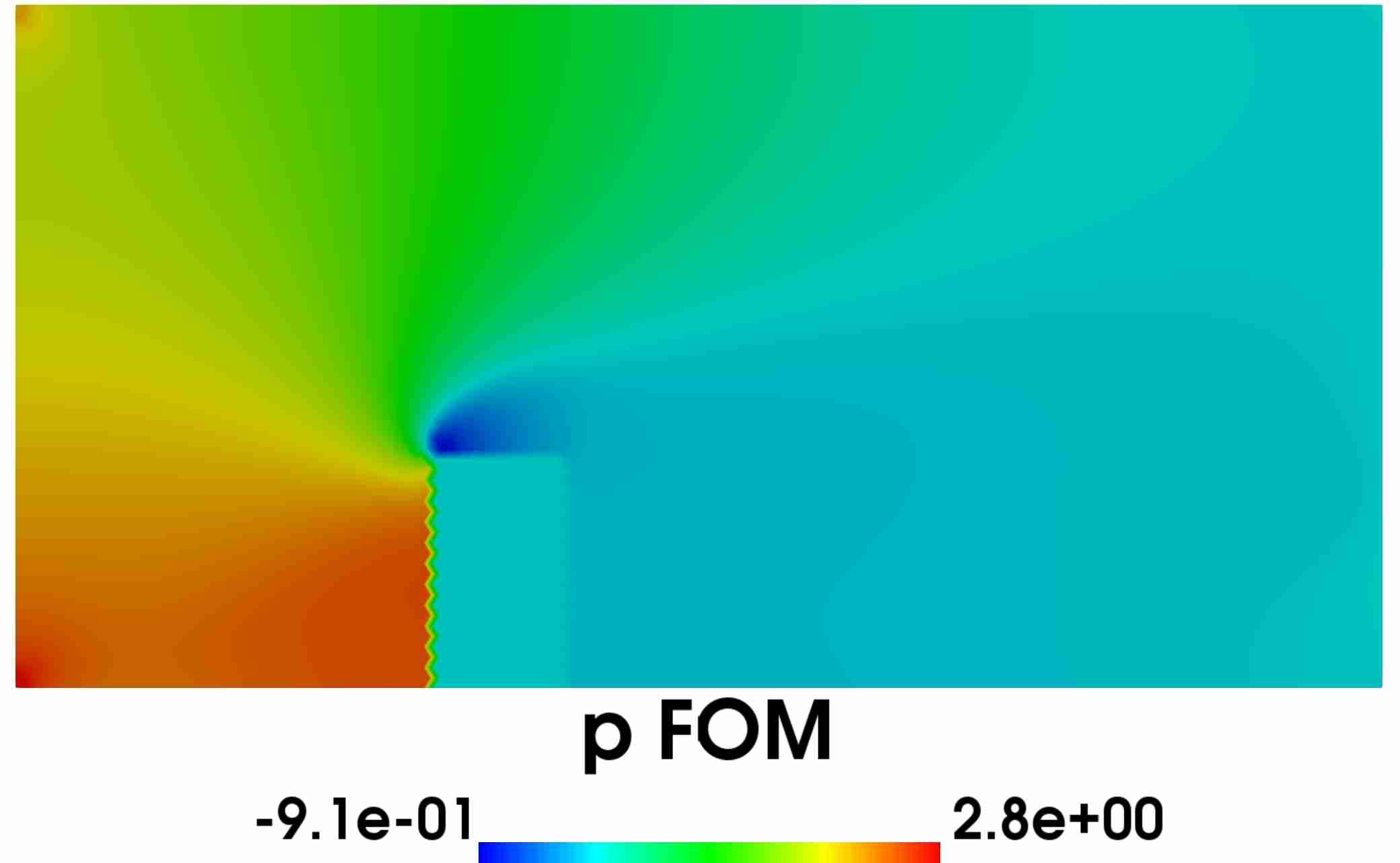}
\end{minipage}
\begin{minipage}{0.24\textwidth}
  \includegraphics[width=\textwidth]{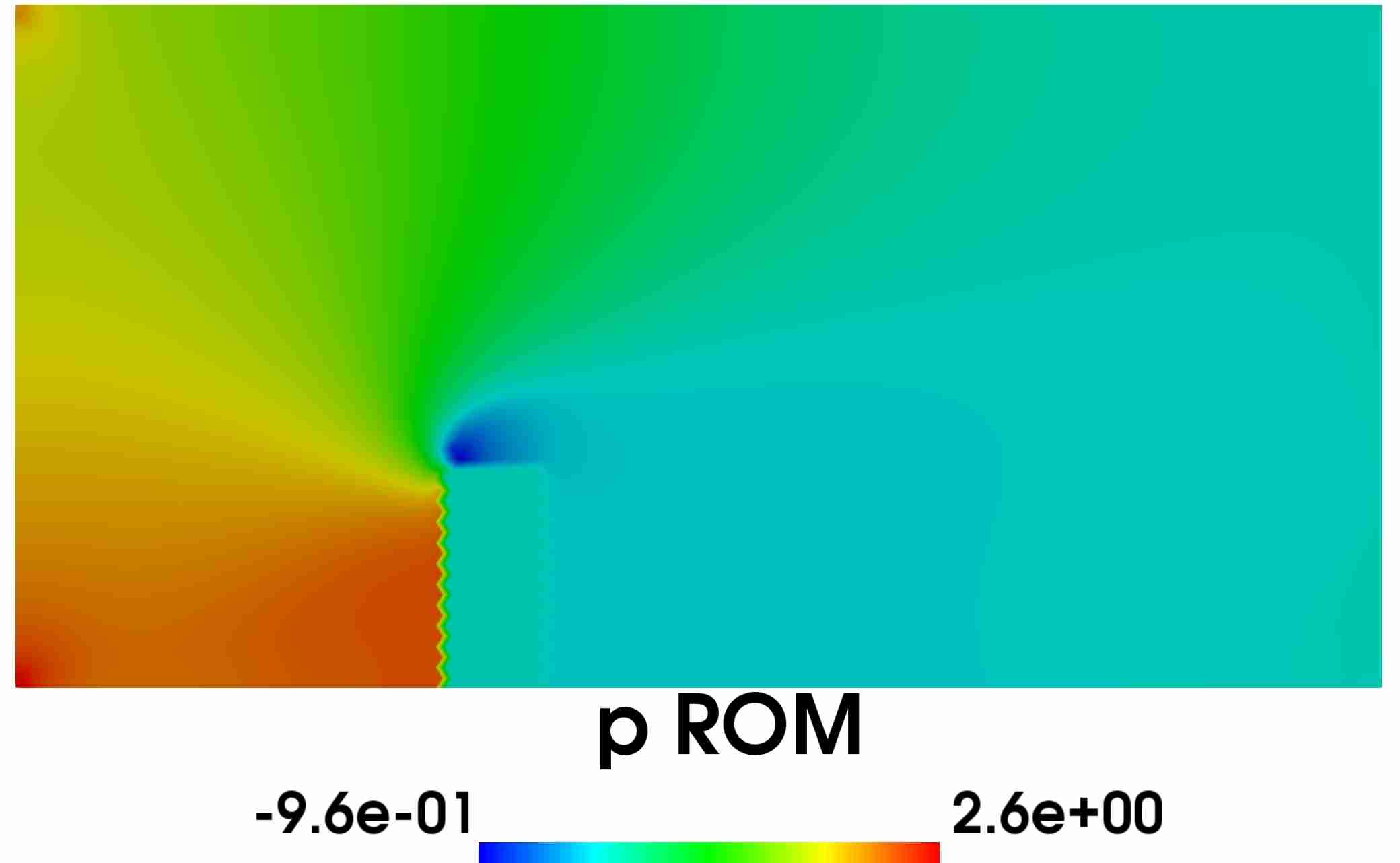} 
\end{minipage}
\begin{minipage}{0.24\textwidth}
  \includegraphics[width=\textwidth]{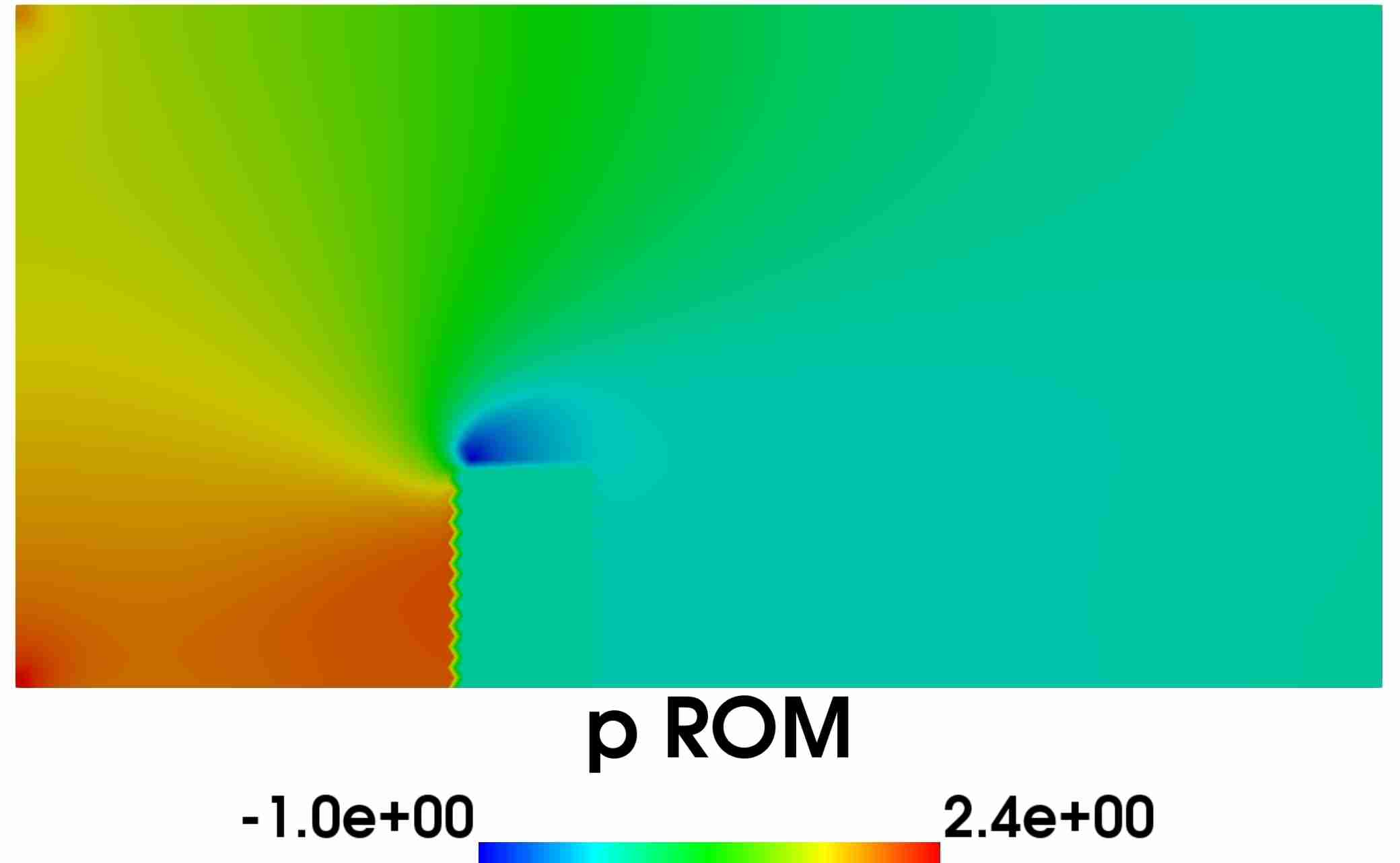}
\end{minipage}
\begin{minipage}{0.24\textwidth}
  \includegraphics[width=\textwidth]{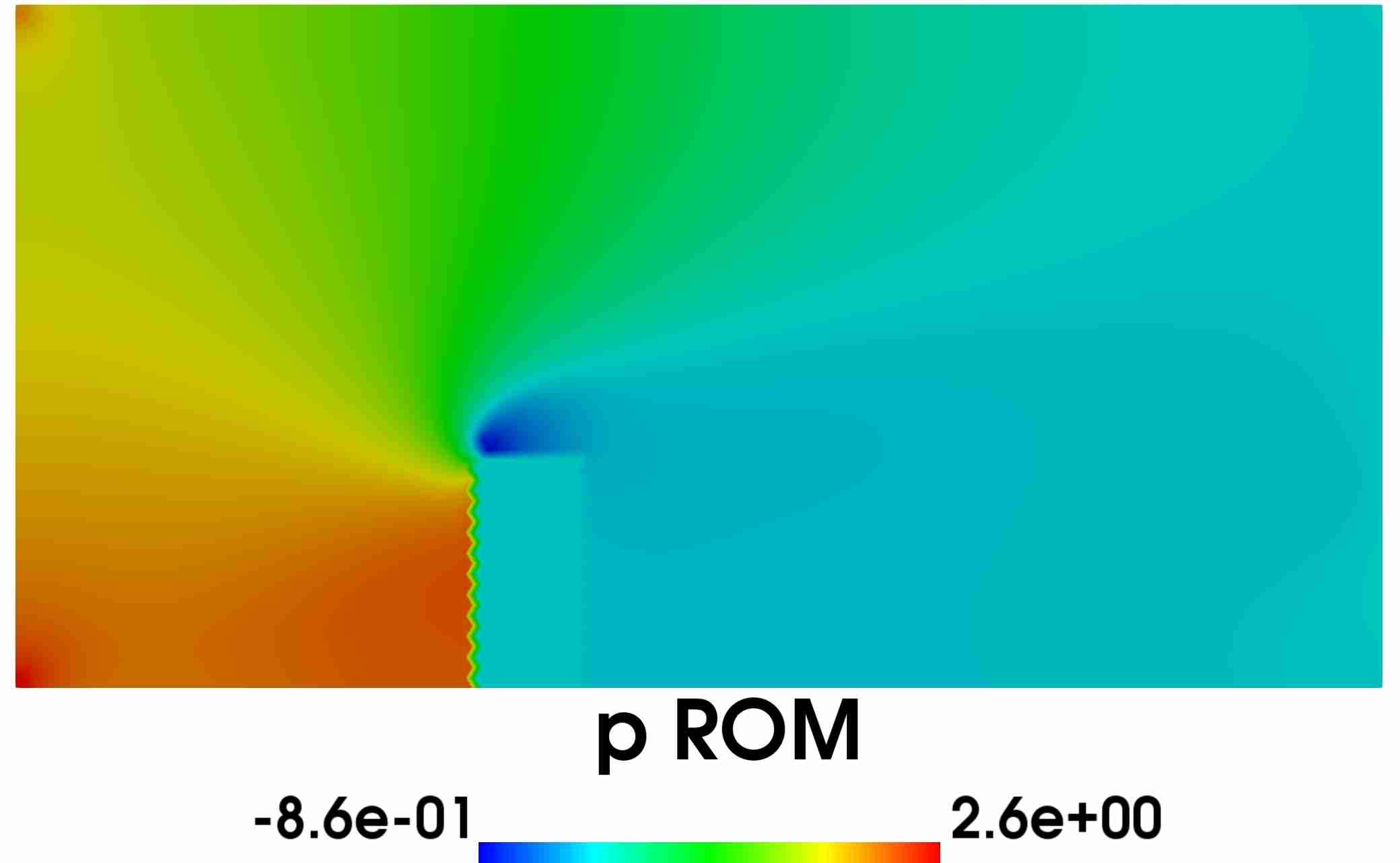}
\end{minipage}
\begin{minipage}{0.24\textwidth}
  \includegraphics[width=\textwidth]{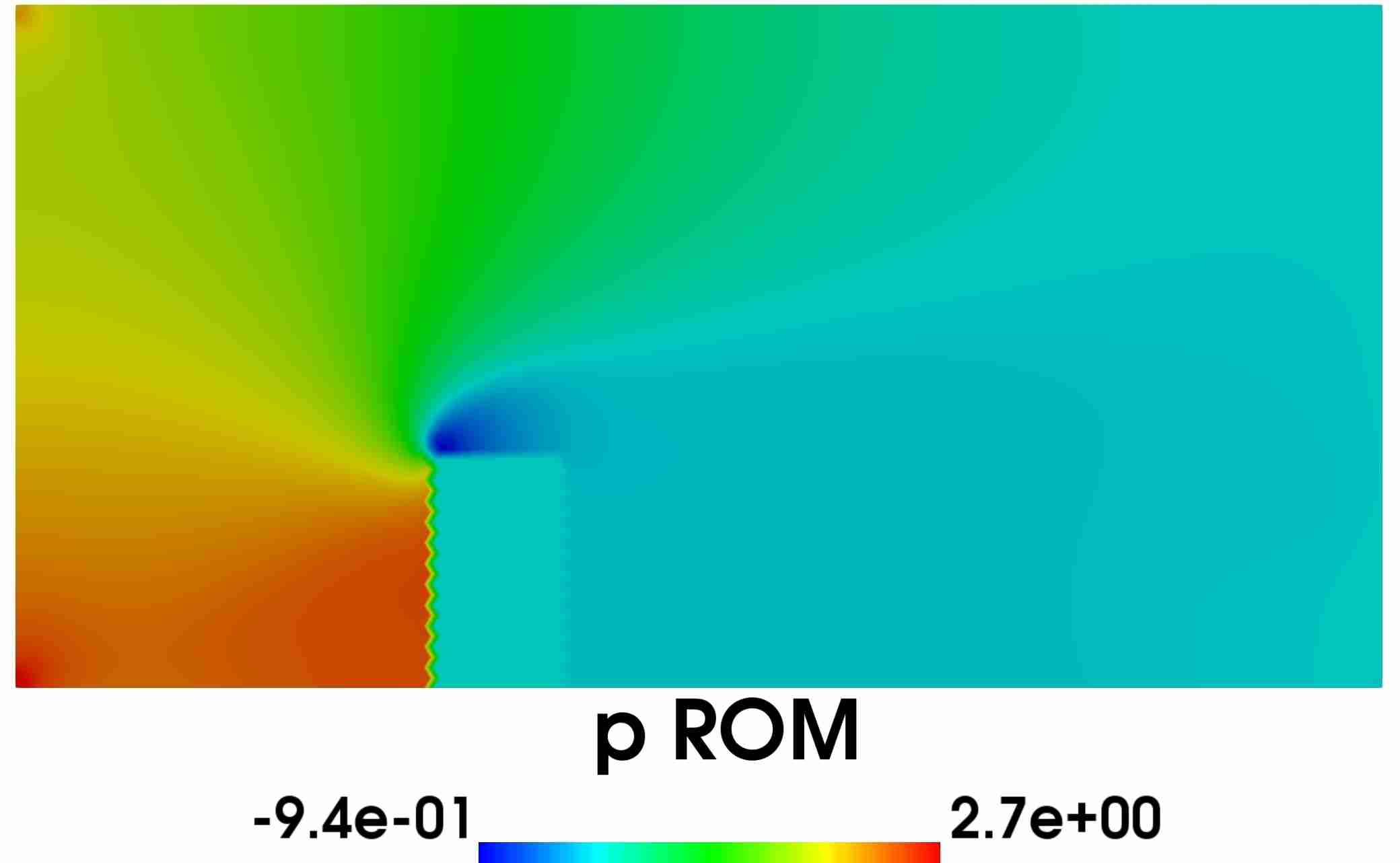}
\end{minipage}
\begin{minipage}{0.24\textwidth}
  \includegraphics[width=\textwidth]{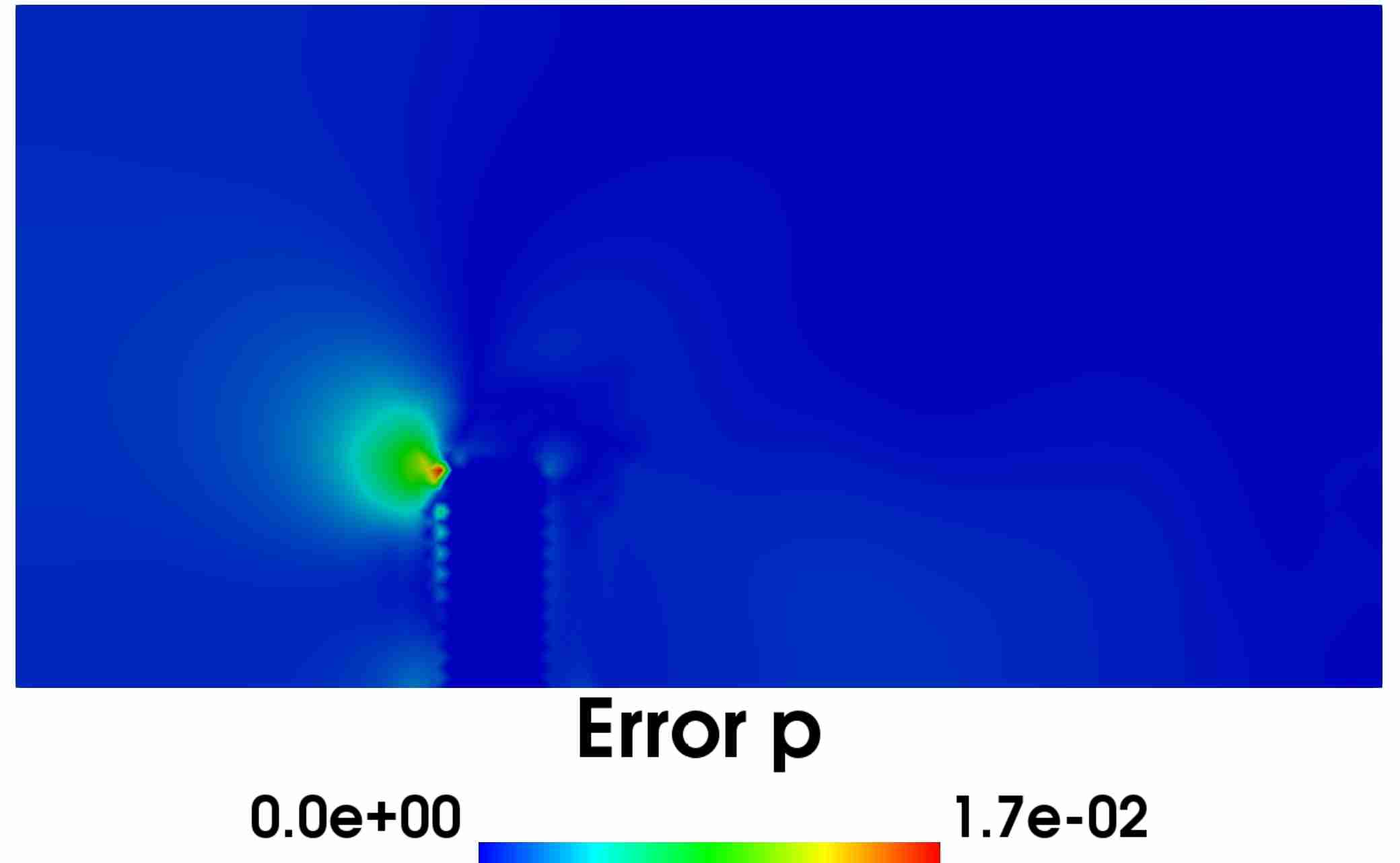} 
\end{minipage}
\begin{minipage}{0.24\textwidth}
  \includegraphics[width=\textwidth]{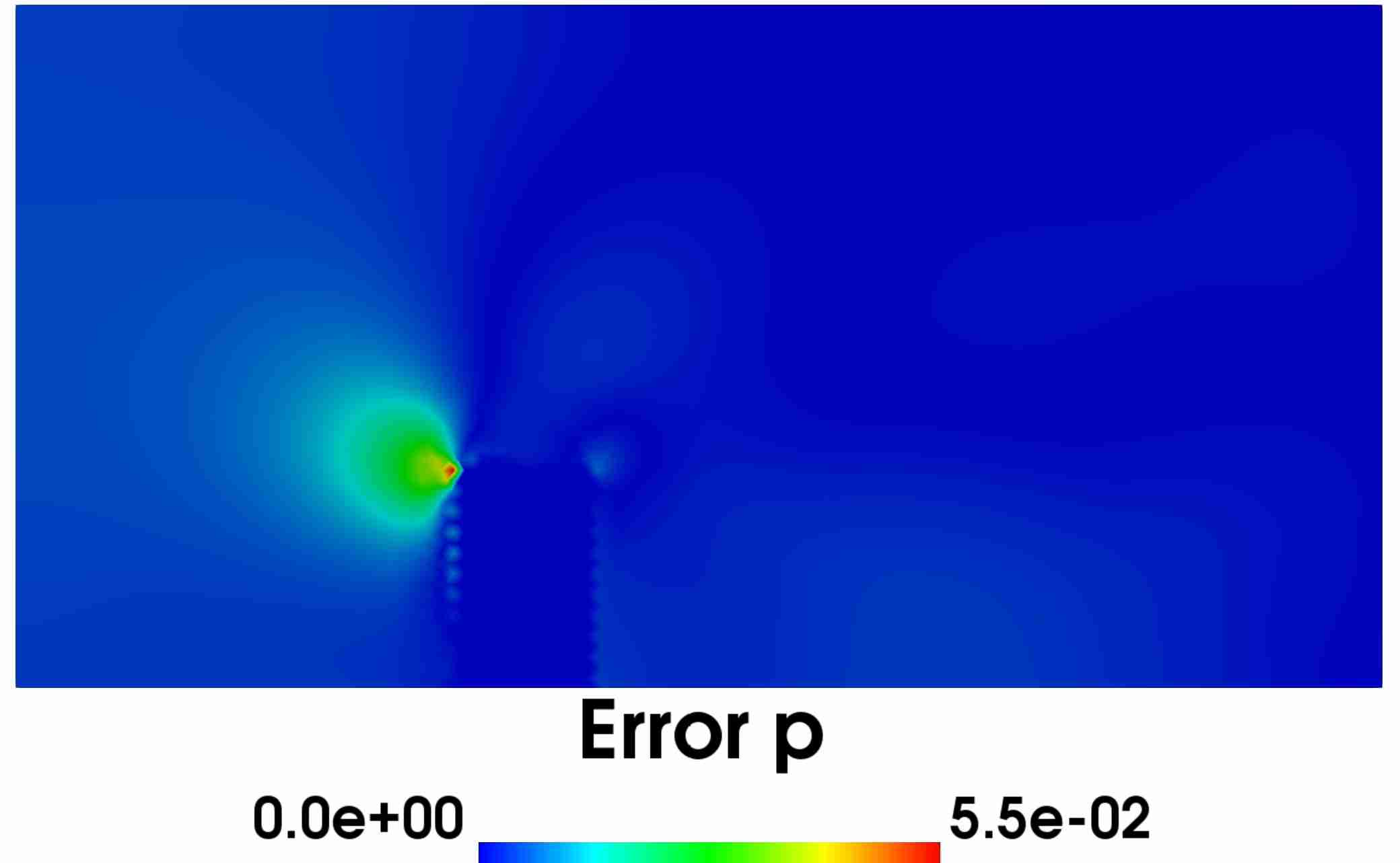}
\end{minipage}
\begin{minipage}{0.24\textwidth}
  \includegraphics[width=\textwidth]{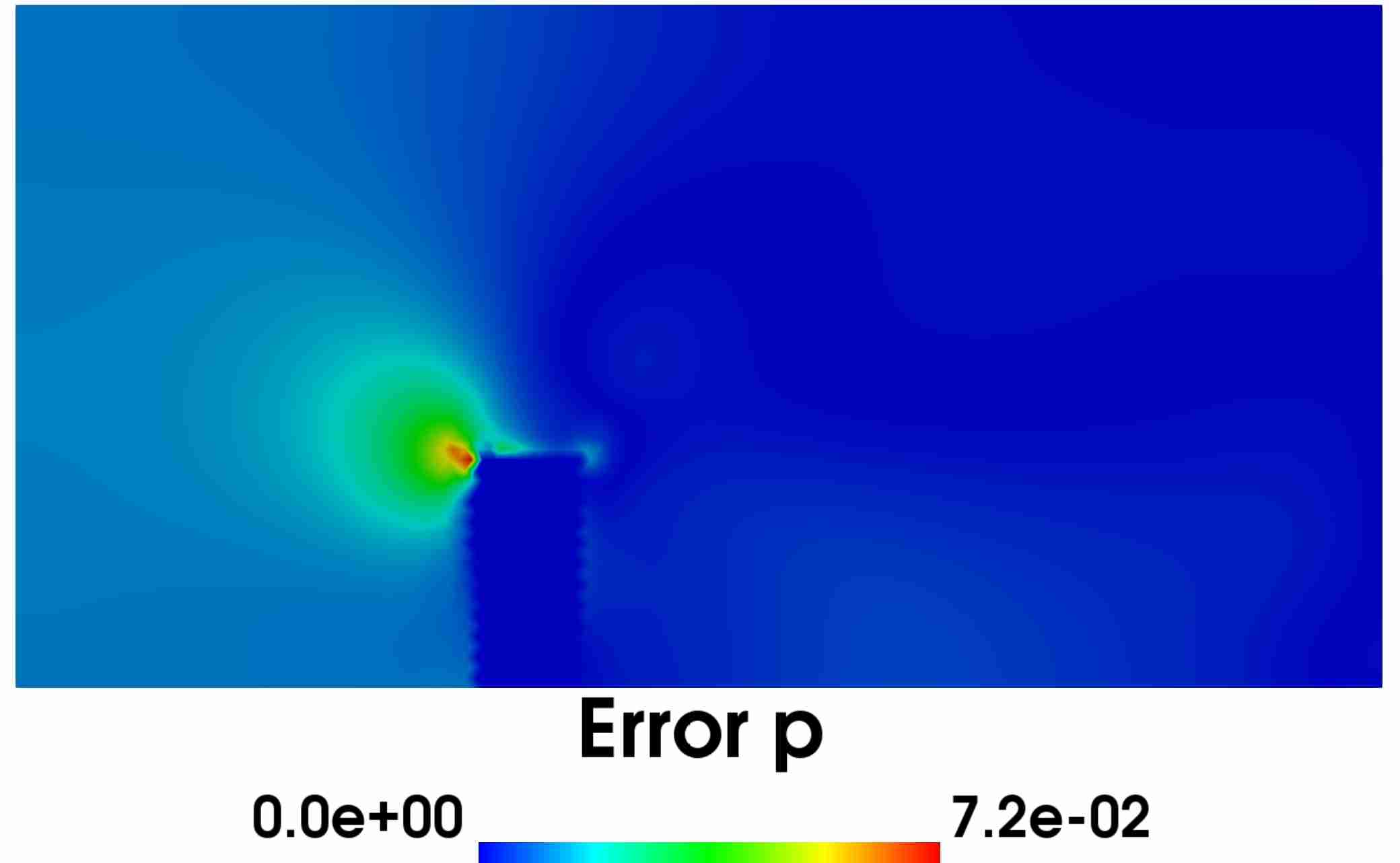}
\end{minipage}
\begin{minipage}{0.24\textwidth}
  \includegraphics[width=\textwidth]{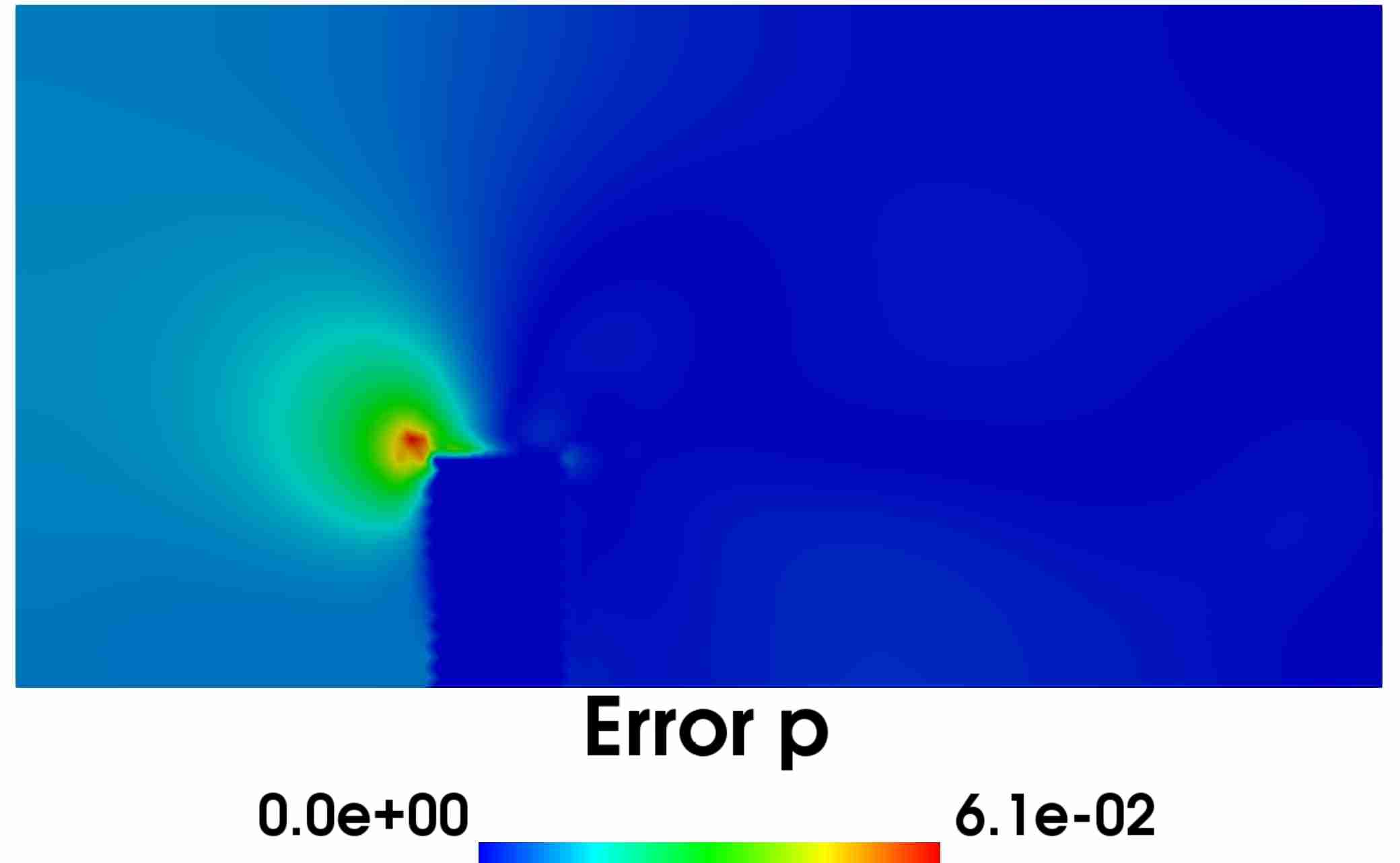}
\end{minipage}
\begin{minipage}{0.24\textwidth}
  \includegraphics[width=\textwidth]{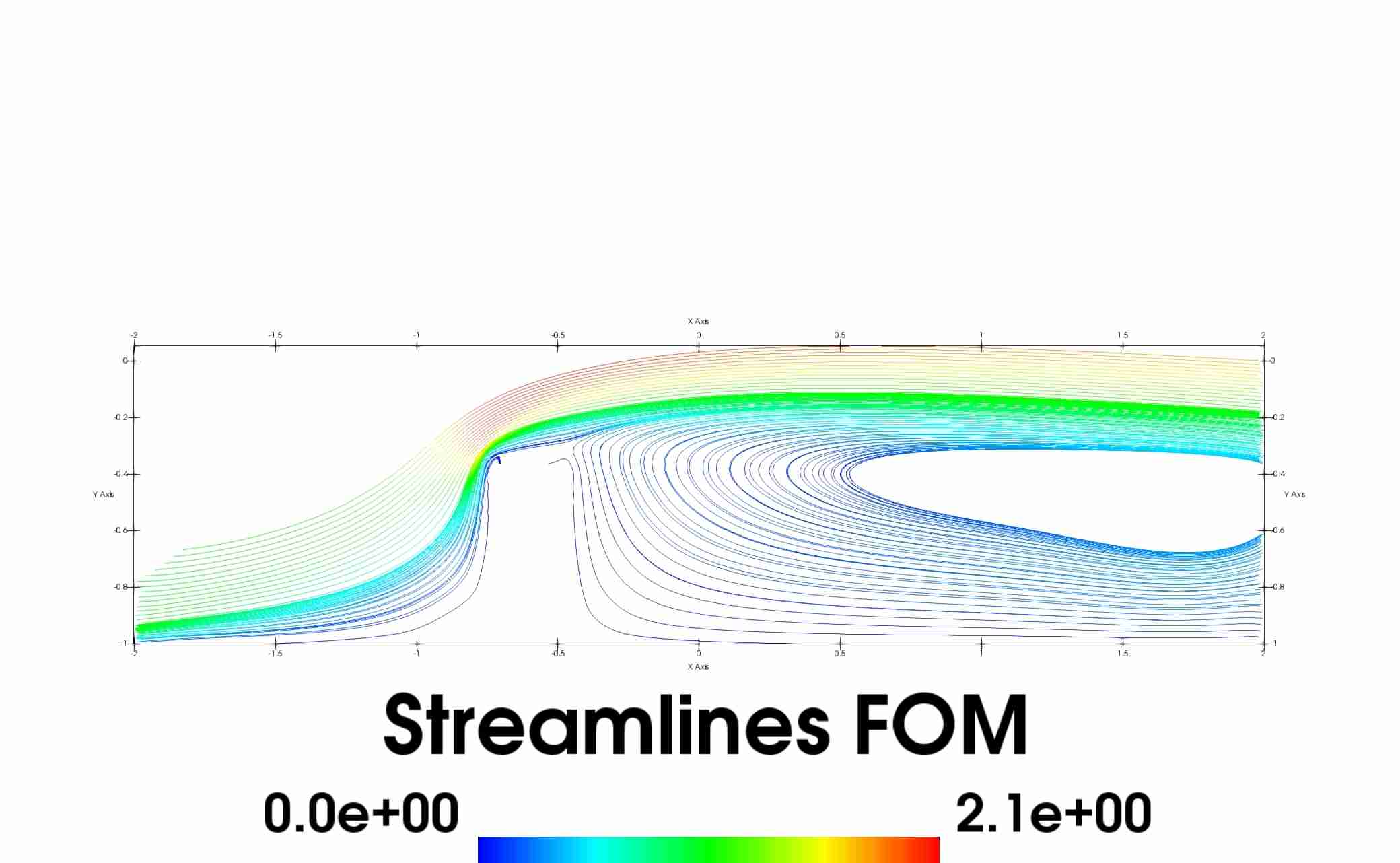} 
\end{minipage}
\begin{minipage}{0.24\textwidth}
  \includegraphics[width=\textwidth]{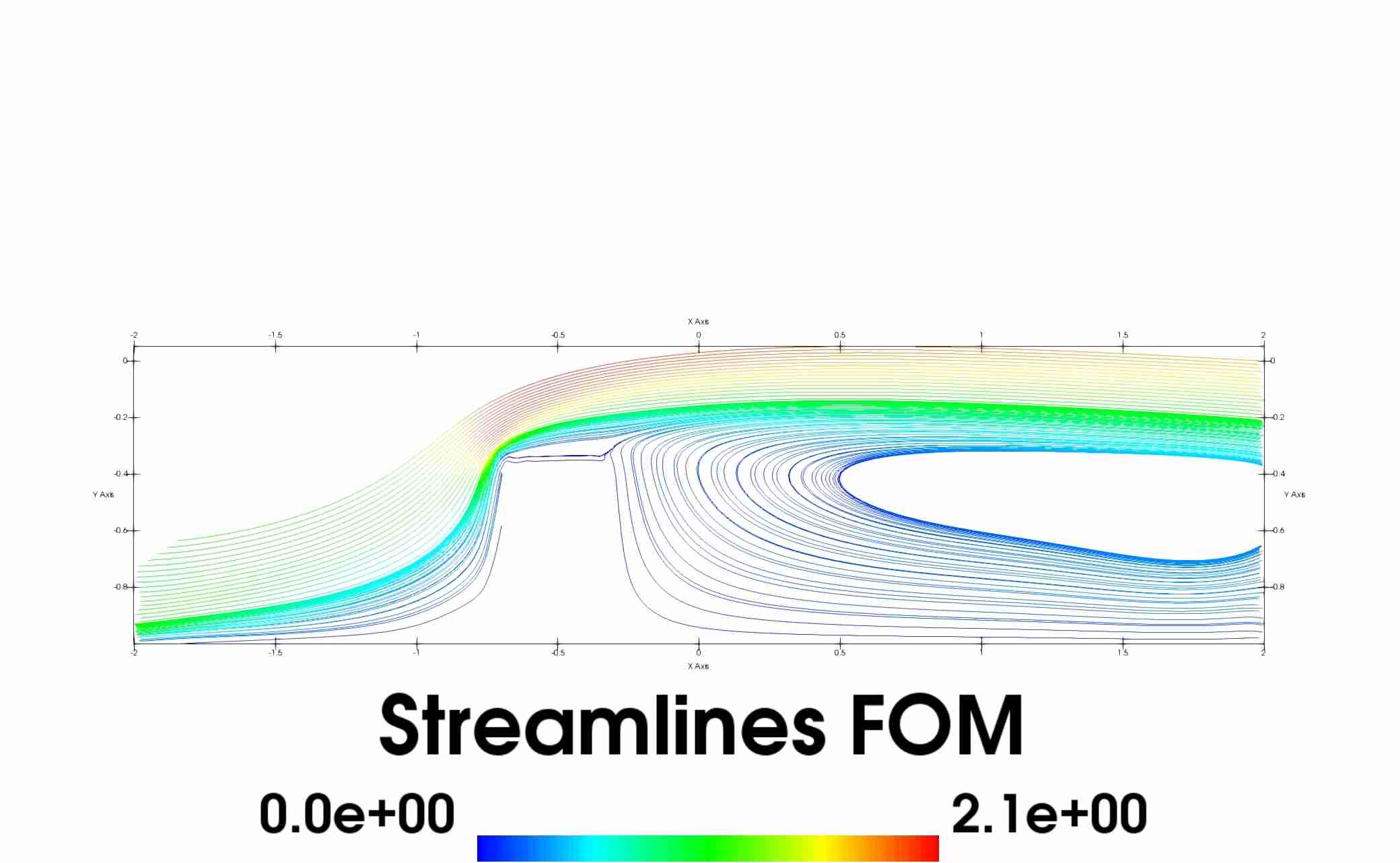}
\end{minipage}
\begin{minipage}{0.24\textwidth}
  \includegraphics[width=\textwidth]{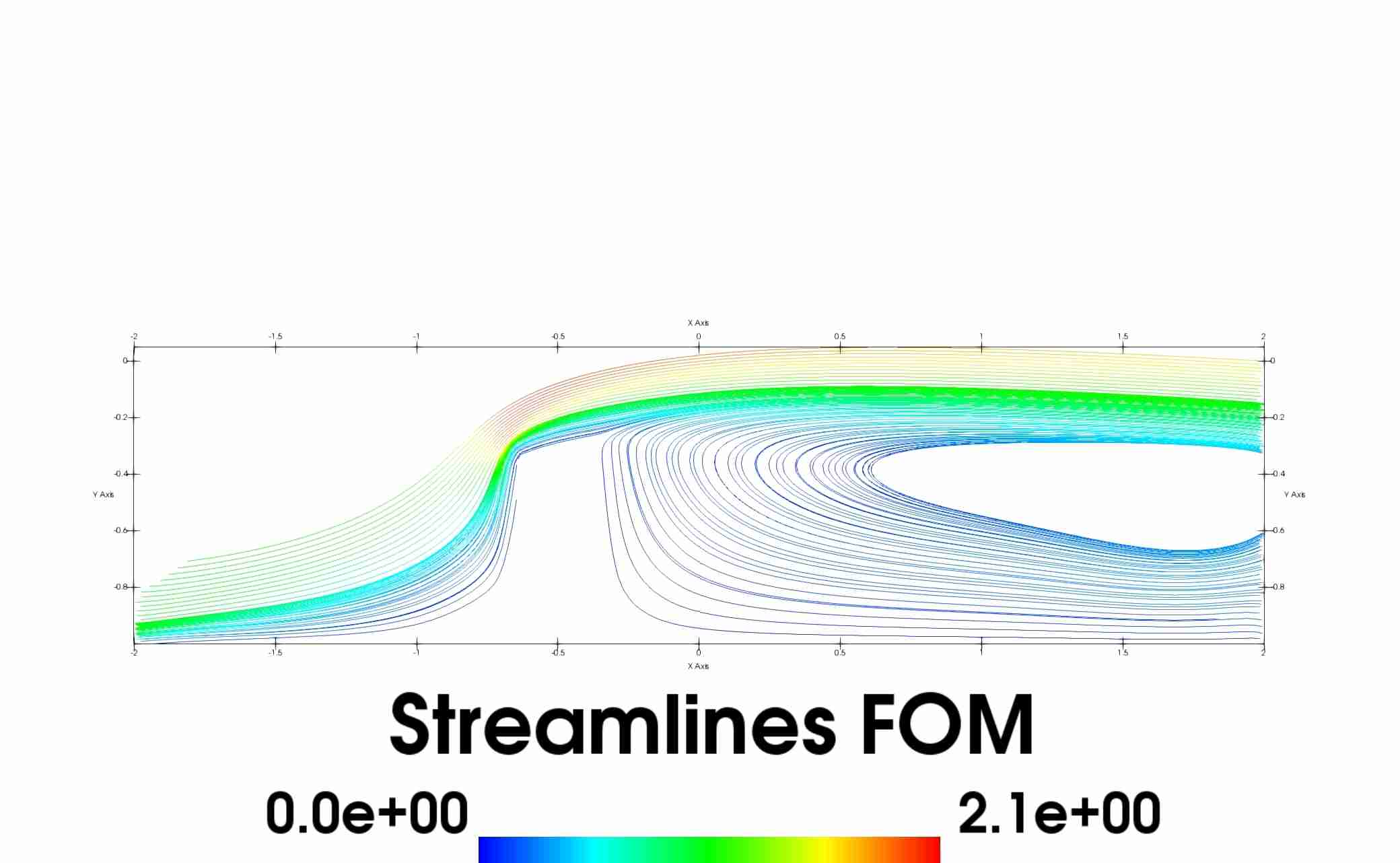}
\end{minipage}
\begin{minipage}{0.24\textwidth}
  \includegraphics[width=\textwidth]{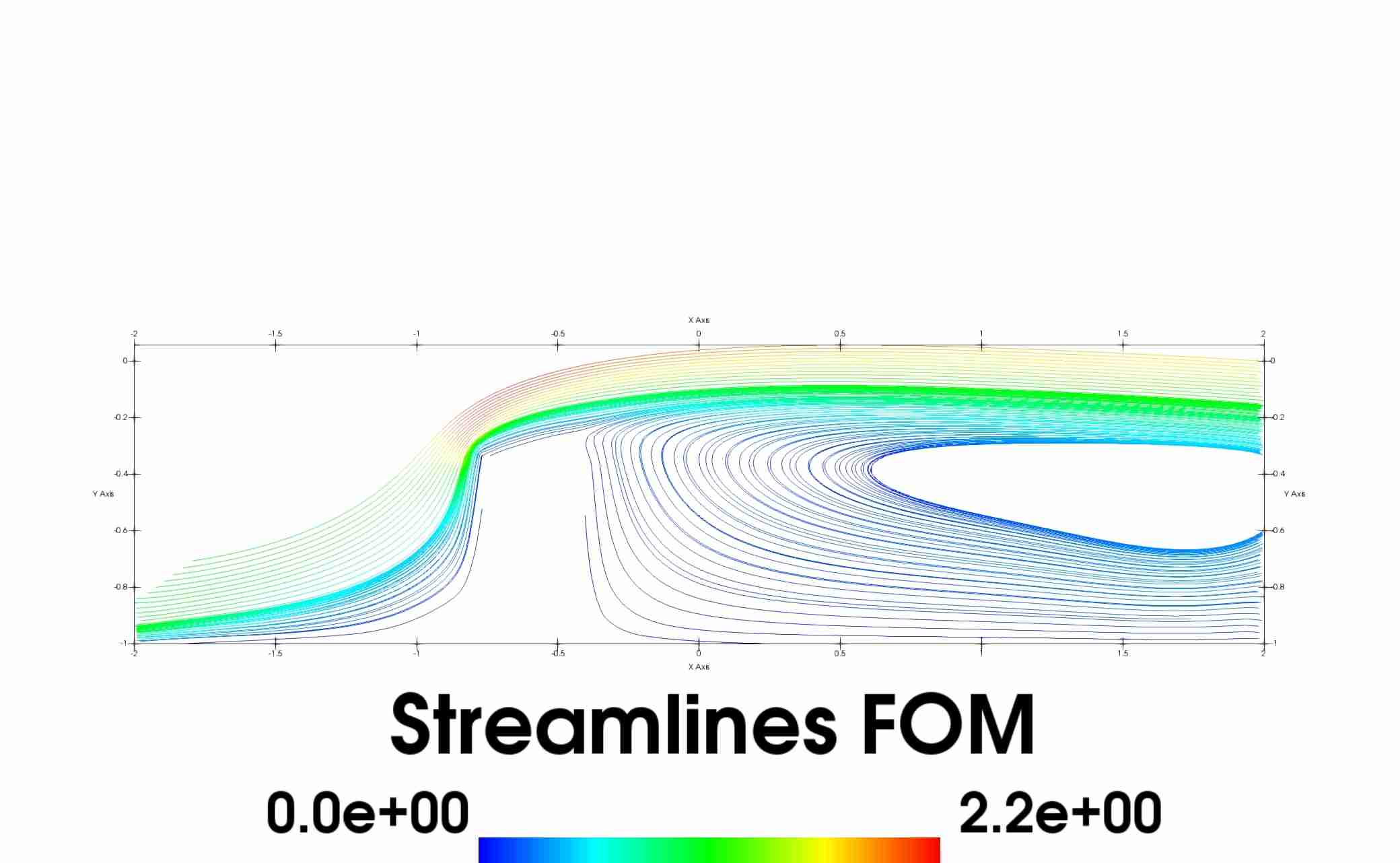}
\end{minipage}
\begin{minipage}{0.24\textwidth}
  \includegraphics[width=\textwidth]{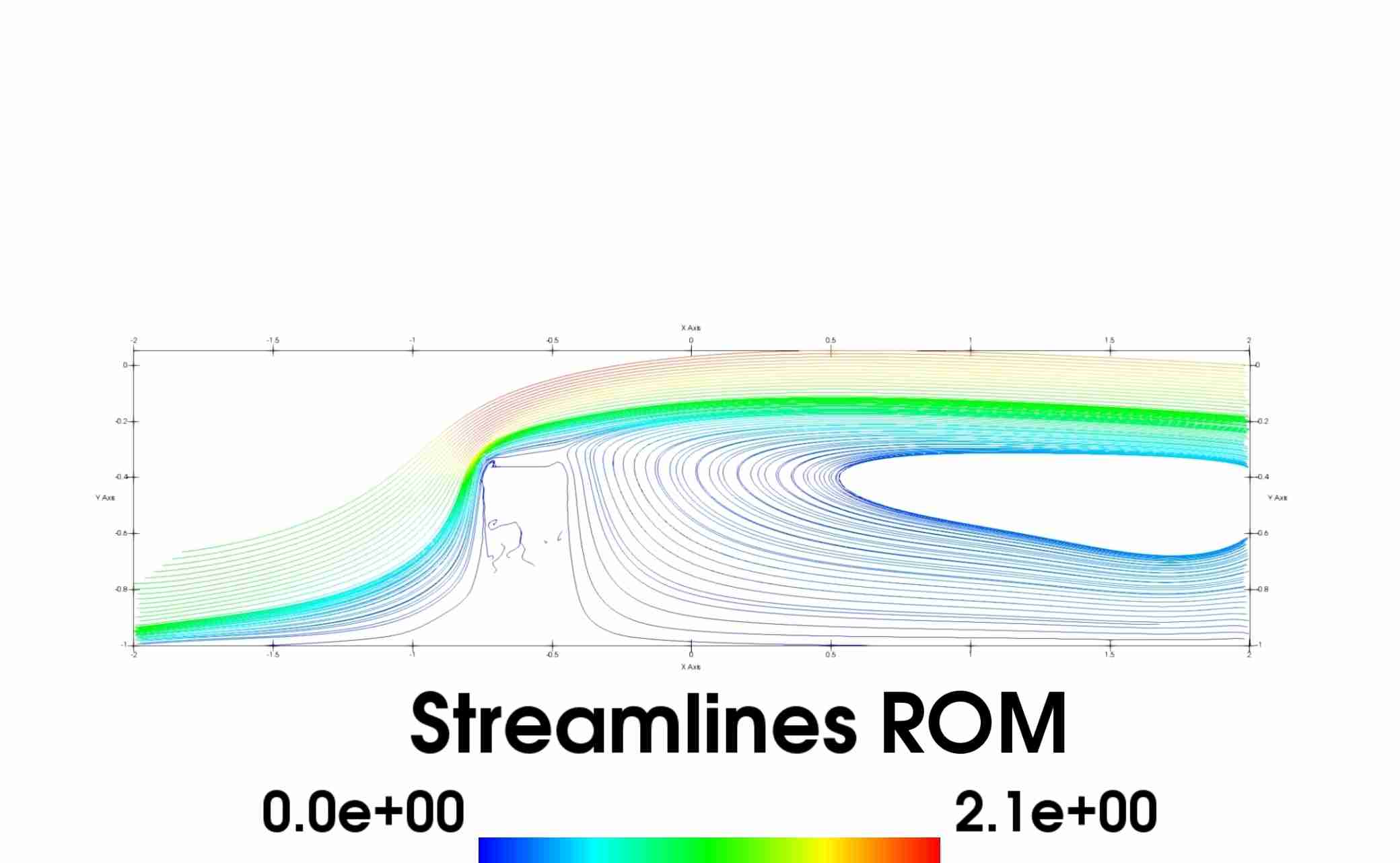} 
\end{minipage}
\begin{minipage}{0.24\textwidth}
  \includegraphics[width=\textwidth]{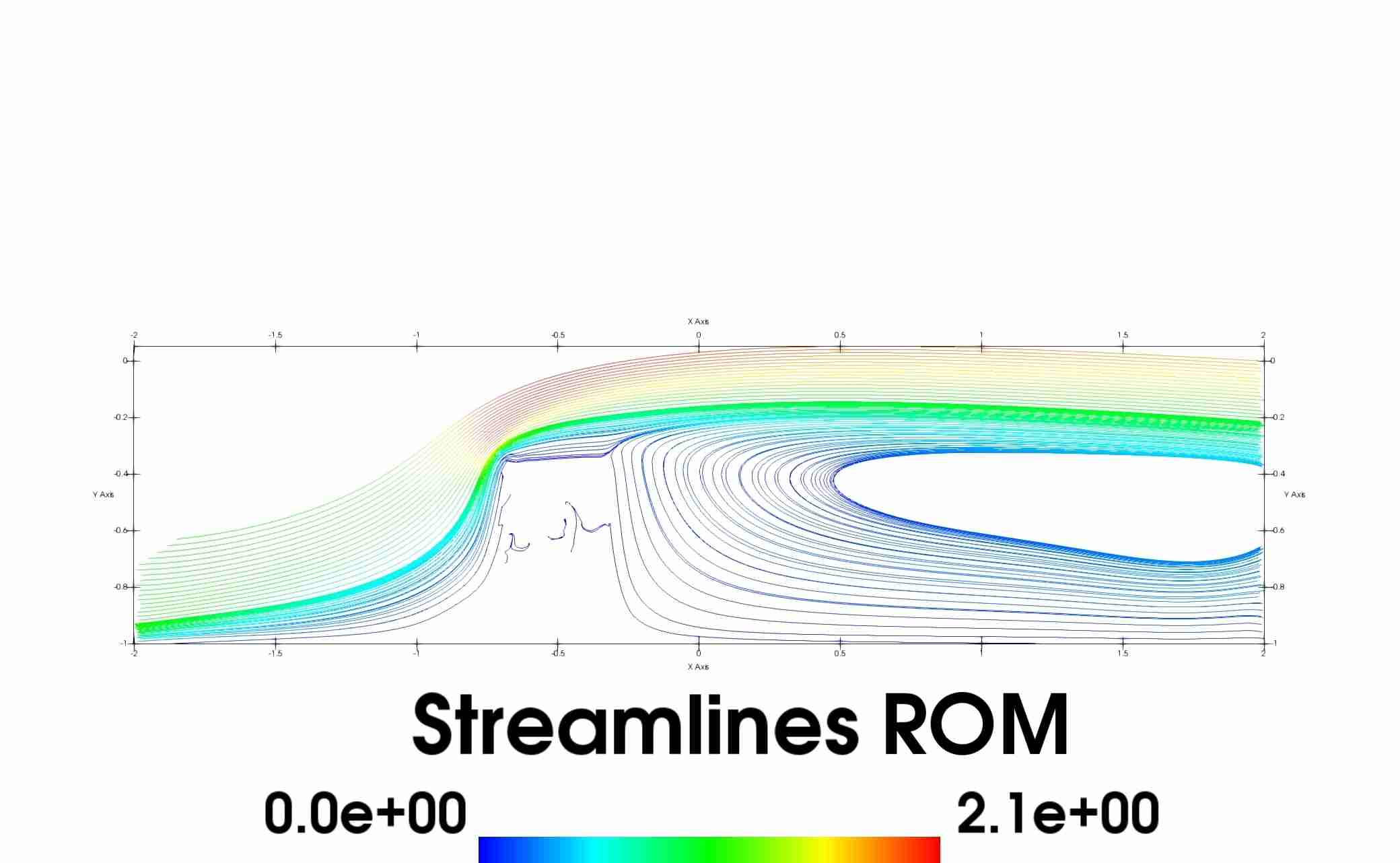}
\end{minipage}
\begin{minipage}{0.24\textwidth}
  \includegraphics[width=\textwidth]{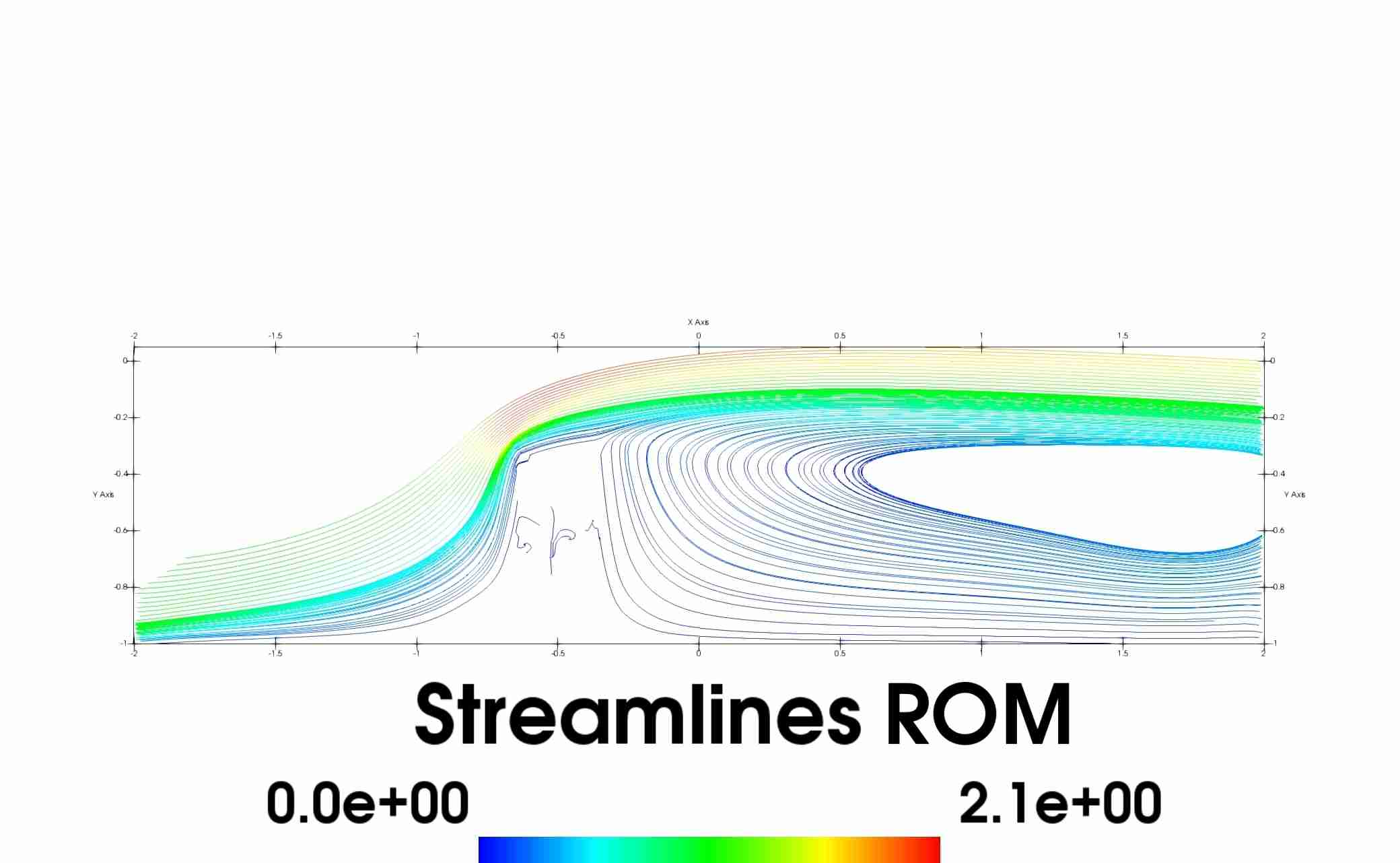}
\end{minipage}
\begin{minipage}{0.24\textwidth}
  \includegraphics[width=\textwidth]{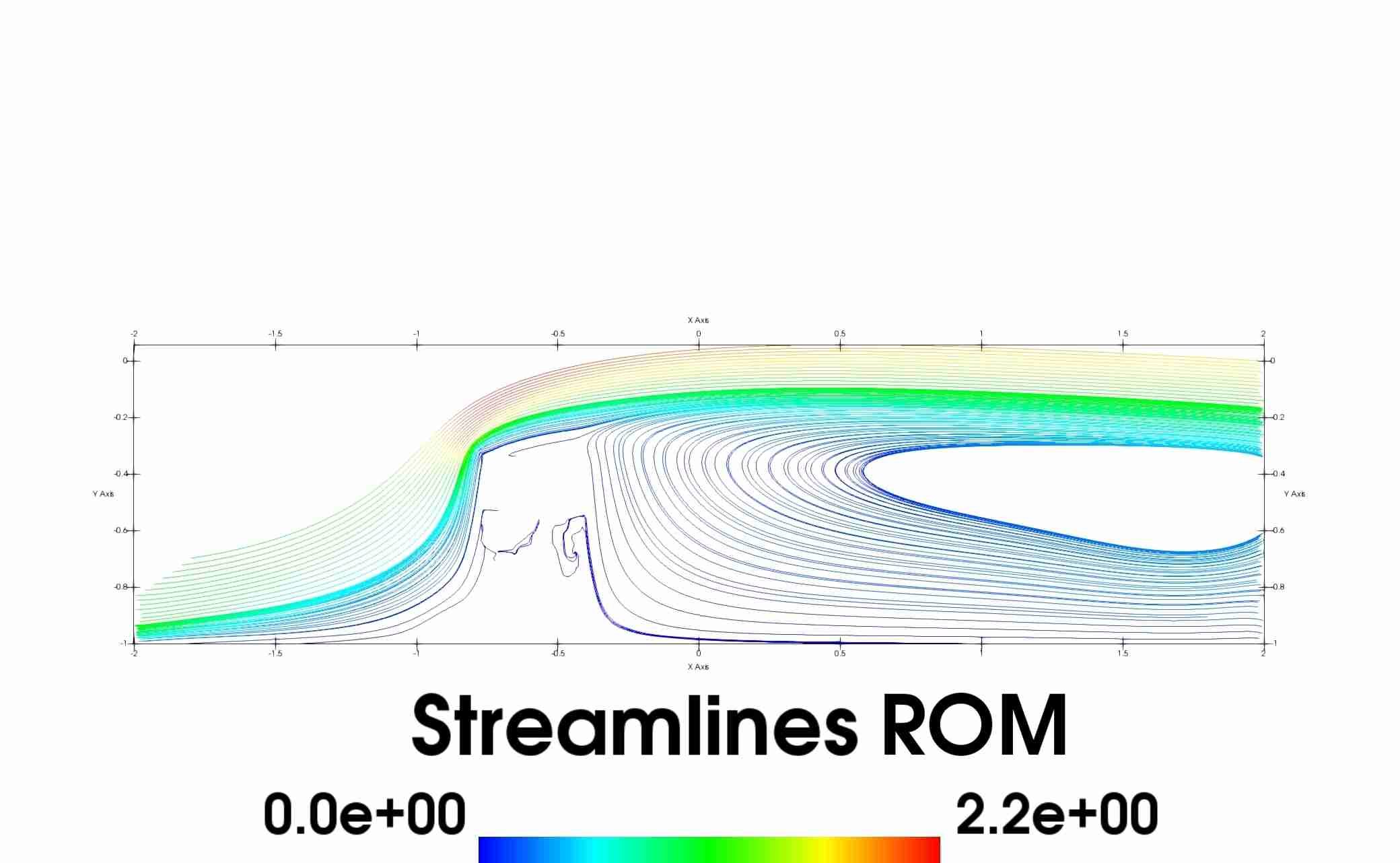}
\end{minipage}
\end{minipage}
\caption{Results for the geometrical parametrization using a three-dimensional parameter space with $(\mu_0,\mu_1, \mu_3) \in  [-0.6, -0.5]\times[0.3, 0.4]\times[1.3, 1.4]$.
In rows $1-3$ we report the full-order solution, the reduced order solution and the absolute error plots for the velocity field while in rows $4-6$ we report the same quantities for the pressure field. The different columns are for four different values of the input parameter $\mu=[(-0.5365, 0.3452, 1.3214)$,$(-0.5313, 0.3164, 1.3822)$,$(-0.5927, 0.3496, 1.3305)$, $(-0.5303, 0.3941, 1.3071)]$. Finally in rows $7-8$ the streamlines for the velocity in the full and in the reduced level are visualized.}
\label{fig:FOM_ROM_ERROR_3D}
\end{figure}
\begin{figure} 
\centering
\includegraphics[width=0.49\textwidth]{{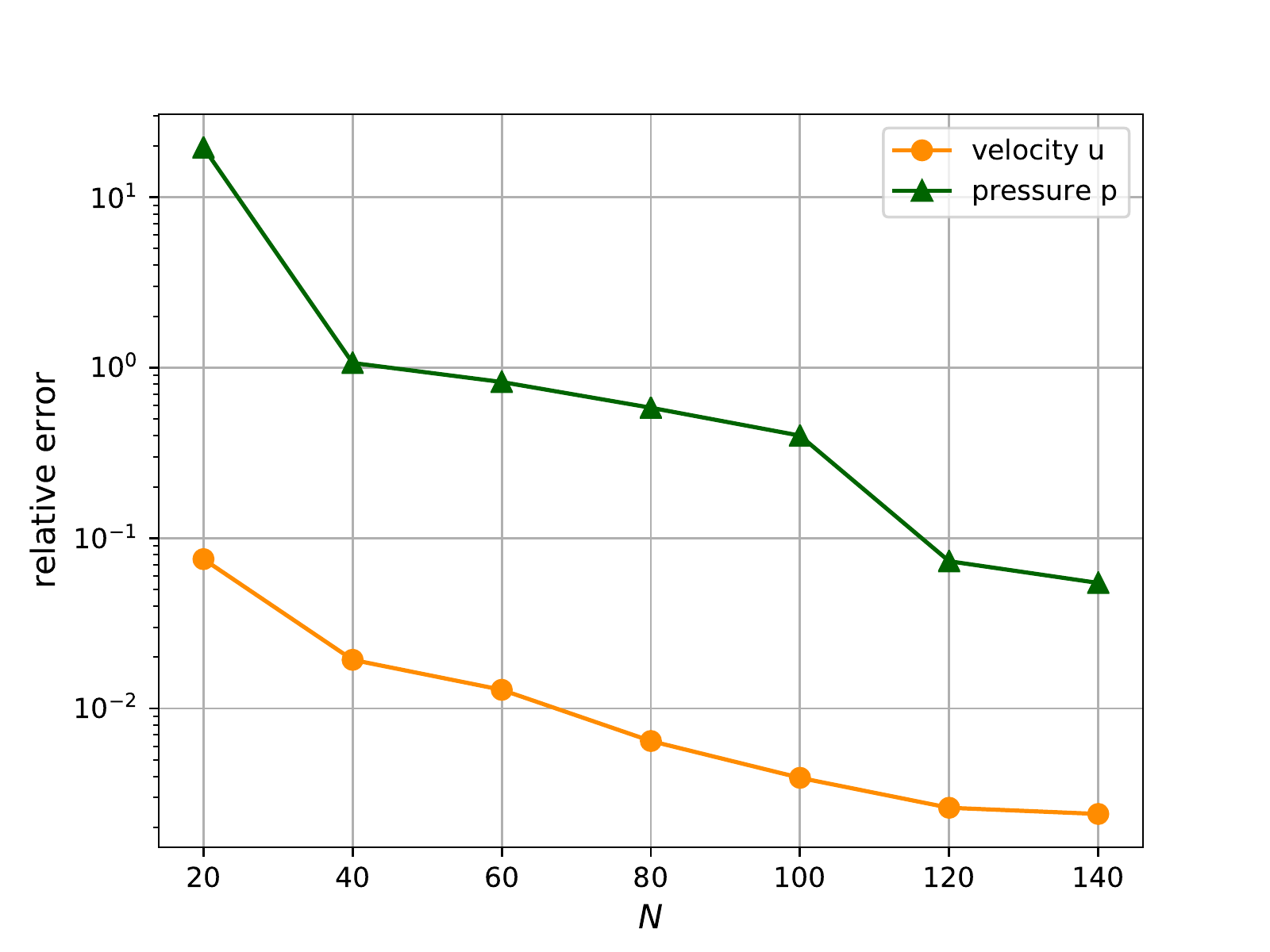}}
\caption{Visualization of the relative errors for velocity and pressure for the case with  geometrical parametrization with a one-dimensional parameter space for various number of modes.}
     \label{fig:errors1D}
\end{figure}

 
In this first experiment, the embedded domain consists of a rectangle of size $\mu_1\times\frac{\mu_2}{2} = 2\times \frac{\mu_2}{2}$ where its aspect ratio inside the domain is parametrized with a geometrical parameter describing the size of the rectangle embedded domain with respect to its y-length with  $\textit{aspect ratio}= \frac{4}{\mu_2}$ as in Figure \ref{background_mesh}. The horizontal x-length of the size of the box is not parametrized and the box's center is located on the left bottom corner of the domain $[-2,-1]$. The ROM has been trained with $600$ samples and tested onto $10$ samples, for a  parameter range  $\mu_2 \in [1, 2]$  chosen randomly inside the parameter space. To test the accuracy of the ROM we compared its results on $10$ additional samples that were not used to create the ROM and were selected randomly within the same range. Under these considerations we record in~\autoref{1DStokes_Pressure_and_Velocity_Components_Modes} the first six modes for the velocity magnitude and pressure, while in~\autoref{table:1Dparametrization} the respective relative errors $||{\bf u}_{FOM}-{\bf{u}}_{ROM}||_{L^2({{\mathcal{D}}})}/||{\bf{u}}_{FOM}||_{L^2({{\mathcal{D}}})}$ and $||p_{FOM} - p_{ROM}||_{L^2({{\mathcal{D}}})}/||p_{FOM}||_{L^2({{\mathcal{D}}})}$ are summarized and they are depicted in~\autoref{fig:errors1D}. In~\autoref{fig:FOM_ROM_ERROR_1D} we depict the FOM and ROM solutions together with the relative error for both velocity magnitude and pressure, as well as the FOM and ROM streamlines for the velocity vector, for four different values of the input parameter.
\begin{table} \centering
  \begin{tabular}{|c||cc|}
   \hline
    Snapshots:& \multicolumn{2}{c|}{600}   \\ \hline \hline
    Modes & \multicolumn{2}{c|}{relative error}  \\
    $N$ & $\bm u$  & $p$  \\
    \hline
     20  & 0.0752959 &  19.590245  \\
     40  & 0.0193090  & 1.0666177  \\
     60  & 0.0128775  & 0.8244898  \\
     80  & 0.0064458  & 0.5823617  \\
     100 & 0.0039197  & 0.3994657  \\
     120 & 0.0026156  & 0.0731744  \\
     140 & 0.0024058  & 0.0545712  \\
    \hline
  \end{tabular}   
  \caption{Relative error between the full-order solution and the reduced basis solution for velocity and pressure in the case of the geometrical parametrization with a one-dimensional parameter space. Results are reported for different dimensions of the reduced basis spaces.}
  \label{table:1Dparametrization}
\end{table} 

\subsection{A geometrical parametrization study with a three-dimensional parameter space}\label{Sec:3d_parametrisation}
\begin{figure} \centering
\begin{minipage}{\textwidth}
\centering
\begin{minipage}{0.48\textwidth}
(I)
\includegraphics[width=\textwidth]{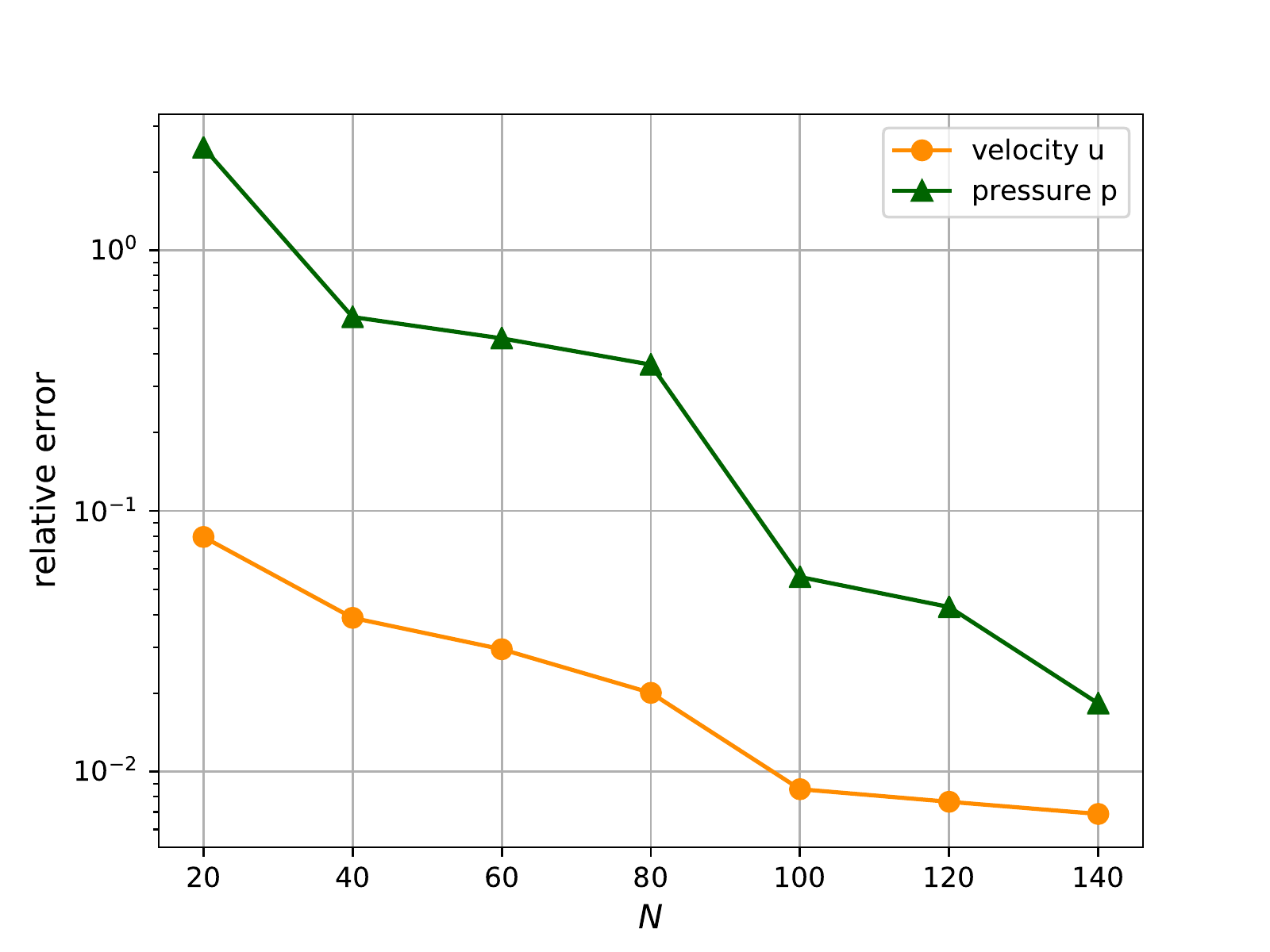}
\end{minipage}
\begin{minipage}{0.48\textwidth}
(II)
\includegraphics[width=\textwidth]{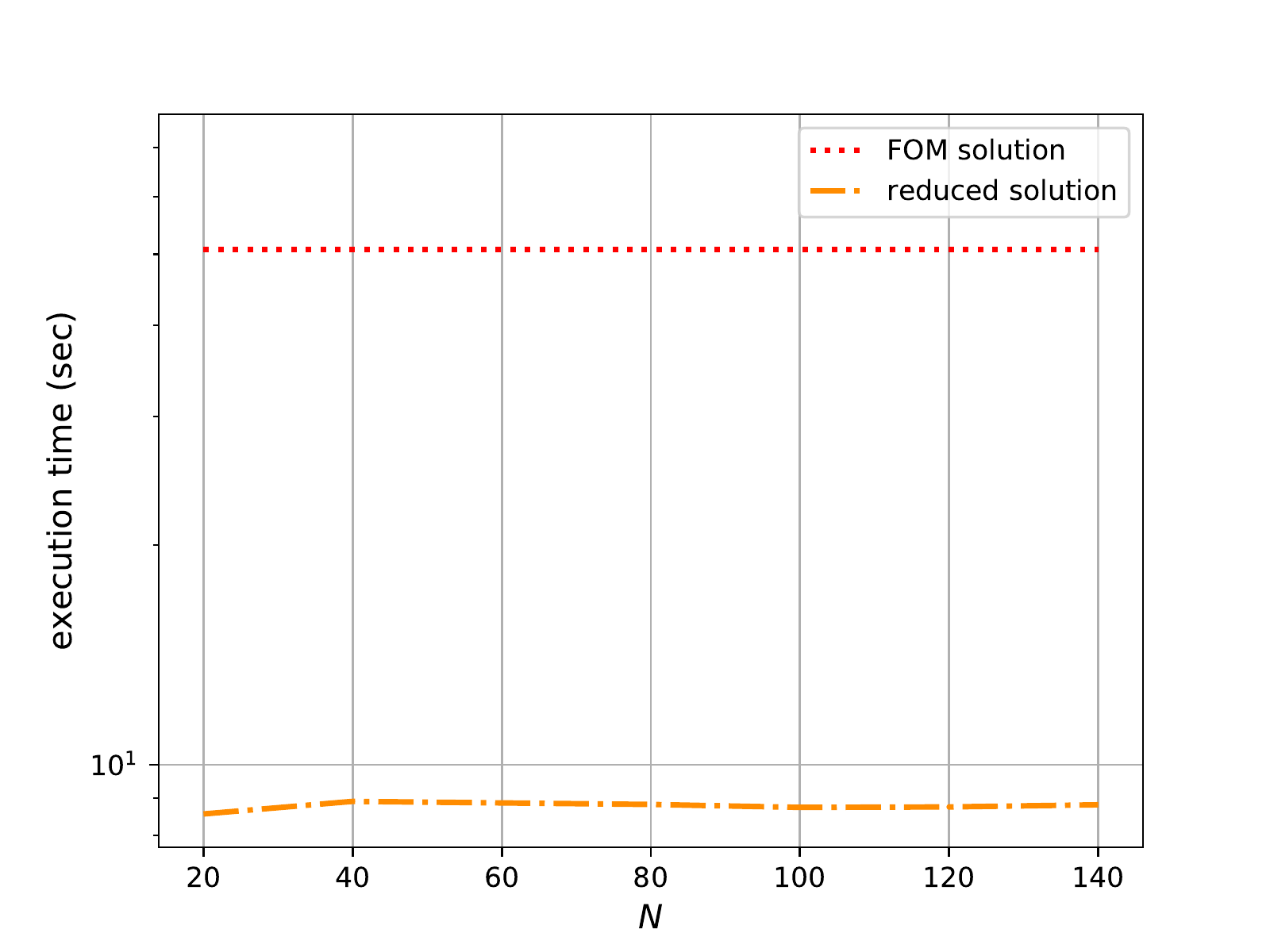}
\end{minipage}
\end{minipage}
\caption{Visualization of the results for the case of geometrical parametrization with a three-dimensional parameter space. On the left plot are depicted the relative errors for velocity and pressure. 
On the right plot we report the execution times of the reduced order problem for one random parameter value. 
In both plots, results are reported for various number of modes.}
\label{fig:3Dparametrization-IErrors-IItimers}
\end{figure}

The second case considers a geometrical parametrization with a three-dimensional parameter space in the range  $(\mu_0,\mu_1, \mu_2) \in  [-0.6, -0.5]\times[0.3, 0.4]\times[1.3, 1.4]$. We perform this test to examine the performances of the methodology on a more complex scenario where the box of the previous numerical test is parametrized onto the $x$ position of its center ($\mu_0$) and its aspect ratio ($\frac{\mu_2}{\mu_1}$) for both the $x$ and $y$ direction, namely width and height. The results are reported in~\autoref{3DStokes_Pressure_and_Velocity_Components_Modes} (the first six modes for the velocity magnitude and pressure), in~\autoref{table:3Dparametrization-IErrors-IItimers} I ($||{\bf u}_{FOM}-{\bf{u}}_{ROM}||_{L^2({{\mathcal{D}}})}/||{\bf{u}}_{FOM}||_{L^2({{\mathcal{D}}})}$ and $||p_{FOM} - p_{ROM}||_{L^2({{\mathcal{D}}})}/||p_{FOM}||_{L^2({{\mathcal{D}}})}$ relative errors report) and they are illustrated in~\autoref{fig:3Dparametrization-IErrors-IItimers} I. In~\autoref{fig:FOM_ROM_ERROR_3D} we visualize the FOM and ROM solutions together with the relative errors for both velocity magnitude and pressure, as well as the FOM and ROM streamlines for the velocity, for four different values of the input parameters. Finally in~\autoref{table:3Dparametrization-IErrors-IItimers} II and ~\autoref{fig:3Dparametrization-IErrors-IItimers} II the time savings are reported and visualized. 
%
%
\begin{table} \centering
  \begin{tabular}{|c||cc||c|c|
}
   \hline
     Snapshots: & \multicolumn{3}{c|}{900} \\\hline\hline
       Modes & \multicolumn{2}{c||}{(I) relative error} & (II) execution time \\
     $N$ & $\bm u$ & $p$ & (sec)\\
    \hline
     20 &   0.0794779  & 2.4790682  & 8.5617990  \\
     40 &    0.0388933  & 0.5546928 & 8.9083982  \\
     60 &    0.0294700  & 0.4594517 & 8.8653870  \\
     80 &    0.0200465  & 0.3642104 & 8.8223756  \\
     100 &   0.0085595  & 0.0558062 & 8.7412506  \\
     120 &   0.0076601  & 0.0428119 & 8.7539559  \\
     140 &   0.0068780  & 0.0182534 & 8.8129689  \\
  \hline
  \end{tabular}
  \caption{(I) The relative error between full-order solution and reduced basis solution for velocity and pressure for $900$ snapshots in the case with a geometrical parametrization using three-dimensional parameter space. (II) Execution time at the reduced order level. The computation time includes the assembling of the full-order matrices, their projection and the resolution of reduced problem. Times are for the resolution of one random value of the input parameter. The time execution at full-order level is equal to $\approx 50.70$ sec.}
  \label{table:3Dparametrization-IErrors-IItimers}
\end{table} 
\subsection{Some Comments} 
Additional testing using more snapshots and with and without supremizers  yielded very similar results, in which supremizers delivered slightly worse velocity results and slightly improved pressure results. We found that the supremizers did not substantially improved the errors. In our opinion, this is a consequence of the incremental iterative formulation in the offline solver, which preserves the effects of the SUPG and PSPG stabilization in the reduced model, which allows an inf-sup stable reduced basis. Some basics related to these stability issues, and  numerical results can be found in {Appendix} A.
Both  the  full-order  and  the  ROM  simulation  were  run  in  serial  on  an  Intel\textsuperscript{\textregistered} Core\textsuperscript{TM}  i7-4770HQ 3.70GHz CPU.
{
\section{Conclusions and future developments}\label{sec:conclusions}
We have introduced a geometrically parametrized ROM model emulator of the two-dimensional Navier--Stokes equations, of much reduced computational cost. The ROM model evaluation through numerical tests shows good convergence properties and low errors. Comparing with high-fidelity solutions, numerical ROM errors improve with the increase of the size of training data and of the number of basis components. As future perspectives, we indicate applications to the time dependent Navier--Stokes systems, fluid structure interaction problems, and shallow water flows. Additionally, from the model reduction point of view, we will pursue further developments in hyper-reduction techniques \cite{Xiao20141,BARRAULT2004667,Carlberg2013623, stabile_geo_}  and transportation methodologies \cite{KaBaRO18}.
}
\section*{Acknowledgments}
This research has been supported by the U.S. Department of Energy, Office of Science, Advanced Scientific Computing Research under Early Career Research Program Grant SC0012169 and the Army Research Office (ARO) under Grant W911NF-18-1-0308, Hellenic Foundation for Research and Innovation (HFRI) and  the  General  Secretariat  for  Research  and  Technology (GSRT),  
under  grant agreement No 1115, the National Infrastructures for Research and Technology S.A. (GRNET S.A.) under project ID pa190902, the European Research Council Executive Agency by means of the H2020 ERC Consolidator Grant project AROMA-CFD ``Advanced Reduced Order Methods  with  Applications  in  Computational  Fluid  Dynamics'' - GA  681447, (PI: Prof. G. Rozza), INdAM-GNCS-2019, and by project FSE - European Social Fund - HEaD "Higher Education and Development" SISSA operazione 1, Regione Autonoma Friuli - Venezia Giulia.
%
\section*{Appendix A}
In this appendix we explore the experiment as described in Section \ref{Sec:3d_parametrisation}. We will justify why in all previous experiments we did not use any additional RB stabilization, verifying numerically that the SUPG and PSPG stabilization which is applied on the high fidelity solver is strongly propagating through the reduced basis construction procedure to the reduced level. In contrast, we refer the interested reader to the classical works of~\cite{ballarin2015supremizer,RoVe07, stabile_stabilized, StaHiMoLoRo17} where the reduced solution stabilization is deemed necessary. 
Subsequently, we introduce the basics related to the  stability issues that could appear in the offline and online stage. In the full order method stage it is well known that the spaces have to  satisfy the, also parametric in our case, Ladyzhenskaya-Brezzi-Babuska ``inf-sup''  condition see e.g.  \cite{boffi_mixed}. In particular, it is required that there should exist a constant $\beta > 0 $, independent to the discretization parameter $h$, such that $$ \inf_{0\neq q \in {Q_h}} \sup_{0\neq \bm{v} \in {\bf V}_h} \frac{\langle \nabla \cdot \bm{v}, q \rangle}{||{\nabla \bm{v}||} \, ||{q}|| }\ge \beta > 0.$$ In the present work this condition is fulfilled for the high fidelity solution through the SUPG and PSPG stabilization. 
Even if, the offline stage and the snapshots are realized and computed by a stabilized numerical method, it is not guaranteed that this stability is preserved onto the reduced basis spaces  \cite{RoVe07,Gerner2012,ballarin2015supremizer}.
Next we briefly introduce the ``inf-sup'' condition enforcement in the reduced level using supremizers, we illustrate the relative errors results and we are comparing them with the case without supremizers enrichment.    
%
Within this approach, the velocity supremizer basis functions $\bm{L}_{\text{sup}} =[\bm{\eta}_1, \dots, \bm{\eta}_{N_{\text{sup}}^r}] \in \mathbb{R}^{N_u^h \times N_{\text{sup}}^r}$ are constructed and added to the reduced velocity space (see Section \ref{subsec_POD_theory}) which is finally transformed into 
$$ 
 \tilde{\bm{L}}_u = [\bm{\varphi}_1,\dots,\bm{\varphi}_{N_u^r},\bm{\eta}_1,\dots,\bm{\eta}_{N_{\text{sup}}^r}] \in \mathbb{R}^{N_u^h \times (N_u^r+N_{\text{sup}}^r)}.$$
To obtain the latter enrichment for each pressure basis function $\chi_i$ the auxiliary ``supremizer" problem: 
$$\Delta {\bm{s}_i} = - {\nabla} \chi_i\mbox{ in } \mathcal D(\mu ^i), \quad{\bm{s}_i}=\bm{0}\mbox{ on } \Gamma(\mu ^i)$$ 
is solved with an SBM Poisson solver starting from the parameter value $\mu^i$. 
For each pressure basis function the corresponding supremizer element can be found and the solution ${\bm{s}_i}$ permits the realization of the ``inf-sup'' condition. 

We emphasize that the above supremizer basis functions do not depend on the particular pressure basis functions and on the geometrical parameters,  they are computed during the offline phase, and their calculation is based on the pressure snapshots.
\begin{table} \centering
  \begin{tabular}{|c||cc|}
   \hline
    Snapshots:& \multicolumn{2}{c|}{900}   \\ \hline \hline
    Modes & \multicolumn{2}{c|}{relative error}  \\
    $N$ & $\bm u$  & $p$  \\
    \hline
     20  & 5.2407257 & 46.628631  \\
     40  & 0.0699833 & 0.3109475 \\
     60  & 0.0537041 & 0.2613186 \\
     80  & 0.0374247 & 0.2116895 \\
     100 & 0.0107190 & 0.0365011 \\
     120 & 0.0093516 & 0.0237161 \\
     140 & 0.0062227 & 0.0132555 \\
    \hline
  \end{tabular}   
  \caption{Supremizers ``inf-sup'' condition  enrichment: Relative error between the full-order solution and the reduced basis solution for velocity and pressure in the case of the geometrical parametrization with a three-dimensional parameter space. Results are reported for different dimensions of the reduced basis spaces.}
  \label{table:S31Dparametrization}
\end{table} 
\begin{figure} 
\centering
\includegraphics[width=0.49\textwidth]{{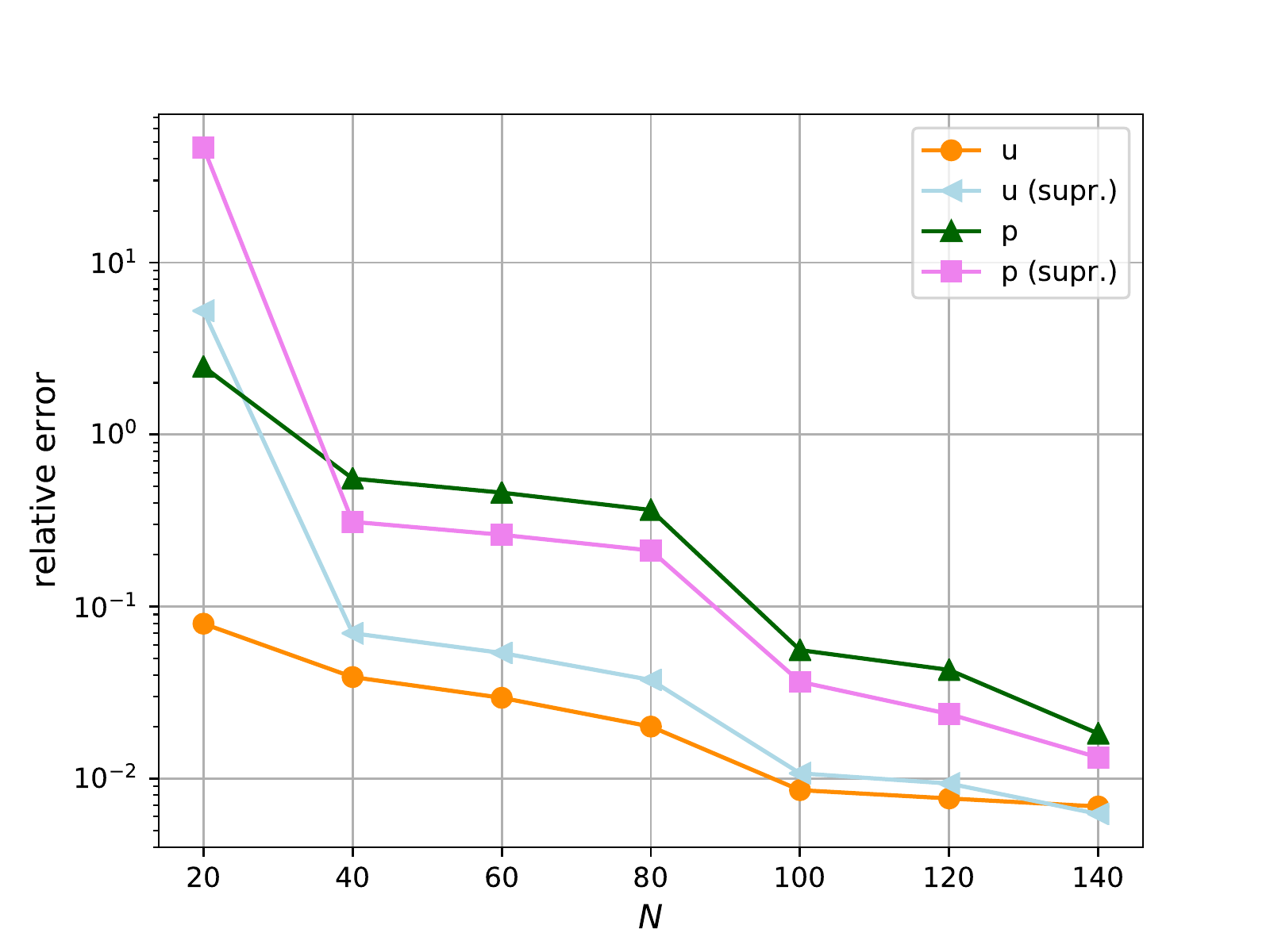}}
\caption{Supremizers ``inf-sup'' condition  enrichment: Visualization of the relative errors for velocity and pressure for the case with  geometrical parametrization and a three-dimensional parameter space for various number of modes.}
     \label{fig:errorsS3D}
\end{figure}
%
%
Obviously, if someone  compares the~\autoref{table:S31Dparametrization} with~\autoref{table:3Dparametrization-IErrors-IItimers} and examines their visualization in \autoref{fig:errorsS3D},  supremizers drove the reduced solution to slightly worse velocity results and slightly improved pressure results. So, the supremizers did not substantially improved the errors and this is the reason that we avoided their application
. In our opinion, this phenomenon is a consequence of the incremental iterative formulation in the offline solver, which preserves the effects of the SUPG and PSPG stabilization in the reduced model and allows an inf-sup stable reduced basis. 
%
\bibliography{bibfile_sissa}
\end{document}